\documentclass{lmcs} 
\pdfoutput=1
\usepackage[utf8]{inputenc}

\usepackage{lastpage}
\lmcsdoi{22}{1}{23}
\lmcsheading{}{\pageref{LastPage}}{}{}%
{Feb.~21,~2025}{Mar.~10,~2026}{}

\keywords{Differential Linear Logic, Differential Categories, Codereliction, Additive Enrichment, Coalgebra Modality, Bialgebra Modality, Hopf Coalgebra Modality}

\usepackage{xspace}
\usepackage{amsmath,amssymb,wasysym}
\usepackage[all]{xy}
\usepackage{proof}
\usepackage{mathtools}
\usepackage{latexsym}
\usepackage{dsfont}
\usepackage{multicol}
\usepackage{cmll}
\usepackage{multirow}
\usepackage{longtable}
\usepackage{array}

\usetikzlibrary{shapes.geometric, arrows}
\usetikzlibrary{positioning,arrows}
\usetikzlibrary{arrows.meta,calc}
\pgfdeclarelayer{edgelayer}
\pgfdeclarelayer{nodelayer}
\pgfsetlayers{edgelayer,nodelayer,main}

\tikzstyle{none}=[inner sep=0pt]
\tikzstyle{port}=[inner sep=0pt]
\tikzstyle{component}=[circle,fill=white,draw=black, inner sep=2.5pt]
\tikzstyle{integral}=[inner sep=0pt]
\tikzstyle{differential}=[inner sep=0pt]
\tikzstyle{codifferential}=[inner sep=0pt]
\tikzstyle{function}=[regular polygon,regular polygon sides=4,fill=white,draw=black]
\tikzstyle{function2}=[regular polygon,regular polygon sides=4,fill=white,draw=black, inner sep=2pt]
\tikzstyle{function3}=[regular polygon,regular polygon sides=4,fill=white,draw=black, inner sep=-2pt]
\tikzstyle{duplicate}=[circle,fill=white,draw=black, inner sep=1pt]
\tikzstyle{wire}=[-,draw=black,line width=1.000]
\tikzstyle{exwire}=[-,draw=black,dashed,line width=1.000]
\tikzstyle{object}=[inner sep=2pt]
\tikzset{every picture/.append style={scale=.65}, transform shape}

\def\eg{{\em e.g.}}
\def\cf{{\em cf.}}

\allowdisplaybreaks
\begin{document}

\title[Additive Enrichment from Coderelictions]{Additive Enrichment from Coderelictions}

\author[J.-S. P. ~Lemay]{Jean-Simon Pacaud Lemay\lmcsorcid{0000-0003-4124-3722}}

\address{Macquarie University, Sydney, New South Whales, Australia}	
\email{js.lemay@mq.edu.au}  






\begin{abstract}
  \noindent Differential linear categories provide the categorical semantics of the multiplicative and exponential fragments of Differential Linear Logic. Briefly, a differential linear category is a symmetric monoidal category that is enriched over commutative monoids (called additive enrichment) and has a monoidal coalgebra modality that is equipped with a codereliction. The codereliction is what captures the ability of differentiating non-linear proofs via linearization in Differential Linear Logic. The additive enrichment plays an important role since it allows one to express the famous Leibniz rule. However, the axioms of a codereliction can be expressed without any sums or zeros. Therefore, it is natural to ask if one can consider a possible non-additive enriched version of differential linear categories. In this paper, we show that even if a codereliction can technically be defined in a non-additive setting, it nevertheless induces an additive enrichment via bialgebra convolution. Thus, we obtain a novel characterization of a differential linear category as a symmetric monoidal category with a monoidal bialgebra modality equipped with a codereliction. Moreover, we also show that coderelictions are, in fact, unique. We also introduce monoidal Hopf coalgebra modalities and discuss how antipodes relate to enrichment over Abelian groups. 
\end{abstract}

\maketitle

\section{Introduction}

Differential categories were introduced by Blute, Cockett, and Seely in \cite{blute2006differential} and provide a categorical framework for the algebraic foundations of differentiation. Differential categories are now a well-established area of research with a rich literature and have been successful in formalizing numerous aspects of differential calculus. The original motivation for developing differential categories was to provide the categorical semantics of Differential Linear Logic, as introduced by Ehrhard and Regnier in \cite{ehrhard2006differential}. In particular, the categorical semantics of the multiplicative and exponential fragments of Differential Linear Logic are captured by what is called a differential \emph{linear} category. 

As is well known, the underlying category for the multiplicative fragment of Linear Logic is a symmetric monoidal category. However, for a differential linear category, the underlying category is, in fact, an \emph{additive} symmetric monoidal category, which is a symmetric monoidal category enriched over commutative monoids. In particular, this means that we can take sums of maps and have zero maps. That additive enrichment plays a role in the foundations of differentiation is to be expected, as it allows one to express the famous Leibniz rule (also called the product rule), $(f \cdot g)^\prime = f^\prime \cdot g + f \cdot g^\prime$, as well as expressing that the derivative of a constant map is equal to zero. Moreover, having (countable) sums (along with invertibility of non-zero natural numbers) is also necessary to properly formalize Taylor series expansion in differential categories \cite{ehrhard2017introduction,entics:14670,kerjean2023taylor}.

On the other hand, the exponential fragment of Linear Logic is captured by the notion of a monoidal coalgebra modality (also sometimes referred to as a linear exponential modality). Briefly, a monoidal coalgebra modality is a symmetric monoidal comonad $\oc$ on a symmetric monoidal category such that for every object $A$, $\oc(A)$ is naturally a cocommutative comonoid, whose comonoid structure is compatible with the symmetric monoidal comonad structure. For a monoidal coalgebra modality on an additive symmetric monoidal category, $\oc(A)$ also has a canonical commutative monoid structure, which in fact, makes $\oc(A)$ a bicommutative bimonoid. Then a differential linear category is an additive symmetric monoidal category with a monoidal coalgebra modality that comes equipped with a \emph{codereliction}, which intuitively captures differentiation via linearization. 

Initially, in \cite[Def 4.11]{blute2006differential}, the codereliction was axiomatized via four axioms: the constant rule, the product rule, the linear rule, and the chain rule. These four rules involved the comonad and bimonoid structural maps for the linear rule and the chain rule, as well as the additive enrichment for the product rule and the constant rule. Fiore proposed an alternative axiomatization of a codereliction (which he called a creation map) in \cite[Def 4.3]{fiore2007differential} by replacing the chain rule with an alternative version plus the monoidal rule (also called the strength rule). This latter additional rule now also involved the monoidal structural maps. For a while, it was unclear whether these two approaches to the axioms of a codereliction were, in fact, the same. Fortunately, as it turns out, they are indeed the same, as was shown by Blute et al. when they revisited the various axiomatization of differential categories in \cite{Blute2019}. Recently, Fiore's version of the axioms of a codereliction is often the definition that is used, especially amongst the Linear Logic community (usually because Fiore's version of the chain rule is easier to state and check). So we may say that a codereliction is equivalently axiomatized via five axioms: the constant rule, the product rule, the linear rule, (Fiore's version of) the chain rule, and the monoidal rule. 

With all that said, it turns out that two of the five axioms of a codereliction are redundant. Indeed, the constant rule and the product rule can in fact be proved from naturality \cite[Lemma 6]{Blute2019} and the linear rule \cite[Prop 4]{Blute2019} respectively. As such, a codereliction can therefore be axiomatized simply in terms of three rules: the linear rule, the chain rule, and the monoidal rule \cite[Cor 5]{Blute2019}. However, none of these three rules involve sums or zeroes; the additive enrichment was only needed for the product rule and the constant rule, which we can remove from the definition of a codereliction. As such, if additive enrichment is not necessary for the axioms of a codereliction, this leads to the natural question of whether we can consider differential categories without additive enrichment. 

Before we answer this question, we need to address the bimonoid structural maps in the non-additive setting. Indeed, while the linear rule and the monoidal rule only need the structural maps of a monoidal coalgebra modality, the chain rule does require the (bi)monoid structural maps. However, while the bimonoid structure for a monoidal coalgebra modality comes for free in an additive setting, it may not necessarily exist in the non-additive setting. So if one wishes to consider coderelictions in a non-additive setting, one still has to ask for bimonoid structure. Thus, in this paper, we introduce the notion of a \emph{monoidal bialgebra modality} (Sec \ref{sec:mon-bialg-mod}), which can be defined for any symmetric monoidal category and is briefly a monoidal coalgebra modality with compatible (bi)monoid structure. Thus one could consider the definition of a non-additive version of a differential linear category to be a symmetric monoidal category with a monoidal bialgebra modality that comes equipped with a codereliction, which is axiomatized by the linear rule, chain rule, and monoidal rule. However, it turns out that this gets us back to additive enrichment anyway. 

The main result of this paper is showing that a monoidal bialgebra modality with a codereliction induces an additive enrichment on the base symmetric monoidal category using \emph{bialgebra convolution} \cite[Pg. 72]{sweedler1969hopf}. This further cements and justifies the central role of additive enrichment in the categorical foundations of differentiation. Moreover, this also provides a novel equivalent axiomatization of a differential linear category as a symmetric monoidal category with a monoidal bialgebra modality that comes equipped with a codereliction (Thm \ref{thm:diff-lin}). 

In fact, the chain rule is not required for the main result. So we introduce the intermediate notion of a \emph{pre-codereliction} (Sec \ref{sec:pre-coder}), which need only satisfy the linear rule and the monoidal rule, and thus can be defined for just a monoidal coalgebra modality. Thus, our main result states that a monoidal bialgebra modality with a pre-codereliction induces an additive enrichment on the base symmetric monoidal category (Thm \ref{thm:additive}). Moreover, we also show that pre-coderelictions are unique (Prop \ref{prop:precoderunique}). As a consequence, we obtain the important observation that coderelictions are, in fact, unique\footnote{It was previously known that coderelictions were unique for \emph{free exponential modalities} \cite[Thm 21]{lemay2021coderelictions}. Here, we extend this fact to all monoidal coalgebra modalities.} as well (Thm \ref{thm:coder-unique}). Thus in a categorical model of Differential Linear Logic, there really is only one way to differentiate non-linear (also called smooth) maps. 

Additionally, if we are able to obtain sums and zeroes, i.e. enrichment over commutative monoids, it is natural to ask if we can also obtain negatives, i.e. enrichment over Abelian groups. Since our additive enrichment was built using bialgebra convolution, it follows that asking for Abelian group enrichment corresponds to having inverses for bialgebra convolution. However, it is well known that having inverses for bialgebra convolution corresponds precisely to a bialgebra having an \emph{antipode}, or in other words, that the bialgebra is in fact a \emph{Hopf algebra} \cite[Pg. 71]{sweedler1969hopf}. Applications of Hopf algebras for Linear Logic were considered by Blute and Scott in \cite{blute1996hopf,blute2004category}, who built models of Linear Logic using modules of Hopf algebras, and by the author in \cite{lemay2019lifting}, who built models of (Differential) Linear Logic using Eilenberg-Moore categories of Hopf monads. It turns out that for an additive bialgebra modality, it has a natural antipode if and only if the base category has negatives \cite[Prop 7.6]{lemay2019lifting}. In light of this, in this paper we introduce the notion of a \emph{(monoidal) Hopf coalgebra modality}, which adds natural Hopf monoid structure to a (monoidal) coalgebra modality. Then we show that a monoidal Hopf coalgebra modality with a pre-codereliction induces an Abelian group enrichment on the base category (Prop \ref{prop:neg}).  

\textbf{Outline:} We have written this paper so that it is as self-contained as possible. As such, we have taken the pain of writing out all the definitions in detail. To help write down and better understand definitions and calculations, we make use of string diagrams, which is a graphical calculus for symmetric monoidal categories. That said, we provide all commutative diagrams for definitions and certain statements in Appendix \ref{sec:appendix}. Section \ref{sec:mon-coalg} is mostly a background section where we set up our string diagram conventions and review monoidal coalgebra modalities. In Section \ref{sec:pre-coder} we introduce pre-coderelictions, and in Section \ref{sec:mon-bialg-mod} we introduce monoidal bialgebra modalities. In Section \ref{sec:add-en} we prove the main result of this paper that a monoidal bialgebra modality with a pre-codereliction induces additive enrichment. In Section \ref{sec:add-bialg-mon}, we revisit the notion of \emph{additive bialgebra modalities} and show that, in the additive setting, they are equivalent to monoidal bialgebra modalities. In Section \ref{sec:hopf-coalg-mod} we take a slight detour where we introduce \emph{Hopf coalgebra modalities} and show how they are connected to being enriched over Abelian groups. Lastly, we conclude with Section \ref{sec:diff-lin} where we revisit differential linear categories. 

\section{Monoidal Coalgebra Modalities}\label{sec:mon-coalg}

To set up terminology, notation, and conventions, in this background section, we quickly review monoidal coalgebra modalities -- including some basic definitions and the graphical calculus we will use in this paper. We will mostly follow the same conventions and style used in \cite{Blute2019}. We assume that the reader is familiar with the basics of monoidal category theory and is also somewhat familiar with string diagrams for symmetric monoidal categories. For an in-depth introduction to the graphical calculus of monoidal categories, see \cite{selinger2010survey}, and for an in-depth introduction to the categorical semantics of Linear Logic, see \cite{mellies_categorical_2008,schalk2004categorical,blute2004category,bierman1995categorical}.  

Let us begin simply with an arbitrary category. For an arbitrary category $\mathbb{X}$, we will denote objects using capital letters $A$, $B$, $C$ etc., homsets will be denoted as $\mathbb{X}(A,B)$ and maps will be denoted by minuscule letters $f,g,h, etc. \in \mathbb{X}(A,B)$. Arbitrary maps will be denoted using an arrow $f: A \to B$, identity maps as $1_A: A \to A$, and for composition we will use \emph{diagrammatic} notation, that is, the composition of $f: A \to B$ followed by $g: B \to C$ is denoted as $f;g: A \to C$. String diagrams are to be read from top to bottom. So a map $f: A \to B$ will be drawn out as follows: 
\[ \begin{tikzpicture}
	\begin{pgfonlayer}{nodelayer}
		\node [style=object] (0) at (0, 2) {$A$};
		\node [style=object] (1) at (0, 0) {$B$};
		\node [style=component] (2) at (0, 1) {$f$};
	\end{pgfonlayer}
	\begin{pgfonlayer}{edgelayer}
		\draw [style=wire] (0) to (2);
		\draw [style=wire] (2) to (1);
	\end{pgfonlayer}
\end{tikzpicture}
 \]
Identity maps $1_A: A \to A$ will be drawn as a simple wire as on the left below, while composition will be drawn out sequentially as on the right below: 
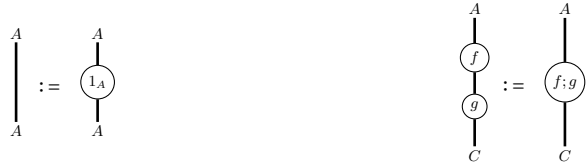
\begin{align*}
\begin{array}[c]{c}
\begin{tikzpicture}
	\begin{pgfonlayer}{nodelayer}
		\node [style=object] (0) at (0, 2) {$A$};
		\node [style=object] (1) at (0, 0) {$A$};
	\end{pgfonlayer}
	\begin{pgfonlayer}{edgelayer}
		\draw [style=wire] (0) to (1);
	\end{pgfonlayer}
\end{tikzpicture}  
\end{array} := \begin{array}[c]{c}
 \begin{tikzpicture}
	\begin{pgfonlayer}{nodelayer}
		\node [style=object] (0) at (0, 2) {$A$};
		\node [style=object] (1) at (0, 0) {$A$};
		\node [style=component] (2) at (0, 1) {$1_A$};
	\end{pgfonlayer}
	\begin{pgfonlayer}{edgelayer}
		\draw [style=wire] (0) to (2);
		\draw [style=wire] (2) to (1);
	\end{pgfonlayer}
\end{tikzpicture}  
\end{array} && \begin{array}[c]{c} \begin{tikzpicture}
	\begin{pgfonlayer}{nodelayer}
		\node [style=object] (5) at (2, 2) {$A$};
		\node [style=object] (6) at (2, -1) {$C$};
		\node [style=component] (7) at (2, 1) {$f$};
		\node [style=component] (12) at (2, 0) {$g$};
	\end{pgfonlayer}
	\begin{pgfonlayer}{edgelayer}
		\draw [style=wire] (5) to (7);
		\draw [style=wire] (7) to (12);
		\draw [style=wire] (12) to (6);
	\end{pgfonlayer}
\end{tikzpicture}
 \end{array} := \begin{array}[c]{c} \begin{tikzpicture}
	\begin{pgfonlayer}{nodelayer}
		\node [style=object] (2) at (1, 2.5) {$A$};
		\node [style=object] (3) at (1, -0.5) {$C$};
		\node [style=component] (4) at (1, 1) {$f;g$};
	\end{pgfonlayer}
	\begin{pgfonlayer}{edgelayer}
		\draw [style=wire] (2) to (4);
		\draw [style=wire] (4) to (3);
	\end{pgfonlayer}
\end{tikzpicture}
 \end{array} 
\end{align*}
When there is no confusion, we will often omit labelling the inputs and outputs. 

The underlying structure of a monodial coalgebra modality is that of a \emph{comonad}, which in particular involves an endofunctor. In our graphical calculus, we will represent applications of our endofunctor using \emph{functor boxes}. So given an endofunctor $\oc: \mathbb{X} \to \mathbb{X}$, for a map $f: A \to B$, we will draw $\oc(f): \oc(A) \to \oc(B)$ as follows: 
 \[ \begin{array}[c]{c} 
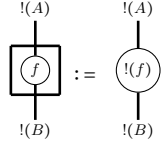
\begin{tikzpicture}
	\begin{pgfonlayer}{nodelayer}
		\node [style=object] (0) at (0, 2.25) {$\oc(A)$};
		\node [style=object] (1) at (0, -0.25) {$\oc(B)$};
		\node [style=component] (2) at (0, 1) {$f$};
		\node [style=none] (3) at (-0.5, 1.5) {};
		\node [style=none] (4) at (0.5, 1.5) {};
		\node [style=none] (5) at (-0.5, 0.5) {};
		\node [style=none] (6) at (0.5, 0.5) {};
	\end{pgfonlayer}
	\begin{pgfonlayer}{edgelayer}
		\draw [style=wire] (0) to (2);
		\draw [style=wire] (2) to (1);
		\draw [style=wire] (3.center) to (4.center);
		\draw [style=wire] (4.center) to (6.center);
		\draw [style=wire] (6.center) to (5.center);
		\draw [style=wire] (5.center) to (3.center);
	\end{pgfonlayer}
\end{tikzpicture} \end{array} := \begin{array}[c]{c}
 \begin{tikzpicture}
	\begin{pgfonlayer}{nodelayer}
		\node [style=object] (0) at (0, 2.25) {$\oc(A)$};
		\node [style=object] (1) at (0, -0.25) {$\oc(B)$};
		\node [style=component] (2) at (0, 1) {$\oc(f)$};
	\end{pgfonlayer}
	\begin{pgfonlayer}{edgelayer}
		\draw [style=wire] (0) to (2);
		\draw [style=wire] (2) to (1);
	\end{pgfonlayer}
\end{tikzpicture}  
\end{array} \]
Then the functoriality of $\oc$ is drawn out as follows: 
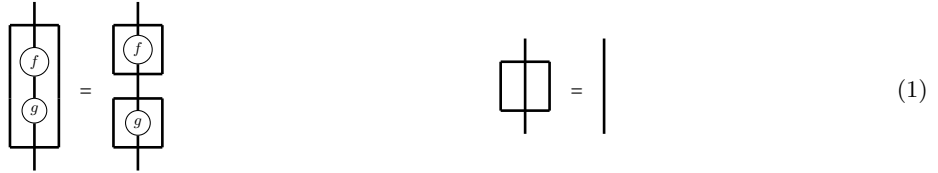
\begin{align}\label{string:functor}
 \begin{array}[c]{c} \begin{tikzpicture}
	\begin{pgfonlayer}{nodelayer}
		\node [style=port] (0) at (0, 2) {};
		\node [style=port] (1) at (0, -1.5) {};
		\node [style=component] (2) at (0, 0.75) {$f$};
		\node [style=none] (3) at (-0.5, 1.5) {};
		\node [style=none] (4) at (0.5, 1.5) {};
		\node [style=none] (5) at (-0.5, 0.5) {};
		\node [style=none] (6) at (0.5, 0.5) {};
		\node [style=component] (7) at (0, -0.25) {$g$};
		\node [style=none] (8) at (-0.5, 0) {};
		\node [style=none] (9) at (0.5, 0) {};
		\node [style=none] (10) at (-0.5, -1) {};
		\node [style=none] (11) at (0.5, -1) {};
	\end{pgfonlayer}
	\begin{pgfonlayer}{edgelayer}
		\draw [style=wire] (0) to (2);
		\draw [style=wire] (3.center) to (4.center);
		\draw [style=wire] (4.center) to (6.center);
		\draw [style=wire] (5.center) to (3.center);
		\draw [style=wire] (9.center) to (11.center);
		\draw [style=wire] (11.center) to (10.center);
		\draw [style=wire] (10.center) to (8.center);
		\draw [style=wire] (2) to (7);
		\draw [style=wire] (7) to (1);
		\draw [style=wire] (5.center) to (8.center);
		\draw [style=wire] (6.center) to (9.center);
	\end{pgfonlayer}
\end{tikzpicture}
\end{array} =  \begin{array}[c]{c} \begin{tikzpicture}
	\begin{pgfonlayer}{nodelayer}
		\node [style=port] (0) at (0, 2) {};
		\node [style=port] (1) at (0, -1.5) {};
		\node [style=component] (2) at (0, 1) {$f$};
		\node [style=none] (3) at (-0.5, 1.5) {};
		\node [style=none] (4) at (0.5, 1.5) {};
		\node [style=none] (5) at (-0.5, 0.5) {};
		\node [style=none] (6) at (0.5, 0.5) {};
		\node [style=component] (7) at (0, -0.5) {$g$};
		\node [style=none] (8) at (-0.5, 0) {};
		\node [style=none] (9) at (0.5, 0) {};
		\node [style=none] (10) at (-0.5, -1) {};
		\node [style=none] (11) at (0.5, -1) {};
	\end{pgfonlayer}
	\begin{pgfonlayer}{edgelayer}
		\draw [style=wire] (0) to (2);
		\draw [style=wire] (3.center) to (4.center);
		\draw [style=wire] (4.center) to (6.center);
		\draw [style=wire] (6.center) to (5.center);
		\draw [style=wire] (5.center) to (3.center);
		\draw [style=wire] (8.center) to (9.center);
		\draw [style=wire] (9.center) to (11.center);
		\draw [style=wire] (11.center) to (10.center);
		\draw [style=wire] (10.center) to (8.center);
		\draw [style=wire] (2) to (7);
		\draw [style=wire] (7) to (1);
	\end{pgfonlayer}
\end{tikzpicture}
\end{array}   &&   \begin{array}[c]{c} \begin{tikzpicture}
	\begin{pgfonlayer}{nodelayer}
		\node [style=port] (0) at (0, 2) {};
		\node [style=port] (1) at (0, 0) {};
		\node [style=none] (3) at (-0.5, 1.5) {};
		\node [style=none] (4) at (0.5, 1.5) {};
		\node [style=none] (10) at (-0.5, 0.5) {};
		\node [style=none] (11) at (0.5, 0.5) {};
	\end{pgfonlayer}
	\begin{pgfonlayer}{edgelayer}
		\draw [style=wire] (3.center) to (4.center);
		\draw [style=wire] (11.center) to (10.center);
		\draw [style=wire] (0) to (1);
		\draw [style=wire] (3.center) to (10.center);
		\draw [style=wire] (4.center) to (11.center);
	\end{pgfonlayer}
\end{tikzpicture}
\end{array} =   \begin{array}[c]{c} \begin{tikzpicture}
	\begin{pgfonlayer}{nodelayer}
		\node [style=port] (0) at (0, 2) {};
		\node [style=port] (1) at (0, 0) {};
	\end{pgfonlayer}
	\begin{pgfonlayer}{edgelayer}
		\draw [style=wire] (0) to (1);
	\end{pgfonlayer}
\end{tikzpicture}
\end{array}
\end{align}
Then a \textbf{comonad} on a category $\mathbb{X}$ is a triple $(\oc, \delta, \varepsilon)$ consisting of an endofunctor $\oc: \mathbb{X} \to \mathbb{X}$ and two natural transformations $\delta_A: \oc(A) \to \oc\oc(A)$ and $\varepsilon_A: \oc(A) \to A$, such that the diagrams in Appendix \ref{sec:diagramsforcomonad} commute. In Linear Logic terminology, the natural transformation $\delta$ is called the \textbf{digging} and the natural transformation $\varepsilon$ is called the \textbf{dereliction}. Respectively $\delta$ and $\varepsilon$ are drawn as follows: 
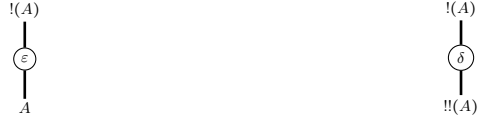
\begin{align*} \begin{array}[c]{c} \begin{tikzpicture}
	\begin{pgfonlayer}{nodelayer}
		\node [style=object] (2) at (1, 2) {$\oc(A)$};
		\node [style=object] (3) at (1, 0) {$A$};
		\node [style=component] (4) at (1, 1) {$\varepsilon$};
	\end{pgfonlayer}
	\begin{pgfonlayer}{edgelayer}
		\draw [style=wire] (2) to (4);
		\draw [style=wire] (4) to (3);
	\end{pgfonlayer}
\end{tikzpicture}
 \end{array} && \begin{array}[c]{c} \begin{tikzpicture}
	\begin{pgfonlayer}{nodelayer}
		\node [style=object] (2) at (1, 2) {$\oc(A)$};
		\node [style=object] (3) at (1, 0) {$\oc\oc(A)$};
		\node [style=component] (4) at (1, 1) {$\delta$};
	\end{pgfonlayer}
	\begin{pgfonlayer}{edgelayer}
		\draw [style=wire] (2) to (4);
		\draw [style=wire] (4) to (3);
	\end{pgfonlayer}
\end{tikzpicture}
 \end{array} 
 \end{align*}
 with their naturality expressed respectively as follows: 
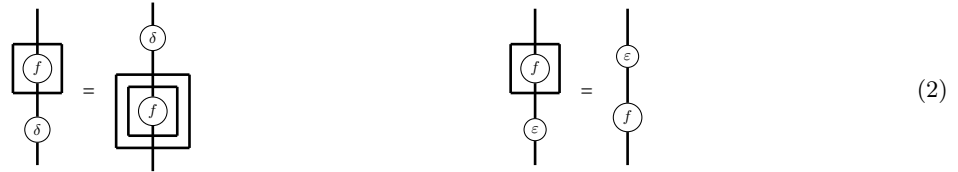
\begin{align}\label{string:comonad-nat}
 \begin{array}[c]{c} \begin{tikzpicture}
	\begin{pgfonlayer}{nodelayer}
		\node [style=port] (16) at (-2.25, 0.5) {};
		\node [style=component] (17) at (-2.25, 1.25) {$\delta$};
		\node [style=port] (18) at (-2.25, 3.75) {};
		\node [style=component] (19) at (-2.25, 2.5) {$f$};
		\node [style=none] (20) at (-2.75, 3) {};
		\node [style=none] (21) at (-1.75, 3) {};
		\node [style=none] (22) at (-2.75, 2) {};
		\node [style=none] (23) at (-1.75, 2) {};
	\end{pgfonlayer}
	\begin{pgfonlayer}{edgelayer}
		\draw [style=wire] (17) to (16);
		\draw [style=wire] (18) to (19);
		\draw [style=wire] (20.center) to (21.center);
		\draw [style=wire] (21.center) to (23.center);
		\draw [style=wire] (23.center) to (22.center);
		\draw [style=wire] (22.center) to (20.center);
		\draw [style=wire] (19) to (17);
	\end{pgfonlayer}
\end{tikzpicture}
\end{array} =  \begin{array}[c]{c}\begin{tikzpicture}
	\begin{pgfonlayer}{nodelayer}
		\node [style=port] (3) at (0, 4) {};
		\node [style=component] (4) at (0, 3.25) {$\delta$};
		\node [style=port] (5) at (0, 0.5) {};
		\node [style=component] (7) at (0, 1.75) {$f$};
		\node [style=none] (8) at (-0.5, 2.25) {};
		\node [style=none] (9) at (0.5, 2.25) {};
		\node [style=none] (10) at (-0.5, 1.25) {};
		\node [style=none] (11) at (0.5, 1.25) {};
		\node [style=none] (12) at (-0.75, 2.5) {};
		\node [style=none] (13) at (0.75, 2.5) {};
		\node [style=none] (14) at (-0.75, 1) {};
		\node [style=none] (15) at (0.75, 1) {};
	\end{pgfonlayer}
	\begin{pgfonlayer}{edgelayer}
		\draw [style=wire] (4) to (3);
		\draw [style=wire] (5) to (7);
		\draw [style=wire] (8.center) to (9.center);
		\draw [style=wire] (9.center) to (11.center);
		\draw [style=wire] (11.center) to (10.center);
		\draw [style=wire] (10.center) to (8.center);
		\draw [style=wire] (12.center) to (13.center);
		\draw [style=wire] (13.center) to (15.center);
		\draw [style=wire] (15.center) to (14.center);
		\draw [style=wire] (14.center) to (12.center);
		\draw [style=wire] (4) to (7);
	\end{pgfonlayer}
\end{tikzpicture}
\end{array}   &&   \begin{array}[c]{c} \begin{tikzpicture}
	\begin{pgfonlayer}{nodelayer}
		\node [style=port] (16) at (-2.25, 0.5) {};
		\node [style=component] (17) at (-2.25, 1.25) {$\varepsilon$};
		\node [style=port] (18) at (-2.25, 3.75) {};
		\node [style=component] (19) at (-2.25, 2.5) {$f$};
		\node [style=none] (20) at (-2.75, 3) {};
		\node [style=none] (21) at (-1.75, 3) {};
		\node [style=none] (22) at (-2.75, 2) {};
		\node [style=none] (23) at (-1.75, 2) {};
	\end{pgfonlayer}
	\begin{pgfonlayer}{edgelayer}
		\draw [style=wire] (17) to (16);
		\draw [style=wire] (18) to (19);
		\draw [style=wire] (20.center) to (21.center);
		\draw [style=wire] (21.center) to (23.center);
		\draw [style=wire] (23.center) to (22.center);
		\draw [style=wire] (22.center) to (20.center);
		\draw [style=wire] (19) to (17);
	\end{pgfonlayer}
\end{tikzpicture}
\end{array} =   \begin{array}[c]{c} \begin{tikzpicture}
	\begin{pgfonlayer}{nodelayer}
		\node [style=port] (3) at (0, 3.75) {};
		\node [style=component] (4) at (0, 2.75) {$\varepsilon$};
		\node [style=port] (5) at (0, 0.5) {};
		\node [style=component] (7) at (0, 1.5) {$f$};
	\end{pgfonlayer}
	\begin{pgfonlayer}{edgelayer}
		\draw [style=wire] (4) to (3);
		\draw [style=wire] (5) to (7);
		\draw [style=wire] (4) to (7);
	\end{pgfonlayer}
\end{tikzpicture}
\end{array}
\end{align}
The comonad identities (\ref{diag:comonad}) are drawn as follows: 
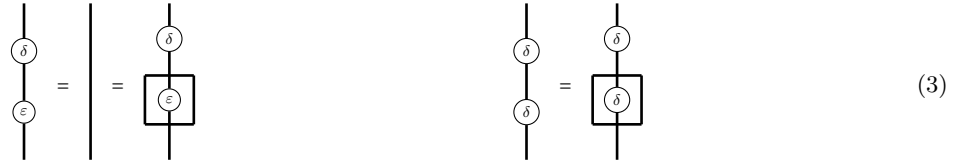
\begin{align}\label{string:comonad}
 \begin{array}[c]{c}\begin{tikzpicture}
	\begin{pgfonlayer}{nodelayer}
		\node [style=port] (16) at (-2.25, 0.5) {};
		\node [style=component] (17) at (-2.25, 1.5) {$\varepsilon$};
		\node [style=port] (18) at (-2.25, 3.75) {};
		\node [style=component] (19) at (-2.25, 2.75) {$\delta$};
	\end{pgfonlayer}
	\begin{pgfonlayer}{edgelayer}
		\draw [style=wire] (17) to (16);
		\draw [style=wire] (18) to (19);
		\draw [style=wire] (19) to (17);
	\end{pgfonlayer}
\end{tikzpicture}
\end{array} =  \begin{array}[c]{c}\begin{tikzpicture}
	\begin{pgfonlayer}{nodelayer}
		\node [style=port] (16) at (-2.25, 0.5) {};
		\node [style=port] (18) at (-2.25, 3.75) {};
	\end{pgfonlayer}
	\begin{pgfonlayer}{edgelayer}
		\draw [style=wire] (18) to (16);
	\end{pgfonlayer}
\end{tikzpicture}
\end{array} =  \begin{array}[c]{c}\begin{tikzpicture}
	\begin{pgfonlayer}{nodelayer}
		\node [style=port] (3) at (0, 3.75) {};
		\node [style=component] (4) at (0, 3) {$\delta$};
		\node [style=port] (5) at (0, 0.5) {};
		\node [style=component] (7) at (0, 1.75) {$\varepsilon$};
		\node [style=none] (8) at (-0.5, 2.25) {};
		\node [style=none] (9) at (0.5, 2.25) {};
		\node [style=none] (10) at (-0.5, 1.25) {};
		\node [style=none] (11) at (0.5, 1.25) {};
	\end{pgfonlayer}
	\begin{pgfonlayer}{edgelayer}
		\draw [style=wire] (4) to (3);
		\draw [style=wire] (5) to (7);
		\draw [style=wire] (8.center) to (9.center);
		\draw [style=wire] (9.center) to (11.center);
		\draw [style=wire] (11.center) to (10.center);
		\draw [style=wire] (10.center) to (8.center);
		\draw [style=wire] (4) to (7);
	\end{pgfonlayer}
\end{tikzpicture}
\end{array}   &&   \begin{array}[c]{c}\begin{tikzpicture}
	\begin{pgfonlayer}{nodelayer}
		\node [style=port] (16) at (-2.25, 0.5) {};
		\node [style=component] (17) at (-2.25, 1.5) {$\delta$};
		\node [style=port] (18) at (-2.25, 3.75) {};
		\node [style=component] (19) at (-2.25, 2.75) {$\delta$};
	\end{pgfonlayer}
	\begin{pgfonlayer}{edgelayer}
		\draw [style=wire] (17) to (16);
		\draw [style=wire] (18) to (19);
		\draw [style=wire] (19) to (17);
	\end{pgfonlayer}
\end{tikzpicture}
\end{array} =   \begin{array}[c]{c} \begin{tikzpicture}
	\begin{pgfonlayer}{nodelayer}
		\node [style=port] (3) at (0, 3.75) {};
		\node [style=component] (4) at (0, 3) {$\delta$};
		\node [style=port] (5) at (0, 0.5) {};
		\node [style=component] (7) at (0, 1.75) {$\delta$};
		\node [style=none] (8) at (-0.5, 2.25) {};
		\node [style=none] (9) at (0.5, 2.25) {};
		\node [style=none] (10) at (-0.5, 1.25) {};
		\node [style=none] (11) at (0.5, 1.25) {};
	\end{pgfonlayer}
	\begin{pgfonlayer}{edgelayer}
		\draw [style=wire] (4) to (3);
		\draw [style=wire] (5) to (7);
		\draw [style=wire] (8.center) to (9.center);
		\draw [style=wire] (9.center) to (11.center);
		\draw [style=wire] (11.center) to (10.center);
		\draw [style=wire] (10.center) to (8.center);
		\draw [style=wire] (4) to (7);
	\end{pgfonlayer}
\end{tikzpicture}
\end{array}
\end{align}

We now, and for the rest of the paper, move to working in a symmetric monoidal category. Following the convention used in most of the literature on differential categories, we will work in a symmetric \emph{strict} monoidal category, so the associativity and unit isomorphisms for the monoidal product are equalities. For an arbitrary symmetric (strict) monoidal category $\mathbb{X}$, we denote its monoidal product as $\otimes$, the monoidal unit as $I$, and the natural symmetry isomorphism as $\sigma_{A, B}: A \otimes B \xrightarrow{\cong} B \otimes A$. Strictness allows us to write $A_1 \otimes A_2 \otimes \hdots \otimes A_n$ and $A \otimes I = A = I \otimes A$. An arbitrary map $f: A_1 \otimes \hdots \otimes A_n \to B_1 \otimes \hdots \otimes B_m$ will be drawn with $n$ input wires from the top, one wire for each $A_i$, and outputs $m$ output wires from the bottom, with one wire for each $B_j$, where the monoidal product of objects is to be read from left to right:
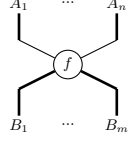
\[ \begin{tikzpicture}
	\begin{pgfonlayer}{nodelayer}
		\node [style=object] (1) at (0, 1.75) {$B_1$};
		\node [style=component] (2) at (1, 3) {$f$};
		\node [style=none] (3) at (0, 2.5) {};
		\node [style=none] (4) at (2, 2.5) {};
		\node [style=object] (5) at (2, 1.75) {$B_m$};
		\node [style=object] (6) at (1, 1.75) {$\hdots$};
		\node [style=object] (7) at (0, 4.25) {$A_1$};
		\node [style=none] (8) at (0, 3.5) {};
		\node [style=none] (9) at (2, 3.5) {};
		\node [style=object] (10) at (2, 4.25) {$A_n$};
		\node [style=object] (11) at (1, 4.25) {$\hdots$};
	\end{pgfonlayer}
	\begin{pgfonlayer}{edgelayer}
		\draw [style=wire] (2) to (3.center);
		\draw [style=wire] (2) to (4.center);
		\draw [style=wire] (4.center) to (5);
		\draw [style=wire] (3.center) to (1);
		\draw [style=wire] (9.center) to (10);
		\draw [style=wire] (8.center) to (7);
		\draw [style=wire] (8.center) to (2);
		\draw [style=wire] (2) to (9.center);
	\end{pgfonlayer}
\end{tikzpicture}
 \]
In the special case that one of the inputs or outputs is the monoidal unit $I$, we do not draw a wire representing it (which lines up with the idea that we are working in the strict case). So for example, this is how we draw maps of type $g: I \to B$ and $h: A \to I$: 
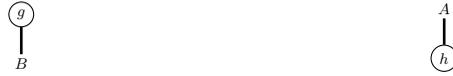
\begin{align*}
\begin{tikzpicture}
	\begin{pgfonlayer}{nodelayer}
		\node [style=object] (13) at (4, 0) {$B$};
		\node [style=component] (14) at (4, 1) {$g$};
	\end{pgfonlayer}
	\begin{pgfonlayer}{edgelayer}
		\draw [style=wire] (14) to (13);
	\end{pgfonlayer}
\end{tikzpicture}
&& \begin{tikzpicture}
	\begin{pgfonlayer}{nodelayer}
		\node [style=object] (13) at (4, 1) {$A$};
		\node [style=component] (14) at (4, 0) {$h$};
	\end{pgfonlayer}
	\begin{pgfonlayer}{edgelayer}
		\draw [style=wire] (14) to (13);
	\end{pgfonlayer}
\end{tikzpicture}
\end{align*}
The monoidal product of maps is drawn simply as drawing them in parallel, so for maps $f: A \to B$ and $g: C \to D$, their monoidal product $f \otimes g: A \otimes C \to B \otimes D$ is drawn as follows: 
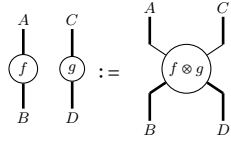
\[  \begin{array}[c]{c}\begin{tikzpicture}
	\begin{pgfonlayer}{nodelayer}
		\node [style=object] (12) at (4, 2) {$A$};
		\node [style=object] (13) at (4, 0) {$B$};
		\node [style=component] (14) at (4, 1) {$f$};
		\node [style=object] (15) at (5, 2) {$C$};
		\node [style=object] (16) at (5, 0) {$D$};
		\node [style=component] (17) at (5, 1) {$g$};
	\end{pgfonlayer}
	\begin{pgfonlayer}{edgelayer}
		\draw [style=wire] (12) to (14);
		\draw [style=wire] (14) to (13);
		\draw [style=wire] (15) to (17);
		\draw [style=wire] (17) to (16);
	\end{pgfonlayer}
\end{tikzpicture}
\end{array} :=   \begin{array}[c]{c} \begin{tikzpicture}
	\begin{pgfonlayer}{nodelayer}
		\node [style=object] (1) at (0.25, 1.75) {$B$};
		\node [style=component] (2) at (1, 3) {$f \otimes g$};
		\node [style=none] (3) at (0.25, 2.5) {};
		\node [style=none] (4) at (1.75, 2.5) {};
		\node [style=object] (5) at (1.75, 1.75) {$D$};
		\node [style=object] (7) at (0.25, 4.25) {$A$};
		\node [style=none] (8) at (0.25, 3.5) {};
		\node [style=none] (9) at (1.75, 3.5) {};
		\node [style=object] (10) at (1.75, 4.25) {$C$};
	\end{pgfonlayer}
	\begin{pgfonlayer}{edgelayer}
		\draw [style=wire] (2) to (3.center);
		\draw [style=wire] (2) to (4.center);
		\draw [style=wire] (4.center) to (5);
		\draw [style=wire] (3.center) to (1);
		\draw [style=wire] (9.center) to (10);
		\draw [style=wire] (8.center) to (7);
		\draw [style=wire] (8.center) to (2);
		\draw [style=wire] (2) to (9.center);
	\end{pgfonlayer}
\end{tikzpicture}
\end{array} \]
The symmetry isomorphism $\sigma_{A,B}: A \otimes B \to B \otimes A$ is drawn simply as twisting the wires, where recall that in the symmetric setting we do not have to worry about which wire passes on top of the other wire -- so we may think of the twisting of wires as the wires passing through each other uninterrupted. So $\sigma_{A,B}$ is drawn as follows:
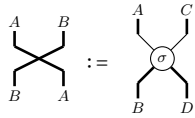
\begin{align*}
\begin{array}[c]{c}\begin{tikzpicture}
	\begin{pgfonlayer}{nodelayer}
		\node [style=object] (19) at (12, 1.25) {$B$};
		\node [style=object] (21) at (12, 2.75) {$A$};
		\node [style=object] (22) at (13, 2.75) {$B$};
		\node [style=object] (23) at (13, 1.25) {$A$};
		\node [style=none] (24) at (12, 2.25) {};
		\node [style=none] (25) at (13, 2.25) {};
		\node [style=none] (26) at (12, 1.75) {};
		\node [style=none] (27) at (13, 1.75) {};
	\end{pgfonlayer}
	\begin{pgfonlayer}{edgelayer}
		\draw [style=wire] (21) to (24.center);
		\draw [style=wire] (22) to (25.center);
		\draw [style=wire] (26.center) to (19);
		\draw [style=wire] (27.center) to (23);
		\draw [style=wire] (24.center) to (27.center);
		\draw [style=wire] (25.center) to (26.center);
	\end{pgfonlayer}
\end{tikzpicture}
\end{array} :=   \begin{array}[c]{c}\begin{tikzpicture}
	\begin{pgfonlayer}{nodelayer}
		\node [style=object] (1) at (9.75, 1) {$B$};
		\node [style=component] (2) at (10.25, 2) {$\sigma$};
		\node [style=none] (3) at (9.75, 1.5) {};
		\node [style=none] (4) at (10.75, 1.5) {};
		\node [style=object] (5) at (10.75, 1) {$D$};
		\node [style=object] (7) at (9.75, 3) {$A$};
		\node [style=none] (8) at (9.75, 2.5) {};
		\node [style=none] (9) at (10.75, 2.5) {};
		\node [style=object] (10) at (10.75, 3) {$C$};
	\end{pgfonlayer}
	\begin{pgfonlayer}{edgelayer}
		\draw [style=wire] (2) to (3.center);
		\draw [style=wire] (2) to (4.center);
		\draw [style=wire] (4.center) to (5);
		\draw [style=wire] (3.center) to (1);
		\draw [style=wire] (9.center) to (10);
		\draw [style=wire] (8.center) to (7);
		\draw [style=wire] (8.center) to (2);
		\draw [style=wire] (2) to (9.center);
	\end{pgfonlayer}
\end{tikzpicture}
\end{array} 
\end{align*}
and its naturality is drawn as: 
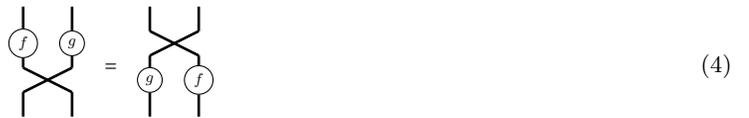
\begin{align}\label{string:nat-sym} \begin{array}[c]{c}\begin{tikzpicture}
	\begin{pgfonlayer}{nodelayer}
		\node [style=none] (12) at (4, 1.25) {};
		\node [style=none] (13) at (4, 0) {};
		\node [style=component] (14) at (4, 0.5) {$f$};
		\node [style=none] (15) at (5, 1.25) {};
		\node [style=none] (16) at (5, 0) {};
		\node [style=component] (17) at (5, 0.5) {$g$};
		\node [style=none] (19) at (4, -1) {};
		\node [style=none] (23) at (5, -1) {};
		\node [style=none] (24) at (4, 0) {};
		\node [style=none] (25) at (5, 0) {};
		\node [style=none] (26) at (4, -0.5) {};
		\node [style=none] (27) at (5, -0.5) {};
	\end{pgfonlayer}
	\begin{pgfonlayer}{edgelayer}
		\draw [style=wire] (12.center) to (14);
		\draw [style=wire] (14) to (13.center);
		\draw [style=wire] (15.center) to (17);
		\draw [style=wire] (17) to (16.center);
		\draw [style=wire] (26.center) to (19.center);
		\draw [style=wire] (27.center) to (23.center);
		\draw [style=wire] (24.center) to (27.center);
		\draw [style=wire] (25.center) to (26.center);
	\end{pgfonlayer}
\end{tikzpicture}
\end{array} = \begin{array}[c]{c}\begin{tikzpicture}
	\begin{pgfonlayer}{nodelayer}
		\node [style=none] (28) at (1, -1) {};
		\node [style=none] (29) at (1, 0.25) {};
		\node [style=component] (30) at (1, -0.25) {$f$};
		\node [style=none] (31) at (0, -1) {};
		\node [style=none] (32) at (0, 0.25) {};
		\node [style=component] (33) at (0, -0.25) {$g$};
		\node [style=none] (34) at (1, 1.25) {};
		\node [style=none] (35) at (0, 1.25) {};
		\node [style=none] (36) at (1, 0.25) {};
		\node [style=none] (37) at (0, 0.25) {};
		\node [style=none] (38) at (1, 0.75) {};
		\node [style=none] (39) at (0, 0.75) {};
	\end{pgfonlayer}
	\begin{pgfonlayer}{edgelayer}
		\draw [style=wire] (28.center) to (30);
		\draw [style=wire] (30) to (29.center);
		\draw [style=wire] (31.center) to (33);
		\draw [style=wire] (33) to (32.center);
		\draw [style=wire] (38.center) to (34.center);
		\draw [style=wire] (39.center) to (35.center);
		\draw [style=wire] (36.center) to (39.center);
		\draw [style=wire] (37.center) to (38.center);
	\end{pgfonlayer}
\end{tikzpicture}
\end{array}
\end{align}

We are now in a position to discuss monoidal coalgebra modalities. Following \cite{blute2006differential,Blute2019}, it will be useful to first review the intermediate notion of simply a \emph{coalgebra modality}, which is a comonad such that each $\oc(A)$ is naturally a cocommutative comonoid. More concretely, a \textbf{coalgebra modality} \cite[Def 1]{Blute2019} on a symmetric monoidal category is a quintuple $(\oc, \delta, \varepsilon, \Delta, \mathsf{e})$ consisting of a comonad $(\oc, \delta, \varepsilon)$ and two natural transformations $\Delta_A: \oc(A) \to \oc(A) \otimes \oc(A)$ and $\mathsf{e}_A: \oc(A) \to I$ such that the diagrams in Appendix \ref{sec:diagramsforcoalgmod} commute. The natural transformation $\Delta$ is called the \textbf{comultiplication} or \textbf{contraction}, and the natural transformation $\mathsf{e}$ is called the \textbf{counit} or \textbf{weakening}, and they are drawn as follows: 
 \begin{align*}\begin{array}[c]{c} 
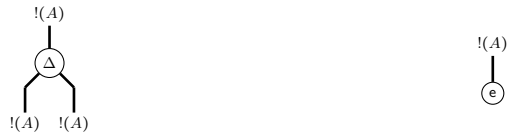
\begin{tikzpicture}
	\begin{pgfonlayer}{nodelayer}
		\node [style=object] (19) at (-7, 4) {$\oc(A)$};
		\node [style=object] (20) at (-7.5, 1.75) {$\oc(A)$};
		\node [style=component] (21) at (-7, 3) {$\Delta$};
		\node [style=none] (22) at (-7.5, 2.5) {};
		\node [style=none] (23) at (-6.5, 2.5) {};
		\node [style=object] (24) at (-6.5, 1.75) {$\oc(A)$};
	\end{pgfonlayer}
	\begin{pgfonlayer}{edgelayer}
		\draw [style=wire] (19) to (21);
		\draw [style=wire] (21) to (22.center);
		\draw [style=wire] (21) to (23.center);
		\draw [style=wire] (23.center) to (24);
		\draw [style=wire] (22.center) to (20);
	\end{pgfonlayer}
\end{tikzpicture}
 \end{array} && \begin{array}[c]{c} \begin{tikzpicture}
	\begin{pgfonlayer}{nodelayer}
		\node [style=object] (25) at (-5, 4) {$\oc(A)$};
		\node [style=component] (27) at (-5, 3) {$\mathsf{e}$};
	\end{pgfonlayer}
	\begin{pgfonlayer}{edgelayer}
		\draw [style=wire] (25) to (27);
	\end{pgfonlayer}
\end{tikzpicture}
 \end{array} 
 \end{align*}
 with their naturality expressed respectively as follows: 
 \begin{align}\label{string:comonoid-nat}\begin{array}[c]{c} \begin{tikzpicture}
	\begin{pgfonlayer}{nodelayer}
		\node [style=component] (7) at (-3.5, 3) {$f$};
		\node [style=none] (8) at (-4, 3.5) {};
		\node [style=none] (9) at (-3, 3.5) {};
		\node [style=none] (10) at (-4, 2.5) {};
		\node [style=none] (11) at (-3, 2.5) {};
		\node [style=port] (20) at (-4, 0.5) {};
		\node [style=component] (21) at (-3.5, 1.75) {$\Delta$};
		\node [style=none] (22) at (-4, 1.25) {};
		\node [style=none] (23) at (-3, 1.25) {};
		\node [style=port] (34) at (-3, 0.5) {};
		\node [style=port] (45) at (-3.5, 4.25) {};
	\end{pgfonlayer}
	\begin{pgfonlayer}{edgelayer}
		\draw [style=wire] (8.center) to (9.center);
		\draw [style=wire] (9.center) to (11.center);
		\draw [style=wire] (11.center) to (10.center);
		\draw [style=wire] (10.center) to (8.center);
		\draw [style=wire] (21) to (22.center);
		\draw [style=wire] (21) to (23.center);
		\draw [style=wire] (22.center) to (20);
		\draw [style=wire] (23.center) to (34);
		\draw [style=wire] (45) to (7);
		\draw [style=wire] (7) to (21);
	\end{pgfonlayer}
\end{tikzpicture}
 \end{array} = \begin{array}[c]{c}\begin{tikzpicture}
	\begin{pgfonlayer}{nodelayer}
		\node [style=component] (46) at (-1.5, 1.5) {$f$};
		\node [style=none] (47) at (-2, 2) {};
		\node [style=none] (48) at (-1, 2) {};
		\node [style=none] (49) at (-2, 1) {};
		\node [style=none] (50) at (-1, 1) {};
		\node [style=port] (51) at (-1.5, 0.5) {};
		\node [style=component] (52) at (-0.75, 3.25) {$\Delta$};
		\node [style=none] (53) at (-1.5, 2.5) {};
		\node [style=none] (54) at (0, 2.5) {};
		\node [style=port] (55) at (0, 0.5) {};
		\node [style=component] (56) at (0, 1.5) {$f$};
		\node [style=none] (57) at (-0.5, 2) {};
		\node [style=none] (58) at (0.5, 2) {};
		\node [style=none] (59) at (-0.5, 1) {};
		\node [style=none] (60) at (0.5, 1) {};
		\node [style=port] (61) at (-0.75, 4.25) {};
	\end{pgfonlayer}
	\begin{pgfonlayer}{edgelayer}
		\draw [style=wire] (47.center) to (48.center);
		\draw [style=wire] (48.center) to (50.center);
		\draw [style=wire] (50.center) to (49.center);
		\draw [style=wire] (49.center) to (47.center);
		\draw [style=wire] (52) to (53.center);
		\draw [style=wire] (52) to (54.center);
		\draw [style=wire] (57.center) to (58.center);
		\draw [style=wire] (58.center) to (60.center);
		\draw [style=wire] (60.center) to (59.center);
		\draw [style=wire] (59.center) to (57.center);
		\draw [style=wire] (46) to (51);
		\draw [style=wire] (56) to (55);
		\draw [style=wire] (53.center) to (46);
		\draw [style=wire] (54.center) to (56);
		\draw [style=wire] (61) to (52);
	\end{pgfonlayer}
\end{tikzpicture}
 \end{array} && \begin{array}[c]{c}\begin{tikzpicture}
	\begin{pgfonlayer}{nodelayer}
		\node [style=component] (7) at (-3.5, 3) {$f$};
		\node [style=none] (8) at (-4, 3.5) {};
		\node [style=none] (9) at (-3, 3.5) {};
		\node [style=none] (10) at (-4, 2.5) {};
		\node [style=none] (11) at (-3, 2.5) {};
		\node [style=component] (21) at (-3.5, 1.75) {$\mathsf{e}$};
		\node [style=port] (45) at (-3.5, 4.25) {};
	\end{pgfonlayer}
	\begin{pgfonlayer}{edgelayer}
		\draw [style=wire] (8.center) to (9.center);
		\draw [style=wire] (9.center) to (11.center);
		\draw [style=wire] (11.center) to (10.center);
		\draw [style=wire] (10.center) to (8.center);
		\draw [style=wire] (45) to (7);
		\draw [style=wire] (7) to (21);
	\end{pgfonlayer}
\end{tikzpicture}
 \end{array} = \begin{array}[c]{c}\begin{tikzpicture}
	\begin{pgfonlayer}{nodelayer}
		\node [style=component] (52) at (-0.75, 3.25) {$\mathsf{e}$};
		\node [style=port] (61) at (-0.75, 4.25) {};
	\end{pgfonlayer}
	\begin{pgfonlayer}{edgelayer}
		\draw [style=wire] (61) to (52);
	\end{pgfonlayer}
\end{tikzpicture}
 \end{array}
 \end{align}
 The requirement that for each object $A$, $(\oc(A), \Delta_A, \mathsf{e}_A)$ is a cocommutative comonoid (\ref{diag:comonoid}) is drawn as follows: 
  \begin{align}\label{string:comonoid}\begin{array}[c]{c} \begin{tikzpicture}
	\begin{pgfonlayer}{nodelayer}
		\node [style=port] (19) at (5, 4) {};
		\node [style=port] (20) at (4.5, 0.5) {};
		\node [style=component] (21) at (5, 3) {$\Delta$};
		\node [style=none] (22) at (4.5, 2.5) {};
		\node [style=none] (23) at (5.5, 2.5) {};
		\node [style=port] (29) at (5, 0.5) {};
		\node [style=component] (30) at (5.5, 1.75) {$\Delta$};
		\node [style=none] (31) at (5, 1.25) {};
		\node [style=none] (32) at (6, 1.25) {};
		\node [style=port] (33) at (6, 0.5) {};
	\end{pgfonlayer}
	\begin{pgfonlayer}{edgelayer}
		\draw [style=wire] (19) to (21);
		\draw [style=wire] (21) to (22.center);
		\draw [style=wire] (21) to (23.center);
		\draw [style=wire] (22.center) to (20);
		\draw [style=wire] (30) to (31.center);
		\draw [style=wire] (30) to (32.center);
		\draw [style=wire] (32.center) to (33);
		\draw [style=wire] (31.center) to (29);
		\draw [style=wire] (23.center) to (30);
	\end{pgfonlayer}
\end{tikzpicture}
 \end{array} = \begin{array}[c]{c} \begin{tikzpicture}
	\begin{pgfonlayer}{nodelayer}
		\node [style=port] (19) at (5.25, 3.975) {};
		\node [style=port] (20) at (5.75, 0.475) {};
		\node [style=component] (21) at (5.25, 2.975) {$\Delta$};
		\node [style=none] (22) at (5.75, 2.475) {};
		\node [style=none] (23) at (4.75, 2.475) {};
		\node [style=port] (29) at (5.25, 0.475) {};
		\node [style=component] (30) at (4.75, 1.725) {$\Delta$};
		\node [style=none] (31) at (5.25, 1.225) {};
		\node [style=none] (32) at (4.25, 1.225) {};
		\node [style=port] (33) at (4.25, 0.475) {};
	\end{pgfonlayer}
	\begin{pgfonlayer}{edgelayer}
		\draw [style=wire] (19) to (21);
		\draw [style=wire] (21) to (22.center);
		\draw [style=wire] (21) to (23.center);
		\draw [style=wire] (22.center) to (20);
		\draw [style=wire] (30) to (31.center);
		\draw [style=wire] (30) to (32.center);
		\draw [style=wire] (32.center) to (33);
		\draw [style=wire] (31.center) to (29);
		\draw [style=wire] (23.center) to (30);
	\end{pgfonlayer}
\end{tikzpicture}
 \end{array} && \begin{array}[c]{c} \begin{tikzpicture}
	\begin{pgfonlayer}{nodelayer}
		\node [style=port] (19) at (5, 4) {};
		\node [style=port] (20) at (4.5, 1.5) {};
		\node [style=component] (21) at (5, 3) {$\Delta$};
		\node [style=none] (22) at (4.5, 2.5) {};
		\node [style=none] (23) at (5.5, 2.5) {};
		\node [style=component] (27) at (5.5, 2) {$\mathsf{e}$};
	\end{pgfonlayer}
	\begin{pgfonlayer}{edgelayer}
		\draw [style=wire] (19) to (21);
		\draw [style=wire] (21) to (22.center);
		\draw [style=wire] (21) to (23.center);
		\draw [style=wire] (22.center) to (20);
		\draw [style=wire] (23.center) to (27);
	\end{pgfonlayer}
\end{tikzpicture}
 \end{array} = \begin{array}[c]{c} \begin{tikzpicture}
	\begin{pgfonlayer}{nodelayer}
		\node [style=port] (16) at (6.5, 1.5) {};
		\node [style=port] (18) at (6.5, 4) {};
	\end{pgfonlayer}
	\begin{pgfonlayer}{edgelayer}
		\draw [style=wire] (18) to (16);
	\end{pgfonlayer}
\end{tikzpicture}
 \end{array} = \begin{array}[c]{c} \begin{tikzpicture}
	\begin{pgfonlayer}{nodelayer}
		\node [style=port] (19) at (5, 4) {};
		\node [style=port] (20) at (5.5, 1.5) {};
		\node [style=component] (21) at (5, 3) {$\Delta$};
		\node [style=none] (22) at (5.5, 2.5) {};
		\node [style=none] (23) at (4.5, 2.5) {};
		\node [style=component] (27) at (4.5, 2) {$\mathsf{e}$};
	\end{pgfonlayer}
	\begin{pgfonlayer}{edgelayer}
		\draw [style=wire] (19) to (21);
		\draw [style=wire] (21) to (22.center);
		\draw [style=wire] (21) to (23.center);
		\draw [style=wire] (22.center) to (20);
		\draw [style=wire] (23.center) to (27);
	\end{pgfonlayer}
\end{tikzpicture}
 \end{array} && \begin{array}[c]{c}\begin{tikzpicture}
	\begin{pgfonlayer}{nodelayer}
		\node [style=port] (35) at (12.5, 4) {};
		\node [style=port] (36) at (12, 1.25) {};
		\node [style=component] (37) at (12.5, 3) {$\Delta$};
		\node [style=none] (38) at (12, 2.5) {};
		\node [style=none] (39) at (13, 2.5) {};
		\node [style=port] (40) at (13, 1.25) {};
		\node [style=none] (41) at (12, 2.25) {};
		\node [style=none] (42) at (13, 2.25) {};
		\node [style=none] (43) at (12, 1.75) {};
		\node [style=none] (44) at (13, 1.75) {};
	\end{pgfonlayer}
	\begin{pgfonlayer}{edgelayer}
		\draw [style=wire] (35) to (37);
		\draw [style=wire] (37) to (38.center);
		\draw [style=wire] (37) to (39.center);
		\draw [style=wire] (38.center) to (41.center);
		\draw [style=wire] (39.center) to (42.center);
		\draw [style=wire] (43.center) to (36);
		\draw [style=wire] (44.center) to (40);
		\draw [style=wire] (41.center) to (44.center);
		\draw [style=wire] (42.center) to (43.center);
	\end{pgfonlayer}
\end{tikzpicture}
 \end{array} = \begin{array}[c]{c} \begin{tikzpicture}
	\begin{pgfonlayer}{nodelayer}
		\node [style=port] (19) at (10.5, 4) {};
		\node [style=port] (20) at (10, 1.25) {};
		\node [style=component] (21) at (10.5, 3) {$\Delta$};
		\node [style=none] (22) at (10, 2.5) {};
		\node [style=none] (23) at (11, 2.5) {};
		\node [style=port] (34) at (11, 1.25) {};
	\end{pgfonlayer}
	\begin{pgfonlayer}{edgelayer}
		\draw [style=wire] (19) to (21);
		\draw [style=wire] (21) to (22.center);
		\draw [style=wire] (21) to (23.center);
		\draw [style=wire] (22.center) to (20);
		\draw [style=wire] (23.center) to (34);
	\end{pgfonlayer}
\end{tikzpicture}
 \end{array}
 \end{align}
 where the first equality is called the coassociativity of the comultiplication, the second is the counit identity, and the third is called the cocommutativity of the comultiplication. It is important to note that naturality of $\Delta$ and $\mathsf{e}$ is precisely the statement that for every map $f$, $\oc(f)$ is a comonoid morphism. The last requirement is that $\delta$ is also a comonoid morphism (\ref{diag:deltacomonoid}), which is drawn as follows: 
  \begin{align}\label{string:delta-comonoid-map}\begin{array}[c]{c}\begin{tikzpicture}
	\begin{pgfonlayer}{nodelayer}
		\node [style=component] (7) at (-3.5, 3.25) {$\delta$};
		\node [style=port] (20) at (-4, 1.25) {};
		\node [style=component] (21) at (-3.5, 2.25) {$\Delta$};
		\node [style=none] (22) at (-4, 1.75) {};
		\node [style=none] (23) at (-3, 1.75) {};
		\node [style=port] (34) at (-3, 1.25) {};
		\node [style=port] (45) at (-3.5, 4) {};
	\end{pgfonlayer}
	\begin{pgfonlayer}{edgelayer}
		\draw [style=wire] (21) to (22.center);
		\draw [style=wire] (21) to (23.center);
		\draw [style=wire] (22.center) to (20);
		\draw [style=wire] (23.center) to (34);
		\draw [style=wire] (45) to (7);
		\draw [style=wire] (7) to (21);
	\end{pgfonlayer}
\end{tikzpicture}
 \end{array} = \begin{array}[c]{c}\begin{tikzpicture}
	\begin{pgfonlayer}{nodelayer}
		\node [style=component] (46) at (-1.5, 1.75) {$\delta$};
		\node [style=port] (51) at (-1.5, 1) {};
		\node [style=component] (52) at (-1, 3) {$\Delta$};
		\node [style=none] (53) at (-1.5, 2.5) {};
		\node [style=none] (54) at (-0.5, 2.5) {};
		\node [style=port] (55) at (-0.5, 1) {};
		\node [style=component] (56) at (-0.5, 1.75) {$\delta$};
		\node [style=port] (61) at (-1, 3.75) {};
	\end{pgfonlayer}
	\begin{pgfonlayer}{edgelayer}
		\draw [style=wire] (52) to (53.center);
		\draw [style=wire] (52) to (54.center);
		\draw [style=wire] (46) to (51);
		\draw [style=wire] (56) to (55);
		\draw [style=wire] (53.center) to (46);
		\draw [style=wire] (54.center) to (56);
		\draw [style=wire] (61) to (52);
	\end{pgfonlayer}
\end{tikzpicture}
 \end{array} && \begin{array}[c]{c}\begin{tikzpicture}
	\begin{pgfonlayer}{nodelayer}
		\node [style=component] (7) at (-3.5, 3.25) {$\delta$};
		\node [style=component] (21) at (-3.5, 2.25) {$\mathsf{e}$};
		\node [style=port] (45) at (-3.5, 4) {};
	\end{pgfonlayer}
	\begin{pgfonlayer}{edgelayer}
		\draw [style=wire] (45) to (7);
		\draw [style=wire] (7) to (21);
	\end{pgfonlayer}
\end{tikzpicture}
 \end{array} = \begin{array}[c]{c}\begin{tikzpicture}
	\begin{pgfonlayer}{nodelayer}
		\node [style=component] (52) at (-0.75, 3.25) {$\mathsf{e}$};
		\node [style=port] (61) at (-0.75, 4.25) {};
	\end{pgfonlayer}
	\begin{pgfonlayer}{edgelayer}
		\draw [style=wire] (61) to (52);
	\end{pgfonlayer}
\end{tikzpicture}
 \end{array}
 \end{align}

Finally, a \emph{monoidal} coalgebra modality (also sometimes called a \emph{linear exponential modality} \cite[Def 10]{schalk2004categorical}) is a coalgebra modality whose underlying endofunctor $\oc$ is symmetric monoidal, all the natural transformations are monoidal, and the comultiplication and counit are also $\oc$-coalgebra morphisms. Let us break this down. First recall that a \textbf{symmetric monoidal endofunctor} on a symmetric monoidal category $\mathbb{X}$ is a triple $(\oc, \mathsf{m}, \mathsf{m}_I)$ consisting of an endofunctor ${\oc: \mathbb{X} \to \mathbb{X}}$, a natural transformation ${\mathsf{m}_{A,B}: \oc(A) \otimes \oc(B) \to \oc(A \otimes B)}$, and a map $\mathsf{m}_I: I \to \oc(I)$, such that the diagrams in Appendix \ref{sec:diagramsforsmc} commute. The natural transformation $\mathsf{m}$ and the map $\mathsf{m}_I$ are drawn as follows\footnote{We note that we are drawing $\mathsf{m}_{A,B}$ differently then how it was drawn in \cite{Blute2019} to help avoid some confusion.}: 
 \begin{align*}\begin{array}[c]{c} \begin{tikzpicture}
	\begin{pgfonlayer}{nodelayer}
		\node [style=object] (62) at (0, -2.25) {$\oc(A\otimes B)$};
		\node [style=object] (63) at (-0.5, 0) {$\oc(A)$};
		\node [style=component] (64) at (0, -1.25) {$\mathsf{m}$};
		\node [style=none] (65) at (-0.5, -0.75) {};
		\node [style=none] (66) at (0.5, -0.75) {};
		\node [style=object] (67) at (0.5, 0) {$\oc(B)$};
	\end{pgfonlayer}
	\begin{pgfonlayer}{edgelayer}
		\draw [style=wire] (62) to (64);
		\draw [style=wire] (64) to (65.center);
		\draw [style=wire] (64) to (66.center);
		\draw [style=wire] (66.center) to (67);
		\draw [style=wire] (65.center) to (63);
	\end{pgfonlayer}
\end{tikzpicture}
 \end{array} && \begin{array}[c]{c}\begin{tikzpicture}
	\begin{pgfonlayer}{nodelayer}
		\node [style=component] (52) at (2, -1) {$\mathsf{m}$};
		\node [style=object] (61) at (2, -2.25) {$\oc(I)$};
	\end{pgfonlayer}
	\begin{pgfonlayer}{edgelayer}
		\draw [style=wire] (61) to (52);
	\end{pgfonlayer}
\end{tikzpicture}
 \end{array} 
 \end{align*}
with the naturality of $\mathsf{m}$ drawn as follows:
 \begin{align}\label{string:monoidal-nat}\begin{array}[c]{c}\begin{tikzpicture}
	\begin{pgfonlayer}{nodelayer}
		\node [style=component] (134) at (22.25, -3.25) {$f$};
		\node [style=none] (135) at (21.75, -3.75) {};
		\node [style=none] (136) at (22.75, -3.75) {};
		\node [style=none] (137) at (21.75, -2.75) {};
		\node [style=none] (138) at (22.75, -2.75) {};
		\node [style=port] (139) at (22.25, -1.75) {};
		\node [style=component] (140) at (23, -5) {$\mathsf{m}$};
		\node [style=none] (141) at (22.25, -4.25) {};
		\node [style=none] (142) at (23.75, -4.25) {};
		\node [style=port] (143) at (23.75, -1.75) {};
		\node [style=component] (144) at (23.75, -3.25) {$g$};
		\node [style=none] (145) at (23.25, -3.75) {};
		\node [style=none] (146) at (24.25, -3.75) {};
		\node [style=none] (147) at (23.25, -2.75) {};
		\node [style=none] (148) at (24.25, -2.75) {};
		\node [style=port] (149) at (23, -6) {};
	\end{pgfonlayer}
	\begin{pgfonlayer}{edgelayer}
		\draw [style=wire] (135.center) to (136.center);
		\draw [style=wire] (136.center) to (138.center);
		\draw [style=wire] (138.center) to (137.center);
		\draw [style=wire] (137.center) to (135.center);
		\draw [style=wire] (140) to (141.center);
		\draw [style=wire] (140) to (142.center);
		\draw [style=wire] (145.center) to (146.center);
		\draw [style=wire] (146.center) to (148.center);
		\draw [style=wire] (148.center) to (147.center);
		\draw [style=wire] (147.center) to (145.center);
		\draw [style=wire] (134) to (139);
		\draw [style=wire] (144) to (143);
		\draw [style=wire] (141.center) to (134);
		\draw [style=wire] (142.center) to (144);
		\draw [style=wire] (149) to (140);
	\end{pgfonlayer}
\end{tikzpicture}
 \end{array} = \begin{array}[c]{c}\begin{tikzpicture}
	\begin{pgfonlayer}{nodelayer}
		\node [style=component] (123) at (19.75, -4.5) {$f$};
		\node [style=none] (124) at (19.25, -5.25) {};
		\node [style=none] (125) at (21.25, -5.25) {};
		\node [style=none] (126) at (19.25, -3.75) {};
		\node [style=none] (127) at (21.25, -3.75) {};
		\node [style=port] (128) at (19.75, -1.75) {};
		\node [style=component] (129) at (20.25, -3) {$\mathsf{m}$};
		\node [style=none] (130) at (19.75, -2.5) {};
		\node [style=none] (131) at (20.75, -2.5) {};
		\node [style=port] (132) at (20.75, -1.75) {};
		\node [style=port] (133) at (20.25, -6) {};
		\node [style=component] (150) at (20.75, -4.5) {$g$};
		\node [style=none] (151) at (19.75, -3.75) {};
		\node [style=none] (152) at (20.75, -3.75) {};
		\node [style=none] (153) at (19.75, -5.25) {};
		\node [style=none] (154) at (20.75, -5.25) {};
		\node [style=none] (155) at (20.25, -3.75) {};
		\node [style=none] (156) at (20.25, -5.25) {};
	\end{pgfonlayer}
	\begin{pgfonlayer}{edgelayer}
		\draw [style=wire] (124.center) to (125.center);
		\draw [style=wire] (125.center) to (127.center);
		\draw [style=wire] (127.center) to (126.center);
		\draw [style=wire] (126.center) to (124.center);
		\draw [style=wire] (129) to (130.center);
		\draw [style=wire] (129) to (131.center);
		\draw [style=wire] (130.center) to (128);
		\draw [style=wire] (131.center) to (132);
		\draw [style=wire] (123) to (153.center);
		\draw [style=wire] (150) to (154.center);
		\draw [style=wire] (152.center) to (150);
		\draw [style=wire] (151.center) to (123);
		\draw [style=wire] (156.center) to (133);
		\draw [style=wire] (129) to (155.center);
	\end{pgfonlayer}
\end{tikzpicture}
 \end{array}
 \end{align}
The axioms for symmetric monoidal endofunctor (\ref{diag:smendo}) are drawn as follows: 
  \begin{align}\label{strings:monoidal-functor}\begin{array}[c]{c} \begin{tikzpicture}
	\begin{pgfonlayer}{nodelayer}
		\node [style=port] (68) at (16, 0) {};
		\node [style=port] (69) at (15.5, 3.5) {};
		\node [style=component] (70) at (16, 1) {$\mathsf{m}$};
		\node [style=none] (71) at (15.5, 1.5) {};
		\node [style=none] (72) at (16.5, 1.5) {};
		\node [style=port] (73) at (16, 3.5) {};
		\node [style=component] (74) at (16.5, 2.25) {$\mathsf{m}$};
		\node [style=none] (75) at (16, 2.75) {};
		\node [style=none] (76) at (17, 2.75) {};
		\node [style=port] (77) at (17, 3.5) {};
	\end{pgfonlayer}
	\begin{pgfonlayer}{edgelayer}
		\draw [style=wire] (68) to (70);
		\draw [style=wire] (70) to (71.center);
		\draw [style=wire] (70) to (72.center);
		\draw [style=wire] (71.center) to (69);
		\draw [style=wire] (74) to (75.center);
		\draw [style=wire] (74) to (76.center);
		\draw [style=wire] (76.center) to (77);
		\draw [style=wire] (75.center) to (73);
		\draw [style=wire] (72.center) to (74);
	\end{pgfonlayer}
\end{tikzpicture}
 \end{array} = \begin{array}[c]{c} \begin{tikzpicture}
	\begin{pgfonlayer}{nodelayer}
		\node [style=port] (68) at (16.5, 0) {};
		\node [style=port] (69) at (17, 3.5) {};
		\node [style=component] (70) at (16.5, 1) {$\mathsf{m}$};
		\node [style=none] (71) at (17, 1.5) {};
		\node [style=none] (72) at (16, 1.5) {};
		\node [style=port] (73) at (16.5, 3.5) {};
		\node [style=component] (74) at (16, 2.25) {$\mathsf{m}$};
		\node [style=none] (75) at (16.5, 2.75) {};
		\node [style=none] (76) at (15.5, 2.75) {};
		\node [style=port] (77) at (15.5, 3.5) {};
	\end{pgfonlayer}
	\begin{pgfonlayer}{edgelayer}
		\draw [style=wire] (68) to (70);
		\draw [style=wire] (70) to (71.center);
		\draw [style=wire] (70) to (72.center);
		\draw [style=wire] (71.center) to (69);
		\draw [style=wire] (74) to (75.center);
		\draw [style=wire] (74) to (76.center);
		\draw [style=wire] (76.center) to (77);
		\draw [style=wire] (75.center) to (73);
		\draw [style=wire] (72.center) to (74);
	\end{pgfonlayer}
\end{tikzpicture}
 \end{array} && \begin{array}[c]{c} \begin{tikzpicture}
	\begin{pgfonlayer}{nodelayer}
		\node [style=port] (78) at (18.5, 1.5) {};
		\node [style=port] (79) at (18, 4) {};
		\node [style=component] (80) at (18.5, 2.5) {$\mathsf{m}$};
		\node [style=none] (81) at (18, 3) {};
		\node [style=none] (82) at (19, 3) {};
		\node [style=component] (83) at (19, 3.5) {$\mathsf{m}$};
	\end{pgfonlayer}
	\begin{pgfonlayer}{edgelayer}
		\draw [style=wire] (78) to (80);
		\draw [style=wire] (80) to (81.center);
		\draw [style=wire] (80) to (82.center);
		\draw [style=wire] (81.center) to (79);
		\draw [style=wire] (82.center) to (83);
	\end{pgfonlayer}
\end{tikzpicture}
 \end{array} = \begin{array}[c]{c} \begin{tikzpicture}
	\begin{pgfonlayer}{nodelayer}
		\node [style=port] (16) at (6.5, 1.5) {};
		\node [style=port] (18) at (6.5, 4) {};
	\end{pgfonlayer}
	\begin{pgfonlayer}{edgelayer}
		\draw [style=wire] (18) to (16);
	\end{pgfonlayer}
\end{tikzpicture}
 \end{array} = \begin{array}[c]{c} \begin{tikzpicture}
	\begin{pgfonlayer}{nodelayer}
		\node [style=port] (78) at (18.5, 1.5) {};
		\node [style=port] (79) at (19, 4) {};
		\node [style=component] (80) at (18.5, 2.5) {$\mathsf{m}$};
		\node [style=none] (81) at (19, 3) {};
		\node [style=none] (82) at (18, 3) {};
		\node [style=component] (83) at (18, 3.5) {$\mathsf{m}$};
	\end{pgfonlayer}
	\begin{pgfonlayer}{edgelayer}
		\draw [style=wire] (78) to (80);
		\draw [style=wire] (80) to (81.center);
		\draw [style=wire] (80) to (82.center);
		\draw [style=wire] (81.center) to (79);
		\draw [style=wire] (82.center) to (83);
	\end{pgfonlayer}
\end{tikzpicture}
 \end{array} && \begin{array}[c]{c}\begin{tikzpicture}
	\begin{pgfonlayer}{nodelayer}
		\node [style=port] (106) at (11.5, -5.25) {};
		\node [style=port] (107) at (11, -1) {};
		\node [style=component] (108) at (11.5, -4.25) {$\mathsf{m}$};
		\node [style=none] (109) at (11, -3.5) {};
		\node [style=none] (110) at (12, -3.5) {};
		\node [style=port] (111) at (12, -1) {};
		\node [style=none] (112) at (11, -2.75) {};
		\node [style=none] (113) at (12, -2.75) {};
		\node [style=none] (114) at (11, -2) {};
		\node [style=none] (115) at (12, -2) {};
	\end{pgfonlayer}
	\begin{pgfonlayer}{edgelayer}
		\draw [style=wire] (106) to (108);
		\draw [style=wire] (108) to (109.center);
		\draw [style=wire] (108) to (110.center);
		\draw [style=wire] (109.center) to (112.center);
		\draw [style=wire] (110.center) to (113.center);
		\draw [style=wire] (114.center) to (107);
		\draw [style=wire] (115.center) to (111);
		\draw [style=wire] (112.center) to (115.center);
		\draw [style=wire] (113.center) to (114.center);
	\end{pgfonlayer}
\end{tikzpicture}
 \end{array} = \begin{array}[c]{c} \begin{tikzpicture}
	\begin{pgfonlayer}{nodelayer}
		\node [style=port] (85) at (13.25, -3) {};
		\node [style=port] (87) at (13.25, -4.5) {};
		\node [style=port] (88) at (14.25, -4.5) {};
		\node [style=port] (89) at (14.25, -3) {};
		\node [style=none] (90) at (13.25, -4) {};
		\node [style=none] (91) at (14.25, -4) {};
		\node [style=none] (92) at (13.25, -3.5) {};
		\node [style=none] (93) at (14.25, -3.5) {};
		\node [style=port] (94) at (13.75, -5.25) {};
		\node [style=port] (95) at (13.25, -1) {};
		\node [style=component] (96) at (13.75, -2.25) {$\mathsf{m}$};
		\node [style=none] (97) at (13.25, -1.75) {};
		\node [style=none] (98) at (14.25, -1.75) {};
		\node [style=port] (99) at (14.25, -1) {};
		\node [style=none] (100) at (13, -3) {};
		\node [style=none] (101) at (13, -4.5) {};
		\node [style=none] (102) at (14.5, -4.5) {};
		\node [style=none] (103) at (14.5, -3) {};
		\node [style=none] (104) at (13.75, -4.5) {};
		\node [style=none] (105) at (13.75, -3) {};
	\end{pgfonlayer}
	\begin{pgfonlayer}{edgelayer}
		\draw [style=wire] (87) to (90.center);
		\draw [style=wire] (88) to (91.center);
		\draw [style=wire] (92.center) to (85);
		\draw [style=wire] (93.center) to (89);
		\draw [style=wire] (90.center) to (93.center);
		\draw [style=wire] (91.center) to (92.center);
		\draw [style=wire] (96) to (97.center);
		\draw [style=wire] (96) to (98.center);
		\draw [style=wire] (97.center) to (95);
		\draw [style=wire] (98.center) to (99);
		\draw [style=wire] (100.center) to (103.center);
		\draw [style=wire] (103.center) to (102.center);
		\draw [style=wire] (102.center) to (101.center);
		\draw [style=wire] (101.center) to (100.center);
		\draw [style=wire] (104.center) to (94);
		\draw [style=wire] (96) to (105.center);
	\end{pgfonlayer}
\end{tikzpicture}
 \end{array}
 \end{align}
A \textbf{monoidal coalgebra modality} \cite[Def 2]{Blute2019} is a septuple $(\oc, \delta, \varepsilon, \Delta, \mathsf{e}, \mathsf{m}, \mathsf{m}_I)$ consisting of a coalgebra modality $(\oc, \delta, \varepsilon, \Delta, \mathsf{e})$ and a symmetric monoidal endofunctor $(\oc, \mathsf{m}, \mathsf{m}_I)$ such that the diagrams in Appendix \ref{sec:diagramsformcoalg} commute. The requirement that $\delta$ and $\varepsilon$ are monoidal transformations (\ref{diag:monoidalcomonad}), which gives us a monoidal comonad, is drawn follows: 
  \begin{align}\label{string:comonad-monoidal}\begin{array}[c]{c}\begin{tikzpicture}
	\begin{pgfonlayer}{nodelayer}
		\node [style=component] (116) at (22.25, -1.25) {$\delta$};
		\node [style=port] (117) at (21.75, 2.5) {};
		\node [style=component] (118) at (22.25, 0) {$\mathsf{m}$};
		\node [style=none] (119) at (21.75, 1) {};
		\node [style=none] (120) at (22.75, 1) {};
		\node [style=port] (121) at (22.75, 2.5) {};
		\node [style=port] (122) at (22.25, -2.25) {};
	\end{pgfonlayer}
	\begin{pgfonlayer}{edgelayer}
		\draw [style=wire] (118) to (119.center);
		\draw [style=wire] (118) to (120.center);
		\draw [style=wire] (119.center) to (117);
		\draw [style=wire] (120.center) to (121);
		\draw [style=wire] (122) to (116);
		\draw [style=wire] (116) to (118);
	\end{pgfonlayer}
\end{tikzpicture}
 \end{array} = \begin{array}[c]{c}\begin{tikzpicture}
	\begin{pgfonlayer}{nodelayer}
		\node [style=none] (124) at (23.5, -1.5) {};
		\node [style=none] (125) at (25.5, -1.5) {};
		\node [style=none] (126) at (23.5, 0) {};
		\node [style=none] (127) at (25.5, 0) {};
		\node [style=port] (128) at (24, 2.5) {};
		\node [style=component] (129) at (24.5, 0.75) {$\mathsf{m}$};
		\node [style=none] (130) at (24, 1.25) {};
		\node [style=none] (131) at (25, 1.25) {};
		\node [style=port] (132) at (25, 2.5) {};
		\node [style=port] (133) at (24.5, -2.25) {};
		\node [style=none] (151) at (24, 0) {};
		\node [style=none] (152) at (25, 0) {};
		\node [style=none] (153) at (24, -1.5) {};
		\node [style=none] (154) at (25, -1.5) {};
		\node [style=none] (155) at (24.5, 0) {};
		\node [style=none] (156) at (24.5, -1.5) {};
		\node [style=component] (157) at (24.5, -1) {$\mathsf{m}$};
		\node [style=none] (158) at (24, -0.5) {};
		\node [style=none] (159) at (25, -0.5) {};
		\node [style=component] (160) at (24, 1.75) {$\delta$};
		\node [style=component] (161) at (25, 1.75) {$\delta$};
	\end{pgfonlayer}
	\begin{pgfonlayer}{edgelayer}
		\draw [style=wire] (124.center) to (125.center);
		\draw [style=wire] (125.center) to (127.center);
		\draw [style=wire] (127.center) to (126.center);
		\draw [style=wire] (126.center) to (124.center);
		\draw [style=wire] (129) to (130.center);
		\draw [style=wire] (129) to (131.center);
		\draw [style=wire] (156.center) to (133);
		\draw [style=wire] (129) to (155.center);
		\draw [style=wire] (157) to (158.center);
		\draw [style=wire] (157) to (159.center);
		\draw [style=wire] (151.center) to (158.center);
		\draw [style=wire] (152.center) to (159.center);
		\draw [style=wire] (157) to (156.center);
		\draw [style=wire] (128) to (160);
		\draw [style=wire] (132) to (161);
		\draw [style=wire] (161) to (131.center);
		\draw [style=wire] (160) to (130.center);
	\end{pgfonlayer}
\end{tikzpicture}
 \end{array} && \begin{array}[c]{c}\begin{tikzpicture}
	\begin{pgfonlayer}{nodelayer}
		\node [style=component] (356) at (65, -1) {$\delta$};
		\node [style=port] (357) at (65, -2) {};
		\node [style=component] (358) at (65, 0) {$\mathsf{m}$};
	\end{pgfonlayer}
	\begin{pgfonlayer}{edgelayer}
		\draw [style=wire] (357) to (356);
		\draw [style=wire] (358) to (356);
	\end{pgfonlayer}
\end{tikzpicture}
 \end{array} = \begin{array}[c]{c}\begin{tikzpicture}
	\begin{pgfonlayer}{nodelayer}
		\node [style=component] (359) at (67, 0) {$\mathsf{m}$};
		\node [style=component] (360) at (67, -1) {$\mathsf{m}$};
		\node [style=none] (361) at (66.5, -1.5) {};
		\node [style=none] (362) at (67.5, -1.5) {};
		\node [style=none] (363) at (66.5, -0.5) {};
		\node [style=none] (364) at (67.5, -0.5) {};
		\node [style=port] (365) at (67, -2) {};
		\node [style=none] (366) at (67, -0.5) {};
	\end{pgfonlayer}
	\begin{pgfonlayer}{edgelayer}
		\draw [style=wire] (361.center) to (362.center);
		\draw [style=wire] (362.center) to (364.center);
		\draw [style=wire] (364.center) to (363.center);
		\draw [style=wire] (363.center) to (361.center);
		\draw [style=wire] (365) to (360);
		\draw [style=wire] (359) to (366.center);
	\end{pgfonlayer}
\end{tikzpicture}
 \end{array} && \begin{array}[c]{c}\begin{tikzpicture}
	\begin{pgfonlayer}{nodelayer}
		\node [style=port] (1065) at (189.25, -3.25) {};
		\node [style=port] (1066) at (190.25, -3.25) {};
		\node [style=component] (1067) at (189.75, -1) {$\mathsf{m}$};
		\node [style=component] (1068) at (189.75, -2) {$\varepsilon$};
		\node [style=none] (1069) at (189.25, -2.5) {};
		\node [style=none] (1070) at (190.25, -2.5) {};
		\node [style=port] (1071) at (189.25, 0) {};
		\node [style=port] (1072) at (190.25, 0) {};
		\node [style=none] (1073) at (189.25, -0.75) {};
		\node [style=none] (1074) at (190.25, -0.75) {};
	\end{pgfonlayer}
	\begin{pgfonlayer}{edgelayer}
		\draw [style=wire] (1067) to (1068);
		\draw [style=wire] (1068) to (1069.center);
		\draw [style=wire] (1068) to (1070.center);
		\draw [style=wire] (1069.center) to (1065);
		\draw [style=wire] (1070.center) to (1066);
		\draw [style=wire] (1073.center) to (1067);
		\draw [style=wire] (1067) to (1074.center);
		\draw [style=wire] (1071) to (1073.center);
		\draw [style=wire] (1074.center) to (1072);
	\end{pgfonlayer}
\end{tikzpicture}
 \end{array} = \begin{array}[c]{c}\begin{tikzpicture}
	\begin{pgfonlayer}{nodelayer}
		\node [style=port] (128) at (24, 2) {};
		\node [style=port] (130) at (24, -1) {};
		\node [style=port] (131) at (25, -1) {};
		\node [style=port] (132) at (25, 2) {};
		\node [style=component] (160) at (24, 0.5) {$\varepsilon$};
		\node [style=component] (161) at (25, 0.5) {$\varepsilon$};
	\end{pgfonlayer}
	\begin{pgfonlayer}{edgelayer}
		\draw [style=wire] (128) to (160);
		\draw [style=wire] (132) to (161);
		\draw [style=wire] (161) to (131);
		\draw [style=wire] (160) to (130);
	\end{pgfonlayer}
\end{tikzpicture}
 \end{array} && \begin{array}[c]{c}\begin{tikzpicture}
	\begin{pgfonlayer}{nodelayer}
		\node [style=component] (215) at (33, -0.5) {$\varepsilon$};
		\node [style=component] (216) at (33, 0.5) {$\mathsf{m}$};
	\end{pgfonlayer}
	\begin{pgfonlayer}{edgelayer}
		\draw [style=wire] (216) to (215);
	\end{pgfonlayer}
\end{tikzpicture}
 \end{array} = 
 \end{align}
 That $\Delta$ and $\mathsf{e}$ are monoidal transformations (\ref{diag:monoidalcomonoid}) is drawn as follows:
  \begin{align}\label{string:monoidalcomonoid}\begin{array}[c]{c}\begin{tikzpicture}
	\begin{pgfonlayer}{nodelayer}
		\node [style=port] (29) at (17, -5.5) {};
		\node [style=component] (30) at (17.5, -4.25) {$\Delta$};
		\node [style=none] (31) at (17, -4.75) {};
		\node [style=none] (32) at (18, -4.75) {};
		\node [style=port] (33) at (18, -5.5) {};
		\node [style=port] (162) at (17, -1.75) {};
		\node [style=component] (163) at (17.5, -3) {$\mathsf{m}$};
		\node [style=none] (164) at (17, -2.5) {};
		\node [style=none] (165) at (18, -2.5) {};
		\node [style=port] (166) at (18, -1.75) {};
	\end{pgfonlayer}
	\begin{pgfonlayer}{edgelayer}
		\draw [style=wire] (30) to (31.center);
		\draw [style=wire] (30) to (32.center);
		\draw [style=wire] (32.center) to (33);
		\draw [style=wire] (31.center) to (29);
		\draw [style=wire] (163) to (164.center);
		\draw [style=wire] (163) to (165.center);
		\draw [style=wire] (164.center) to (162);
		\draw [style=wire] (165.center) to (166);
		\draw [style=wire] (163) to (30);
	\end{pgfonlayer}
\end{tikzpicture}
 \end{array} = \begin{array}[c]{c}\begin{tikzpicture}
	\begin{pgfonlayer}{nodelayer}
		\node [style=component] (167) at (20.25, -2.5) {$\Delta$};
		\node [style=none] (168) at (19.75, -3) {};
		\node [style=none] (169) at (20.75, -3) {};
		\node [style=port] (170) at (20.25, -1.75) {};
		\node [style=port] (171) at (21.75, -1.75) {};
		\node [style=component] (172) at (21.75, -2.5) {$\Delta$};
		\node [style=none] (173) at (21.25, -3) {};
		\node [style=none] (174) at (22.25, -3) {};
		\node [style=none] (179) at (20.75, -4) {};
		\node [style=none] (180) at (21.25, -4) {};
		\node [style=none] (181) at (20.75, -3.25) {};
		\node [style=none] (182) at (21.25, -3.25) {};
		\node [style=component] (183) at (20.25, -4.75) {$\mathsf{m}$};
		\node [style=none] (184) at (19.75, -4.25) {};
		\node [style=none] (185) at (20.75, -4.25) {};
		\node [style=component] (186) at (21.75, -4.75) {$\mathsf{m}$};
		\node [style=none] (187) at (21.25, -4.25) {};
		\node [style=none] (188) at (22.25, -4.25) {};
		\node [style=port] (189) at (20.25, -5.5) {};
		\node [style=port] (190) at (21.75, -5.5) {};
	\end{pgfonlayer}
	\begin{pgfonlayer}{edgelayer}
		\draw [style=wire] (167) to (168.center);
		\draw [style=wire] (167) to (169.center);
		\draw [style=wire] (172) to (173.center);
		\draw [style=wire] (172) to (174.center);
		\draw [style=wire] (170) to (167);
		\draw [style=wire] (171) to (172);
		\draw [style=wire] (179.center) to (182.center);
		\draw [style=wire] (180.center) to (181.center);
		\draw [style=wire] (169.center) to (181.center);
		\draw [style=wire] (173.center) to (182.center);
		\draw [style=wire] (183) to (184.center);
		\draw [style=wire] (183) to (185.center);
		\draw [style=wire] (186) to (187.center);
		\draw [style=wire] (186) to (188.center);
		\draw [style=wire] (168.center) to (184.center);
		\draw [style=wire] (174.center) to (188.center);
		\draw [style=wire] (179.center) to (185.center);
		\draw [style=wire] (180.center) to (187.center);
		\draw [style=wire] (183) to (189);
		\draw [style=wire] (186) to (190);
	\end{pgfonlayer}
\end{tikzpicture}
 \end{array} && \begin{array}[c]{c}\begin{tikzpicture}
	\begin{pgfonlayer}{nodelayer}
		\node [style=component] (116) at (22.25, -0.25) {$\mathsf{e}$};
		\node [style=port] (117) at (21.75, 2) {};
		\node [style=component] (118) at (22.25, 0.75) {$\mathsf{m}$};
		\node [style=none] (119) at (21.75, 1.25) {};
		\node [style=none] (120) at (22.75, 1.25) {};
		\node [style=port] (121) at (22.75, 2) {};
	\end{pgfonlayer}
	\begin{pgfonlayer}{edgelayer}
		\draw [style=wire] (118) to (119.center);
		\draw [style=wire] (118) to (120.center);
		\draw [style=wire] (119.center) to (117);
		\draw [style=wire] (120.center) to (121);
		\draw [style=wire] (116) to (118);
	\end{pgfonlayer}
\end{tikzpicture}
 \end{array} = \begin{array}[c]{c}\begin{tikzpicture}
	\begin{pgfonlayer}{nodelayer}
		\node [style=port] (128) at (24, 2) {};
		\node [style=port] (132) at (25, 2) {};
		\node [style=component] (160) at (24, 0.5) {$\mathsf{e}$};
		\node [style=component] (161) at (25, 0.5) {$\mathsf{e}$};
	\end{pgfonlayer}
	\begin{pgfonlayer}{edgelayer}
		\draw [style=wire] (128) to (160);
		\draw [style=wire] (132) to (161);
	\end{pgfonlayer}
\end{tikzpicture}
 \end{array} && \begin{array}[c]{c}\begin{tikzpicture}
	\begin{pgfonlayer}{nodelayer}
		\node [style=component] (116) at (22.25, 2) {$\mathsf{m}$};
		\node [style=port] (117) at (21.75, -0.25) {};
		\node [style=component] (118) at (22.25, 1) {$\Delta$};
		\node [style=none] (119) at (21.75, 0.5) {};
		\node [style=none] (120) at (22.75, 0.5) {};
		\node [style=port] (121) at (22.75, -0.25) {};
	\end{pgfonlayer}
	\begin{pgfonlayer}{edgelayer}
		\draw [style=wire] (118) to (119.center);
		\draw [style=wire] (118) to (120.center);
		\draw [style=wire] (119.center) to (117);
		\draw [style=wire] (120.center) to (121);
		\draw [style=wire] (116) to (118);
	\end{pgfonlayer}
\end{tikzpicture}
 \end{array} = \begin{array}[c]{c}\begin{tikzpicture}
	\begin{pgfonlayer}{nodelayer}
		\node [style=port] (128) at (24, 0.5) {};
		\node [style=port] (132) at (25, 0.5) {};
		\node [style=component] (160) at (24, 2) {$\mathsf{m}$};
		\node [style=component] (161) at (25, 2) {$\mathsf{m}$};
	\end{pgfonlayer}
	\begin{pgfonlayer}{edgelayer}
		\draw [style=wire] (128) to (160);
		\draw [style=wire] (132) to (161);
	\end{pgfonlayer}
\end{tikzpicture}
 \end{array} && \begin{array}[c]{c}\begin{tikzpicture}
	\begin{pgfonlayer}{nodelayer}
		\node [style=component] (215) at (33, -0.5) {$\mathsf{e}$};
		\node [style=component] (216) at (33, 0.5) {$\mathsf{m}$};
	\end{pgfonlayer}
	\begin{pgfonlayer}{edgelayer}
		\draw [style=wire] (216) to (215);
	\end{pgfonlayer}
\end{tikzpicture}
 \end{array} = 
 \end{align}
 We note that saying that $\Delta$ and $\mathsf{e}$ are monoidal is equivalent to saying that $\mathsf{m}$ and $\mathsf{m}_I$ are comonoid morphisms. Lastly, that $\Delta$ and $\mathsf{e}$ are $\oc$-coalgebra morphisms (\ref{diag:!coalgcomonoid}) is drawn out as follows:
  \begin{align}\label{string:comonoid-!coalg}\begin{array}[c]{c}\begin{tikzpicture}
	\begin{pgfonlayer}{nodelayer}
		\node [style=none] (191) at (30.5, -1) {};
		\node [style=none] (192) at (32.5, -1) {};
		\node [style=none] (193) at (30.5, -2.75) {};
		\node [style=none] (194) at (32.5, -2.75) {};
		\node [style=port] (200) at (31.5, -3.5) {};
		\node [style=none] (201) at (31, -2.75) {};
		\node [style=none] (202) at (32, -2.75) {};
		\node [style=none] (203) at (31, -1) {};
		\node [style=none] (204) at (32, -1) {};
		\node [style=none] (205) at (31.5, -2.75) {};
		\node [style=none] (206) at (31.5, -1) {};
		\node [style=component] (207) at (31.5, -1.75) {$\Delta$};
		\node [style=none] (208) at (31, -2.25) {};
		\node [style=none] (209) at (32, -2.25) {};
		\node [style=port] (212) at (31.5, 0.5) {};
		\node [style=component] (213) at (31.5, -0.25) {$\delta$};
	\end{pgfonlayer}
	\begin{pgfonlayer}{edgelayer}
		\draw [style=wire] (191.center) to (192.center);
		\draw [style=wire] (192.center) to (194.center);
		\draw [style=wire] (194.center) to (193.center);
		\draw [style=wire] (193.center) to (191.center);
		\draw [style=wire] (207) to (208.center);
		\draw [style=wire] (207) to (209.center);
		\draw [style=wire] (201.center) to (208.center);
		\draw [style=wire] (202.center) to (209.center);
		\draw [style=wire] (207) to (206.center);
		\draw [style=wire] (205.center) to (200);
		\draw [style=wire] (212) to (213);
		\draw [style=wire] (213) to (206.center);
	\end{pgfonlayer}
\end{tikzpicture}
 \end{array} = \begin{array}[c]{c}\begin{tikzpicture}
	\begin{pgfonlayer}{nodelayer}
		\node [style=component] (196) at (34, -2.75) {$\mathsf{m}$};
		\node [style=none] (197) at (33.5, -2.25) {};
		\node [style=none] (198) at (34.5, -2.25) {};
		\node [style=component] (210) at (33.5, -1.5) {$\delta$};
		\node [style=component] (211) at (34.5, -1.5) {$\delta$};
		\node [style=component] (214) at (34, -0.25) {$\Delta$};
		\node [style=none] (215) at (33.5, -0.75) {};
		\node [style=none] (216) at (34.5, -0.75) {};
		\node [style=port] (217) at (34, 0.5) {};
		\node [style=port] (218) at (34, -3.5) {};
	\end{pgfonlayer}
	\begin{pgfonlayer}{edgelayer}
		\draw [style=wire] (196) to (197.center);
		\draw [style=wire] (196) to (198.center);
		\draw [style=wire] (211) to (198.center);
		\draw [style=wire] (210) to (197.center);
		\draw [style=wire] (214) to (215.center);
		\draw [style=wire] (214) to (216.center);
		\draw [style=wire] (215.center) to (210);
		\draw [style=wire] (216.center) to (211);
		\draw [style=wire] (217) to (214);
		\draw [style=wire] (196) to (218);
	\end{pgfonlayer}
\end{tikzpicture}
 \end{array} && \begin{array}[c]{c}\begin{tikzpicture}
	\begin{pgfonlayer}{nodelayer}
		\node [style=none] (191) at (30.75, -1) {};
		\node [style=none] (192) at (31.75, -1) {};
		\node [style=none] (193) at (30.75, -2) {};
		\node [style=none] (194) at (31.75, -2) {};
		\node [style=port] (200) at (31.25, -2.75) {};
		\node [style=none] (205) at (31.25, -2) {};
		\node [style=none] (206) at (31.25, -1) {};
		\node [style=component] (207) at (31.25, -1.5) {$\mathsf{e}$};
		\node [style=port] (212) at (31.25, 0.5) {};
		\node [style=component] (213) at (31.25, -0.25) {$\delta$};
	\end{pgfonlayer}
	\begin{pgfonlayer}{edgelayer}
		\draw [style=wire] (191.center) to (192.center);
		\draw [style=wire] (192.center) to (194.center);
		\draw [style=wire] (194.center) to (193.center);
		\draw [style=wire] (193.center) to (191.center);
		\draw [style=wire] (207) to (206.center);
		\draw [style=wire] (205.center) to (200);
		\draw [style=wire] (212) to (213);
		\draw [style=wire] (213) to (206.center);
	\end{pgfonlayer}
\end{tikzpicture}
 \end{array} = \begin{array}[c]{c}\begin{tikzpicture}
	\begin{pgfonlayer}{nodelayer}
		\node [style=component] (215) at (33, -0.5) {$\mathsf{e}$};
		\node [style=port] (216) at (33, 0.5) {};
		\node [style=component] (217) at (33, -1.75) {$\mathsf{m}$};
		\node [style=port] (218) at (33, -2.75) {};
	\end{pgfonlayer}
	\begin{pgfonlayer}{edgelayer}
		\draw [style=wire] (216) to (215);
		\draw [style=wire] (218) to (217);
	\end{pgfonlayer}
\end{tikzpicture}
 \end{array}
 \end{align}

 Many examples of (monoidal) coalgebra modalities can be found in \cite[Sec 9]{Blute2019}, \cite[Sec 2.4]{hyland2003glueing}, and \cite[Ex 4.6]{lemay2019lifting}. An important class of monoidal coalgebra modalities worth mentioning are the \emph{free exponential modalities} \cite{mellies2009explicit}, which are monoidal coalgebra modalities where $\oc(A)$ is the cofree cocommutative comonoid over $A$. There are also many other equivalent alternative ways of defining a monoidal coalgebra modality, we invite the curious reader to see \cite{schalk2004categorical} for some of them. One in particular that is very well known and worth mentioning is that in the presence of finite products (with binary product $\times$ and terminal object $\top$), a monoidal coalgebra modality is equivalently a coalgebra modality that has the \emph{Seely isomorphisms}, that is, that the canonical maps $\oc(A \times B) \to \oc(A) \otimes \oc(B)$ and $\oc(\top) \to I$ are isomorphisms \cite{bierman1995categorical}. However, products are not necessary for the story of this paper, and so we will not discuss the Seely isomorphisms further. We invite the curious reader to see \cite[Sec 7]{Blute2019} for a detailed review on the Seely isomorphisms. 
 
 We conclude this section with the observation that for a coalgebra modality, being a monoidal coalgebra modality is a property rather than extra structure. In other words, if it has a \emph{monoidal structure}, then it is unique. By a \textbf{monoidal structure} for a  coalgebra modality $(\oc, \delta, \varepsilon, \Delta, \mathsf{e})$, we mean a natural transformation $\mathsf{m}$ and a map $\mathsf{m}_I$ such that $(\oc, \delta, \varepsilon, \Delta, \mathsf{e}, \mathsf{m}, \mathsf{m}_I)$ is a monoidal coalgebra modality.
 
\begin{prop} If there exists a monoidal structure for a coalgebra modality, then it is unique.
\end{prop}
\begin{proof} Let $(\oc, \delta, \varepsilon, \Delta, \mathsf{e})$ be a coalgebra modality, and suppose that it has two possible monoidal structures, that is, we have that $(\oc, \delta, \varepsilon, \Delta, \mathsf{e}, \mathsf{m}, \mathsf{m}_I)$ and $(\oc, \delta, \varepsilon, \Delta, \mathsf{e}, \mathsf{m}^\prime, \mathsf{m}^\prime_I)$ are monoidal coalgebra modalities. Our goal is to show that $\mathsf{m} = \mathsf{m}^\prime$ and $\mathsf{m}_I = \mathsf{m}^\prime_I$. Starting with the latter, we compute that: 
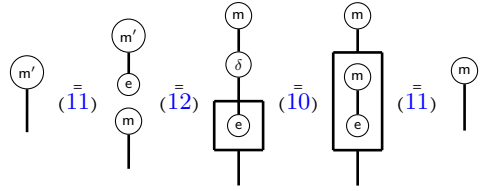
\begin{gather*}
\begin{array}[c]{c} \begin{tikzpicture}
	\begin{pgfonlayer}{nodelayer}
		\node [style=component] (64) at (45, 0.75) {$\mathsf{m}^\prime$};
		\node [style=port] (65) at (45, -0.5) {};
	\end{pgfonlayer}
	\begin{pgfonlayer}{edgelayer}
		\draw [style=wire] (64) to (65);
	\end{pgfonlayer}
\end{tikzpicture}
 \end{array} \substack{=\\(\ref{string:monoidalcomonoid})} \begin{array}[c]{c}\begin{tikzpicture}
	\begin{pgfonlayer}{nodelayer}
		\node [style=component] (11) at (69, -0.75) {$\mathsf{e}$};
		\node [style=component] (12) at (69, 0.25) {$\mathsf{m}$};
		\node [style=component] (14) at (69, -1.5) {$\mathsf{m}^\prime$};
		\node [style=port] (17) at (69, -2.5) {};
	\end{pgfonlayer}
	\begin{pgfonlayer}{edgelayer}
		\draw [style=wire] (17) to (14);
		\draw [style=wire] (12) to (11);
	\end{pgfonlayer}
\end{tikzpicture}
 \end{array} \substack{=\\(\ref{string:comonoid-!coalg})}
 \begin{array}[c]{c}
 \begin{tikzpicture}
	\begin{pgfonlayer}{nodelayer}
		\node [style=none] (45) at (73.25, -1) {};
		\node [style=none] (46) at (74.25, -1) {};
		\node [style=none] (47) at (73.25, -2) {};
		\node [style=none] (48) at (74.25, -2) {};
		\node [style=port] (49) at (73.75, -2.75) {};
		\node [style=none] (50) at (73.75, -2) {};
		\node [style=none] (51) at (73.75, -1) {};
		\node [style=component] (52) at (73.75, -1.5) {$\mathsf{e}$};
		\node [style=component] (54) at (73.75, -0.25) {$\delta$};
		\node [style=component] (55) at (73.75, 0.75) {$\mathsf{m}$};
	\end{pgfonlayer}
	\begin{pgfonlayer}{edgelayer}
		\draw [style=wire] (45.center) to (46.center);
		\draw [style=wire] (46.center) to (48.center);
		\draw [style=wire] (48.center) to (47.center);
		\draw [style=wire] (47.center) to (45.center);
		\draw [style=wire] (52) to (51.center);
		\draw [style=wire] (50.center) to (49);
		\draw [style=wire] (54) to (51.center);
		\draw [style=wire] (55) to (54);
	\end{pgfonlayer}
\end{tikzpicture}
 \end{array} \substack{=\\(\ref{string:comonad-monoidal})}  \begin{array}[c]{c} \begin{tikzpicture}
	\begin{pgfonlayer}{nodelayer}
		\node [style=component] (11) at (55.5, -5) {$\mathsf{e}$};
		\node [style=component] (12) at (55.5, -4) {$\mathsf{m}$};
		\node [style=port] (17) at (55.5, -6.25) {};
		\node [style=none] (22) at (55, -5.5) {};
		\node [style=none] (23) at (56, -5.5) {};
		\node [style=none] (24) at (55, -3.5) {};
		\node [style=none] (25) at (56, -3.5) {};
		\node [style=none] (56) at (55.5, -5.5) {};
		\node [style=none] (57) at (55.5, -3.5) {};
		\node [style=component] (58) at (55.5, -2.75) {$\mathsf{m}$};
	\end{pgfonlayer}
	\begin{pgfonlayer}{edgelayer}
		\draw [style=wire] (12) to (11);
		\draw [style=wire] (22.center) to (23.center);
		\draw [style=wire] (23.center) to (25.center);
		\draw [style=wire] (25.center) to (24.center);
		\draw [style=wire] (24.center) to (22.center);
		\draw [style=wire] (58) to (57.center);
		\draw [style=wire] (56.center) to (17);
	\end{pgfonlayer}
\end{tikzpicture}
 \end{array} \substack{=\\(\ref{string:monoidalcomonoid})} \begin{array}[c]{c} \begin{tikzpicture}
	\begin{pgfonlayer}{nodelayer}
		\node [style=component] (64) at (45, 0.75) {$\mathsf{m}$};
		\node [style=port] (65) at (45, -0.5) {};
	\end{pgfonlayer}
	\begin{pgfonlayer}{edgelayer}
		\draw [style=wire] (64) to (65);
	\end{pgfonlayer}
\end{tikzpicture}
 \end{array}
\end{gather*}
So $\mathsf{m}_I = \mathsf{m}^\prime_I$. The other desired equality requires some more work. So we compute that: 
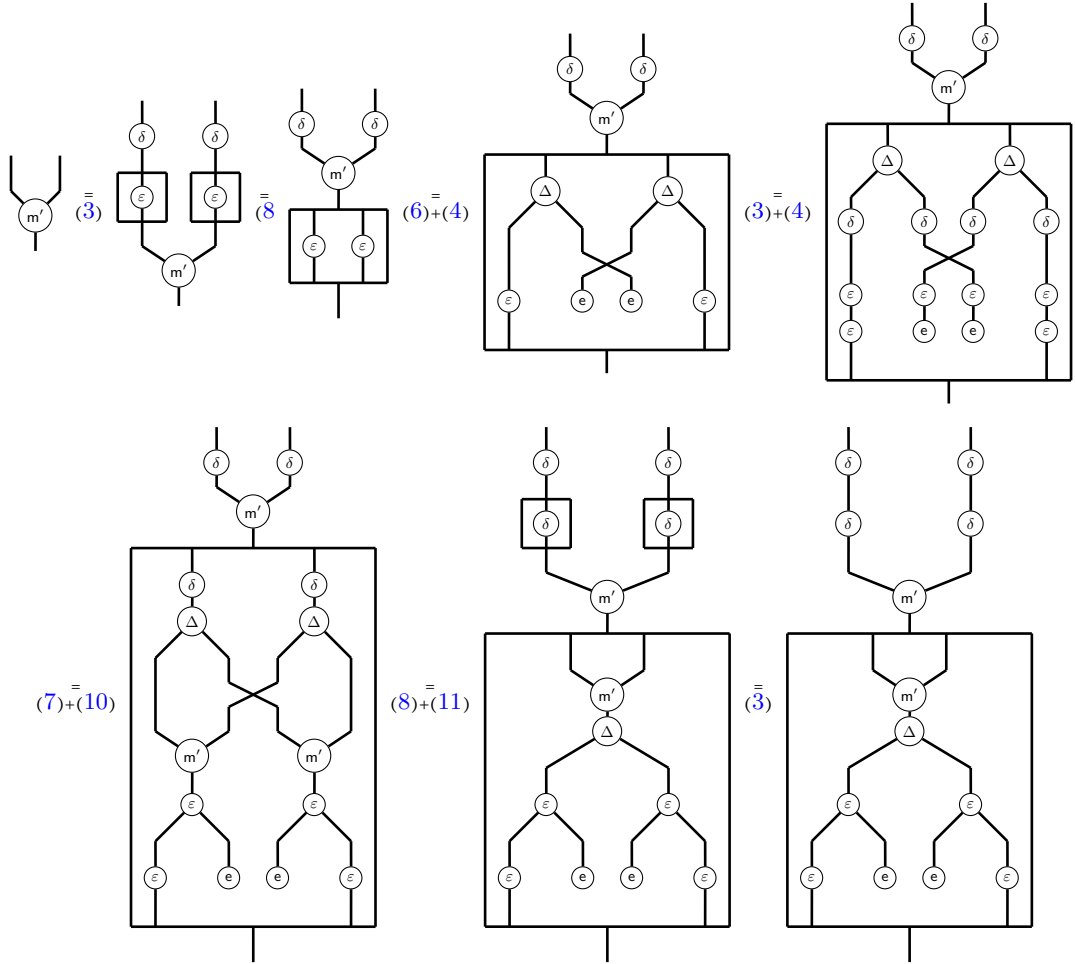
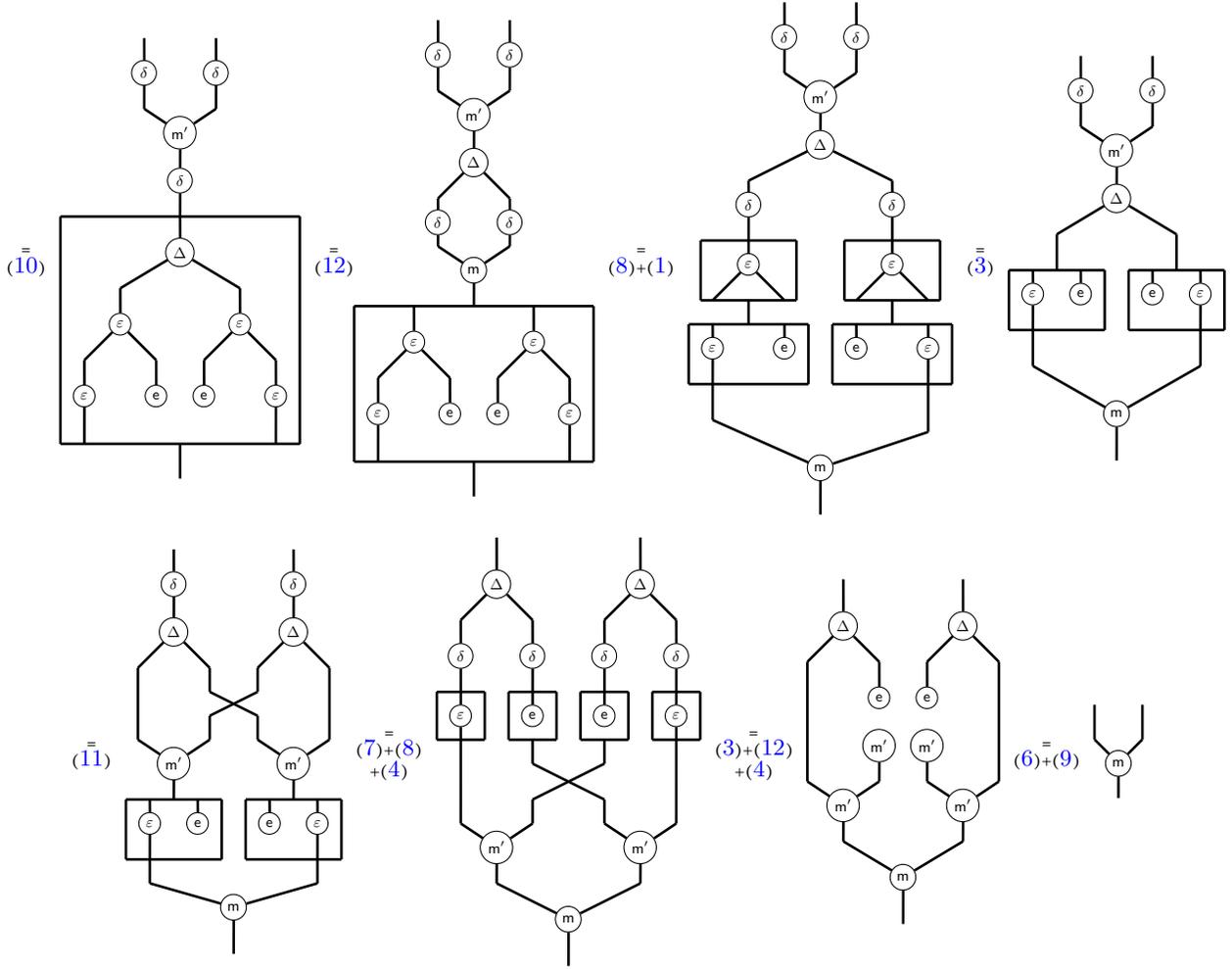
\begin{gather*}
\begin{array}[c]{c}\begin{tikzpicture}
	\begin{pgfonlayer}{nodelayer}
		\node [style=component] (59) at (73.25, 5) {$\mathsf{m}^\prime$};
		\node [style=none] (60) at (72.75, 5.5) {};
		\node [style=none] (61) at (73.75, 5.5) {};
		\node [style=port] (62) at (72.75, 6.25) {};
		\node [style=port] (63) at (73.75, 6.25) {};
		\node [style=port] (64) at (73.25, 4.25) {};
	\end{pgfonlayer}
	\begin{pgfonlayer}{edgelayer}
		\draw [style=wire] (59) to (60.center);
		\draw [style=wire] (59) to (61.center);
		\draw [style=wire] (62) to (60.center);
		\draw [style=wire] (63) to (61.center);
		\draw [style=wire] (59) to (64);
	\end{pgfonlayer}
\end{tikzpicture}
 \end{array} \substack{=\\(\ref{string:comonad})} \begin{array}[c]{c} \begin{tikzpicture}
	\begin{pgfonlayer}{nodelayer}
		\node [style=component] (59) at (78.5, -0.75) {$\mathsf{m}^\prime$};
		\node [style=none] (60) at (77.75, -0.25) {};
		\node [style=none] (61) at (79.25, -0.25) {};
		\node [style=port] (64) at (78.5, -1.5) {};
		\node [style=port] (65) at (77.75, 2.75) {};
		\node [style=component] (66) at (77.75, 2) {$\delta$};
		\node [style=component] (68) at (77.75, 0.75) {$\varepsilon$};
		\node [style=none] (69) at (77.25, 1.25) {};
		\node [style=none] (70) at (78.25, 1.25) {};
		\node [style=none] (71) at (77.25, 0.25) {};
		\node [style=none] (72) at (78.25, 0.25) {};
		\node [style=port] (73) at (79.25, 2.75) {};
		\node [style=component] (74) at (79.25, 2) {$\delta$};
		\node [style=component] (75) at (79.25, 0.75) {$\varepsilon$};
		\node [style=none] (76) at (78.75, 1.25) {};
		\node [style=none] (77) at (79.75, 1.25) {};
		\node [style=none] (78) at (78.75, 0.25) {};
		\node [style=none] (79) at (79.75, 0.25) {};
	\end{pgfonlayer}
	\begin{pgfonlayer}{edgelayer}
		\draw [style=wire] (59) to (60.center);
		\draw [style=wire] (59) to (61.center);
		\draw [style=wire] (59) to (64);
		\draw [style=wire] (66) to (65);
		\draw [style=wire] (69.center) to (70.center);
		\draw [style=wire] (70.center) to (72.center);
		\draw [style=wire] (72.center) to (71.center);
		\draw [style=wire] (71.center) to (69.center);
		\draw [style=wire] (66) to (68);
		\draw [style=wire] (74) to (73);
		\draw [style=wire] (76.center) to (77.center);
		\draw [style=wire] (77.center) to (79.center);
		\draw [style=wire] (79.center) to (78.center);
		\draw [style=wire] (78.center) to (76.center);
		\draw [style=wire] (74) to (75);
		\draw [style=wire] (68) to (60.center);
		\draw [style=wire] (75) to (61.center);
	\end{pgfonlayer}
\end{tikzpicture}
  \end{array} \substack{=\\(\ref{string:monoidal-nat})} \begin{array}[c]{c} \begin{tikzpicture}
	\begin{pgfonlayer}{nodelayer}
		\node [style=component] (59) at (78.25, -1) {$\mathsf{m}^\prime$};
		\node [style=none] (60) at (77.5, -0.5) {};
		\node [style=none] (61) at (79, -0.5) {};
		\node [style=port] (64) at (78.25, -4) {};
		\node [style=port] (65) at (77.5, 0.75) {};
		\node [style=component] (66) at (77.5, 0) {$\delta$};
		\node [style=component] (68) at (77.75, -2.5) {$\varepsilon$};
		\node [style=none] (69) at (77.25, -1.75) {};
		\node [style=none] (71) at (77.25, -3.25) {};
		\node [style=port] (73) at (79, 0.75) {};
		\node [style=component] (74) at (79, 0) {$\delta$};
		\node [style=component] (75) at (78.75, -2.5) {$\varepsilon$};
		\node [style=none] (77) at (79.25, -1.75) {};
		\node [style=none] (79) at (79.25, -3.25) {};
		\node [style=none] (80) at (77.75, -1.75) {};
		\node [style=none] (81) at (78.75, -1.75) {};
		\node [style=none] (82) at (77.75, -3.25) {};
		\node [style=none] (83) at (78.75, -3.25) {};
		\node [style=none] (84) at (78.25, -1.75) {};
		\node [style=none] (85) at (78.25, -3.25) {};
	\end{pgfonlayer}
	\begin{pgfonlayer}{edgelayer}
		\draw [style=wire] (59) to (60.center);
		\draw [style=wire] (59) to (61.center);
		\draw [style=wire] (66) to (65);
		\draw [style=wire] (71.center) to (69.center);
		\draw [style=wire] (74) to (73);
		\draw [style=wire] (77.center) to (79.center);
		\draw [style=wire] (71.center) to (79.center);
		\draw [style=wire] (69.center) to (77.center);
		\draw [style=wire] (80.center) to (68);
		\draw [style=wire] (68) to (82.center);
		\draw [style=wire] (81.center) to (75);
		\draw [style=wire] (75) to (83.center);
		\draw [style=wire] (59) to (84.center);
		\draw [style=wire] (85.center) to (64);
		\draw [style=wire] (66) to (60.center);
		\draw [style=wire] (74) to (61.center);
	\end{pgfonlayer}
\end{tikzpicture}
  \end{array} \substack{=\\(\ref{string:comonoid})+(\ref{strings:monoidal-functor})} \begin{array}[c]{c}\begin{tikzpicture}
	\begin{pgfonlayer}{nodelayer}
		\node [style=component] (59) at (73.25, 5) {$\mathsf{m}$};
		\node [style=none] (60) at (72.75, 5.5) {};
		\node [style=none] (61) at (73.75, 5.5) {};
		\node [style=port] (62) at (72.75, 6.25) {};
		\node [style=port] (63) at (73.75, 6.25) {};
		\node [style=port] (64) at (73.25, 4.25) {};
	\end{pgfonlayer}
	\begin{pgfonlayer}{edgelayer}
		\draw [style=wire] (59) to (60.center);
		\draw [style=wire] (59) to (61.center);
		\draw [style=wire] (62) to (60.center);
		\draw [style=wire] (63) to (61.center);
		\draw [style=wire] (59) to (64);
	\end{pgfonlayer}
\end{tikzpicture}
 \end{array}
\end{gather*}
So $\mathsf{m} = \mathsf{m}^\prime$. Therefore, we conclude that monoidal structure for a coalgebra modality is unique. 
\end{proof}

 \section{Pre-Coderelictions}\label{sec:pre-coder}

As explained in the introduction, for the main result of this paper (Thm \ref{thm:additive}) we do not require the full version of a codereliction. As such, in this section we introduce the notion of a \emph{pre-codereliction}, which, as the name suggests, is a slightly weaker version of a codereliction. 

For a monoidal coalgebra modality $(\oc, \delta, \varepsilon, \Delta, \mathsf{e}, \mathsf{m}, \mathsf{m}_I)$ on a symmetric monoidal category, a \textbf{pre-codereliction} is a natural transformation ${\eta_A: A \to \oc(A)}$ such that diagrams in Appendix \ref{sec:precoder} commute. In string diagrams, a pre-codereliction $\eta$ is drawn as follows: 
\begin{align*} \begin{array}[c]{c} \begin{tikzpicture}
	\begin{pgfonlayer}{nodelayer}
		\node [style=object] (2) at (1, 2) {$A$};
		\node [style=object] (3) at (1, 0) {$\oc(A)$};
		\node [style=component] (4) at (1, 1) {$\eta$};
	\end{pgfonlayer}
	\begin{pgfonlayer}{edgelayer}
		\draw [style=wire] (2) to (4);
		\draw [style=wire] (4) to (3);
	\end{pgfonlayer}
\end{tikzpicture}
 \end{array} 
 \end{align*}
with its naturality drawn as follows: 

\begin{align}\label{string:coder-nat}
   \begin{array}[c]{c} \begin{tikzpicture}
	\begin{pgfonlayer}{nodelayer}
		\node [style=port] (416) at (83, 0.5) {};
		\node [style=component] (417) at (83, 1.5) {$\eta$};
		\node [style=port] (418) at (83, 3.75) {};
		\node [style=component] (419) at (83, 2.75) {$f$};
	\end{pgfonlayer}
	\begin{pgfonlayer}{edgelayer}
		\draw [style=wire] (417) to (416);
		\draw [style=wire] (418) to (419);
		\draw [style=wire] (417) to (419);
	\end{pgfonlayer}
\end{tikzpicture}
\end{array} =   \begin{array}[c]{c} \begin{tikzpicture}
	\begin{pgfonlayer}{nodelayer}
		\node [style=port] (408) at (81.5, 3.75) {};
		\node [style=component] (409) at (81.5, 3) {$\eta$};
		\node [style=port] (410) at (81.5, 0.5) {};
		\node [style=component] (411) at (81.5, 1.75) {$f$};
		\node [style=none] (412) at (81, 1.25) {};
		\node [style=none] (413) at (82, 1.25) {};
		\node [style=none] (414) at (81, 2.25) {};
		\node [style=none] (415) at (82, 2.25) {};
	\end{pgfonlayer}
	\begin{pgfonlayer}{edgelayer}
		\draw [style=wire] (409) to (408);
		\draw [style=wire] (410) to (411);
		\draw [style=wire] (412.center) to (413.center);
		\draw [style=wire] (413.center) to (415.center);
		\draw [style=wire] (415.center) to (414.center);
		\draw [style=wire] (414.center) to (412.center);
		\draw [style=wire] (411) to (409);
	\end{pgfonlayer}
\end{tikzpicture}
\end{array}
\end{align}
The axioms for a pre-codereliction (\ref{diag:precoder}) are drawn as: 
\begin{align}\label{string:coder}
 \begin{array}[c]{c}\begin{tikzpicture}
	\begin{pgfonlayer}{nodelayer}
		\node [style=port] (16) at (-2.25, 0.5) {};
		\node [style=component] (17) at (-2.25, 1.5) {$\varepsilon$};
		\node [style=port] (18) at (-2.25, 3.75) {};
		\node [style=component] (19) at (-2.25, 2.75) {$\eta$};
	\end{pgfonlayer}
	\begin{pgfonlayer}{edgelayer}
		\draw [style=wire] (17) to (16);
		\draw [style=wire] (18) to (19);
		\draw [style=wire] (19) to (17);
	\end{pgfonlayer}
\end{tikzpicture}
\end{array} \substack{=\\\text{\textbf{[cd.3]}}}  \begin{array}[c]{c}\begin{tikzpicture}
	\begin{pgfonlayer}{nodelayer}
		\node [style=port] (16) at (-2.25, 0.5) {};
		\node [style=port] (18) at (-2.25, 3.75) {};
	\end{pgfonlayer}
	\begin{pgfonlayer}{edgelayer}
		\draw [style=wire] (18) to (16);
	\end{pgfonlayer}
\end{tikzpicture}
\end{array} && \begin{array}[c]{c}\begin{tikzpicture}
	\begin{pgfonlayer}{nodelayer}
		\node [style=port] (420) at (84.5, -7) {};
		\node [style=port] (421) at (85, -4.25) {};
		\node [style=component] (422) at (84.5, -6) {$\mathsf{m}$};
		\node [style=none] (423) at (85, -5.5) {};
		\node [style=none] (424) at (84, -5.5) {};
		\node [style=component] (425) at (84, -5) {$\eta$};
		\node [style=port] (426) at (84, -4.25) {};
	\end{pgfonlayer}
	\begin{pgfonlayer}{edgelayer}
		\draw [style=wire] (420) to (422);
		\draw [style=wire] (422) to (423.center);
		\draw [style=wire] (422) to (424.center);
		\draw [style=wire] (423.center) to (421);
		\draw [style=wire] (424.center) to (425);
		\draw [style=wire] (426) to (425);
	\end{pgfonlayer}
\end{tikzpicture}
 \end{array} \substack{=\\\text{\textbf{[cd.m.l]}}} \begin{array}[c]{c}\begin{tikzpicture}
	\begin{pgfonlayer}{nodelayer}
		\node [style=port] (420) at (84.5, -7) {};
		\node [style=port] (421) at (84, -4.25) {};
		\node [style=component] (422) at (84.5, -6) {$\eta$};
		\node [style=none] (423) at (84, -5.5) {};
		\node [style=none] (424) at (85, -5.5) {};
		\node [style=component] (425) at (85, -5) {$\varepsilon$};
		\node [style=port] (426) at (85, -4.25) {};
	\end{pgfonlayer}
	\begin{pgfonlayer}{edgelayer}
		\draw [style=wire] (420) to (422);
		\draw [style=wire] (422) to (423.center);
		\draw [style=wire] (422) to (424.center);
		\draw [style=wire] (423.center) to (421);
		\draw [style=wire] (424.center) to (425);
		\draw [style=wire] (426) to (425);
	\end{pgfonlayer}
\end{tikzpicture} \end{array} && \begin{array}[c]{c}\begin{tikzpicture}
	\begin{pgfonlayer}{nodelayer}
		\node [style=port] (420) at (84.5, -7) {};
		\node [style=port] (421) at (84, -4.25) {};
		\node [style=component] (422) at (84.5, -6) {$\mathsf{m}$};
		\node [style=none] (423) at (84, -5.5) {};
		\node [style=none] (424) at (85, -5.5) {};
		\node [style=component] (425) at (85, -5) {$\eta$};
		\node [style=port] (426) at (85, -4.25) {};
	\end{pgfonlayer}
	\begin{pgfonlayer}{edgelayer}
		\draw [style=wire] (420) to (422);
		\draw [style=wire] (422) to (423.center);
		\draw [style=wire] (422) to (424.center);
		\draw [style=wire] (423.center) to (421);
		\draw [style=wire] (424.center) to (425);
		\draw [style=wire] (426) to (425);
	\end{pgfonlayer}
\end{tikzpicture}
 \end{array} \substack{=\\\text{\textbf{[cd.m.r]}}} \begin{array}[c]{c}\begin{tikzpicture}
	\begin{pgfonlayer}{nodelayer}
		\node [style=port] (420) at (84.5, -7) {};
		\node [style=port] (421) at (85, -4.25) {};
		\node [style=component] (422) at (84.5, -6) {$\eta$};
		\node [style=none] (423) at (85, -5.5) {};
		\node [style=none] (424) at (84, -5.5) {};
		\node [style=component] (425) at (84, -5) {$\varepsilon$};
		\node [style=port] (426) at (84, -4.25) {};
	\end{pgfonlayer}
	\begin{pgfonlayer}{edgelayer}
		\draw [style=wire] (420) to (422);
		\draw [style=wire] (422) to (423.center);
		\draw [style=wire] (422) to (424.center);
		\draw [style=wire] (423.center) to (421);
		\draw [style=wire] (424.center) to (425);
		\draw [style=wire] (426) to (425);
	\end{pgfonlayer}
\end{tikzpicture}\end{array} 
\end{align}
Using the same terminology from \cite{Blute2019}, \textbf{[cd.3]} is called the linear rule and \textbf{[cd.m.l]} (resp. \textbf{[cd.m.r]}) is called the left (resp. right) monoidal rule. There is a bit of redundancy since \textbf{[cd.m.l]} and \textbf{[cd.m.r]} are equivalent. 

\begin{lem}\label{lemma:precoder-left=right} For a monoidal coalgebra modality $(\oc, \delta, \varepsilon, \Delta, \mathsf{e}, \mathsf{m}, \mathsf{m}_I)$, a natural transformation ${\eta_A: A \to \oc(A)}$ satisfies \textbf{[cd.m.l]} if and only if $\eta$ satisfies \textbf{[cd.m.r]}. 
\end{lem}
\begin{proof} Suppose that $\eta$ satisfies \textbf{[cd.m.l]}. Then we compute that: 
\begin{gather*}
\begin{array}[c]{c}\begin{tikzpicture}
	\begin{pgfonlayer}{nodelayer}
		\node [style=port] (420) at (84.5, -7) {};
		\node [style=port] (421) at (84, -4.25) {};
		\node [style=component] (422) at (84.5, -6) {$\mathsf{m}$};
		\node [style=none] (423) at (84, -5.5) {};
		\node [style=none] (424) at (85, -5.5) {};
		\node [style=component] (425) at (85, -5) {$\eta$};
		\node [style=port] (426) at (85, -4.25) {};
	\end{pgfonlayer}
	\begin{pgfonlayer}{edgelayer}
		\draw [style=wire] (420) to (422);
		\draw [style=wire] (422) to (423.center);
		\draw [style=wire] (422) to (424.center);
		\draw [style=wire] (423.center) to (421);
		\draw [style=wire] (424.center) to (425);
		\draw [style=wire] (426) to (425);
	\end{pgfonlayer}
\end{tikzpicture}
 \end{array} \substack{=\\\text{sym.}}  \begin{array}[c]{c} \begin{tikzpicture}
	\begin{pgfonlayer}{nodelayer}
		\node [style=port] (0) at (77, -9.5) {};
		\node [style=port] (1) at (76.5, -3.75) {};
		\node [style=component] (2) at (77, -8.5) {$\mathsf{m}$};
		\node [style=none] (3) at (76.5, -8) {};
		\node [style=none] (4) at (77.5, -8) {};
		\node [style=component] (5) at (77.5, -4.5) {$\eta$};
		\node [style=port] (6) at (77.5, -3.75) {};
		\node [style=none] (14) at (76.5, -7.5) {};
		\node [style=none] (15) at (77.5, -7.5) {};
		\node [style=none] (16) at (77.5, -5) {};
		\node [style=none] (17) at (76.5, -5) {};
		\node [style=none] (18) at (76.5, -6.75) {};
		\node [style=none] (19) at (77.5, -6.75) {};
		\node [style=none] (20) at (76.5, -6.25) {};
		\node [style=none] (21) at (77.5, -6.25) {};
		\node [style=none] (22) at (76.5, -5.5) {};
		\node [style=none] (23) at (77.5, -5.5) {};
	\end{pgfonlayer}
	\begin{pgfonlayer}{edgelayer}
		\draw [style=wire] (0) to (2);
		\draw [style=wire] (2) to (3.center);
		\draw [style=wire] (2) to (4.center);
		\draw [style=wire] (6) to (5);
		\draw [style=wire] (3.center) to (14.center);
		\draw [style=wire] (15.center) to (4.center);
		\draw [style=wire] (14.center) to (19.center);
		\draw [style=wire] (18.center) to (15.center);
		\draw [style=wire] (20.center) to (18.center);
		\draw [style=wire] (19.center) to (21.center);
		\draw [style=wire] (22.center) to (21.center);
		\draw [style=wire] (20.center) to (23.center);
		\draw [style=wire] (23.center) to (16.center);
		\draw [style=wire] (17.center) to (22.center);
		\draw [style=wire] (17.center) to (1);
		\draw [style=wire] (5) to (16.center);
	\end{pgfonlayer}
\end{tikzpicture}
 \end{array} \substack{=\\ (\ref{string:coder-nat})}  \begin{array}[c]{c}\begin{tikzpicture}
	\begin{pgfonlayer}{nodelayer}
		\node [style=port] (44) at (72.5, -10.5) {};
		\node [style=port] (45) at (72, -4.25) {};
		\node [style=component] (46) at (72.5, -9.5) {$\eta$};
		\node [style=none] (47) at (72, -9) {};
		\node [style=none] (48) at (73, -9) {};
		\node [style=component] (49) at (73, -7) {$\varepsilon$};
		\node [style=port] (50) at (73, -4.25) {};
		\node [style=none] (51) at (72, -8.5) {};
		\node [style=none] (52) at (73, -8.5) {};
		\node [style=none] (53) at (73, -5) {};
		\node [style=none] (54) at (72, -5) {};
		\node [style=none] (55) at (72, -7.75) {};
		\node [style=none] (56) at (73, -7.75) {};
		\node [style=none] (57) at (72, -6.25) {};
		\node [style=none] (58) at (73, -6.25) {};
		\node [style=none] (59) at (72, -5.5) {};
		\node [style=none] (60) at (73, -5.5) {};
	\end{pgfonlayer}
	\begin{pgfonlayer}{edgelayer}
		\draw [style=wire] (44) to (46);
		\draw [style=wire] (46) to (47.center);
		\draw [style=wire] (46) to (48.center);
		\draw [style=wire] (47.center) to (51.center);
		\draw [style=wire] (52.center) to (48.center);
		\draw [style=wire] (51.center) to (56.center);
		\draw [style=wire] (55.center) to (52.center);
		\draw [style=wire] (57.center) to (55.center);
		\draw [style=wire] (59.center) to (58.center);
		\draw [style=wire] (57.center) to (60.center);
		\draw [style=wire] (60.center) to (53.center);
		\draw [style=wire] (54.center) to (59.center);
		\draw [style=wire] (54.center) to (45);
		\draw [style=wire] (58.center) to (49);
		\draw [style=wire] (49) to (56.center);
		\draw [style=wire] (53.center) to (50);
	\end{pgfonlayer}
\end{tikzpicture}
 \end{array} \substack{=\\(\ref{string:nat-sym})}  \begin{array}[c]{c}\begin{tikzpicture}
	\begin{pgfonlayer}{nodelayer}
		\node [style=port] (0) at (77, -9.5) {};
		\node [style=port] (1) at (77.5, -3.75) {};
		\node [style=component] (2) at (77, -8.5) {$\eta$};
		\node [style=none] (3) at (77.5, -8) {};
		\node [style=none] (4) at (76.5, -8) {};
		\node [style=component] (5) at (76.5, -4.5) {$\varepsilon$};
		\node [style=port] (6) at (76.5, -3.75) {};
		\node [style=none] (14) at (77.5, -7.5) {};
		\node [style=none] (15) at (76.5, -7.5) {};
		\node [style=none] (16) at (76.5, -5) {};
		\node [style=none] (17) at (77.5, -5) {};
		\node [style=none] (18) at (77.5, -6.75) {};
		\node [style=none] (19) at (76.5, -6.75) {};
		\node [style=none] (20) at (77.5, -6.25) {};
		\node [style=none] (21) at (76.5, -6.25) {};
		\node [style=none] (22) at (77.5, -5.5) {};
		\node [style=none] (23) at (76.5, -5.5) {};
	\end{pgfonlayer}
	\begin{pgfonlayer}{edgelayer}
		\draw [style=wire] (0) to (2);
		\draw [style=wire] (2) to (3.center);
		\draw [style=wire] (2) to (4.center);
		\draw [style=wire] (6) to (5);
		\draw [style=wire] (3.center) to (14.center);
		\draw [style=wire] (15.center) to (4.center);
		\draw [style=wire] (14.center) to (19.center);
		\draw [style=wire] (18.center) to (15.center);
		\draw [style=wire] (20.center) to (18.center);
		\draw [style=wire] (19.center) to (21.center);
		\draw [style=wire] (22.center) to (21.center);
		\draw [style=wire] (20.center) to (23.center);
		\draw [style=wire] (23.center) to (16.center);
		\draw [style=wire] (17.center) to (22.center);
		\draw [style=wire] (17.center) to (1);
		\draw [style=wire] (5) to (16.center);
	\end{pgfonlayer}
\end{tikzpicture}
 \end{array} \substack{=\\\text{sym.}} \begin{array}[c]{c}\begin{tikzpicture}
	\begin{pgfonlayer}{nodelayer}
		\node [style=port] (420) at (84.5, -7) {};
		\node [style=port] (421) at (85, -4.25) {};
		\node [style=component] (422) at (84.5, -6) {$\eta$};
		\node [style=none] (423) at (85, -5.5) {};
		\node [style=none] (424) at (84, -5.5) {};
		\node [style=component] (425) at (84, -5) {$\varepsilon$};
		\node [style=port] (426) at (84, -4.25) {};
	\end{pgfonlayer}
	\begin{pgfonlayer}{edgelayer}
		\draw [style=wire] (420) to (422);
		\draw [style=wire] (422) to (423.center);
		\draw [style=wire] (422) to (424.center);
		\draw [style=wire] (423.center) to (421);
		\draw [style=wire] (424.center) to (425);
		\draw [style=wire] (426) to (425);
	\end{pgfonlayer}
\end{tikzpicture}\end{array}
\end{gather*}
So $\eta$ satisfies \textbf{[cd.m.r]}. The converse is shown using similar calculations. 
\end{proof}

Here are two useful identities for a pre-codereliction. 

\begin{lem}\label{lemma:precoder} For a monoidal coalgebra modality $(\oc, \delta, \varepsilon, \Delta, \mathsf{e}, \mathsf{m}, \mathsf{m}_I)$ with a pre-codereliction $\eta$, the diagrams in Appendix \ref{sec:lemmaprecoder} commute, that is, the following equalities hold: 
\begin{align}\label{string:eta-lemma} 
\begin{array}[c]{c}\begin{tikzpicture}
	\begin{pgfonlayer}{nodelayer}
		\node [style=port] (420) at (84.5, -7) {};
		\node [style=port] (421) at (85, -4.25) {};
		\node [style=component] (422) at (84.5, -6) {$\mathsf{m}$};
		\node [style=none] (423) at (85, -5.5) {};
		\node [style=none] (424) at (84, -5.5) {};
		\node [style=component] (425) at (84, -5) {$\eta$};
	\end{pgfonlayer}
	\begin{pgfonlayer}{edgelayer}
		\draw [style=wire] (420) to (422);
		\draw [style=wire] (422) to (423.center);
		\draw [style=wire] (422) to (424.center);
		\draw [style=wire] (423.center) to (421);
		\draw [style=wire] (424.center) to (425);
	\end{pgfonlayer}
\end{tikzpicture}
\end{array} = 
 \begin{array}[c]{c}\begin{tikzpicture}
	\begin{pgfonlayer}{nodelayer}
		\node [style=port] (16) at (-2.25, 0.5) {};
		\node [style=component] (17) at (-2.25, 1.5) {$\eta$};
		\node [style=port] (18) at (-2.25, 3.75) {};
		\node [style=component] (19) at (-2.25, 2.75) {$\varepsilon$};
	\end{pgfonlayer}
	\begin{pgfonlayer}{edgelayer}
		\draw [style=wire] (17) to (16);
		\draw [style=wire] (18) to (19);
		\draw [style=wire] (19) to (17);
	\end{pgfonlayer}
\end{tikzpicture}
\end{array} =  \begin{array}[c]{c}\begin{tikzpicture}
	\begin{pgfonlayer}{nodelayer}
		\node [style=port] (420) at (84.5, -7) {};
		\node [style=port] (421) at (84, -4.25) {};
		\node [style=component] (422) at (84.5, -6) {$\mathsf{m}$};
		\node [style=none] (423) at (84, -5.5) {};
		\node [style=none] (424) at (85, -5.5) {};
		\node [style=component] (425) at (85, -5) {$\eta$};
	\end{pgfonlayer}
	\begin{pgfonlayer}{edgelayer}
		\draw [style=wire] (420) to (422);
		\draw [style=wire] (422) to (423.center);
		\draw [style=wire] (422) to (424.center);
		\draw [style=wire] (423.center) to (421);
		\draw [style=wire] (424.center) to (425);
	\end{pgfonlayer}
\end{tikzpicture}
\end{array} && \begin{array}[c]{c}\begin{tikzpicture}
	\begin{pgfonlayer}{nodelayer}
		\node [style=component] (1251) at (213.25, -1.75) {$\eta$};
		\node [style=port] (1252) at (213.25, -1) {};
		\node [style=component] (1253) at (214.25, -1.75) {$\eta$};
		\node [style=port] (1254) at (214.25, -1) {};
		\node [style=component] (1255) at (213.75, -2.75) {$\mathsf{m}$};
		\node [style=none] (1256) at (213.25, -2.25) {};
		\node [style=none] (1257) at (214.25, -2.25) {};
		\node [style=port] (1258) at (213.75, -3.5) {};
	\end{pgfonlayer}
	\begin{pgfonlayer}{edgelayer}
		\draw [style=wire] (1252) to (1251);
		\draw [style=wire] (1254) to (1253);
		\draw [style=wire] (1255) to (1256.center);
		\draw [style=wire] (1255) to (1257.center);
		\draw [style=wire] (1251) to (1256.center);
		\draw [style=wire] (1253) to (1257.center);
		\draw [style=wire] (1255) to (1258);
	\end{pgfonlayer}
\end{tikzpicture}
\end{array} = 
 \begin{array}[c]{c}\begin{tikzpicture}
	\begin{pgfonlayer}{nodelayer}
		\node [style=port] (1260) at (215.25, -1) {};
		\node [style=port] (1262) at (216.25, -1) {};
		\node [style=component] (1263) at (215.75, -2.75) {$\eta$};
		\node [style=none] (1264) at (215.25, -2.25) {};
		\node [style=none] (1265) at (216.25, -2.25) {};
		\node [style=port] (1266) at (215.75, -3.5) {};
	\end{pgfonlayer}
	\begin{pgfonlayer}{edgelayer}
		\draw [style=wire] (1263) to (1264.center);
		\draw [style=wire] (1263) to (1265.center);
		\draw [style=wire] (1263) to (1266);
		\draw [style=wire] (1260) to (1264.center);
		\draw [style=wire] (1262) to (1265.center);
	\end{pgfonlayer}
\end{tikzpicture}
\end{array} 
\end{align}
\end{lem}
\begin{proof} Observe that the left most equalities is simply \textbf{[cd.m.l]} (resp. \textbf{[cd.m.r]}) in the special case that $A=I$ (resp. $B=I$). On the other hand for the equality on the right, we easily compute that: 
\begin{gather*}
\begin{array}[c]{c}\begin{tikzpicture}
	\begin{pgfonlayer}{nodelayer}
		\node [style=component] (1251) at (213.25, -1.75) {$\eta$};
		\node [style=port] (1252) at (213.25, -1) {};
		\node [style=component] (1253) at (214.25, -1.75) {$\eta$};
		\node [style=port] (1254) at (214.25, -1) {};
		\node [style=component] (1255) at (213.75, -2.75) {$\mathsf{m}$};
		\node [style=none] (1256) at (213.25, -2.25) {};
		\node [style=none] (1257) at (214.25, -2.25) {};
		\node [style=port] (1258) at (213.75, -3.5) {};
	\end{pgfonlayer}
	\begin{pgfonlayer}{edgelayer}
		\draw [style=wire] (1252) to (1251);
		\draw [style=wire] (1254) to (1253);
		\draw [style=wire] (1255) to (1256.center);
		\draw [style=wire] (1255) to (1257.center);
		\draw [style=wire] (1251) to (1256.center);
		\draw [style=wire] (1253) to (1257.center);
		\draw [style=wire] (1255) to (1258);
	\end{pgfonlayer}
\end{tikzpicture}
\end{array} \substack{=\\(\ref{string:coder})} \begin{array}[c]{c}\begin{tikzpicture}
	\begin{pgfonlayer}{nodelayer}
		\node [style=component] (1604) at (269, -1.75) {$\eta$};
		\node [style=port] (1605) at (269, -1) {};
		\node [style=port] (1607) at (270, -1) {};
		\node [style=component] (1608) at (269.5, -4) {$\eta$};
		\node [style=none] (1609) at (269, -3.25) {};
		\node [style=none] (1610) at (270, -3.25) {};
		\node [style=port] (1611) at (269.5, -4.75) {};
		\node [style=component] (1612) at (269, -2.5) {$\varepsilon$};
	\end{pgfonlayer}
	\begin{pgfonlayer}{edgelayer}
		\draw [style=wire] (1605) to (1604);
		\draw [style=wire] (1608) to (1609.center);
		\draw [style=wire] (1608) to (1610.center);
		\draw [style=wire] (1608) to (1611);
		\draw [style=wire] (1604) to (1612);
		\draw [style=wire] (1612) to (1609.center);
		\draw [style=wire] (1607) to (1610.center);
	\end{pgfonlayer}
\end{tikzpicture}
 \end{array}  \substack{=\\(\ref{string:coder})}  \begin{array}[c]{c}\begin{tikzpicture}
	\begin{pgfonlayer}{nodelayer}
		\node [style=port] (1260) at (215.25, -1) {};
		\node [style=port] (1262) at (216.25, -1) {};
		\node [style=component] (1263) at (215.75, -2.75) {$\eta$};
		\node [style=none] (1264) at (215.25, -2.25) {};
		\node [style=none] (1265) at (216.25, -2.25) {};
		\node [style=port] (1266) at (215.75, -3.5) {};
	\end{pgfonlayer}
	\begin{pgfonlayer}{edgelayer}
		\draw [style=wire] (1263) to (1264.center);
		\draw [style=wire] (1263) to (1265.center);
		\draw [style=wire] (1263) to (1266);
		\draw [style=wire] (1260) to (1264.center);
		\draw [style=wire] (1262) to (1265.center);
	\end{pgfonlayer}
\end{tikzpicture}
\end{array}
\end{gather*}
So the desired equality holds. 
\end{proof}

It is important to note that the right identity of (\ref{string:eta-lemma}) does not say that a pre-codereliction is monoidal. For a pre-codereliction to be monoidal, we would also need that $\eta_I$ be equal to $\mathsf{m}_I$, however in general this is not necessarily true. On the other hand, we can use the left identity of (\ref{string:eta-lemma}) to show that pre-coderelictions are in fact unique!

\begin{prop}\label{prop:precoderunique} For a monoidal coalgebra modality, if a pre-codereliction exists, then it is unique. 
\end{prop}
\begin{proof} Suppose that we have two pre-coderelictions $\eta$ and $\eta^\prime$. Then we compute: 
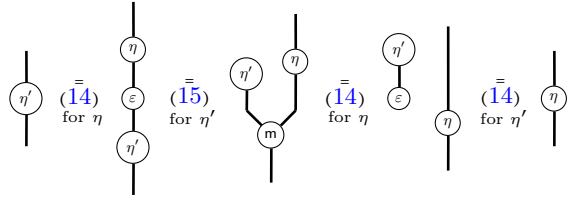
\begin{gather*}
\begin{array}[c]{c}\begin{tikzpicture}
	\begin{pgfonlayer}{nodelayer}
		\node [style=port] (1591) at (263, 0) {};
		\node [style=port] (1592) at (263, -2) {};
		\node [style=component] (1593) at (263, -1) {$\eta^\prime$};
	\end{pgfonlayer}
	\begin{pgfonlayer}{edgelayer}
		\draw [style=wire] (1591) to (1593);
		\draw [style=wire] (1593) to (1592);
	\end{pgfonlayer}
\end{tikzpicture}
 \end{array} \substack{=\\(\ref{string:coder})\\\text{ for }\eta} \begin{array}[c]{c}\begin{tikzpicture}
	\begin{pgfonlayer}{nodelayer}
		\node [style=port] (1594) at (265, 0) {};
		\node [style=port] (1595) at (265, -4) {};
		\node [style=component] (1596) at (265, -3) {$\eta^\prime$};
		\node [style=component] (1597) at (265, -1) {$\eta$};
		\node [style=component] (1598) at (265, -2) {$\varepsilon$};
	\end{pgfonlayer}
	\begin{pgfonlayer}{edgelayer}
		\draw [style=wire] (1596) to (1595);
		\draw [style=wire] (1594) to (1597);
		\draw [style=wire] (1597) to (1598);
		\draw [style=wire] (1598) to (1596);
	\end{pgfonlayer}
\end{tikzpicture}
 \end{array}  \substack{=\\(\ref{string:eta-lemma})\\\text{ for }\eta^\prime} \begin{array}[c]{c}\begin{tikzpicture}
	\begin{pgfonlayer}{nodelayer}
		\node [style=port] (1594) at (266, 0) {};
		\node [style=port] (1595) at (265.5, -3.5) {};
		\node [style=component] (1597) at (266, -1) {$\eta$};
		\node [style=component] (1599) at (265.5, -2.5) {$\mathsf{m}$};
		\node [style=none] (1600) at (266, -2) {};
		\node [style=none] (1601) at (265, -2) {};
		\node [style=component] (1602) at (265, -1.25) {$\eta^\prime$};
	\end{pgfonlayer}
	\begin{pgfonlayer}{edgelayer}
		\draw [style=wire] (1594) to (1597);
		\draw [style=wire] (1599) to (1600.center);
		\draw [style=wire] (1599) to (1601.center);
		\draw [style=wire] (1601.center) to (1602);
		\draw [style=wire] (1597) to (1600.center);
		\draw [style=wire] (1599) to (1595);
	\end{pgfonlayer}
\end{tikzpicture}
 \end{array}  \substack{=\\(\ref{string:coder})\\\text{ for }\eta} \begin{array}[c]{c}\begin{tikzpicture}
	\begin{pgfonlayer}{nodelayer}
		\node [style=port] (1594) at (266, 0) {};
		\node [style=port] (1595) at (266, -3) {};
		\node [style=component] (1597) at (266, -2) {$\eta$};
		\node [style=component] (1601) at (265, -1.5) {$\varepsilon$};
		\node [style=component] (1602) at (265, -0.5) {$\eta^\prime$};
	\end{pgfonlayer}
	\begin{pgfonlayer}{edgelayer}
		\draw [style=wire] (1594) to (1597);
		\draw [style=wire] (1601) to (1602);
		\draw [style=wire] (1597) to (1595);
	\end{pgfonlayer}
\end{tikzpicture}
 \end{array} \substack{=\\(\ref{string:coder})\\\text{ for }\eta^\prime} \begin{array}[c]{c}\begin{tikzpicture}
	\begin{pgfonlayer}{nodelayer}
		\node [style=port] (1591) at (263, 0) {};
		\node [style=port] (1592) at (263, -2) {};
		\node [style=component] (1593) at (263, -1) {$\eta$};
	\end{pgfonlayer}
	\begin{pgfonlayer}{edgelayer}
		\draw [style=wire] (1591) to (1593);
		\draw [style=wire] (1593) to (1592);
	\end{pgfonlayer}
\end{tikzpicture}
 \end{array} 
\end{gather*}
So $\eta = \eta^\prime$, and therefore we conclude that a pre-codereliction must be unique. 
\end{proof}

\section{Monoidal Bialgebra Modalities}\label{sec:mon-bialg-mod}

As explained in the introduction, if one wishes to define a codereliction in a non-additive setting, one still requires the (bi)monoid structural maps to express the chain rule. Thus in this section we introduce the concept of a \emph{monoidal bialgebra modality}, which is a monoidal coalgebra modality where $\oc(A)$ also has a natural (bi)monoid structure which is compatible with the symmetric monoidal comonad structure. The axioms for a monoidal bialgebra modality are based on the derivable identities one obtains for the canonical (bi)monoid structure of a monoidal coalgebra modality on an additive symmetric monoidal category, specifically those of \cite[Thm 3.1]{fiore2007differential} and \cite[Prop 2]{Blute2019}. 

We propose a slight change in terminology regarding the term \emph{bialgebra modality} in comparison to previous papers. Blute, Cockett, and Seely first used the term bialgebra modality in \cite[Def 4.8]{blute2006differential} for a coalgebra modality on an additive symmetric monoidal category with additional natural (bi)monoid structural maps which also satisfied extra compatibilities with the dereliction which involved the additive enrichment. Here, however, we find it more natural to refer to a bialgebra modality for the natural notion the name suggests: as simply a coalgebra modality with natural bimonoid structure and no extra requirements, which can thus be defined in any symmetric monoidal category. We rename the Blute, Cockett, and Seely version as a \emph{pre-additive} bialgebra modality, which we revisit in Sec \ref{sec:add-bialg-mon}. 

So in this paper, a \textbf{bialgebra modality} on a symmetric monoidal category $\mathbb{X}$ is a septuple $(\oc, \delta, \varepsilon, \Delta, \mathsf{e}, \nabla, \mathsf{u})$ consisting of a coalgebra modality $(\oc, \delta, \varepsilon, \Delta, \mathsf{e})$ and two natural transformations $\nabla_A: \oc(A) \otimes \oc(A) \to \oc(A)$ and ${\mathsf{u}_A: I \to \oc(A)}$, such that the diagrams in Appendix \ref{sec:diagbialgmod} commute. The natural transformation $\nabla$ is called the \textbf{multiplication} or \textbf{cocontraction} and the natural transformation $\mathsf{u}$ is called the \textbf{unit} or \textbf{coweakening}, and are drawn as follows: 
 \begin{align*}\begin{array}[c]{c} 
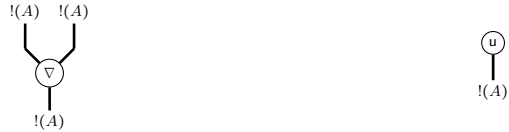
\begin{tikzpicture}
	\begin{pgfonlayer}{nodelayer}
		\node [style=object] (219) at (34.5, 1.75) {$\oc(A)$};
		\node [style=object] (220) at (34, 4) {$\oc(A)$};
		\node [style=component] (221) at (34.5, 2.75) {$\nabla$};
		\node [style=none] (222) at (34, 3.25) {};
		\node [style=none] (223) at (35, 3.25) {};
		\node [style=object] (224) at (35, 4) {$\oc(A)$};
	\end{pgfonlayer}
	\begin{pgfonlayer}{edgelayer}
		\draw [style=wire] (219) to (221);
		\draw [style=wire] (221) to (222.center);
		\draw [style=wire] (221) to (223.center);
		\draw [style=wire] (223.center) to (224);
		\draw [style=wire] (222.center) to (220);
	\end{pgfonlayer}
\end{tikzpicture}
 \end{array} && \begin{array}[c]{c} \begin{tikzpicture}
	\begin{pgfonlayer}{nodelayer}
		\node [style=object] (225) at (37, 1.75) {$\oc(A)$};
		\node [style=component] (226) at (37, 2.75) {$\mathsf{u}$};
	\end{pgfonlayer}
	\begin{pgfonlayer}{edgelayer}
		\draw [style=wire] (225) to (226);
	\end{pgfonlayer}
\end{tikzpicture}
 \end{array} 
 \end{align*}
with their naturality drawn as follows: 
 \begin{align}\label{string:monoid-nat}\begin{array}[c]{c} 
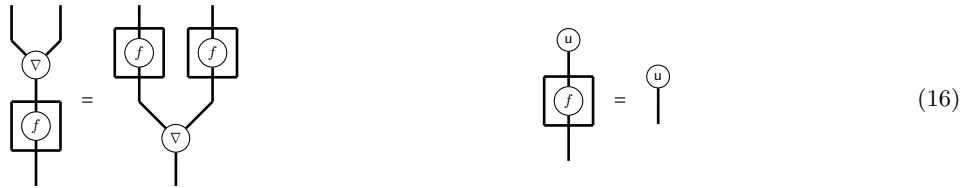
\begin{tikzpicture}
	\begin{pgfonlayer}{nodelayer}
		\node [style=component] (259) at (46.75, 1.75) {$f$};
		\node [style=none] (260) at (46.25, 1.25) {};
		\node [style=none] (261) at (47.25, 1.25) {};
		\node [style=none] (262) at (46.25, 2.25) {};
		\node [style=none] (263) at (47.25, 2.25) {};
		\node [style=port] (264) at (46.25, 4.25) {};
		\node [style=component] (265) at (46.75, 3) {$\nabla$};
		\node [style=none] (266) at (46.25, 3.5) {};
		\node [style=none] (267) at (47.25, 3.5) {};
		\node [style=port] (268) at (47.25, 4.25) {};
		\node [style=port] (269) at (46.75, 0.5) {};
	\end{pgfonlayer}
	\begin{pgfonlayer}{edgelayer}
		\draw [style=wire] (260.center) to (261.center);
		\draw [style=wire] (261.center) to (263.center);
		\draw [style=wire] (263.center) to (262.center);
		\draw [style=wire] (262.center) to (260.center);
		\draw [style=wire] (265) to (266.center);
		\draw [style=wire] (265) to (267.center);
		\draw [style=wire] (266.center) to (264);
		\draw [style=wire] (267.center) to (268);
		\draw [style=wire] (269) to (259);
		\draw [style=wire] (259) to (265);
	\end{pgfonlayer}
\end{tikzpicture}
 \end{array} = \begin{array}[c]{c}\begin{tikzpicture}
	\begin{pgfonlayer}{nodelayer}
		\node [style=component] (270) at (48.75, 3.25) {$f$};
		\node [style=none] (271) at (48.25, 2.75) {};
		\node [style=none] (272) at (49.25, 2.75) {};
		\node [style=none] (273) at (48.25, 3.75) {};
		\node [style=none] (274) at (49.25, 3.75) {};
		\node [style=port] (275) at (48.75, 4.25) {};
		\node [style=component] (276) at (49.5, 1.5) {$\nabla$};
		\node [style=none] (277) at (48.75, 2.25) {};
		\node [style=none] (278) at (50.25, 2.25) {};
		\node [style=port] (279) at (50.25, 4.25) {};
		\node [style=component] (280) at (50.25, 3.25) {$f$};
		\node [style=none] (281) at (49.75, 2.75) {};
		\node [style=none] (282) at (50.75, 2.75) {};
		\node [style=none] (283) at (49.75, 3.75) {};
		\node [style=none] (284) at (50.75, 3.75) {};
		\node [style=port] (285) at (49.5, 0.5) {};
	\end{pgfonlayer}
	\begin{pgfonlayer}{edgelayer}
		\draw [style=wire] (271.center) to (272.center);
		\draw [style=wire] (272.center) to (274.center);
		\draw [style=wire] (274.center) to (273.center);
		\draw [style=wire] (273.center) to (271.center);
		\draw [style=wire] (276) to (277.center);
		\draw [style=wire] (276) to (278.center);
		\draw [style=wire] (281.center) to (282.center);
		\draw [style=wire] (282.center) to (284.center);
		\draw [style=wire] (284.center) to (283.center);
		\draw [style=wire] (283.center) to (281.center);
		\draw [style=wire] (270) to (275);
		\draw [style=wire] (280) to (279);
		\draw [style=wire] (277.center) to (270);
		\draw [style=wire] (278.center) to (280);
		\draw [style=wire] (285) to (276);
	\end{pgfonlayer}
\end{tikzpicture}
 \end{array} && \begin{array}[c]{c}\begin{tikzpicture}
	\begin{pgfonlayer}{nodelayer}
		\node [style=component] (286) at (52.25, 3) {$f$};
		\node [style=none] (287) at (51.75, 2.5) {};
		\node [style=none] (288) at (52.75, 2.5) {};
		\node [style=none] (289) at (51.75, 3.5) {};
		\node [style=none] (290) at (52.75, 3.5) {};
		\node [style=component] (291) at (52.25, 4.25) {$\mathsf{u}$};
		\node [style=port] (292) at (52.25, 1.75) {};
	\end{pgfonlayer}
	\begin{pgfonlayer}{edgelayer}
		\draw [style=wire] (287.center) to (288.center);
		\draw [style=wire] (288.center) to (290.center);
		\draw [style=wire] (290.center) to (289.center);
		\draw [style=wire] (289.center) to (287.center);
		\draw [style=wire] (292) to (286);
		\draw [style=wire] (286) to (291);
	\end{pgfonlayer}
\end{tikzpicture}
 \end{array} = \begin{array}[c]{c}\begin{tikzpicture}
	\begin{pgfonlayer}{nodelayer}
		\node [style=component] (293) at (53.75, 4.25) {$\mathsf{u}$};
		\node [style=port] (294) at (53.75, 3.25) {};
	\end{pgfonlayer}
	\begin{pgfonlayer}{edgelayer}
		\draw [style=wire] (294) to (293);
	\end{pgfonlayer}
\end{tikzpicture}
 \end{array}
 \end{align}
  The requirement that for each object $A$, $(\oc(A), \nabla_A, \mathsf{u}_A)$ is a commutative monoid (\ref{diag:monoid}) is drawn as follows: 
  \begin{align}\label{string:monoid}\begin{array}[c]{c}\begin{tikzpicture}
	\begin{pgfonlayer}{nodelayer}
		\node [style=port] (227) at (38.25, -3) {};
		\node [style=port] (228) at (37.75, 0.5) {};
		\node [style=component] (229) at (38.25, -2) {$\nabla$};
		\node [style=none] (230) at (37.75, -1.5) {};
		\node [style=none] (231) at (38.75, -1.5) {};
		\node [style=port] (232) at (38.25, 0.5) {};
		\node [style=component] (233) at (38.75, -0.75) {$\nabla$};
		\node [style=none] (234) at (38.25, -0.25) {};
		\node [style=none] (235) at (39.25, -0.25) {};
		\node [style=port] (236) at (39.25, 0.5) {};
	\end{pgfonlayer}
	\begin{pgfonlayer}{edgelayer}
		\draw [style=wire] (227) to (229);
		\draw [style=wire] (229) to (230.center);
		\draw [style=wire] (229) to (231.center);
		\draw [style=wire] (230.center) to (228);
		\draw [style=wire] (233) to (234.center);
		\draw [style=wire] (233) to (235.center);
		\draw [style=wire] (235.center) to (236);
		\draw [style=wire] (234.center) to (232);
		\draw [style=wire] (231.center) to (233);
	\end{pgfonlayer}
\end{tikzpicture}
 \end{array} = \begin{array}[c]{c}\begin{tikzpicture}
	\begin{pgfonlayer}{nodelayer}
		\node [style=port] (227) at (38.75, -3) {};
		\node [style=port] (228) at (39.25, 0.5) {};
		\node [style=component] (229) at (38.75, -2) {$\nabla$};
		\node [style=none] (230) at (39.25, -1.5) {};
		\node [style=none] (231) at (38.25, -1.5) {};
		\node [style=port] (232) at (38.75, 0.5) {};
		\node [style=component] (233) at (38.25, -0.75) {$\nabla$};
		\node [style=none] (234) at (38.75, -0.25) {};
		\node [style=none] (235) at (37.75, -0.25) {};
		\node [style=port] (236) at (37.75, 0.5) {};
	\end{pgfonlayer}
	\begin{pgfonlayer}{edgelayer}
		\draw [style=wire] (227) to (229);
		\draw [style=wire] (229) to (230.center);
		\draw [style=wire] (229) to (231.center);
		\draw [style=wire] (230.center) to (228);
		\draw [style=wire] (233) to (234.center);
		\draw [style=wire] (233) to (235.center);
		\draw [style=wire] (235.center) to (236);
		\draw [style=wire] (234.center) to (232);
		\draw [style=wire] (231.center) to (233);
	\end{pgfonlayer}
\end{tikzpicture}
 \end{array} && \begin{array}[c]{c} \begin{tikzpicture}
	\begin{pgfonlayer}{nodelayer}
		\node [style=port] (237) at (40.75, 1.5) {};
		\node [style=port] (238) at (40.25, 4) {};
		\node [style=component] (239) at (40.75, 2.5) {$\nabla$};
		\node [style=none] (240) at (40.25, 3) {};
		\node [style=none] (241) at (41.25, 3) {};
		\node [style=component] (242) at (41.25, 3.5) {$\mathsf{u}$};
	\end{pgfonlayer}
	\begin{pgfonlayer}{edgelayer}
		\draw [style=wire] (237) to (239);
		\draw [style=wire] (239) to (240.center);
		\draw [style=wire] (239) to (241.center);
		\draw [style=wire] (240.center) to (238);
		\draw [style=wire] (241.center) to (242);
	\end{pgfonlayer}
\end{tikzpicture}
 \end{array} = \begin{array}[c]{c} \begin{tikzpicture}
	\begin{pgfonlayer}{nodelayer}
		\node [style=port] (16) at (6.5, 1.5) {};
		\node [style=port] (18) at (6.5, 4) {};
	\end{pgfonlayer}
	\begin{pgfonlayer}{edgelayer}
		\draw [style=wire] (18) to (16);
	\end{pgfonlayer}
\end{tikzpicture}
 \end{array} = \begin{array}[c]{c}\begin{tikzpicture}
	\begin{pgfonlayer}{nodelayer}
		\node [style=port] (237) at (40.75, 1.5) {};
		\node [style=port] (238) at (41.25, 4) {};
		\node [style=component] (239) at (40.75, 2.5) {$\nabla$};
		\node [style=none] (240) at (41.25, 3) {};
		\node [style=none] (241) at (40.25, 3) {};
		\node [style=component] (242) at (40.25, 3.5) {$\mathsf{u}$};
	\end{pgfonlayer}
	\begin{pgfonlayer}{edgelayer}
		\draw [style=wire] (237) to (239);
		\draw [style=wire] (239) to (240.center);
		\draw [style=wire] (239) to (241.center);
		\draw [style=wire] (240.center) to (238);
		\draw [style=wire] (241.center) to (242);
	\end{pgfonlayer}
\end{tikzpicture}
 \end{array} && \begin{array}[c]{c}\begin{tikzpicture}
	\begin{pgfonlayer}{nodelayer}
		\node [style=port] (243) at (42.75, 1.25) {};
		\node [style=port] (244) at (42.25, 4) {};
		\node [style=component] (245) at (42.75, 2.25) {$\nabla$};
		\node [style=none] (246) at (42.25, 2.75) {};
		\node [style=none] (247) at (43.25, 2.75) {};
		\node [style=port] (248) at (43.25, 4) {};
		\node [style=none] (249) at (42.25, 3) {};
		\node [style=none] (250) at (43.25, 3) {};
		\node [style=none] (251) at (42.25, 3.5) {};
		\node [style=none] (252) at (43.25, 3.5) {};
	\end{pgfonlayer}
	\begin{pgfonlayer}{edgelayer}
		\draw [style=wire] (243) to (245);
		\draw [style=wire] (245) to (246.center);
		\draw [style=wire] (245) to (247.center);
		\draw [style=wire] (246.center) to (249.center);
		\draw [style=wire] (247.center) to (250.center);
		\draw [style=wire] (251.center) to (244);
		\draw [style=wire] (252.center) to (248);
		\draw [style=wire] (249.center) to (252.center);
		\draw [style=wire] (250.center) to (251.center);
	\end{pgfonlayer}
\end{tikzpicture}
 \end{array} = \begin{array}[c]{c} \begin{tikzpicture}
	\begin{pgfonlayer}{nodelayer}
		\node [style=port] (253) at (44.75, 1.25) {};
		\node [style=port] (254) at (44.25, 4) {};
		\node [style=component] (255) at (44.75, 2.25) {$\nabla$};
		\node [style=none] (256) at (44.25, 2.75) {};
		\node [style=none] (257) at (45.25, 2.75) {};
		\node [style=port] (258) at (45.25, 4) {};
	\end{pgfonlayer}
	\begin{pgfonlayer}{edgelayer}
		\draw [style=wire] (253) to (255);
		\draw [style=wire] (255) to (256.center);
		\draw [style=wire] (255) to (257.center);
		\draw [style=wire] (256.center) to (254);
		\draw [style=wire] (257.center) to (258);
	\end{pgfonlayer}
\end{tikzpicture}
 \end{array}
 \end{align}
 where the first equality is called the associativity of the multiplication, the second is the unit axioms, and the third is called the commutativity of the multiplication. As before, note that naturality of $\nabla$ and $\mathsf{u}$ are together precisely the statement that for every map $f$, $\oc(f)$ is a monoid morphism. We then also ask that for each object $A$, $(\oc(A), \nabla_A, \mathsf{u}_A, \Delta_A, \mathsf{e}_A)$ is a bimonoid (\ref{diag:bimonoid}), whose necessary axioms are drawn out as follows: 
  \begin{align}\label{string:bimonoid}\begin{array}[c]{c}\begin{tikzpicture}
	\begin{pgfonlayer}{nodelayer}
		\node [style=port] (29) at (17, -5.5) {};
		\node [style=component] (30) at (17.5, -4.25) {$\Delta$};
		\node [style=none] (31) at (17, -4.75) {};
		\node [style=none] (32) at (18, -4.75) {};
		\node [style=port] (33) at (18, -5.5) {};
		\node [style=port] (162) at (17, -1.75) {};
		\node [style=component] (163) at (17.5, -3) {$\nabla$};
		\node [style=none] (164) at (17, -2.5) {};
		\node [style=none] (165) at (18, -2.5) {};
		\node [style=port] (166) at (18, -1.75) {};
	\end{pgfonlayer}
	\begin{pgfonlayer}{edgelayer}
		\draw [style=wire] (30) to (31.center);
		\draw [style=wire] (30) to (32.center);
		\draw [style=wire] (32.center) to (33);
		\draw [style=wire] (31.center) to (29);
		\draw [style=wire] (163) to (164.center);
		\draw [style=wire] (163) to (165.center);
		\draw [style=wire] (164.center) to (162);
		\draw [style=wire] (165.center) to (166);
		\draw [style=wire] (163) to (30);
	\end{pgfonlayer}
\end{tikzpicture}
 \end{array} = \begin{array}[c]{c}\begin{tikzpicture}
	\begin{pgfonlayer}{nodelayer}
		\node [style=component] (167) at (20.25, -2.5) {$\Delta$};
		\node [style=none] (168) at (19.75, -3) {};
		\node [style=none] (169) at (20.75, -3) {};
		\node [style=port] (170) at (20.25, -1.75) {};
		\node [style=port] (171) at (21.75, -1.75) {};
		\node [style=component] (172) at (21.75, -2.5) {$\Delta$};
		\node [style=none] (173) at (21.25, -3) {};
		\node [style=none] (174) at (22.25, -3) {};
		\node [style=none] (179) at (20.75, -4) {};
		\node [style=none] (180) at (21.25, -4) {};
		\node [style=none] (181) at (20.75, -3.25) {};
		\node [style=none] (182) at (21.25, -3.25) {};
		\node [style=component] (183) at (20.25, -4.75) {$\nabla$};
		\node [style=none] (184) at (19.75, -4.25) {};
		\node [style=none] (185) at (20.75, -4.25) {};
		\node [style=component] (186) at (21.75, -4.75) {$\nabla$};
		\node [style=none] (187) at (21.25, -4.25) {};
		\node [style=none] (188) at (22.25, -4.25) {};
		\node [style=port] (189) at (20.25, -5.5) {};
		\node [style=port] (190) at (21.75, -5.5) {};
	\end{pgfonlayer}
	\begin{pgfonlayer}{edgelayer}
		\draw [style=wire] (167) to (168.center);
		\draw [style=wire] (167) to (169.center);
		\draw [style=wire] (172) to (173.center);
		\draw [style=wire] (172) to (174.center);
		\draw [style=wire] (170) to (167);
		\draw [style=wire] (171) to (172);
		\draw [style=wire] (179.center) to (182.center);
		\draw [style=wire] (180.center) to (181.center);
		\draw [style=wire] (169.center) to (181.center);
		\draw [style=wire] (173.center) to (182.center);
		\draw [style=wire] (183) to (184.center);
		\draw [style=wire] (183) to (185.center);
		\draw [style=wire] (186) to (187.center);
		\draw [style=wire] (186) to (188.center);
		\draw [style=wire] (168.center) to (184.center);
		\draw [style=wire] (174.center) to (188.center);
		\draw [style=wire] (179.center) to (185.center);
		\draw [style=wire] (180.center) to (187.center);
		\draw [style=wire] (183) to (189);
		\draw [style=wire] (186) to (190);
	\end{pgfonlayer}
\end{tikzpicture}
 \end{array} && \begin{array}[c]{c}\begin{tikzpicture}
	\begin{pgfonlayer}{nodelayer}
		\node [style=component] (116) at (22.25, -0.25) {$\mathsf{e}$};
		\node [style=port] (117) at (21.75, 2) {};
		\node [style=component] (118) at (22.25, 0.75) {$\nabla$};
		\node [style=none] (119) at (21.75, 1.25) {};
		\node [style=none] (120) at (22.75, 1.25) {};
		\node [style=port] (121) at (22.75, 2) {};
	\end{pgfonlayer}
	\begin{pgfonlayer}{edgelayer}
		\draw [style=wire] (118) to (119.center);
		\draw [style=wire] (118) to (120.center);
		\draw [style=wire] (119.center) to (117);
		\draw [style=wire] (120.center) to (121);
		\draw [style=wire] (116) to (118);
	\end{pgfonlayer}
\end{tikzpicture}
 \end{array} = \begin{array}[c]{c}\begin{tikzpicture}
	\begin{pgfonlayer}{nodelayer}
		\node [style=port] (128) at (24, 2) {};
		\node [style=port] (132) at (25, 2) {};
		\node [style=component] (160) at (24, 0.5) {$\mathsf{e}$};
		\node [style=component] (161) at (25, 0.5) {$\mathsf{e}$};
	\end{pgfonlayer}
	\begin{pgfonlayer}{edgelayer}
		\draw [style=wire] (128) to (160);
		\draw [style=wire] (132) to (161);
	\end{pgfonlayer}
\end{tikzpicture}
 \end{array} && \begin{array}[c]{c}\begin{tikzpicture}
	\begin{pgfonlayer}{nodelayer}
		\node [style=component] (116) at (22.25, 2) {$\mathsf{u}$};
		\node [style=port] (117) at (21.75, -0.25) {};
		\node [style=component] (118) at (22.25, 1) {$\Delta$};
		\node [style=none] (119) at (21.75, 0.5) {};
		\node [style=none] (120) at (22.75, 0.5) {};
		\node [style=port] (121) at (22.75, -0.25) {};
	\end{pgfonlayer}
	\begin{pgfonlayer}{edgelayer}
		\draw [style=wire] (118) to (119.center);
		\draw [style=wire] (118) to (120.center);
		\draw [style=wire] (119.center) to (117);
		\draw [style=wire] (120.center) to (121);
		\draw [style=wire] (116) to (118);
	\end{pgfonlayer}
\end{tikzpicture}
 \end{array} = \begin{array}[c]{c}\begin{tikzpicture}
	\begin{pgfonlayer}{nodelayer}
		\node [style=port] (128) at (24, 0.5) {};
		\node [style=port] (132) at (25, 0.5) {};
		\node [style=component] (160) at (24, 2) {$\mathsf{u}$};
		\node [style=component] (161) at (25, 2) {$\mathsf{u}$};
	\end{pgfonlayer}
	\begin{pgfonlayer}{edgelayer}
		\draw [style=wire] (128) to (160);
		\draw [style=wire] (132) to (161);
	\end{pgfonlayer}
\end{tikzpicture}
 \end{array} && \begin{array}[c]{c}\begin{tikzpicture}
	\begin{pgfonlayer}{nodelayer}
		\node [style=component] (215) at (33, -0.5) {$\mathsf{e}$};
		\node [style=component] (216) at (33, 0.5) {$\mathsf{u}$};
	\end{pgfonlayer}
	\begin{pgfonlayer}{edgelayer}
		\draw [style=wire] (216) to (215);
	\end{pgfonlayer}
\end{tikzpicture}
 \end{array} = 
 \end{align}
 It is worth recalling that the axioms of a bimonoid are equivalent to asking that $\Delta$ and $\mathsf{e}$ be monoid morphism, or equivalently that $\nabla$ and $\mathsf{u}$ be comonoid morphisms. 

Then a \textbf{monoidal bialgebra modality} on a symmetric monoidal category $\mathbb{X}$ is a nonuple $(\oc, \delta, \varepsilon,  \Delta, \mathsf{e}, \mathsf{m}, \mathsf{m}_I, \nabla, \mathsf{u})$ consisting of a monoidal coalgebra modality $(\oc, \delta, \varepsilon,  \Delta, \mathsf{e}, \mathsf{m}, \mathsf{m}_I)$ and a bialgebra modality $(\oc, \delta, \varepsilon, \Delta, \mathsf{e}, \nabla, \mathsf{u})$, such that the diagrams in Appendix \ref{sec:diagmonbialgmod} commute. The first requirement is that the monoid structure is compatible with $\mathsf{m}$ (\ref{diag:nablamonspecial}) in the sense that following equalities hold: 
   \begin{align}\label{string:m-monoid}\begin{array}[c]{c}\begin{tikzpicture}
	\begin{pgfonlayer}{nodelayer}
		\node [style=port] (237) at (79, -7) {};
		\node [style=port] (238) at (78.5, -1.75) {};
		\node [style=component] (239) at (79, -6.25) {$\mathsf{m}$};
		\node [style=none] (240) at (78.5, -5.5) {};
		\node [style=none] (241) at (79.5, -5.5) {};
		\node [style=port] (367) at (80, -1.75) {};
		\node [style=component] (368) at (79.5, -4.25) {$\nabla$};
		\node [style=none] (369) at (80, -3.25) {};
		\node [style=none] (370) at (79, -3.25) {};
		\node [style=port] (371) at (79, -1.75) {};
	\end{pgfonlayer}
	\begin{pgfonlayer}{edgelayer}
		\draw [style=wire] (237) to (239);
		\draw [style=wire] (239) to (240.center);
		\draw [style=wire] (239) to (241.center);
		\draw [style=wire] (240.center) to (238);
		\draw [style=wire] (368) to (369.center);
		\draw [style=wire] (368) to (370.center);
		\draw [style=wire] (370.center) to (371);
		\draw [style=wire] (369.center) to (367);
		\draw [style=wire] (368) to (241.center);
	\end{pgfonlayer}
\end{tikzpicture}
 \end{array} = \begin{array}[c]{c}\begin{tikzpicture}
	\begin{pgfonlayer}{nodelayer}
		\node [style=component] (381) at (76.25, -6.25) {$\nabla$};
		\node [style=none] (382) at (77, -5.5) {};
		\node [style=none] (383) at (75.5, -5.5) {};
		\node [style=component] (388) at (75.5, -2.5) {$\Delta$};
		\node [style=none] (389) at (75, -3) {};
		\node [style=none] (390) at (76, -3) {};
		\node [style=port] (391) at (75.5, -1.75) {};
		\node [style=port] (392) at (77.5, -1.75) {};
		\node [style=none] (394) at (76.5, -3) {};
		\node [style=none] (395) at (77.5, -3) {};
		\node [style=none] (396) at (76, -4) {};
		\node [style=none] (397) at (76.5, -4) {};
		\node [style=none] (398) at (76, -3.25) {};
		\node [style=none] (399) at (76.5, -3.25) {};
		\node [style=component] (400) at (75.5, -4.75) {$\mathsf{m}$};
		\node [style=none] (401) at (75, -4.25) {};
		\node [style=none] (402) at (76, -4.25) {};
		\node [style=component] (403) at (77, -4.75) {$\mathsf{m}$};
		\node [style=none] (404) at (76.5, -4.25) {};
		\node [style=none] (405) at (77.5, -4.25) {};
		\node [style=port] (406) at (76.5, -1.75) {};
		\node [style=port] (407) at (76.25, -7) {};
	\end{pgfonlayer}
	\begin{pgfonlayer}{edgelayer}
		\draw [style=wire] (381) to (382.center);
		\draw [style=wire] (381) to (383.center);
		\draw [style=wire] (388) to (389.center);
		\draw [style=wire] (388) to (390.center);
		\draw [style=wire] (391) to (388);
		\draw [style=wire] (396.center) to (399.center);
		\draw [style=wire] (397.center) to (398.center);
		\draw [style=wire] (390.center) to (398.center);
		\draw [style=wire] (394.center) to (399.center);
		\draw [style=wire] (400) to (401.center);
		\draw [style=wire] (400) to (402.center);
		\draw [style=wire] (403) to (404.center);
		\draw [style=wire] (403) to (405.center);
		\draw [style=wire] (389.center) to (401.center);
		\draw [style=wire] (395.center) to (405.center);
		\draw [style=wire] (396.center) to (402.center);
		\draw [style=wire] (397.center) to (404.center);
		\draw [style=wire] (400) to (383.center);
		\draw [style=wire] (403) to (382.center);
		\draw [style=wire] (406) to (399.center);
		\draw [style=wire] (392) to (395.center);
		\draw [style=wire] (381) to (407);
	\end{pgfonlayer}
\end{tikzpicture}
 \end{array} && \begin{array}[c]{c}\begin{tikzpicture}
	\begin{pgfonlayer}{nodelayer}
		\node [style=port] (237) at (79.5, -7) {};
		\node [style=port] (238) at (80, -1.75) {};
		\node [style=component] (239) at (79.5, -6.25) {$\mathsf{m}$};
		\node [style=none] (240) at (80, -5.5) {};
		\node [style=none] (241) at (79, -5.5) {};
		\node [style=port] (367) at (78.5, -1.75) {};
		\node [style=component] (368) at (79, -4.25) {$\nabla$};
		\node [style=none] (369) at (78.5, -3.25) {};
		\node [style=none] (370) at (79.5, -3.25) {};
		\node [style=port] (371) at (79.5, -1.75) {};
	\end{pgfonlayer}
	\begin{pgfonlayer}{edgelayer}
		\draw [style=wire] (237) to (239);
		\draw [style=wire] (239) to (240.center);
		\draw [style=wire] (239) to (241.center);
		\draw [style=wire] (240.center) to (238);
		\draw [style=wire] (368) to (369.center);
		\draw [style=wire] (368) to (370.center);
		\draw [style=wire] (370.center) to (371);
		\draw [style=wire] (369.center) to (367);
		\draw [style=wire] (368) to (241.center);
	\end{pgfonlayer}
\end{tikzpicture}
 \end{array} = \begin{array}[c]{c}\begin{tikzpicture}
	\begin{pgfonlayer}{nodelayer}
		\node [style=component] (381) at (76.25, -6.25) {$\nabla$};
		\node [style=none] (382) at (75.5, -5.5) {};
		\node [style=none] (383) at (77, -5.5) {};
		\node [style=component] (388) at (77, -2.5) {$\Delta$};
		\node [style=none] (389) at (77.5, -3) {};
		\node [style=none] (390) at (76.5, -3) {};
		\node [style=port] (391) at (77, -1.75) {};
		\node [style=port] (392) at (75, -1.75) {};
		\node [style=none] (394) at (76, -3) {};
		\node [style=none] (395) at (75, -3) {};
		\node [style=none] (396) at (76.5, -4) {};
		\node [style=none] (397) at (76, -4) {};
		\node [style=none] (398) at (76.5, -3.25) {};
		\node [style=none] (399) at (76, -3.25) {};
		\node [style=component] (400) at (77, -4.75) {$\mathsf{m}$};
		\node [style=none] (401) at (77.5, -4.25) {};
		\node [style=none] (402) at (76.5, -4.25) {};
		\node [style=component] (403) at (75.5, -4.75) {$\mathsf{m}$};
		\node [style=none] (404) at (76, -4.25) {};
		\node [style=none] (405) at (75, -4.25) {};
		\node [style=port] (406) at (76, -1.75) {};
		\node [style=port] (407) at (76.25, -7) {};
	\end{pgfonlayer}
	\begin{pgfonlayer}{edgelayer}
		\draw [style=wire] (381) to (382.center);
		\draw [style=wire] (381) to (383.center);
		\draw [style=wire] (388) to (389.center);
		\draw [style=wire] (388) to (390.center);
		\draw [style=wire] (391) to (388);
		\draw [style=wire] (396.center) to (399.center);
		\draw [style=wire] (397.center) to (398.center);
		\draw [style=wire] (390.center) to (398.center);
		\draw [style=wire] (394.center) to (399.center);
		\draw [style=wire] (400) to (401.center);
		\draw [style=wire] (400) to (402.center);
		\draw [style=wire] (403) to (404.center);
		\draw [style=wire] (403) to (405.center);
		\draw [style=wire] (389.center) to (401.center);
		\draw [style=wire] (395.center) to (405.center);
		\draw [style=wire] (396.center) to (402.center);
		\draw [style=wire] (397.center) to (404.center);
		\draw [style=wire] (400) to (383.center);
		\draw [style=wire] (403) to (382.center);
		\draw [style=wire] (406) to (399.center);
		\draw [style=wire] (392) to (395.center);
		\draw [style=wire] (381) to (407);
	\end{pgfonlayer}
\end{tikzpicture}
 \end{array} && \begin{array}[c]{c}\begin{tikzpicture}
	\begin{pgfonlayer}{nodelayer}
		\node [style=port] (237) at (35.5, -7) {};
		\node [style=port] (238) at (36, -4.5) {};
		\node [style=component] (239) at (35.5, -6) {$\mathsf{m}$};
		\node [style=none] (240) at (36, -5.5) {};
		\node [style=none] (241) at (35, -5.5) {};
		\node [style=component] (242) at (35, -5) {$\mathsf{u}$};
	\end{pgfonlayer}
	\begin{pgfonlayer}{edgelayer}
		\draw [style=wire] (237) to (239);
		\draw [style=wire] (239) to (240.center);
		\draw [style=wire] (239) to (241.center);
		\draw [style=wire] (240.center) to (238);
		\draw [style=wire] (241.center) to (242);
	\end{pgfonlayer}
\end{tikzpicture}
 \end{array} = \begin{array}[c]{c}\begin{tikzpicture}
	\begin{pgfonlayer}{nodelayer}
		\node [style=component] (215) at (37.25, -5.25) {$\mathsf{e}$};
		\node [style=port] (216) at (37.25, -4.5) {};
		\node [style=component] (217) at (37.25, -6.25) {$\mathsf{u}$};
		\node [style=port] (218) at (37.25, -7) {};
	\end{pgfonlayer}
	\begin{pgfonlayer}{edgelayer}
		\draw [style=wire] (216) to (215);
		\draw [style=wire] (218) to (217);
	\end{pgfonlayer}
\end{tikzpicture} \end{array}
&& \begin{array}[c]{c}\begin{tikzpicture}
	\begin{pgfonlayer}{nodelayer}
		\node [style=port] (237) at (35.5, -7) {};
		\node [style=port] (238) at (35, -4.5) {};
		\node [style=component] (239) at (35.5, -6) {$\mathsf{m}$};
		\node [style=none] (240) at (35, -5.5) {};
		\node [style=none] (241) at (36, -5.5) {};
		\node [style=component] (242) at (36, -5) {$\mathsf{u}$};
	\end{pgfonlayer}
	\begin{pgfonlayer}{edgelayer}
		\draw [style=wire] (237) to (239);
		\draw [style=wire] (239) to (240.center);
		\draw [style=wire] (239) to (241.center);
		\draw [style=wire] (240.center) to (238);
		\draw [style=wire] (241.center) to (242);
	\end{pgfonlayer}
\end{tikzpicture}
 \end{array} = \begin{array}[c]{c}\begin{tikzpicture}
	\begin{pgfonlayer}{nodelayer}
		\node [style=component] (215) at (37.25, -5.25) {$\mathsf{e}$};
		\node [style=port] (216) at (37.25, -4.5) {};
		\node [style=component] (217) at (37.25, -6.25) {$\mathsf{u}$};
		\node [style=port] (218) at (37.25, -7) {};
	\end{pgfonlayer}
	\begin{pgfonlayer}{edgelayer}
		\draw [style=wire] (216) to (215);
		\draw [style=wire] (218) to (217);
	\end{pgfonlayer}
\end{tikzpicture} \end{array}
 \end{align}
 The other requirement is that $\nabla$ and $\mathsf{u}$ are $\oc$-coalgebra morphisms (\ref{diag:!coalgmonoid}), which is drawn as follows: 
  \begin{align}\label{string:monoid-!coalg}\begin{array}[c]{c}\begin{tikzpicture}
	\begin{pgfonlayer}{nodelayer}
		\node [style=component] (259) at (58, -1) {$\delta$};
		\node [style=port] (264) at (57.5, 2.5) {};
		\node [style=component] (265) at (58, 0) {$\nabla$};
		\node [style=none] (266) at (57.5, 1) {};
		\node [style=none] (267) at (58.5, 1) {};
		\node [style=port] (268) at (58.5, 2.5) {};
		\node [style=port] (269) at (58, -2.25) {};
	\end{pgfonlayer}
	\begin{pgfonlayer}{edgelayer}
		\draw [style=wire] (265) to (266.center);
		\draw [style=wire] (265) to (267.center);
		\draw [style=wire] (266.center) to (264);
		\draw [style=wire] (267.center) to (268);
		\draw [style=wire] (269) to (259);
		\draw [style=wire] (259) to (265);
	\end{pgfonlayer}
\end{tikzpicture}
 \end{array} = \begin{array}[c]{c}\begin{tikzpicture}
	\begin{pgfonlayer}{nodelayer}
		\node [style=none] (321) at (54.75, -1.5) {};
		\node [style=none] (322) at (56.75, -1.5) {};
		\node [style=none] (323) at (54.75, 0) {};
		\node [style=none] (324) at (56.75, 0) {};
		\node [style=port] (325) at (55.25, 2.5) {};
		\node [style=component] (326) at (55.75, 0.75) {$\mathsf{m}$};
		\node [style=none] (327) at (55.25, 1.25) {};
		\node [style=none] (328) at (56.25, 1.25) {};
		\node [style=port] (329) at (56.25, 2.5) {};
		\node [style=port] (330) at (55.75, -2.25) {};
		\node [style=none] (331) at (55.25, 0) {};
		\node [style=none] (332) at (56.25, 0) {};
		\node [style=none] (333) at (55.25, -1.5) {};
		\node [style=none] (334) at (56.25, -1.5) {};
		\node [style=none] (335) at (55.75, 0) {};
		\node [style=none] (336) at (55.75, -1.5) {};
		\node [style=component] (337) at (55.75, -1) {$\nabla$};
		\node [style=none] (338) at (55.25, -0.5) {};
		\node [style=none] (339) at (56.25, -0.5) {};
		\node [style=component] (340) at (55.25, 1.75) {$\delta$};
		\node [style=component] (341) at (56.25, 1.75) {$\delta$};
	\end{pgfonlayer}
	\begin{pgfonlayer}{edgelayer}
		\draw [style=wire] (321.center) to (322.center);
		\draw [style=wire] (322.center) to (324.center);
		\draw [style=wire] (324.center) to (323.center);
		\draw [style=wire] (323.center) to (321.center);
		\draw [style=wire] (326) to (327.center);
		\draw [style=wire] (326) to (328.center);
		\draw [style=wire] (336.center) to (330);
		\draw [style=wire] (326) to (335.center);
		\draw [style=wire] (337) to (338.center);
		\draw [style=wire] (337) to (339.center);
		\draw [style=wire] (331.center) to (338.center);
		\draw [style=wire] (332.center) to (339.center);
		\draw [style=wire] (337) to (336.center);
		\draw [style=wire] (325) to (340);
		\draw [style=wire] (329) to (341);
		\draw [style=wire] (341) to (328.center);
		\draw [style=wire] (340) to (327.center);
	\end{pgfonlayer}
\end{tikzpicture}
 \end{array} && \begin{array}[c]{c}\begin{tikzpicture}
	\begin{pgfonlayer}{nodelayer}
		\node [style=component] (356) at (65, -1) {$\delta$};
		\node [style=port] (357) at (65, -2) {};
		\node [style=component] (358) at (65, 0) {$\mathsf{u}$};
	\end{pgfonlayer}
	\begin{pgfonlayer}{edgelayer}
		\draw [style=wire] (357) to (356);
		\draw [style=wire] (358) to (356);
	\end{pgfonlayer}
\end{tikzpicture}
 \end{array} = \begin{array}[c]{c}\begin{tikzpicture}
	\begin{pgfonlayer}{nodelayer}
		\node [style=component] (359) at (67, 0) {$\mathsf{m}$};
		\node [style=component] (360) at (67, -1) {$\mathsf{u}$};
		\node [style=none] (361) at (66.5, -1.5) {};
		\node [style=none] (362) at (67.5, -1.5) {};
		\node [style=none] (363) at (66.5, -0.5) {};
		\node [style=none] (364) at (67.5, -0.5) {};
		\node [style=port] (365) at (67, -2) {};
		\node [style=none] (366) at (67, -0.5) {};
	\end{pgfonlayer}
	\begin{pgfonlayer}{edgelayer}
		\draw [style=wire] (361.center) to (362.center);
		\draw [style=wire] (362.center) to (364.center);
		\draw [style=wire] (364.center) to (363.center);
		\draw [style=wire] (363.center) to (361.center);
		\draw [style=wire] (365) to (360);
		\draw [style=wire] (359) to (366.center);
	\end{pgfonlayer}
\end{tikzpicture}
 \end{array}
 \end{align}

It is worthwhile to mention that there is technically some redundancy in the definition of a monoidal bialgebra modality. For (\ref{string:m-monoid}), the left (resp. right) versions of the axioms implies the right (resp. left) versions. 

 \begin{lem} For a monoidal coalgebra modality $(\oc, \delta, \varepsilon,  \Delta, \mathsf{e}, \mathsf{m}, \mathsf{m}_I)$, a natural transformation $\nabla: \oc(A) \otimes \oc(A) \to \oc(A)$ (resp. $\mathsf{u}: I \to \oc(A)$) satisfies the first (resp. third) equality in (\ref{string:m-monoid}) if and only if it satisfies the second (resp. fourth) equality in (\ref{string:m-monoid}).
 \end{lem}
 \begin{proof} The proof uses similar techniques as in Lemma \ref{lemma:precoder-left=right} using symmetry and naturality, so we leave this as an exercise for the reader. 
 \end{proof}

 On the other hand, it turns out that asking $\nabla$ and $\mathsf{u}$ to be $\oc$-coalgebra morphisms implies the bimonoid axioms. 

  \begin{lem} For a monoidal coalgebra modality $(\oc, \delta, \varepsilon,  \Delta, \mathsf{e}, \mathsf{m}, \mathsf{m}_I)$, if a natural transformation $\nabla: \oc(A) \otimes \oc(A) \to \oc(A)$ (resp. $\mathsf{u}: I \to \oc(A)$) satisfies the first (resp. second) equality in (\ref{string:monoid-!coalg}), then it satisfies the two left (resp. right) most equalities in (\ref{string:bimonoid}). 
 \end{lem}
 \begin{proof} The efficient version of the proof follows from results about $\oc$-coalgebras of monoidal coalgebra modalities. For a monoidal coalgebra modality, every $\oc$-coalgebra has a canonical comonoid structure and every $\oc$-coalgebra morphism is a comonoid morphism \cite[Prop 2]{schalk2004categorical}. For $\oc(A)$, which recall are always $\oc$-coalgebras (sometimes called the cofree $\oc$-coalgebras), the induced canonical comonoid structure is precisely $(\oc(A), \Delta_A, \mathsf{e}_A)$. Thus asking that $\nabla$ and $\mathsf{u}$ be $\oc$-coalgebra morphisms implies that they are also comonoid morphisms, which is precisely (\ref{string:bimonoid}), as desired. 

 However, since we have not reviewed some of these concepts in this paper, it may be worthwhile for the reader to see the direct computations using string diagrams. Let us begin by showing that if $\mathsf{u}$ satisfies (\ref{string:monoid-!coalg}), then it is also a comonoid morphism, that is, the right most equalities of (\ref{string:bimonoid}) hold. Starting with the compatibility between $\mathsf{u}$ and $\mathsf{e}$, we compute that: 
\begin{gather*}
\begin{array}[c]{c}\begin{tikzpicture}
	\begin{pgfonlayer}{nodelayer}
		\node [style=component] (215) at (33, -0.5) {$\mathsf{e}$};
		\node [style=component] (216) at (33, 0.5) {$\mathsf{u}$};
	\end{pgfonlayer}
	\begin{pgfonlayer}{edgelayer}
		\draw [style=wire] (216) to (215);
	\end{pgfonlayer}
\end{tikzpicture}
 \end{array} \substack{=\\(\ref{string:delta-comonoid-map})}\begin{array}[c]{c}\begin{tikzpicture}
	\begin{pgfonlayer}{nodelayer}
		\node [style=component] (369) at (210, -2) {$\mathsf{e}$};
		\node [style=component] (370) at (210, 0) {$\mathsf{u}$};
		\node [style=component] (371) at (210, -1) {$\delta$};
	\end{pgfonlayer}
	\begin{pgfonlayer}{edgelayer}
		\draw [style=wire] (370) to (371);
		\draw [style=wire] (371) to (369);
	\end{pgfonlayer}
\end{tikzpicture}
 \end{array}\substack{=\\(\ref{string:monoid-!coalg})}\begin{array}[c]{c}\begin{tikzpicture}
	\begin{pgfonlayer}{nodelayer}
		\node [style=component] (369) at (211.5, -2) {$\mathsf{e}$};
		\node [style=component] (372) at (211.5, 0) {$\mathsf{m}$};
		\node [style=component] (373) at (211.5, -1) {$\mathsf{u}$};
		\node [style=none] (374) at (211, -1.5) {};
		\node [style=none] (375) at (212, -1.5) {};
		\node [style=none] (376) at (211, -0.5) {};
		\node [style=none] (377) at (212, -0.5) {};
		\node [style=none] (379) at (211.5, -0.5) {};
	\end{pgfonlayer}
	\begin{pgfonlayer}{edgelayer}
		\draw [style=wire] (374.center) to (375.center);
		\draw [style=wire] (375.center) to (377.center);
		\draw [style=wire] (377.center) to (376.center);
		\draw [style=wire] (376.center) to (374.center);
		\draw [style=wire] (372) to (379.center);
		\draw [style=wire] (373) to (369);
	\end{pgfonlayer}
\end{tikzpicture}
 \end{array}\substack{=\\(\ref{string:comonoid-nat})} \begin{array}[c]{c}\begin{tikzpicture}
	\begin{pgfonlayer}{nodelayer}
		\node [style=component] (215) at (33, -0.5) {$\mathsf{e}$};
		\node [style=component] (216) at (33, 0.5) {$\mathsf{m}$};
	\end{pgfonlayer}
	\begin{pgfonlayer}{edgelayer}
		\draw [style=wire] (216) to (215);
	\end{pgfonlayer}
\end{tikzpicture}
 \end{array} \substack{=\\(\ref{string:monoidalcomonoid})}
\end{gather*}
 Next we compute the compatibility between $\mathsf{u}$ and $\Delta$: 
 \begin{gather*}
 \begin{array}[c]{c}\begin{tikzpicture}
	\begin{pgfonlayer}{nodelayer}
		\node [style=component] (116) at (22.25, 2) {$\mathsf{u}$};
		\node [style=port] (117) at (21.75, -0.25) {};
		\node [style=component] (118) at (22.25, 1) {$\Delta$};
		\node [style=none] (119) at (21.75, 0.5) {};
		\node [style=none] (120) at (22.75, 0.5) {};
		\node [style=port] (121) at (22.75, -0.25) {};
	\end{pgfonlayer}
	\begin{pgfonlayer}{edgelayer}
		\draw [style=wire] (118) to (119.center);
		\draw [style=wire] (118) to (120.center);
		\draw [style=wire] (119.center) to (117);
		\draw [style=wire] (120.center) to (121);
		\draw [style=wire] (116) to (118);
	\end{pgfonlayer}
\end{tikzpicture}
 \end{array} \substack{=\\(\ref{string:comonad})}  \begin{array}[c]{c}\begin{tikzpicture}
	\begin{pgfonlayer}{nodelayer}
		\node [style=component] (380) at (210, -4) {$\mathsf{u}$};
		\node [style=port] (381) at (209.5, -7.75) {};
		\node [style=component] (382) at (210, -5) {$\Delta$};
		\node [style=none] (383) at (209.5, -5.5) {};
		\node [style=none] (384) at (210.5, -5.5) {};
		\node [style=port] (385) at (210.5, -7.75) {};
		\node [style=component] (386) at (209.5, -6.25) {$\delta$};
		\node [style=component] (387) at (209.5, -7) {$\varepsilon$};
		\node [style=component] (388) at (210.5, -6.25) {$\delta$};
		\node [style=component] (389) at (210.5, -7) {$\varepsilon$};
	\end{pgfonlayer}
	\begin{pgfonlayer}{edgelayer}
		\draw [style=wire] (382) to (383.center);
		\draw [style=wire] (382) to (384.center);
		\draw [style=wire] (380) to (382);
		\draw [style=wire] (387) to (386);
		\draw [style=wire] (389) to (388);
		\draw [style=wire] (383.center) to (386);
		\draw [style=wire] (387) to (381);
		\draw [style=wire] (389) to (385);
		\draw [style=wire] (384.center) to (388);
	\end{pgfonlayer}
\end{tikzpicture}
 \end{array} \substack{=\\(\ref{string:delta-comonoid-map})}  \begin{array}[c]{c}\begin{tikzpicture}
	\begin{pgfonlayer}{nodelayer}
		\node [style=component] (380) at (210, -4) {$\mathsf{u}$};
		\node [style=port] (381) at (209.5, -7.75) {};
		\node [style=component] (382) at (210, -6) {$\Delta$};
		\node [style=none] (383) at (209.5, -6.5) {};
		\node [style=none] (384) at (210.5, -6.5) {};
		\node [style=port] (385) at (210.5, -7.75) {};
		\node [style=component] (386) at (210, -5) {$\delta$};
		\node [style=component] (387) at (209.5, -7) {$\varepsilon$};
		\node [style=component] (389) at (210.5, -7) {$\varepsilon$};
	\end{pgfonlayer}
	\begin{pgfonlayer}{edgelayer}
		\draw [style=wire] (382) to (383.center);
		\draw [style=wire] (382) to (384.center);
		\draw [style=wire] (387) to (381);
		\draw [style=wire] (389) to (385);
		\draw [style=wire] (383.center) to (387);
		\draw [style=wire] (384.center) to (389);
		\draw [style=wire] (380) to (386);
		\draw [style=wire] (386) to (382);
	\end{pgfonlayer}
\end{tikzpicture}
 \end{array} \substack{=\\(\ref{string:monoid-!coalg})}  \begin{array}[c]{c}\begin{tikzpicture}
	\begin{pgfonlayer}{nodelayer}
		\node [style=component] (372) at (210, -4) {$\mathsf{m}$};
		\node [style=component] (373) at (210, -5) {$\mathsf{u}$};
		\node [style=none] (374) at (209.5, -5.5) {};
		\node [style=none] (375) at (210.5, -5.5) {};
		\node [style=none] (376) at (209.5, -4.5) {};
		\node [style=none] (377) at (210.5, -4.5) {};
		\node [style=none] (379) at (210, -4.5) {};
		\node [style=port] (381) at (209.5, -7.75) {};
		\node [style=component] (382) at (210, -6) {$\Delta$};
		\node [style=none] (383) at (209.5, -6.5) {};
		\node [style=none] (384) at (210.5, -6.5) {};
		\node [style=port] (385) at (210.5, -7.75) {};
		\node [style=component] (387) at (209.5, -7) {$\varepsilon$};
		\node [style=component] (389) at (210.5, -7) {$\varepsilon$};
	\end{pgfonlayer}
	\begin{pgfonlayer}{edgelayer}
		\draw [style=wire] (374.center) to (375.center);
		\draw [style=wire] (375.center) to (377.center);
		\draw [style=wire] (377.center) to (376.center);
		\draw [style=wire] (376.center) to (374.center);
		\draw [style=wire] (372) to (379.center);
		\draw [style=wire] (382) to (383.center);
		\draw [style=wire] (382) to (384.center);
		\draw [style=wire] (387) to (381);
		\draw [style=wire] (389) to (385);
		\draw [style=wire] (383.center) to (387);
		\draw [style=wire] (384.center) to (389);
		\draw [style=wire] (373) to (382);
	\end{pgfonlayer}
\end{tikzpicture}
 \end{array} \substack{=\\(\ref{string:comonoid-nat})+(\ref{string:comonad-nat})}  \begin{array}[c]{c}\begin{tikzpicture}
	\begin{pgfonlayer}{nodelayer}
		\node [style=component] (372) at (210, -4) {$\mathsf{m}$};
		\node [style=component] (382) at (210, -4.75) {$\Delta$};
		\node [style=none] (383) at (209.5, -5.25) {};
		\node [style=none] (384) at (210.5, -5.25) {};
		\node [style=component] (387) at (209.5, -5.75) {$\varepsilon$};
		\node [style=component] (389) at (210.5, -5.75) {$\varepsilon$};
		\node [style=port] (390) at (209.5, -7.25) {};
		\node [style=port] (391) at (210.5, -7.25) {};
		\node [style=component] (392) at (209.5, -6.5) {$\mathsf{u}$};
		\node [style=component] (393) at (210.5, -6.5) {$\mathsf{u}$};
	\end{pgfonlayer}
	\begin{pgfonlayer}{edgelayer}
		\draw [style=wire] (382) to (383.center);
		\draw [style=wire] (382) to (384.center);
		\draw [style=wire] (383.center) to (387);
		\draw [style=wire] (384.center) to (389);
		\draw [style=wire] (372) to (382);
		\draw [style=wire] (390) to (392);
		\draw [style=wire] (391) to (393);
	\end{pgfonlayer}
\end{tikzpicture}
 \end{array} \substack{=\\(\ref{string:monoidalcomonoid})}  \begin{array}[c]{c}\begin{tikzpicture}
	\begin{pgfonlayer}{nodelayer}
		\node [style=component] (372) at (209.5, -5) {$\mathsf{m}$};
		\node [style=component] (387) at (209.5, -5.75) {$\varepsilon$};
		\node [style=component] (389) at (210.5, -5.75) {$\varepsilon$};
		\node [style=port] (390) at (209.5, -7.25) {};
		\node [style=port] (391) at (210.5, -7.25) {};
		\node [style=component] (392) at (209.5, -6.5) {$\mathsf{u}$};
		\node [style=component] (393) at (210.5, -6.5) {$\mathsf{u}$};
		\node [style=component] (394) at (210.5, -5) {$\mathsf{m}$};
	\end{pgfonlayer}
	\begin{pgfonlayer}{edgelayer}
		\draw [style=wire] (390) to (392);
		\draw [style=wire] (391) to (393);
		\draw [style=wire] (372) to (387);
		\draw [style=wire] (394) to (389);
	\end{pgfonlayer}
\end{tikzpicture}
 \end{array} \substack{=\\(\ref{string:comonad-monoidal})} \begin{array}[c]{c}\begin{tikzpicture}
	\begin{pgfonlayer}{nodelayer}
		\node [style=port] (128) at (24, 0.5) {};
		\node [style=port] (132) at (25, 0.5) {};
		\node [style=component] (160) at (24, 2) {$\mathsf{u}$};
		\node [style=component] (161) at (25, 2) {$\mathsf{u}$};
	\end{pgfonlayer}
	\begin{pgfonlayer}{edgelayer}
		\draw [style=wire] (128) to (160);
		\draw [style=wire] (132) to (161);
	\end{pgfonlayer}
\end{tikzpicture}
 \end{array}
 \end{gather*}
 So we conclude that $\mathsf{u}$ is a comonoid morphism as desired. Now suppose that $\nabla$ satisfies (\ref{string:monoid-!coalg}), then we will show that it is a comonoid morphism, that is, the left most equalities of (\ref{string:bimonoid}) hold. Starting with the compatibility between $\nabla$ and $\mathsf{e}$, we compute that:
\begin{gather*}
\begin{array}[c]{c}\begin{tikzpicture}
	\begin{pgfonlayer}{nodelayer}
		\node [style=component] (116) at (22.25, -0.25) {$\mathsf{e}$};
		\node [style=port] (117) at (21.75, 2) {};
		\node [style=component] (118) at (22.25, 0.75) {$\nabla$};
		\node [style=none] (119) at (21.75, 1.25) {};
		\node [style=none] (120) at (22.75, 1.25) {};
		\node [style=port] (121) at (22.75, 2) {};
	\end{pgfonlayer}
	\begin{pgfonlayer}{edgelayer}
		\draw [style=wire] (118) to (119.center);
		\draw [style=wire] (118) to (120.center);
		\draw [style=wire] (119.center) to (117);
		\draw [style=wire] (120.center) to (121);
		\draw [style=wire] (116) to (118);
	\end{pgfonlayer}
\end{tikzpicture}
 \end{array}   \substack{=\\(\ref{string:delta-comonoid-map})}\begin{array}[c]{c}\begin{tikzpicture}
	\begin{pgfonlayer}{nodelayer}
		\node [style=component] (388) at (215.25, 0) {$\delta$};
		\node [style=component] (395) at (215.25, -0.75) {$\mathsf{e}$};
		\node [style=port] (396) at (214.75, 2) {};
		\node [style=component] (397) at (215.25, 0.75) {$\nabla$};
		\node [style=none] (398) at (214.75, 1.25) {};
		\node [style=none] (399) at (215.75, 1.25) {};
		\node [style=port] (400) at (215.75, 2) {};
	\end{pgfonlayer}
	\begin{pgfonlayer}{edgelayer}
		\draw [style=wire] (397) to (398.center);
		\draw [style=wire] (397) to (399.center);
		\draw [style=wire] (398.center) to (396);
		\draw [style=wire] (399.center) to (400);
		\draw [style=wire] (397) to (388);
		\draw [style=wire] (388) to (395);
	\end{pgfonlayer}
\end{tikzpicture}
 \end{array} \substack{=\\(\ref{string:comonoid-nat})} \begin{array}[c]{c} \begin{tikzpicture}
	\begin{pgfonlayer}{nodelayer}
		\node [style=port] (405) at (217.25, 2.5) {};
		\node [style=component] (406) at (217.75, 0.75) {$\mathsf{m}$};
		\node [style=none] (407) at (217.25, 1.25) {};
		\node [style=none] (408) at (218.25, 1.25) {};
		\node [style=port] (409) at (218.25, 2.5) {};
		\node [style=component] (420) at (217.25, 1.75) {$\delta$};
		\node [style=component] (421) at (218.25, 1.75) {$\delta$};
		\node [style=component] (443) at (217.75, 0) {$\mathsf{e}$};
	\end{pgfonlayer}
	\begin{pgfonlayer}{edgelayer}
		\draw [style=wire] (406) to (407.center);
		\draw [style=wire] (406) to (408.center);
		\draw [style=wire] (405) to (420);
		\draw [style=wire] (409) to (421);
		\draw [style=wire] (421) to (408.center);
		\draw [style=wire] (420) to (407.center);
		\draw [style=wire] (406) to (443);
	\end{pgfonlayer}
\end{tikzpicture}
 \end{array} \substack{=\\(\ref{string:monoidalcomonoid})} \begin{array}[c]{c} \begin{tikzpicture}
	\begin{pgfonlayer}{nodelayer}
		\node [style=port] (405) at (217.25, 2.5) {};
		\node [style=port] (409) at (218.25, 2.5) {};
		\node [style=component] (420) at (217.25, 1.75) {$\delta$};
		\node [style=component] (421) at (218.25, 1.75) {$\delta$};
		\node [style=component] (443) at (217.25, 0.75) {$\mathsf{e}$};
		\node [style=component] (444) at (218.25, 0.75) {$\mathsf{e}$};
	\end{pgfonlayer}
	\begin{pgfonlayer}{edgelayer}
		\draw [style=wire] (405) to (420);
		\draw [style=wire] (409) to (421);
		\draw [style=wire] (420) to (443);
		\draw [style=wire] (421) to (444);
	\end{pgfonlayer}
\end{tikzpicture}
 \end{array} \substack{=\\(\ref{string:delta-comonoid-map})} \begin{array}[c]{c}\begin{tikzpicture}
	\begin{pgfonlayer}{nodelayer}
		\node [style=port] (128) at (24, 2) {};
		\node [style=port] (132) at (25, 2) {};
		\node [style=component] (160) at (24, 0.5) {$\mathsf{e}$};
		\node [style=component] (161) at (25, 0.5) {$\mathsf{e}$};
	\end{pgfonlayer}
	\begin{pgfonlayer}{edgelayer}
		\draw [style=wire] (128) to (160);
		\draw [style=wire] (132) to (161);
	\end{pgfonlayer}
\end{tikzpicture}
 \end{array}
\end{gather*}
 Next we compute the compatibility between $\nabla$ and $\Delta$: 
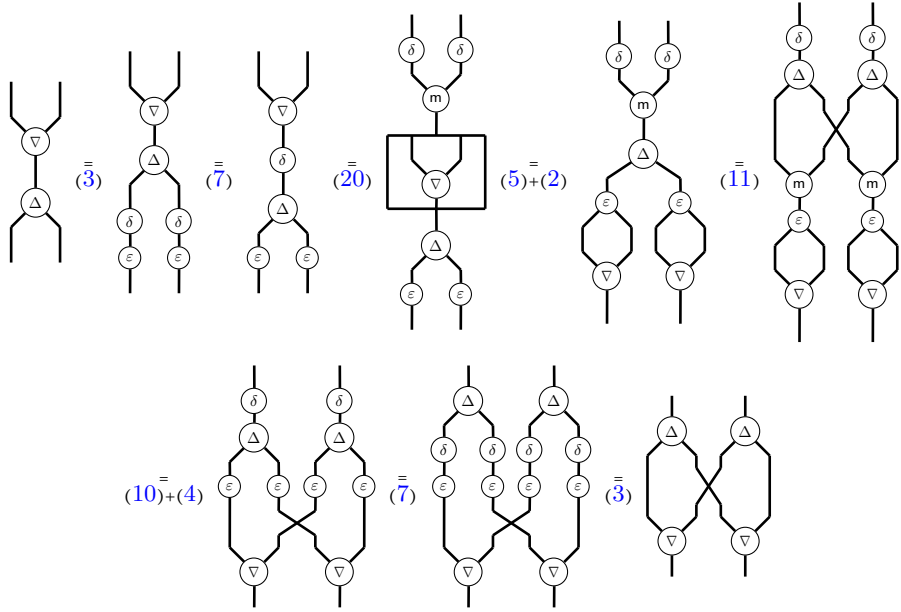
 \begin{gather*}
 \begin{array}[c]{c}\begin{tikzpicture}
	\begin{pgfonlayer}{nodelayer}
		\node [style=port] (29) at (17, -5.5) {};
		\node [style=component] (30) at (17.5, -4.25) {$\Delta$};
		\node [style=none] (31) at (17, -4.75) {};
		\node [style=none] (32) at (18, -4.75) {};
		\node [style=port] (33) at (18, -5.5) {};
		\node [style=port] (162) at (17, -1.75) {};
		\node [style=component] (163) at (17.5, -3) {$\nabla$};
		\node [style=none] (164) at (17, -2.5) {};
		\node [style=none] (165) at (18, -2.5) {};
		\node [style=port] (166) at (18, -1.75) {};
	\end{pgfonlayer}
	\begin{pgfonlayer}{edgelayer}
		\draw [style=wire] (30) to (31.center);
		\draw [style=wire] (30) to (32.center);
		\draw [style=wire] (32.center) to (33);
		\draw [style=wire] (31.center) to (29);
		\draw [style=wire] (163) to (164.center);
		\draw [style=wire] (163) to (165.center);
		\draw [style=wire] (164.center) to (162);
		\draw [style=wire] (165.center) to (166);
		\draw [style=wire] (163) to (30);
	\end{pgfonlayer}
\end{tikzpicture}
 \end{array} \substack{=\\(\ref{string:comonad})}  \begin{array}[c]{c}\begin{tikzpicture}
	\begin{pgfonlayer}{nodelayer}
		\node [style=port] (450) at (221.25, -2.75) {};
		\node [style=component] (451) at (221.75, -4) {$\nabla$};
		\node [style=none] (452) at (221.25, -3.5) {};
		\node [style=none] (453) at (222.25, -3.5) {};
		\node [style=port] (454) at (222.25, -2.75) {};
		\node [style=port] (456) at (221.25, -7.75) {};
		\node [style=component] (457) at (221.75, -5) {$\Delta$};
		\node [style=none] (458) at (221.25, -5.5) {};
		\node [style=none] (459) at (222.25, -5.5) {};
		\node [style=port] (460) at (222.25, -7.75) {};
		\node [style=component] (461) at (221.25, -6.25) {$\delta$};
		\node [style=component] (462) at (221.25, -7) {$\varepsilon$};
		\node [style=component] (463) at (222.25, -6.25) {$\delta$};
		\node [style=component] (464) at (222.25, -7) {$\varepsilon$};
	\end{pgfonlayer}
	\begin{pgfonlayer}{edgelayer}
		\draw [style=wire] (451) to (452.center);
		\draw [style=wire] (451) to (453.center);
		\draw [style=wire] (452.center) to (450);
		\draw [style=wire] (453.center) to (454);
		\draw [style=wire] (457) to (458.center);
		\draw [style=wire] (457) to (459.center);
		\draw [style=wire] (462) to (461);
		\draw [style=wire] (464) to (463);
		\draw [style=wire] (458.center) to (461);
		\draw [style=wire] (462) to (456);
		\draw [style=wire] (464) to (460);
		\draw [style=wire] (459.center) to (463);
		\draw [style=wire] (451) to (457);
	\end{pgfonlayer}
\end{tikzpicture} 
 \end{array} \substack{=\\(\ref{string:delta-comonoid-map})}  \begin{array}[c]{c}\begin{tikzpicture}
	\begin{pgfonlayer}{nodelayer}
		\node [style=port] (450) at (221.25, -2.75) {};
		\node [style=component] (451) at (221.75, -4) {$\nabla$};
		\node [style=none] (452) at (221.25, -3.5) {};
		\node [style=none] (453) at (222.25, -3.5) {};
		\node [style=port] (454) at (222.25, -2.75) {};
		\node [style=port] (466) at (221.25, -7.75) {};
		\node [style=component] (467) at (221.75, -6) {$\Delta$};
		\node [style=none] (468) at (221.25, -6.5) {};
		\node [style=none] (469) at (222.25, -6.5) {};
		\node [style=port] (470) at (222.25, -7.75) {};
		\node [style=component] (471) at (221.75, -5) {$\delta$};
		\node [style=component] (472) at (221.25, -7) {$\varepsilon$};
		\node [style=component] (473) at (222.25, -7) {$\varepsilon$};
	\end{pgfonlayer}
	\begin{pgfonlayer}{edgelayer}
		\draw [style=wire] (451) to (452.center);
		\draw [style=wire] (451) to (453.center);
		\draw [style=wire] (452.center) to (450);
		\draw [style=wire] (453.center) to (454);
		\draw [style=wire] (467) to (468.center);
		\draw [style=wire] (467) to (469.center);
		\draw [style=wire] (472) to (466);
		\draw [style=wire] (473) to (470);
		\draw [style=wire] (468.center) to (472);
		\draw [style=wire] (469.center) to (473);
		\draw [style=wire] (471) to (467);
		\draw [style=wire] (451) to (471);
	\end{pgfonlayer}
\end{tikzpicture}
 \end{array}%
 
 \end{gather*}
 So we conclude that $\nabla$ is a comonoid morphism as desired.
 \end{proof}
 
\section{Additive Enrichment}\label{sec:add-en}

In this section we prove the main objective of this paper: that a monoidal bialgebra modality with a pre-codereliction induces an additive enrichment. In this paper, following \cite{blute2006differential,Blute2019}, by an \emph{additive} category, we mean a category enriched over commutative monoids, which differs from the use of the name in other category theory literature. Explicitly, an \textbf{additive category} \cite[Def 3]{Blute2019} is a category $\mathbb{X}$ such that each homset $\mathbb{X}(A,B)$ is a commutative monoid, with binary operation $+: \mathbb{X}(A,B) \times \mathbb{X}(A,B) \to \mathbb{X}(A,B)$ and zero element $0 \in \mathbb{X}(A,B)$, and such that composition preserves the additive structure, that is: 
 \begin{equation}\begin{gathered} \label{eq:add}
 f;(g+h) = (f;g)+(f;h) \qquad f;0 = 0 \qquad (g+h);k = (g;k)+(h;k) \qquad 0;k = 0
 \end{gathered}\end{equation}
 Then an \textbf{additive symmetric monoidal category} \cite[Def 3]{Blute2019} is a symmetric monoidal category $\mathbb{X}$ such that $\mathbb{X}$ is also an additive category and the monoidal product also preserves the additive structure, that is:
  \begin{equation}\begin{gathered}\label{eq:add-tensor}
 f \otimes (g+h) = (f \otimes g)+(f \otimes h) \qquad f \otimes 0 = 0 \\
 (g+h) \otimes k = (g \otimes k)+(h \otimes k) \qquad 0 \otimes k = 0
 \end{gathered}\end{equation}

\begin{thm}\label{thm:additive} Let $\mathbb{X}$ be a symmetric monoidal category with a monoidal bialgebra modality $(\oc, \delta, \varepsilon,  \Delta, \mathsf{e}, \mathsf{m}, \mathsf{m}_I, \nabla, \mathsf{u})$ which comes equipped with a pre-codereliction $\eta$. Then $\mathbb{X}$ is an additive symmetric monoidal category where:
\begin{enumerate}[{\em (i)}]
\item For parallel maps $f: A \to B$ and $g: A \to B$, their sum $f+g: A \to B$ is defined as the following composite: 
 \begin{equation}\begin{gathered}
  \xymatrixrowsep{1.75pc}\xymatrixcolsep{3pc}\xymatrix{A  \ar[r]^-{\eta_A} & \oc(A) \ar[r]^-{\Delta_A} & \oc(A) \otimes \oc(A) \ar[r]^-{\oc(f) \otimes \oc(g)} & \oc(B) \otimes \oc(B) \ar[r]^-{\nabla_B} & \oc(B) \ar[r]^-{\varepsilon_B} & B    } 
\end{gathered}\end{equation}
which is drawn as follows: 
\begin{align}\label{string:+def}
 \begin{array}[c]{c}\begin{tikzpicture}
	\begin{pgfonlayer}{nodelayer}
		\node [style=port] (427) at (86.5, 1.25) {};
		\node [style=port] (428) at (86.5, -0.75) {};
		\node [style=component] (429) at (86.5, 0.25) {$f+g$};
	\end{pgfonlayer}
	\begin{pgfonlayer}{edgelayer}
		\draw [style=wire] (427) to (429);
		\draw [style=wire] (429) to (428);
	\end{pgfonlayer}
\end{tikzpicture}
 \end{array} = \begin{array}[c]{c} \begin{tikzpicture}
	\begin{pgfonlayer}{nodelayer}
		\node [style=component] (470) at (91.75, -7) {$\nabla$};
		\node [style=none] (471) at (91, -6.5) {};
		\node [style=none] (472) at (92.5, -6.5) {};
		\node [style=component] (473) at (91, -5.75) {$f$};
		\node [style=component] (474) at (92.5, -5.75) {$g$};
		\node [style=component] (475) at (91.75, -4.5) {$\Delta$};
		\node [style=none] (476) at (91, -5) {};
		\node [style=none] (477) at (92.5, -5) {};
		\node [style=component] (478) at (91.75, -3.5) {$\eta$};
		\node [style=component] (479) at (91.75, -8) {$\varepsilon$};
		\node [style=port] (480) at (91.75, -2.75) {};
		\node [style=port] (481) at (91.75, -8.75) {};
		\node [style=none] (482) at (90.5, -6.25) {};
		\node [style=none] (483) at (91.5, -6.25) {};
		\node [style=none] (484) at (90.5, -5.25) {};
		\node [style=none] (485) at (91.5, -5.25) {};
		\node [style=none] (486) at (92, -6.25) {};
		\node [style=none] (487) at (93, -6.25) {};
		\node [style=none] (488) at (92, -5.25) {};
		\node [style=none] (489) at (93, -5.25) {};
	\end{pgfonlayer}
	\begin{pgfonlayer}{edgelayer}
		\draw [style=wire] (470) to (471.center);
		\draw [style=wire] (470) to (472.center);
		\draw [style=wire] (474) to (472.center);
		\draw [style=wire] (473) to (471.center);
		\draw [style=wire] (475) to (476.center);
		\draw [style=wire] (475) to (477.center);
		\draw [style=wire] (476.center) to (473);
		\draw [style=wire] (477.center) to (474);
		\draw [style=wire] (478) to (475);
		\draw [style=wire] (470) to (479);
		\draw [style=wire] (480) to (478);
		\draw [style=wire] (479) to (481);
		\draw [style=wire] (482.center) to (483.center);
		\draw [style=wire] (483.center) to (485.center);
		\draw [style=wire] (485.center) to (484.center);
		\draw [style=wire] (484.center) to (482.center);
		\draw [style=wire] (486.center) to (487.center);
		\draw [style=wire] (487.center) to (489.center);
		\draw [style=wire] (489.center) to (488.center);
		\draw [style=wire] (488.center) to (486.center);
	\end{pgfonlayer}
\end{tikzpicture}
 \end{array}  
\end{align}
\item The zero map $0: A \to B$ is defined as the following composite:  
 \begin{equation}\begin{gathered}
  \xymatrixrowsep{1.75pc}\xymatrixcolsep{5pc}\xymatrix{A  \ar[r]^-{\eta_A} & \oc(A) \ar[r]^-{\mathsf{e}_A} & I \ar[r]^-{\mathsf{u}_B} & \oc(B) \ar[r]^-{\varepsilon_B} & B    } 
\end{gathered}\end{equation}
which is drawn as follows: 
\begin{align}\label{string:0def}
 \begin{array}[c]{c}\begin{tikzpicture}
	\begin{pgfonlayer}{nodelayer}
		\node [style=port] (427) at (86.5, 1.25) {};
		\node [style=port] (428) at (86.5, -0.75) {};
		\node [style=component] (429) at (86.5, 0.25) {$0$};
	\end{pgfonlayer}
	\begin{pgfonlayer}{edgelayer}
		\draw [style=wire] (427) to (429);
		\draw [style=wire] (429) to (428);
	\end{pgfonlayer}
\end{tikzpicture}
 \end{array} = \begin{array}[c]{c}\begin{tikzpicture}
	\begin{pgfonlayer}{nodelayer}
		\node [style=component] (452) at (93.5, 1.5) {$\mathsf{e}$};
		\node [style=component] (453) at (93.5, 0.5) {$\mathsf{u}$};
		\node [style=component] (462) at (93.5, 2.5) {$\eta$};
		\node [style=component] (463) at (93.5, -0.5) {$\varepsilon$};
		\node [style=port] (464) at (93.5, 3.25) {};
		\node [style=port] (465) at (93.5, -1.25) {};
	\end{pgfonlayer}
	\begin{pgfonlayer}{edgelayer}
		\draw [style=wire] (464) to (462);
		\draw [style=wire] (463) to (465);
		\draw [style=wire] (462) to (452);
		\draw [style=wire] (453) to (463);
	\end{pgfonlayer}
\end{tikzpicture}
 \end{array}  
\end{align}
 \end{enumerate}
\end{thm}
\begin{proof} The keen-eyed reader will observe that the formulas for sum and the zero are given by taking \textbf{bialgebra convolution} \cite[Pg. 72]{sweedler1969hopf}, then pre-composing with the pre-codereliction and post-composing with the dereliction. However, while it is well known that bialgebra convolution does indeed give a commutative monoid structure on $\mathbb{X}(\oc(A), \oc(B))$ \cite[Pg. 70]{sweedler1969hopf}, in order to show that this also gives a commutative structure on $\mathbb{X}(A,B)$ we will first need to prove that the following equalities holds: 
\begin{align}\label{string:!+0}
\begin{array}[c]{c}\begin{tikzpicture}
	\begin{pgfonlayer}{nodelayer}
		\node [style=port] (427) at (86.5, 2) {};
		\node [style=port] (428) at (86.5, -1.5) {};
		\node [style=component] (429) at (86.5, 0.25) {$f+g$};
		\node [style=none] (466) at (85.75, 1) {};
		\node [style=none] (467) at (87.25, 1) {};
		\node [style=none] (468) at (85.75, -0.5) {};
		\node [style=none] (469) at (87.25, -0.5) {};
	\end{pgfonlayer}
	\begin{pgfonlayer}{edgelayer}
		\draw [style=wire] (427) to (429);
		\draw [style=wire] (429) to (428);
		\draw [style=wire] (466.center) to (467.center);
		\draw [style=wire] (467.center) to (469.center);
		\draw [style=wire] (469.center) to (468.center);
		\draw [style=wire] (468.center) to (466.center);
	\end{pgfonlayer}
\end{tikzpicture}
 \end{array} = \begin{array}[c]{c}\begin{tikzpicture}
	\begin{pgfonlayer}{nodelayer}
		\node [style=component] (430) at (89, -1) {$\nabla$};
		\node [style=none] (431) at (88.25, -0.5) {};
		\node [style=none] (432) at (89.75, -0.5) {};
		\node [style=component] (433) at (88.25, 0.25) {$f$};
		\node [style=component] (434) at (89.75, 0.25) {$g$};
		\node [style=component] (435) at (89, 1.5) {$\Delta$};
		\node [style=none] (436) at (88.25, 1) {};
		\node [style=none] (437) at (89.75, 1) {};
		\node [style=port] (440) at (89, 2.5) {};
		\node [style=port] (441) at (89, -2) {};
		\node [style=none] (443) at (87.75, -0.25) {};
		\node [style=none] (444) at (88.75, -0.25) {};
		\node [style=none] (445) at (87.75, 0.75) {};
		\node [style=none] (446) at (88.75, 0.75) {};
		\node [style=none] (448) at (89.25, -0.25) {};
		\node [style=none] (449) at (90.25, -0.25) {};
		\node [style=none] (450) at (89.25, 0.75) {};
		\node [style=none] (451) at (90.25, 0.75) {};
	\end{pgfonlayer}
	\begin{pgfonlayer}{edgelayer}
		\draw [style=wire] (430) to (431.center);
		\draw [style=wire] (430) to (432.center);
		\draw [style=wire] (434) to (432.center);
		\draw [style=wire] (433) to (431.center);
		\draw [style=wire] (435) to (436.center);
		\draw [style=wire] (435) to (437.center);
		\draw [style=wire] (436.center) to (433);
		\draw [style=wire] (437.center) to (434);
		\draw [style=wire] (443.center) to (444.center);
		\draw [style=wire] (444.center) to (446.center);
		\draw [style=wire] (446.center) to (445.center);
		\draw [style=wire] (445.center) to (443.center);
		\draw [style=wire] (448.center) to (449.center);
		\draw [style=wire] (449.center) to (451.center);
		\draw [style=wire] (451.center) to (450.center);
		\draw [style=wire] (450.center) to (448.center);
		\draw [style=wire] (440) to (435);
		\draw [style=wire] (441) to (430);
	\end{pgfonlayer}
\end{tikzpicture}
 \end{array}  &&  \begin{array}[c]{c}\begin{tikzpicture}
	\begin{pgfonlayer}{nodelayer}
		\node [style=port] (427) at (86.5, 1.5) {};
		\node [style=port] (428) at (86.5, -1) {};
		\node [style=component] (429) at (86.5, 0.25) {$0$};
		\node [style=none] (466) at (86, 0.75) {};
		\node [style=none] (467) at (87, 0.75) {};
		\node [style=none] (468) at (86, -0.25) {};
		\node [style=none] (469) at (87, -0.25) {};
	\end{pgfonlayer}
	\begin{pgfonlayer}{edgelayer}
		\draw [style=wire] (427) to (429);
		\draw [style=wire] (429) to (428);
		\draw [style=wire] (466.center) to (467.center);
		\draw [style=wire] (467.center) to (469.center);
		\draw [style=wire] (469.center) to (468.center);
		\draw [style=wire] (468.center) to (466.center);
	\end{pgfonlayer}
\end{tikzpicture}
 \end{array} = \begin{array}[c]{c}\begin{tikzpicture}
	\begin{pgfonlayer}{nodelayer}
		\node [style=component] (215) at (37.25, -5.25) {$\mathsf{e}$};
		\node [style=port] (216) at (37.25, -4.5) {};
		\node [style=component] (217) at (37.25, -6.25) {$\mathsf{u}$};
		\node [style=port] (218) at (37.25, -7) {};
	\end{pgfonlayer}
	\begin{pgfonlayer}{edgelayer}
		\draw [style=wire] (216) to (215);
		\draw [style=wire] (218) to (217);
	\end{pgfonlayer}
\end{tikzpicture}
 \end{array} 
\end{align}
For the equality on the left, we compute: 
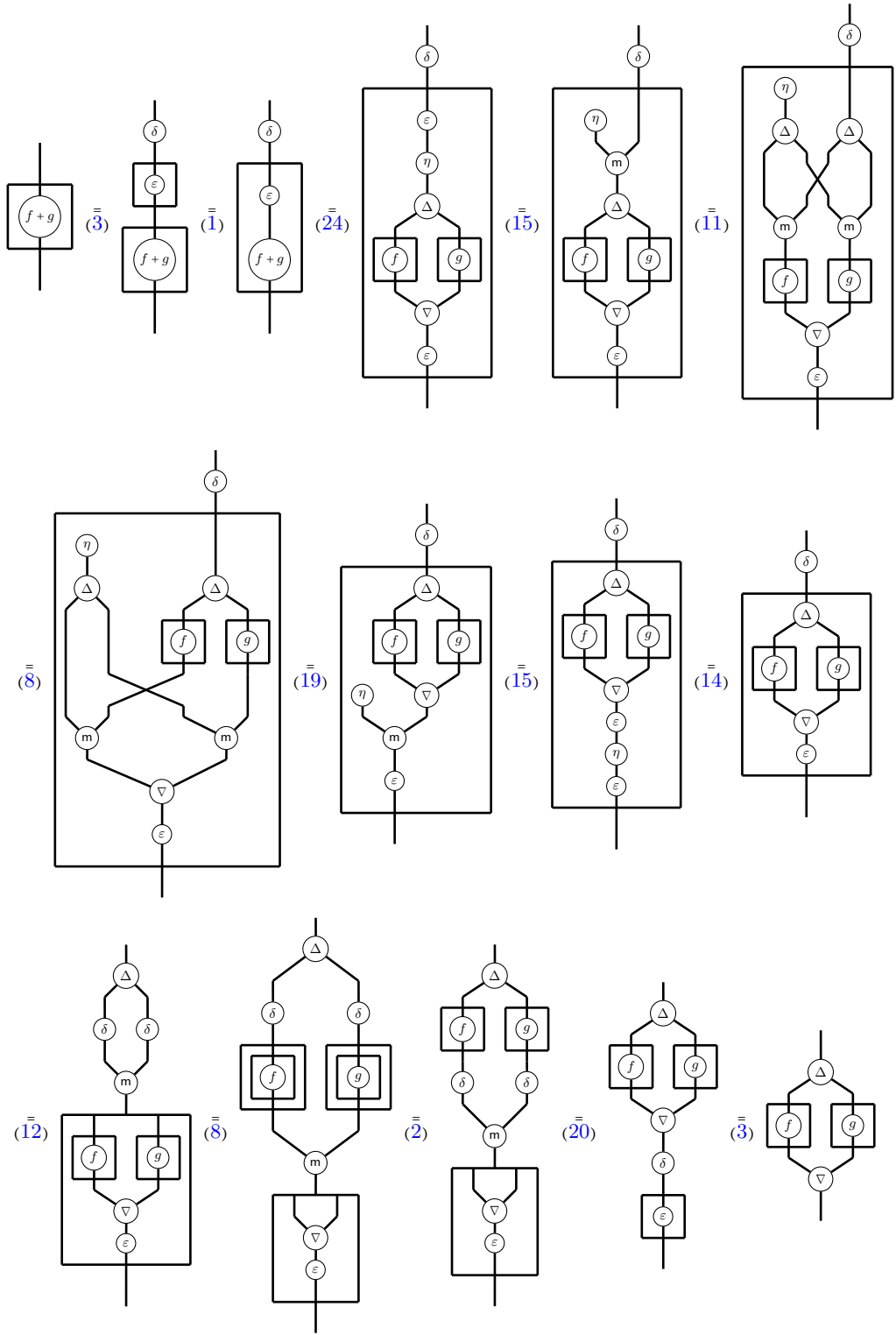
\begin{gather*}
 \begin{array}[c]{c}\begin{tikzpicture}
	\begin{pgfonlayer}{nodelayer}
		\node [style=port] (427) at (86.5, 2) {};
		\node [style=port] (428) at (86.5, -1.5) {};
		\node [style=component] (429) at (86.5, 0.25) {$f+g$};
		\node [style=none] (466) at (85.75, 1) {};
		\node [style=none] (467) at (87.25, 1) {};
		\node [style=none] (468) at (85.75, -0.5) {};
		\node [style=none] (469) at (87.25, -0.5) {};
	\end{pgfonlayer}
	\begin{pgfonlayer}{edgelayer}
		\draw [style=wire] (427) to (429);
		\draw [style=wire] (429) to (428);
		\draw [style=wire] (466.center) to (467.center);
		\draw [style=wire] (467.center) to (469.center);
		\draw [style=wire] (469.center) to (468.center);
		\draw [style=wire] (468.center) to (466.center);
	\end{pgfonlayer}
\end{tikzpicture}
 \end{array} \substack{=\\(\ref{string:comonad})}  \begin{array}[c]{c}\begin{tikzpicture}
	\begin{pgfonlayer}{nodelayer}
		\node [style=port] (428) at (101.25, -5) {};
		\node [style=component] (429) at (101.25, -3.25) {$f+g$};
		\node [style=none] (466) at (100.5, -2.5) {};
		\node [style=none] (467) at (102, -2.5) {};
		\node [style=none] (468) at (100.5, -4) {};
		\node [style=none] (469) at (102, -4) {};
		\node [style=port] (494) at (101.25, 0.5) {};
		\node [style=component] (495) at (101.25, -0.25) {$\delta$};
		\node [style=component] (497) at (101.25, -1.5) {$\varepsilon$};
		\node [style=none] (498) at (100.75, -1) {};
		\node [style=none] (499) at (101.75, -1) {};
		\node [style=none] (500) at (100.75, -2) {};
		\node [style=none] (501) at (101.75, -2) {};
	\end{pgfonlayer}
	\begin{pgfonlayer}{edgelayer}
		\draw [style=wire] (429) to (428);
		\draw [style=wire] (466.center) to (467.center);
		\draw [style=wire] (467.center) to (469.center);
		\draw [style=wire] (469.center) to (468.center);
		\draw [style=wire] (468.center) to (466.center);
		\draw [style=wire] (495) to (494);
		\draw [style=wire] (498.center) to (499.center);
		\draw [style=wire] (499.center) to (501.center);
		\draw [style=wire] (501.center) to (500.center);
		\draw [style=wire] (500.center) to (498.center);
		\draw [style=wire] (495) to (497);
		\draw [style=wire] (497) to (429);
	\end{pgfonlayer}
\end{tikzpicture}
 \end{array} \substack{=\\(\ref{string:functor})}  \begin{array}[c]{c}\begin{tikzpicture}
	\begin{pgfonlayer}{nodelayer}
		\node [style=port] (428) at (101.25, -5) {};
		\node [style=component] (429) at (101.25, -3.25) {$f+g$};
		\node [style=none] (468) at (100.5, -4) {};
		\node [style=none] (469) at (102, -4) {};
		\node [style=port] (494) at (101.25, 0.5) {};
		\node [style=component] (495) at (101.25, -0.25) {$\delta$};
		\node [style=component] (497) at (101.25, -1.75) {$\varepsilon$};
		\node [style=none] (498) at (100.5, -1) {};
		\node [style=none] (499) at (102, -1) {};
	\end{pgfonlayer}
	\begin{pgfonlayer}{edgelayer}
		\draw [style=wire] (429) to (428);
		\draw [style=wire] (469.center) to (468.center);
		\draw [style=wire] (495) to (494);
		\draw [style=wire] (498.center) to (499.center);
		\draw [style=wire] (495) to (497);
		\draw [style=wire] (497) to (429);
		\draw [style=wire] (498.center) to (468.center);
		\draw [style=wire] (499.center) to (469.center);
	\end{pgfonlayer}
\end{tikzpicture}
 \end{array} \substack{=\\(\ref{string:comonad})} \begin{array}[c]{c}\begin{tikzpicture}
	\begin{pgfonlayer}{nodelayer}
		\node [style=component] (430) at (89, -1) {$\nabla$};
		\node [style=none] (431) at (88.25, -0.5) {};
		\node [style=none] (432) at (89.75, -0.5) {};
		\node [style=component] (433) at (88.25, 0.25) {$f$};
		\node [style=component] (434) at (89.75, 0.25) {$g$};
		\node [style=component] (435) at (89, 1.5) {$\Delta$};
		\node [style=none] (436) at (88.25, 1) {};
		\node [style=none] (437) at (89.75, 1) {};
		\node [style=port] (440) at (89, 2.5) {};
		\node [style=port] (441) at (89, -2) {};
		\node [style=none] (443) at (87.75, -0.25) {};
		\node [style=none] (444) at (88.75, -0.25) {};
		\node [style=none] (445) at (87.75, 0.75) {};
		\node [style=none] (446) at (88.75, 0.75) {};
		\node [style=none] (448) at (89.25, -0.25) {};
		\node [style=none] (449) at (90.25, -0.25) {};
		\node [style=none] (450) at (89.25, 0.75) {};
		\node [style=none] (451) at (90.25, 0.75) {};
	\end{pgfonlayer}
	\begin{pgfonlayer}{edgelayer}
		\draw [style=wire] (430) to (431.center);
		\draw [style=wire] (430) to (432.center);
		\draw [style=wire] (434) to (432.center);
		\draw [style=wire] (433) to (431.center);
		\draw [style=wire] (435) to (436.center);
		\draw [style=wire] (435) to (437.center);
		\draw [style=wire] (436.center) to (433);
		\draw [style=wire] (437.center) to (434);
		\draw [style=wire] (443.center) to (444.center);
		\draw [style=wire] (444.center) to (446.center);
		\draw [style=wire] (446.center) to (445.center);
		\draw [style=wire] (445.center) to (443.center);
		\draw [style=wire] (448.center) to (449.center);
		\draw [style=wire] (449.center) to (451.center);
		\draw [style=wire] (451.center) to (450.center);
		\draw [style=wire] (450.center) to (448.center);
		\draw [style=wire] (440) to (435);
		\draw [style=wire] (441) to (430);
	\end{pgfonlayer}
\end{tikzpicture}
 \end{array}
\end{gather*}
For the equality on the right, we compute:
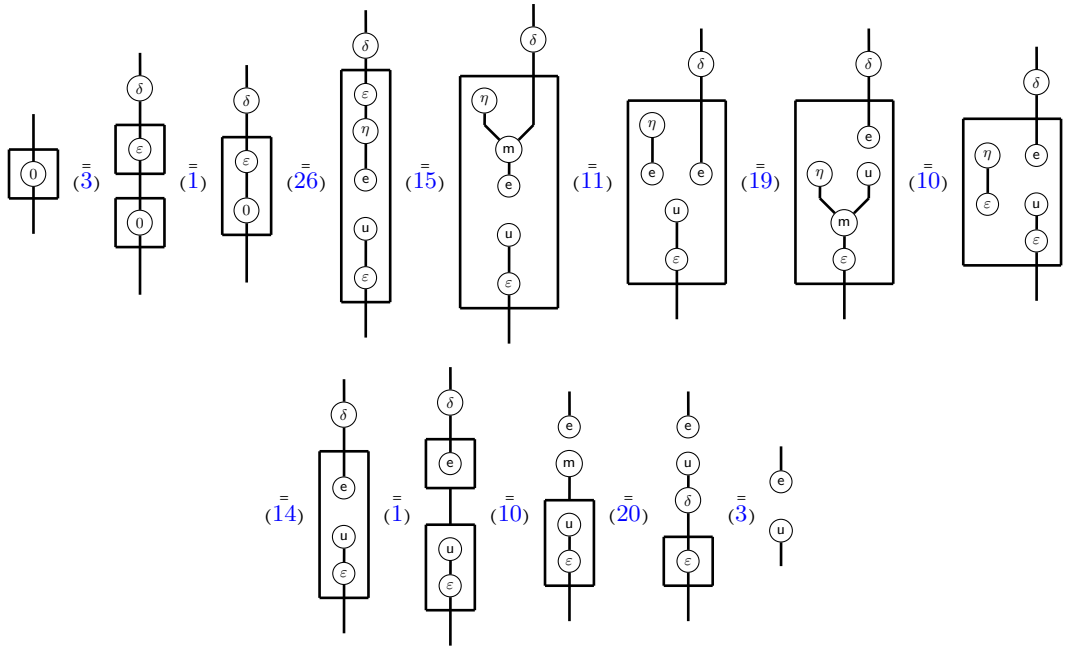
\begin{gather*}
\begin{array}[c]{c}\begin{tikzpicture}
	\begin{pgfonlayer}{nodelayer}
		\node [style=port] (427) at (86.5, 1.5) {};
		\node [style=port] (428) at (86.5, -1) {};
		\node [style=component] (429) at (86.5, 0.25) {$0$};
		\node [style=none] (466) at (86, 0.75) {};
		\node [style=none] (467) at (87, 0.75) {};
		\node [style=none] (468) at (86, -0.25) {};
		\node [style=none] (469) at (87, -0.25) {};
	\end{pgfonlayer}
	\begin{pgfonlayer}{edgelayer}
		\draw [style=wire] (427) to (429);
		\draw [style=wire] (429) to (428);
		\draw [style=wire] (466.center) to (467.center);
		\draw [style=wire] (467.center) to (469.center);
		\draw [style=wire] (469.center) to (468.center);
		\draw [style=wire] (468.center) to (466.center);
	\end{pgfonlayer}
\end{tikzpicture}
 \end{array}  \substack{=\\(\ref{string:comonad})}  \begin{array}[c]{c}\begin{tikzpicture}
	\begin{pgfonlayer}{nodelayer}
		\node [style=port] (686) at (121.75, -4.5) {};
		\node [style=component] (687) at (121.75, -3) {$0$};
		\node [style=none] (688) at (121.25, -2.5) {};
		\node [style=none] (689) at (122.25, -2.5) {};
		\node [style=none] (690) at (121.25, -3.5) {};
		\node [style=none] (691) at (122.25, -3.5) {};
		\node [style=port] (692) at (121.75, 0.5) {};
		\node [style=component] (693) at (121.75, -0.25) {$\delta$};
		\node [style=component] (694) at (121.75, -1.5) {$\varepsilon$};
		\node [style=none] (695) at (121.25, -1) {};
		\node [style=none] (696) at (122.25, -1) {};
		\node [style=none] (697) at (121.25, -2) {};
		\node [style=none] (698) at (122.25, -2) {};
	\end{pgfonlayer}
	\begin{pgfonlayer}{edgelayer}
		\draw [style=wire] (687) to (686);
		\draw [style=wire] (688.center) to (689.center);
		\draw [style=wire] (689.center) to (691.center);
		\draw [style=wire] (691.center) to (690.center);
		\draw [style=wire] (690.center) to (688.center);
		\draw [style=wire] (693) to (692);
		\draw [style=wire] (695.center) to (696.center);
		\draw [style=wire] (696.center) to (698.center);
		\draw [style=wire] (698.center) to (697.center);
		\draw [style=wire] (697.center) to (695.center);
		\draw [style=wire] (693) to (694);
		\draw [style=wire] (694) to (687);
	\end{pgfonlayer}
\end{tikzpicture} \end{array}%
 
\end{gather*} 
Now that we have these identities, we may prove that the homsets are indeed commutative monoids. So we need to check that $+$ is commutative and associative, and that $0$ is a unit for $+$ (once we have shown commutativity, it suffices to show that $f+0=f$, as $0+f=f$ will follow). 

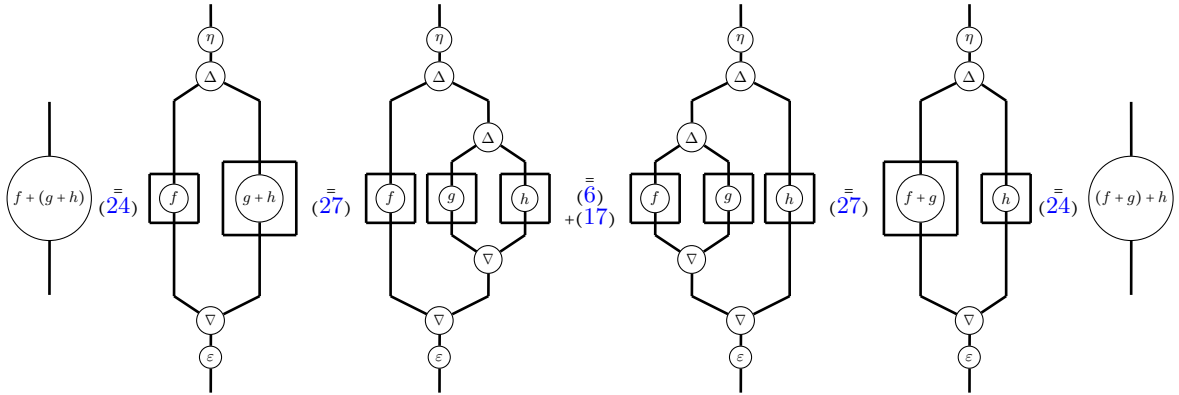
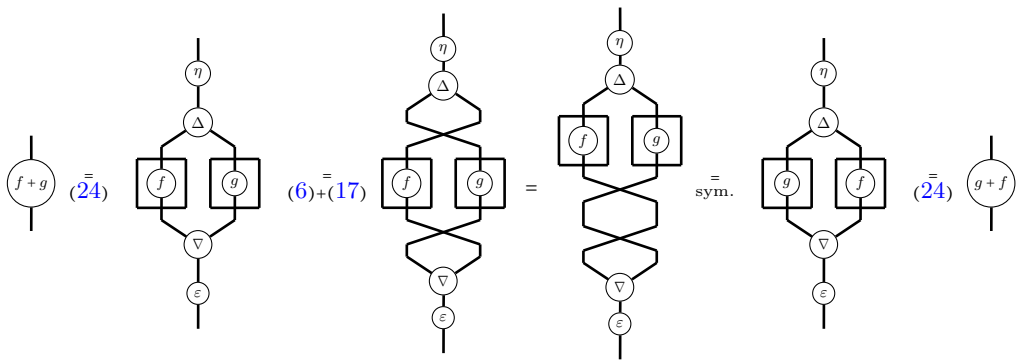
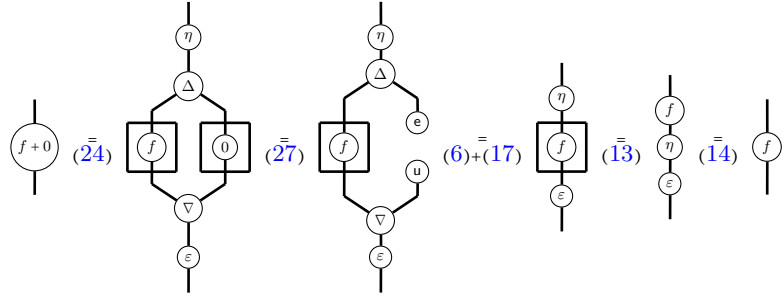
\begin{enumerate}
    \item Associativity: 
    \begin{gather*}
\begin{array}[c]{c}\begin{tikzpicture}
	\begin{pgfonlayer}{nodelayer}
		\node [style=port] (728) at (130, 1) {};
		\node [style=port] (729) at (130, -3) {};
		\node [style=component] (730) at (130, -1) {$f+(g+h)$};
	\end{pgfonlayer}
	\begin{pgfonlayer}{edgelayer}
		\draw [style=wire] (730) to (729);
		\draw [style=wire] (728) to (730);
	\end{pgfonlayer}
\end{tikzpicture}
 \end{array}  \!\!\substack{=\\(\ref{string:+def})} \!\begin{array}[c]{c}\begin{tikzpicture}
	\begin{pgfonlayer}{nodelayer}
		\node [style=port] (728) at (130, 1) {};
		\node [style=port] (729) at (130, -3) {};
		\node [style=component] (730) at (130, -1) {$(f+g)+h$};
	\end{pgfonlayer}
	\begin{pgfonlayer}{edgelayer}
		\draw [style=wire] (730) to (729);
		\draw [style=wire] (728) to (730);
	\end{pgfonlayer}
\end{tikzpicture}
 \end{array}%
 
    \end{gather*}
\end{enumerate}
Next we show that composition is compatible with the additive structure. To do so, we will in fact prove directly that $k;(f+g);h = (k;f;h)+(k;g;h)$ and $f;0;g=0$, which is of course equivalent to (\ref{eq:add}). Starting with the sum, we compute that: 
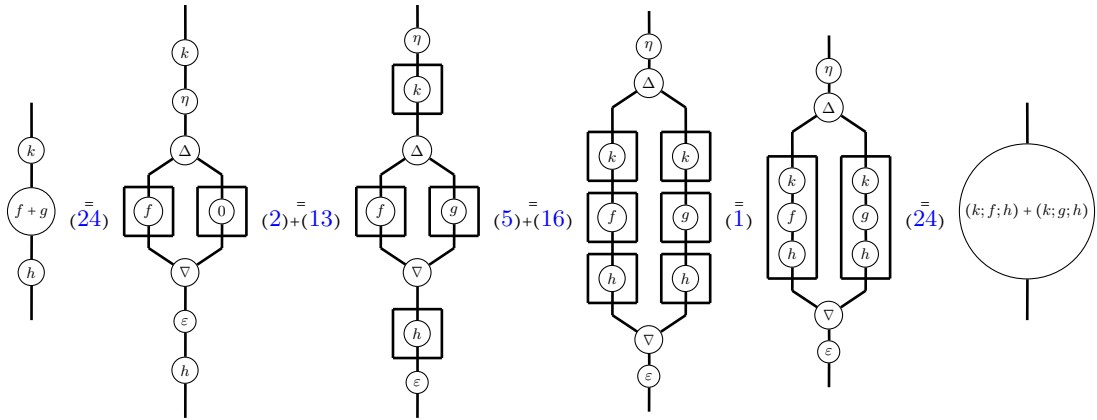
\begin{gather*}
%
  
\end{gather*}
Next we show that composition also preserves the zero: 
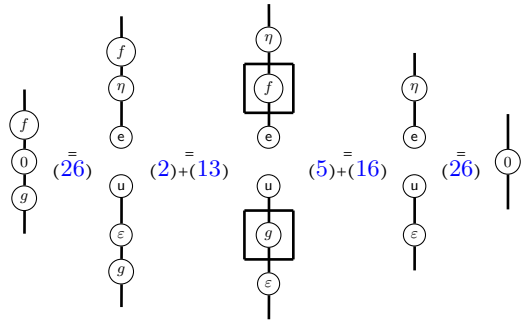
\begin{gather*}
%
\end{gather*}
Lastly, we prove that the monoidal product preserves the additive structure as well. Starting with the sum, we compute that:  
\begin{gather*}  %
 
\end{gather*}
Next we show that the monoidal product also preserves the zero: 
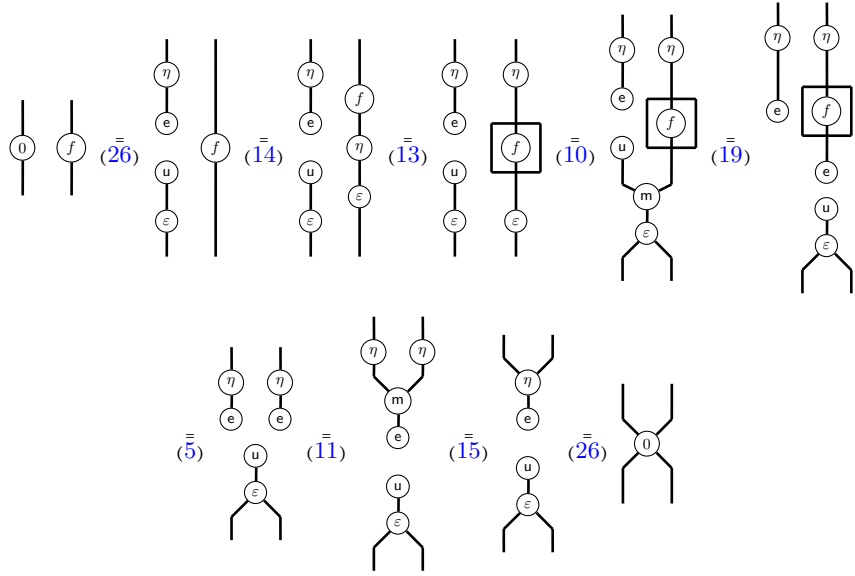
\begin{gather*} %
\end{gather*}
Similarly, we can also show that $k \otimes  (f+g) = (k \otimes f) + (k \otimes g)$ and $0 \otimes f = 0$. So we conclude that we have an additive symmetric monoidal category. 
\end{proof}

\section{Additive Bialgebra Modalities Revisited}\label{sec:add-bialg-mon}

One of the fundamental results of \cite{Blute2019} was showing that for an additive symmetric monoidal category, there was a bijective correspondence between additive bialgebra modalities and monoidal coalgebra modalities \cite[Thm 1]{Blute2019}. In this section, we will show that we can in fact add monoidal bialgebra modalities to this list as well. 

So throughout this section, we now work in an additive symmetric monoidal category $\mathbb{X}$. For convenience, to help draw addition in our string diagrams, following \cite{Blute2019}, we will write sums of maps as follows: 
\begin{align}\label{string:+def2}
 \begin{array}[c]{c}\begin{tikzpicture}
	\begin{pgfonlayer}{nodelayer}
		\node [style=port] (427) at (86.5, 1.25) {};
		\node [style=port] (428) at (86.5, -0.75) {};
		\node [style=component] (429) at (86.5, 0.25) {$f+g$};
	\end{pgfonlayer}
	\begin{pgfonlayer}{edgelayer}
		\draw [style=wire] (427) to (429);
		\draw [style=wire] (429) to (428);
	\end{pgfonlayer}
\end{tikzpicture}
 \end{array} =  \begin{array}[c]{c}\begin{tikzpicture}
	\begin{pgfonlayer}{nodelayer}
		\node [style=port] (427) at (86.5, 1.25) {};
		\node [style=port] (428) at (86.5, -0.75) {};
		\node [style=component] (429) at (86.5, 0.25) {$f$};
	\end{pgfonlayer}
	\begin{pgfonlayer}{edgelayer}
		\draw [style=wire] (427) to (429);
		\draw [style=wire] (429) to (428);
	\end{pgfonlayer}
\end{tikzpicture}
 \end{array} +  \begin{array}[c]{c}\begin{tikzpicture}
	\begin{pgfonlayer}{nodelayer}
		\node [style=port] (427) at (86.5, 1.25) {};
		\node [style=port] (428) at (86.5, -0.75) {};
		\node [style=component] (429) at (86.5, 0.25) {$g$};
	\end{pgfonlayer}
	\begin{pgfonlayer}{edgelayer}
		\draw [style=wire] (427) to (429);
		\draw [style=wire] (429) to (428);
	\end{pgfonlayer}
\end{tikzpicture}
 \end{array}
\end{align}

A \textbf{pre-additive bialgebra modality} \cite[Def 4]{Blute2019} (which as mentioned above is our new name for what was a bialgebra modality in \cite{blute2006differential,Blute2019}) on  an additive symmetric monoidal category is a bialgebra modality $(\oc, \delta, \varepsilon, \Delta, \mathsf{e}, \nabla, \mathsf{u})$ such that the diagrams from Appendix \ref{sec:diagaddbialgmod} commute. The axioms of a pre-additive bialgebra modality state the compatibility between $\varepsilon$ and the monoid structure (\ref{diag:ep-monoid}), which is drawn out as follows: 
\begin{align}\label{string:ep-monoid}
 \begin{array}[c]{c}\begin{tikzpicture}
	\begin{pgfonlayer}{nodelayer}
		\node [style=component] (0) at (22.25, -0.25) {$\varepsilon$};
		\node [style=port] (1) at (21.75, 2) {};
		\node [style=component] (2) at (22.25, 0.75) {$\nabla$};
		\node [style=none] (3) at (21.75, 1.25) {};
		\node [style=none] (4) at (22.75, 1.25) {};
		\node [style=port] (5) at (22.75, 2) {};
		\node [style=port] (6) at (22.25, -1) {};
	\end{pgfonlayer}
	\begin{pgfonlayer}{edgelayer}
		\draw [style=wire] (2) to (3.center);
		\draw [style=wire] (2) to (4.center);
		\draw [style=wire] (3.center) to (1);
		\draw [style=wire] (4.center) to (5);
		\draw [style=wire] (0) to (2);
		\draw [style=wire] (0) to (6);
	\end{pgfonlayer}
\end{tikzpicture}
 \end{array} = \begin{array}[c]{c}\begin{tikzpicture}
	\begin{pgfonlayer}{nodelayer}
		\node [style=port] (7) at (24, 2) {};
		\node [style=port] (8) at (25, 2) {};
		\node [style=component] (9) at (24, 0.5) {$\mathsf{e}$};
		\node [style=component] (10) at (25, 0.5) {$\varepsilon$};
		\node [style=port] (11) at (25, -1) {};
	\end{pgfonlayer}
	\begin{pgfonlayer}{edgelayer}
		\draw [style=wire] (7) to (9);
		\draw [style=wire] (8) to (10);
		\draw [style=wire] (10) to (11);
	\end{pgfonlayer}
\end{tikzpicture}
 \end{array} + \begin{array}[c]{c}\begin{tikzpicture}
	\begin{pgfonlayer}{nodelayer}
		\node [style=port] (7) at (25, 2) {};
		\node [style=port] (8) at (24, 2) {};
		\node [style=component] (9) at (25, 0.5) {$\mathsf{e}$};
		\node [style=component] (10) at (24, 0.5) {$\varepsilon$};
		\node [style=port] (11) at (24, -1) {};
	\end{pgfonlayer}
	\begin{pgfonlayer}{edgelayer}
		\draw [style=wire] (7) to (9);
		\draw [style=wire] (8) to (10);
		\draw [style=wire] (10) to (11);
	\end{pgfonlayer}
\end{tikzpicture}
 \end{array} &&  \begin{array}[c]{c}\begin{tikzpicture}
	\begin{pgfonlayer}{nodelayer}
		\node [style=component] (0) at (22.25, -0.25) {$\varepsilon$};
		\node [style=component] (2) at (22.25, 0.75) {$\mathsf{u}$};
		\node [style=port] (6) at (22.25, -1) {};
	\end{pgfonlayer}
	\begin{pgfonlayer}{edgelayer}
		\draw [style=wire] (0) to (2);
		\draw [style=wire] (0) to (6);
	\end{pgfonlayer}
\end{tikzpicture}
 \end{array} =  \begin{array}[c]{c}\begin{tikzpicture}
	\begin{pgfonlayer}{nodelayer}
		\node [style=component] (10) at (24, 0.5) {$0$};
		\node [style=port] (11) at (24, -1) {};
	\end{pgfonlayer}
	\begin{pgfonlayer}{edgelayer}
		\draw [style=wire] (10) to (11);
	\end{pgfonlayer}
\end{tikzpicture}
 \end{array}
\end{align}
Before we review additive bialgebra modalities, it will also be useful to introduce an intermediate concept. So a \textbf{convolution bialgebra modality} on an additive symmetric monoidal category is a bialgebra modality $(\oc, \delta, \varepsilon, \Delta, \mathsf{e}, \nabla, \mathsf{u})$ such that the diagrams from Appendix \ref{sec:diagaddbialgmod} commute. This states the compatibility between the additive enrichment and the bimonoid structure via convolution (\ref{diag:add-bialg}), which is drawn as follows: 
\begin{align}\label{string:!+02}
 \begin{array}[c]{c}\begin{tikzpicture}
	\begin{pgfonlayer}{nodelayer}
		\node [style=port] (12) at (33, 1.5) {};
		\node [style=port] (13) at (33, -2) {};
		\node [style=component] (14) at (32.5, -0.25) {$f$};
		\node [style=none] (15) at (32, 0.75) {};
		\node [style=none] (16) at (34, 0.75) {};
		\node [style=none] (17) at (32, -1.25) {};
		\node [style=none] (18) at (34, -1.25) {};
		\node [style=component] (19) at (33.5, -0.25) {$g$};
		\node [style=port] (20) at (32.5, 0.75) {};
		\node [style=port] (21) at (33.5, 0.75) {};
		\node [style=port] (22) at (33, -1.25) {};
		\node [style=port] (23) at (32.5, -1.25) {};
		\node [style=port] (24) at (33.5, -1.25) {};
		\node [style=port] (25) at (33, 0.75) {};
		\node [style=port] (26) at (33, -0.25) {$+$};
	\end{pgfonlayer}
	\begin{pgfonlayer}{edgelayer}
		\draw [style=wire] (15.center) to (16.center);
		\draw [style=wire] (16.center) to (18.center);
		\draw [style=wire] (18.center) to (17.center);
		\draw [style=wire] (17.center) to (15.center);
		\draw [style=wire] (20) to (14);
		\draw [style=wire] (14) to (23);
		\draw [style=wire] (21) to (19);
		\draw [style=wire] (19) to (24);
		\draw [style=wire] (22) to (13);
		\draw [style=wire] (12) to (25);
	\end{pgfonlayer}
\end{tikzpicture}
 \end{array} = \begin{array}[c]{c}\begin{tikzpicture}
	\begin{pgfonlayer}{nodelayer}
		\node [style=component] (430) at (89, -1) {$\nabla$};
		\node [style=none] (431) at (88.25, -0.5) {};
		\node [style=none] (432) at (89.75, -0.5) {};
		\node [style=component] (433) at (88.25, 0.25) {$f$};
		\node [style=component] (434) at (89.75, 0.25) {$g$};
		\node [style=component] (435) at (89, 1.5) {$\Delta$};
		\node [style=none] (436) at (88.25, 1) {};
		\node [style=none] (437) at (89.75, 1) {};
		\node [style=port] (440) at (89, 2.5) {};
		\node [style=port] (441) at (89, -2) {};
		\node [style=none] (443) at (87.75, -0.25) {};
		\node [style=none] (444) at (88.75, -0.25) {};
		\node [style=none] (445) at (87.75, 0.75) {};
		\node [style=none] (446) at (88.75, 0.75) {};
		\node [style=none] (448) at (89.25, -0.25) {};
		\node [style=none] (449) at (90.25, -0.25) {};
		\node [style=none] (450) at (89.25, 0.75) {};
		\node [style=none] (451) at (90.25, 0.75) {};
	\end{pgfonlayer}
	\begin{pgfonlayer}{edgelayer}
		\draw [style=wire] (430) to (431.center);
		\draw [style=wire] (430) to (432.center);
		\draw [style=wire] (434) to (432.center);
		\draw [style=wire] (433) to (431.center);
		\draw [style=wire] (435) to (436.center);
		\draw [style=wire] (435) to (437.center);
		\draw [style=wire] (436.center) to (433);
		\draw [style=wire] (437.center) to (434);
		\draw [style=wire] (443.center) to (444.center);
		\draw [style=wire] (444.center) to (446.center);
		\draw [style=wire] (446.center) to (445.center);
		\draw [style=wire] (445.center) to (443.center);
		\draw [style=wire] (448.center) to (449.center);
		\draw [style=wire] (449.center) to (451.center);
		\draw [style=wire] (451.center) to (450.center);
		\draw [style=wire] (450.center) to (448.center);
		\draw [style=wire] (440) to (435);
		\draw [style=wire] (441) to (430);
	\end{pgfonlayer}
\end{tikzpicture}
 \end{array}  &&  \begin{array}[c]{c}\begin{tikzpicture}
	\begin{pgfonlayer}{nodelayer}
		\node [style=port] (427) at (86.5, 1.5) {};
		\node [style=port] (428) at (86.5, -1) {};
		\node [style=component] (429) at (86.5, 0.25) {$0$};
		\node [style=none] (466) at (86, 0.75) {};
		\node [style=none] (467) at (87, 0.75) {};
		\node [style=none] (468) at (86, -0.25) {};
		\node [style=none] (469) at (87, -0.25) {};
	\end{pgfonlayer}
	\begin{pgfonlayer}{edgelayer}
		\draw [style=wire] (427) to (429);
		\draw [style=wire] (429) to (428);
		\draw [style=wire] (466.center) to (467.center);
		\draw [style=wire] (467.center) to (469.center);
		\draw [style=wire] (469.center) to (468.center);
		\draw [style=wire] (468.center) to (466.center);
	\end{pgfonlayer}
\end{tikzpicture}
 \end{array} = \begin{array}[c]{c}\begin{tikzpicture}
	\begin{pgfonlayer}{nodelayer}
		\node [style=component] (215) at (37.25, -5.25) {$\mathsf{e}$};
		\node [style=port] (216) at (37.25, -4.5) {};
		\node [style=component] (217) at (37.25, -6.25) {$\mathsf{u}$};
		\node [style=port] (218) at (37.25, -7) {};
	\end{pgfonlayer}
	\begin{pgfonlayer}{edgelayer}
		\draw [style=wire] (216) to (215);
		\draw [style=wire] (218) to (217);
	\end{pgfonlayer}
\end{tikzpicture}
 \end{array} 
\end{align}
Note that we have already encountered the above identities in the proof of Theorem \ref{thm:additive}. Then using the terminology from this paper, an \textbf{additive bialgebra modality} \cite[Def 5]{Blute2019} on an additive symmetric monoidal category is a convolution bialgebra modality $(\oc, \delta, \varepsilon, \Delta, \mathsf{e}, \nabla, \mathsf{u})$ which is also pre-additive. 

It turns out that for an additive bialgebra modality, its bimonoid structural maps are the unique ones that make its underlying coalgebra modality into a bialgebra modality. In fact, this is even true for a convolution bialgebra modality. Explicitly, by a \textbf{bialgebra structure} for a  coalgebra modality $(\oc, \delta, \varepsilon, \Delta, \mathsf{e})$, we mean natural transformations $\nabla$ and a map $\mathsf{u}$ such that $(\oc, \delta, \varepsilon, \Delta, \mathsf{e}, \nabla, \mathsf{u})$ is a bialgebra modality.

\begin{lem} For a convolution bialgebra modality, its underlying coalgebra modality has a unique bialgebra structure. 
\end{lem}
\begin{proof} Let $(\oc, \delta, \varepsilon, \Delta, \mathsf{e}, \nabla, \mathsf{u})$ be a convolution bialgebra modality. Suppose we have another bialgebra structure on the underlying coalgebra modality, that is, there are natural transformations $\nabla^\prime$ and $\mathsf{u}^\prime$ which also make $(\oc, \delta, \varepsilon, \Delta, \mathsf{e}, \nabla^\prime, \mathsf{u}^\prime)$ a bialgebra modality. We will show that $\nabla = \nabla^\prime$ and $\mathsf{u}= \mathsf{u}^\prime$. We begin by showing that the units are in fact equal: 
\begin{gather*}

\end{gather*}
So we have that $\mathsf{u}=\mathsf{u}^\prime$ as desired. To show that the multiplications are also equal, following \cite[App A]{Blute2019}, it will be useful to consider the natural transformation $\phi_A: \oc(A) \otimes \oc(A) \to A$ which is defined as follows:
\begin{gather}\label{string:phi-def}
%
\end{gather}
We first show the following useful equality\footnote{This equality was shown to hold for additive bialgebra modalities in \cite[Lemma 16]{Blute2019}. But here we check that it also holds for just a convolution bialgebra modality.}: 
\begin{align}\label{string:phi-lemma}
%
\end{align}
So we compute that:
\begin{gather*}
%
 
\end{gather*}
And similarly for the other equality of (\ref{string:phi-lemma}). Then using this we compute the following: 
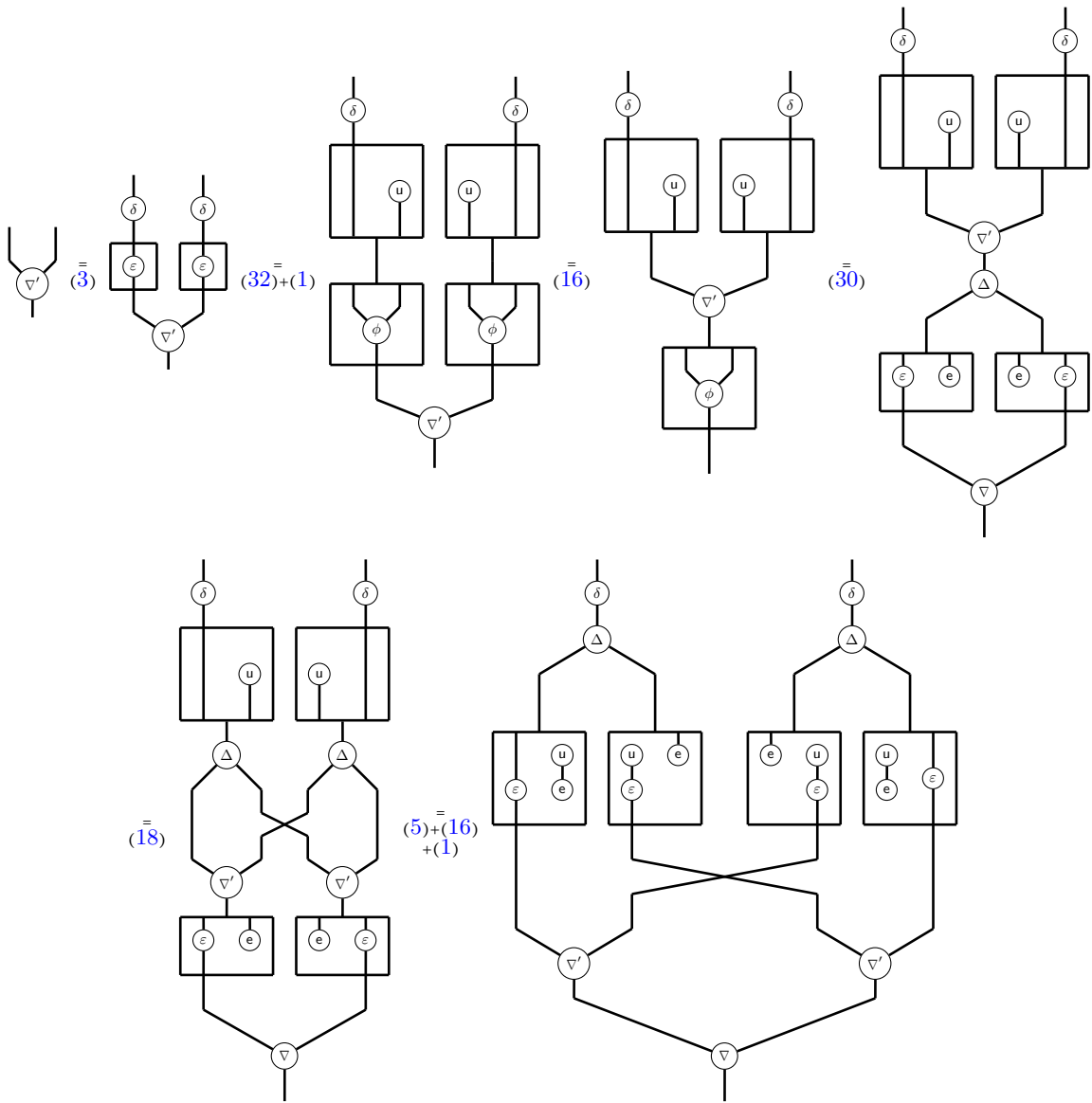
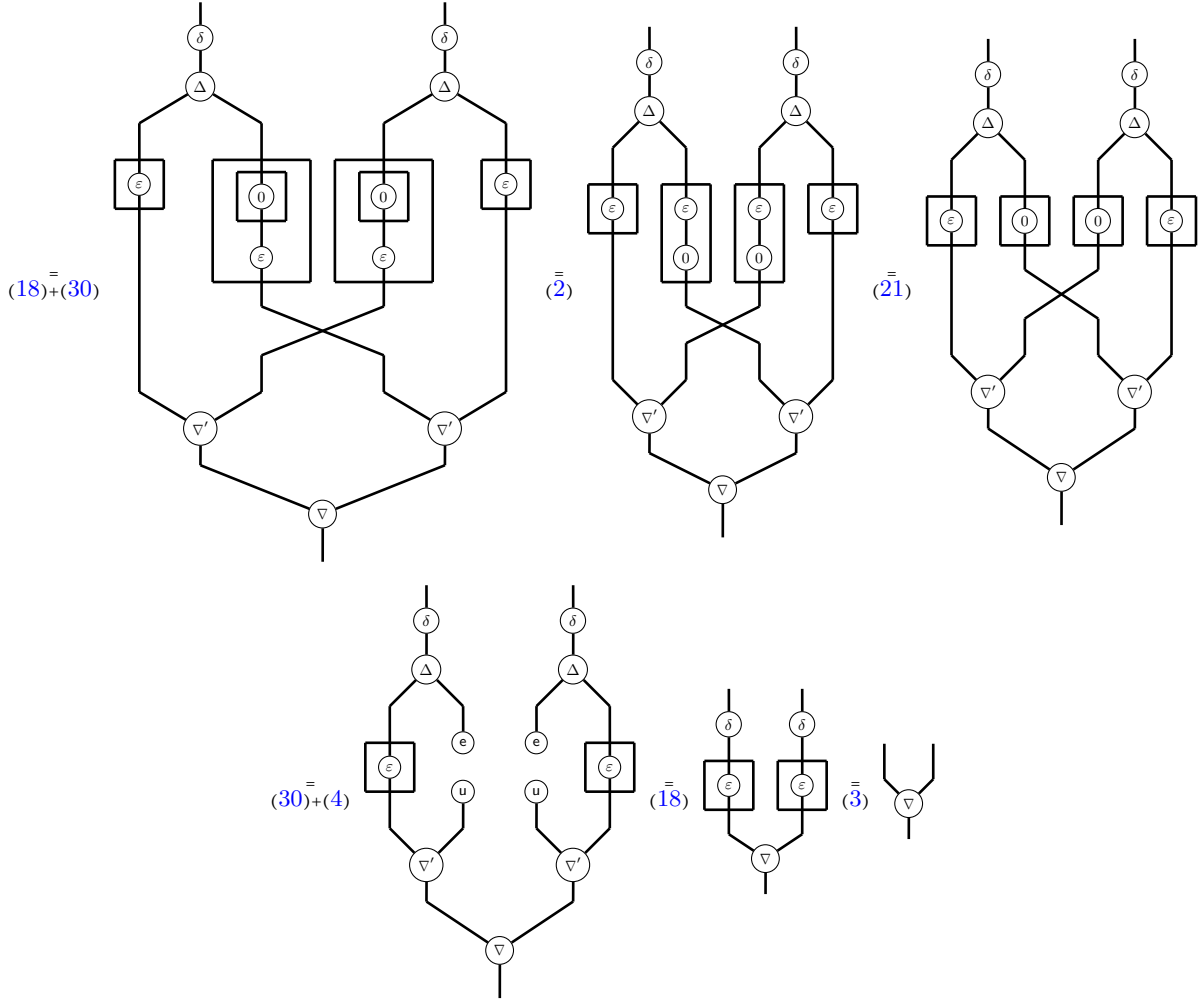
\begin{gather*}
%
 
\end{gather*}
So we have that $\nabla= \nabla^\prime$ as desired. Therefore, we conclude that $(\oc, \delta, \varepsilon, \Delta, \mathsf{e})$ has a unique bialgebra structure. 
\end{proof}

\begin{cor}\label{cor:add-bialg-unique} For an additive bialgebra modality, its underlying coalgebra modality has a unique bialgebra structure. 
\end{cor}

Let us now briefly review how, for an additive symmetric monoidal category, every monoidal coalgebra modality is an additive bialgebra modality, and vice-versa. Starting with a monoidal coalgebra modality $(\oc, \delta, \varepsilon, \Delta, \mathsf{e}, \mathsf{m}, \mathsf{m}_I)$, define $\nabla$ and $\mathsf{u}$ as in Appendix \ref{sec:moncoalg-addbialg}, which are drawn respectively as follows: 
 \begin{align}\label{string:mcm.to.abm} \begin{array}[c]{c}
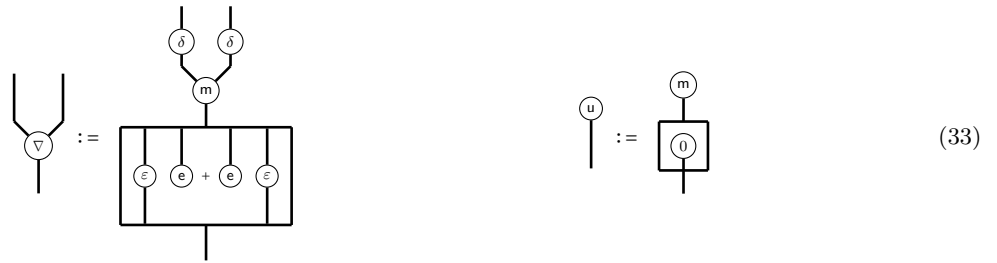
\begin{tikzpicture}
	\begin{pgfonlayer}{nodelayer}
		\node [style=component] (57) at (42.25, 2.75) {$\nabla$};
		\node [style=none] (58) at (41.75, 3.25) {};
		\node [style=none] (59) at (42.75, 3.25) {};
		\node [style=port] (60) at (41.75, 4.25) {};
		\node [style=port] (61) at (42.75, 4.25) {};
		\node [style=port] (62) at (42.25, 1.75) {};
	\end{pgfonlayer}
	\begin{pgfonlayer}{edgelayer}
		\draw [style=wire] (57) to (58.center);
		\draw [style=wire] (57) to (59.center);
		\draw [style=wire] (60) to (58.center);
		\draw [style=wire] (61) to (59.center);
		\draw [style=wire] (57) to (62);
	\end{pgfonlayer}
\end{tikzpicture}
 \end{array} := \begin{array}[c]{c}\begin{tikzpicture}
	\begin{pgfonlayer}{nodelayer}
		\node [style=component] (0) at (37.75, -0.75) {$\varepsilon$};
		\node [style=port] (6) at (37.75, 0.25) {};
		\node [style=component] (10) at (38.5, -0.75) {$\mathsf{e}$};
		\node [style=port] (11) at (38.5, 0.25) {};
		\node [style=none] (27) at (37.25, -1.75) {};
		\node [style=none] (28) at (40.75, -1.75) {};
		\node [style=none] (29) at (37.25, 0.25) {};
		\node [style=none] (30) at (40.75, 0.25) {};
		\node [style=port] (31) at (38.5, 2.75) {};
		\node [style=component] (32) at (39, 1) {$\mathsf{m}$};
		\node [style=none] (33) at (38.5, 1.5) {};
		\node [style=none] (34) at (39.5, 1.5) {};
		\node [style=port] (35) at (39.5, 2.75) {};
		\node [style=port] (36) at (39, -2.5) {};
		\node [style=none] (37) at (38.5, 0.25) {};
		\node [style=none] (38) at (39.5, 0.25) {};
		\node [style=none] (41) at (39, 0.25) {};
		\node [style=none] (42) at (39, -1.75) {};
		\node [style=component] (46) at (38.5, 2) {$\delta$};
		\node [style=component] (47) at (39.5, 2) {$\delta$};
		\node [style=port] (48) at (37.75, -1.75) {};
		\node [style=component] (49) at (40.25, -0.75) {$\varepsilon$};
		\node [style=port] (50) at (40.25, 0.25) {};
		\node [style=component] (51) at (39.5, -0.75) {$\mathsf{e}$};
		\node [style=port] (52) at (39.5, 0.25) {};
		\node [style=port] (53) at (40.25, -1.75) {};
		\node [style=port] (54) at (39, -0.75) {$+$};
	\end{pgfonlayer}
	\begin{pgfonlayer}{edgelayer}
		\draw [style=wire] (0) to (6);
		\draw [style=wire] (10) to (11);
		\draw [style=wire] (27.center) to (28.center);
		\draw [style=wire] (28.center) to (30.center);
		\draw [style=wire] (30.center) to (29.center);
		\draw [style=wire] (29.center) to (27.center);
		\draw [style=wire] (32) to (33.center);
		\draw [style=wire] (32) to (34.center);
		\draw [style=wire] (42.center) to (36);
		\draw [style=wire] (32) to (41.center);
		\draw [style=wire] (31) to (46);
		\draw [style=wire] (35) to (47);
		\draw [style=wire] (47) to (34.center);
		\draw [style=wire] (46) to (33.center);
		\draw [style=wire] (0) to (48);
		\draw [style=wire] (49) to (50);
		\draw [style=wire] (51) to (52);
		\draw [style=wire] (49) to (53);
	\end{pgfonlayer}
\end{tikzpicture}
 \end{array} && \begin{array}[c]{c} \begin{tikzpicture}
	\begin{pgfonlayer}{nodelayer}
		\node [style=component] (64) at (45, 0.75) {$\mathsf{u}$};
		\node [style=port] (65) at (45, -0.5) {};
	\end{pgfonlayer}
	\begin{pgfonlayer}{edgelayer}
		\draw [style=wire] (64) to (65);
	\end{pgfonlayer}
\end{tikzpicture}
 \end{array} := \begin{array}[c]{c}\begin{tikzpicture}
	\begin{pgfonlayer}{nodelayer}
		\node [style=component] (106) at (67, 0.25) {$\mathsf{m}$};
		\node [style=component] (107) at (67, -1) {$0$};
		\node [style=none] (108) at (66.5, -1.5) {};
		\node [style=none] (109) at (67.5, -1.5) {};
		\node [style=none] (110) at (66.5, -0.5) {};
		\node [style=none] (111) at (67.5, -0.5) {};
		\node [style=port] (112) at (67, -2) {};
		\node [style=none] (113) at (67, -0.5) {};
	\end{pgfonlayer}
	\begin{pgfonlayer}{edgelayer}
		\draw [style=wire] (108.center) to (109.center);
		\draw [style=wire] (109.center) to (111.center);
		\draw [style=wire] (111.center) to (110.center);
		\draw [style=wire] (110.center) to (108.center);
		\draw [style=wire] (112) to (107);
		\draw [style=wire] (106) to (113.center);
	\end{pgfonlayer}
\end{tikzpicture}
 \end{array}
 \end{align}
 Then $(\oc, \delta, \varepsilon, \Delta, \mathsf{e}, \nabla, \mathsf{u})$ is an additive bialgebra modality \cite[Prop 1]{Blute2019}. Conversely, starting with an additive bialgebra modality $(\oc, \delta, \varepsilon, \Delta, \mathsf{e}, \nabla, \mathsf{u})$, define $\mathsf{m}$ and $\mathsf{m}_I$ as in Appendix \ref{sec:moncoalg-addbialg}, which are drawn respectively as follows: 
  \begin{align}\label{string:abm.to.mcm} \begin{array}[c]{c}\begin{tikzpicture}
	\begin{pgfonlayer}{nodelayer}
		\node [style=component] (57) at (42.25, 2.75) {$\mathsf{m}$};
		\node [style=none] (58) at (41.75, 3.25) {};
		\node [style=none] (59) at (42.75, 3.25) {};
		\node [style=port] (60) at (41.75, 4.25) {};
		\node [style=port] (61) at (42.75, 4.25) {};
		\node [style=port] (62) at (42.25, 1.75) {};
	\end{pgfonlayer}
	\begin{pgfonlayer}{edgelayer}
		\draw [style=wire] (57) to (58.center);
		\draw [style=wire] (57) to (59.center);
		\draw [style=wire] (60) to (58.center);
		\draw [style=wire] (61) to (59.center);
		\draw [style=wire] (57) to (62);
	\end{pgfonlayer}
\end{tikzpicture}
 \end{array} := \begin{array}[c]{c}\begin{tikzpicture}
	\begin{pgfonlayer}{nodelayer}
		\node [style=component] (268) at (53, -4) {$\varepsilon$};
		\node [style=port] (269) at (53, -3) {};
		\node [style=component] (270) at (54, -4) {$\mathsf{e}$};
		\node [style=port] (271) at (54, -3) {};
		\node [style=component] (272) at (54, 3) {$\mathsf{u}$};
		\node [style=none] (273) at (54.5, 4) {};
		\node [style=none] (274) at (52.5, 4) {};
		\node [style=none] (275) at (54.5, 2) {};
		\node [style=none] (276) at (52.5, 2) {};
		\node [style=port] (277) at (53.5, 2) {};
		\node [style=port] (278) at (54, 2) {};
		\node [style=port] (279) at (53, 2) {};
		\node [style=none] (280) at (52.5, -5) {};
		\node [style=none] (281) at (54.5, -5) {};
		\node [style=none] (282) at (52.5, -3) {};
		\node [style=none] (283) at (54.5, -3) {};
		\node [style=port] (284) at (53, 5.5) {};
		\node [style=component] (285) at (54.75, 0.75) {$\nabla$};
		\node [style=none] (286) at (53.5, 1.5) {};
		\node [style=none] (287) at (56, 1.5) {};
		\node [style=component] (288) at (53, -5.75) {$\varepsilon$};
		\node [style=component] (289) at (53, 4.75) {$\delta$};
		\node [style=component] (290) at (54.75, -0.25) {$\delta$};
		\node [style=component] (291) at (55.5, 3) {$\mathsf{u}$};
		\node [style=none] (292) at (55, 4) {};
		\node [style=none] (293) at (57, 4) {};
		\node [style=none] (294) at (55, 2) {};
		\node [style=none] (295) at (57, 2) {};
		\node [style=port] (296) at (56, 2) {};
		\node [style=port] (297) at (55.5, 2) {};
		\node [style=port] (298) at (56.5, 2) {};
		\node [style=port] (299) at (56.5, 5.5) {};
		\node [style=component] (300) at (56.5, 4.75) {$\delta$};
		\node [style=component] (301) at (54.75, -1.75) {$\Delta$};
		\node [style=none] (302) at (53.5, -2.5) {};
		\node [style=none] (303) at (56, -2.5) {};
		\node [style=port] (304) at (53.5, -3) {};
		\node [style=component] (305) at (56.5, -4) {$\varepsilon$};
		\node [style=port] (306) at (56.5, -3) {};
		\node [style=component] (307) at (55.5, -4) {$\mathsf{e}$};
		\node [style=port] (308) at (55.5, -3) {};
		\node [style=none] (309) at (57, -5) {};
		\node [style=none] (310) at (55, -5) {};
		\node [style=none] (311) at (57, -3) {};
		\node [style=none] (312) at (55, -3) {};
		\node [style=component] (313) at (56.5, -5.75) {$\varepsilon$};
		\node [style=port] (314) at (56, -3) {};
		\node [style=port] (315) at (54.75, -6.75) {};
		\node [style=port] (316) at (54.75, -7.75) {};
		\node [style=none] (317) at (57.5, -1) {};
		\node [style=none] (318) at (52, -1) {};
		\node [style=none] (319) at (52, -6.75) {};
		\node [style=none] (320) at (57.5, -6.75) {};
		\node [style=port] (321) at (53, -6.75) {};
		\node [style=port] (322) at (56.5, -6.75) {};
	\end{pgfonlayer}
	\begin{pgfonlayer}{edgelayer}
		\draw [style=wire] (268) to (269);
		\draw [style=wire] (270) to (271);
		\draw [style=wire] (273.center) to (274.center);
		\draw [style=wire] (274.center) to (276.center);
		\draw [style=wire] (276.center) to (275.center);
		\draw [style=wire] (275.center) to (273.center);
		\draw [style=wire] (272) to (278);
		\draw [style=wire] (280.center) to (281.center);
		\draw [style=wire] (281.center) to (283.center);
		\draw [style=wire] (283.center) to (282.center);
		\draw [style=wire] (282.center) to (280.center);
		\draw [style=wire] (285) to (286.center);
		\draw [style=wire] (285) to (287.center);
		\draw [style=wire] (284) to (289);
		\draw [style=wire] (289) to (279);
		\draw [style=wire] (292.center) to (293.center);
		\draw [style=wire] (293.center) to (295.center);
		\draw [style=wire] (295.center) to (294.center);
		\draw [style=wire] (294.center) to (292.center);
		\draw [style=wire] (291) to (297);
		\draw [style=wire] (299) to (300);
		\draw [style=wire] (300) to (298);
		\draw [style=wire] (285) to (290);
		\draw [style=wire] (301) to (302.center);
		\draw [style=wire] (301) to (303.center);
		\draw [style=wire] (290) to (301);
		\draw [style=wire] (268) to (288);
		\draw [style=wire] (302.center) to (304);
		\draw [style=wire] (305) to (306);
		\draw [style=wire] (307) to (308);
		\draw [style=wire] (309.center) to (310.center);
		\draw [style=wire] (310.center) to (312.center);
		\draw [style=wire] (312.center) to (311.center);
		\draw [style=wire] (311.center) to (309.center);
		\draw [style=wire] (305) to (313);
		\draw [style=wire] (303.center) to (314);
		\draw [style=wire] (315) to (316);
		\draw [style=wire] (317.center) to (320.center);
		\draw [style=wire] (320.center) to (319.center);
		\draw [style=wire] (319.center) to (318.center);
		\draw [style=wire] (318.center) to (317.center);
		\draw [style=wire] (288) to (321);
		\draw [style=wire] (313) to (322);
		\draw [style=wire] (296) to (287.center);
		\draw [style=wire] (277) to (286.center);
	\end{pgfonlayer}
\end{tikzpicture}
 \end{array}  && \begin{array}[c]{c} \begin{tikzpicture}
	\begin{pgfonlayer}{nodelayer}
		\node [style=component] (64) at (45, 0.75) {$\mathsf{m}$};
		\node [style=port] (65) at (45, -0.5) {};
	\end{pgfonlayer}
	\begin{pgfonlayer}{edgelayer}
		\draw [style=wire] (64) to (65);
	\end{pgfonlayer}
\end{tikzpicture}
 \end{array} := \begin{array}[c]{c}\begin{tikzpicture}
	\begin{pgfonlayer}{nodelayer}
		\node [style=none] (96) at (46, -1) {};
		\node [style=none] (97) at (47, -1) {};
		\node [style=none] (98) at (46, -2) {};
		\node [style=none] (99) at (47, -2) {};
		\node [style=port] (100) at (46.5, -2.75) {};
		\node [style=none] (101) at (46.5, -2) {};
		\node [style=none] (102) at (46.5, -1) {};
		\node [style=component] (103) at (46.5, -1.5) {$\mathsf{e}$};
		\node [style=component] (104) at (46.5, 0.75) {$\mathsf{u}$};
		\node [style=component] (105) at (46.5, -0.25) {$\delta$};
	\end{pgfonlayer}
	\begin{pgfonlayer}{edgelayer}
		\draw [style=wire] (96.center) to (97.center);
		\draw [style=wire] (97.center) to (99.center);
		\draw [style=wire] (99.center) to (98.center);
		\draw [style=wire] (98.center) to (96.center);
		\draw [style=wire] (103) to (102.center);
		\draw [style=wire] (101.center) to (100);
		\draw [style=wire] (104) to (105);
		\draw [style=wire] (105) to (102.center);
	\end{pgfonlayer}
\end{tikzpicture}
 \end{array}
 \end{align}
Then $(\oc, \delta, \varepsilon, \Delta, \mathsf{e}, \mathsf{m}, \mathsf{m}_I)$ is a monoidal coalgebra modality \cite[Prop 3]{Blute2019}. Moreover, these constructions are inverses of each other \cite[Thm 1]{Blute2019}. In particular, this means that for an additive bialgebra modality $(\oc, \delta, \varepsilon, \Delta, \mathsf{e}, \nabla, \mathsf{u})$, we have that $\mathsf{m}$ and $\mathsf{m}_I$ as defined above are the unique such natural transformations which make $(\oc, \delta, \varepsilon, \Delta, \mathsf{e}, \mathsf{m}, \mathsf{m}_I)$ a monoidal coalgebra modality, and conversely, that for a monoidal coalgebra modality $(\oc, \delta, \varepsilon, \Delta, \mathsf{e}, \mathsf{m}, \mathsf{m}_I)$, $\nabla$ and $\mathsf{u}$ are the unique such natural transformations which make $(\oc, \delta, \varepsilon, \Delta, \mathsf{e}, \nabla, \mathsf{u})$ an additive bialgebra modality. 

 Let us now bring monoidal bialgebra modalities into the story. It turns out that \cite[Prop 2]{Blute2019} already tells us that an additive bialgebra modality is a monoidal bialgebra modality. 

\begin{prop}  Let let $(\oc, \delta, \varepsilon, \Delta, \mathsf{e}, \nabla, \mathsf{u})$ be an additive bialgebra modality with induced monoidal coalgebra modality $(\oc, \delta, \varepsilon,  \Delta, \mathsf{e}, \mathsf{m}, \mathsf{m}_I)$ (or equivalently let $(\oc, \delta, \varepsilon,  \Delta, \mathsf{e}, \mathsf{m}, \mathsf{m}_I)$ be a monoidal coalgebra modality with induced additive bialgebra modality $(\oc, \delta, \varepsilon, \Delta, \mathsf{e}, \nabla, \mathsf{u})$). Then $(\oc, \delta, \varepsilon,  \Delta, \mathsf{e}, \mathsf{m}, \mathsf{m}_I, \nabla, \mathsf{u})$ is a monoidal bialgebra modality. 
\end{prop}
\begin{proof} All that remains to check for the definition of a monoidal bialgebra modality are the equalities of (\ref{string:m-monoid}) and (\ref{string:monoid-!coalg}). However these equalities are precisely those from \cite[Prop 2]{Blute2019}, which were proved using string diagrams in \cite[App B]{Blute2019}. So we conclude that for an additive symmetric monoidal category, every additive bialgebra modality (or equivalently every monoidal coalgebra modality) is also a monoidal bialgebra modality. 
 \end{proof}

Conversely, let us now explain how, in the additive setting, a monoidal bialgebra modality is also an additive bialgebra modality. 
 
\begin{prop} Let $(\oc, \delta, \varepsilon,  \Delta, \mathsf{e}, \mathsf{m}, \mathsf{m}_I, \nabla, \mathsf{u})$ be a monoidal bialgebra modality on an additive symmetric monoidal category $\mathbb{X}$. Then: 
\begin{enumerate}[{\em (i)}]
\item The underlying bialgebra modality $(\oc, \delta, \varepsilon, \Delta, \mathsf{e}, \nabla, \mathsf{u})$ is an additive bialgebra modality, and the induced monoidal coalgebra modality via the construction of (\ref{string:abm.to.mcm}) is precisely $(\oc, \delta, \varepsilon,  \Delta, \mathsf{e}, \mathsf{m}, \mathsf{m}_I)$. 
\item For the underlying monoidal coalgebra modality $(\oc, \delta, \varepsilon,  \Delta, \mathsf{e}, \mathsf{m}, \mathsf{m}_I)$, the induced additive bialgebra modality via the construction of (\ref{string:mcm.to.abm}) is precisely $(\oc, \delta, \varepsilon, \Delta, \mathsf{e}, \nabla, \mathsf{u})$. 
  \end{enumerate}
\end{prop}
\begin{proof} Given a monoidal bialgebra modality $(\oc, \delta, \varepsilon,  \Delta, \mathsf{e}, \mathsf{m}, \mathsf{m}_I, \nabla, \mathsf{u})$, by definition we also have a monodial coalgebra modality $(\oc, \delta, \varepsilon,  \Delta, \mathsf{e}, \mathsf{m}, \mathsf{m}_I)$. Then by \cite[Thm 1]{Blute2019}, there are natural transformations $\nabla^\prime$ and $\mathsf{u}^\prime$ constructed as in (\ref{string:mcm.to.abm}) such that $(\oc, \delta, \varepsilon, \Delta, \mathsf{e}, \nabla^\prime, \mathsf{u}^\prime)$ is an additive bialgebra modality. However by assumption, $(\oc, \delta, \varepsilon,  \Delta, \mathsf{e}, \nabla, \mathsf{u})$ is also a bialgebra modality. Then by Cor \ref{cor:add-bialg-unique}, it follows that $\nabla=\nabla^\prime$ and $\mathsf{u}=\mathsf{u}^\prime$. So $(\oc, \delta, \varepsilon,  \Delta, \mathsf{e}, \nabla, \mathsf{u})$ is an additive bialgebra modality. Moreover, applying the construction of \cite[Thm 1]{Blute2019} as above gives back $\mathsf{m}$ and $\mathsf{m}_I$ as well. 
\end{proof}

Therefore we may extend \cite[Thm 1]{Blute2019} as follows:
 
\begin{thm} For an additive symmetric monoidal category, there is a bijective correspondence between the following:
\begin{enumerate}[{\em (i)}]
\item Monoidal coalgebra modalities;
\item Additive bialgebra modalities;
\item Monoidal bialgebra modalities. 
  \end{enumerate}
\end{thm}

As an immediate consequence, a monoidal bialgebra modality with a pre-codereliction is an additive bialgebra modality. 

\begin{cor}\label{cor:mon-add-1} Let $(\oc, \delta, \varepsilon,  \Delta, \mathsf{e}, \mathsf{m}, \mathsf{m}_I, \nabla, \mathsf{u})$ be a monoidal bialgebra modality equipped with a pre-codereliction $\eta$ on a symmetric monoidal category $\mathbb{X}$. Then for the induced additive symmetric monoidal category $\mathbb{X}$ as defined in Thm \ref{thm:additive}, $(\oc, \delta, \varepsilon,  \Delta, \mathsf{e}, \nabla, \mathsf{u})$ is an additive bialgebra modality. 
\end{cor}

That said, additive enrichment is not unique for a given category. Therefore, we need to verify that for an additive symmetric monoidal category, the additive enrichment induced from a monoidal bialgebra modality with a pre-codereliction is the same. 

\begin{lem} Let $(\oc, \delta, \varepsilon,  \Delta, \mathsf{e}, \mathsf{m}, \mathsf{m}_I, \nabla, \mathsf{u})$ be a monoidal bialgebra modality equipped with a pre-codereliction $\eta$ on an additive symmetric monoidal category $\mathbb{X}$. Then the induced additive structure from Thm \ref{thm:additive} is precisely the same as the starting additive structure on $\mathbb{X}$. 
\end{lem}
 \begin{proof} We need to show that the formulas from Thm \ref{thm:additive} give us back the starting sum and zero. This essentially follows from the convolution bialgebra modality axioms. So starting with checking the zero, we compute: 
 \begin{gather*} \begin{array}[c]{c}
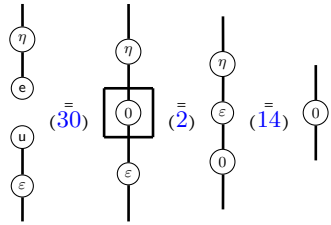
\begin{tikzpicture}
	\begin{pgfonlayer}{nodelayer}
		\node [style=component] (452) at (93.5, 1.5) {$\mathsf{e}$};
		\node [style=component] (453) at (93.5, 0.5) {$\mathsf{u}$};
		\node [style=component] (462) at (93.5, 2.5) {$\eta$};
		\node [style=component] (463) at (93.5, -0.5) {$\varepsilon$};
		\node [style=port] (464) at (93.5, 3.25) {};
		\node [style=port] (465) at (93.5, -1.25) {};
	\end{pgfonlayer}
	\begin{pgfonlayer}{edgelayer}
		\draw [style=wire] (464) to (462);
		\draw [style=wire] (463) to (465);
		\draw [style=wire] (462) to (452);
		\draw [style=wire] (453) to (463);
	\end{pgfonlayer}
\end{tikzpicture}
 \end{array}  \substack{=\\(\ref{string:!+02})}  \begin{array}[c]{c}\begin{tikzpicture}
	\begin{pgfonlayer}{nodelayer}
		\node [style=component] (267) at (117.5, 2.25) {$\varepsilon$};
		\node [style=port] (268) at (117.5, 5.75) {};
		\node [style=port] (269) at (117.5, 1.25) {};
		\node [style=component] (271) at (117.5, 3.5) {$0$};
		\node [style=none] (272) at (117, 3) {};
		\node [style=none] (273) at (118, 3) {};
		\node [style=none] (274) at (117, 4) {};
		\node [style=none] (275) at (118, 4) {};
		\node [style=none] (276) at (117.5, 4) {};
		\node [style=component] (354) at (117.5, 4.75) {$\eta$};
	\end{pgfonlayer}
	\begin{pgfonlayer}{edgelayer}
		\draw [style=wire] (267) to (269);
		\draw [style=wire] (272.center) to (273.center);
		\draw [style=wire] (273.center) to (275.center);
		\draw [style=wire] (275.center) to (274.center);
		\draw [style=wire] (274.center) to (272.center);
		\draw [style=wire] (271) to (267);
		\draw [style=wire] (354) to (271);
		\draw [style=wire] (268) to (354);
	\end{pgfonlayer}
\end{tikzpicture}
 \end{array}  \substack{=\\(\ref{string:comonad-nat})}   \begin{array}[c]{c}\begin{tikzpicture}
	\begin{pgfonlayer}{nodelayer}
		\node [style=port] (1594) at (265, 0) {};
		\node [style=port] (1595) at (265, -4) {};
		\node [style=component] (1596) at (265, -3) {$0$};
		\node [style=component] (1597) at (265, -1) {$\eta$};
		\node [style=component] (1598) at (265, -2) {$\varepsilon$};
	\end{pgfonlayer}
	\begin{pgfonlayer}{edgelayer}
		\draw [style=wire] (1596) to (1595);
		\draw [style=wire] (1594) to (1597);
		\draw [style=wire] (1597) to (1598);
		\draw [style=wire] (1598) to (1596);
	\end{pgfonlayer}
\end{tikzpicture}
 \end{array} \substack{=\\(\ref{string:coder})} \begin{array}[c]{c}\begin{tikzpicture}
	\begin{pgfonlayer}{nodelayer}
		\node [style=port] (427) at (86.5, 1.25) {};
		\node [style=port] (428) at (86.5, -0.75) {};
		\node [style=component] (429) at (86.5, 0.25) {$0$};
	\end{pgfonlayer}
	\begin{pgfonlayer}{edgelayer}
		\draw [style=wire] (427) to (429);
		\draw [style=wire] (429) to (428);
	\end{pgfonlayer}
\end{tikzpicture}
 \end{array} 
 \end{gather*}
Then for the sum we compute: 
 \begin{gather*} 
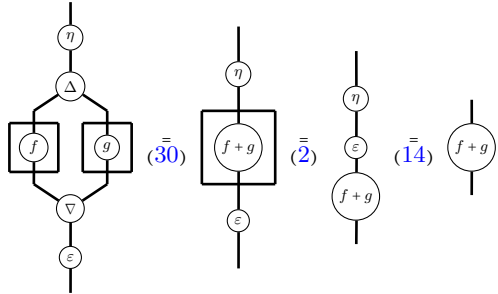
 \begin{array}[c]{c} \begin{tikzpicture}
	\begin{pgfonlayer}{nodelayer}
		\node [style=component] (470) at (91.75, -7) {$\nabla$};
		\node [style=none] (471) at (91, -6.5) {};
		\node [style=none] (472) at (92.5, -6.5) {};
		\node [style=component] (473) at (91, -5.75) {$f$};
		\node [style=component] (474) at (92.5, -5.75) {$g$};
		\node [style=component] (475) at (91.75, -4.5) {$\Delta$};
		\node [style=none] (476) at (91, -5) {};
		\node [style=none] (477) at (92.5, -5) {};
		\node [style=component] (478) at (91.75, -3.5) {$\eta$};
		\node [style=component] (479) at (91.75, -8) {$\varepsilon$};
		\node [style=port] (480) at (91.75, -2.75) {};
		\node [style=port] (481) at (91.75, -8.75) {};
		\node [style=none] (482) at (90.5, -6.25) {};
		\node [style=none] (483) at (91.5, -6.25) {};
		\node [style=none] (484) at (90.5, -5.25) {};
		\node [style=none] (485) at (91.5, -5.25) {};
		\node [style=none] (486) at (92, -6.25) {};
		\node [style=none] (487) at (93, -6.25) {};
		\node [style=none] (488) at (92, -5.25) {};
		\node [style=none] (489) at (93, -5.25) {};
	\end{pgfonlayer}
	\begin{pgfonlayer}{edgelayer}
		\draw [style=wire] (470) to (471.center);
		\draw [style=wire] (470) to (472.center);
		\draw [style=wire] (474) to (472.center);
		\draw [style=wire] (473) to (471.center);
		\draw [style=wire] (475) to (476.center);
		\draw [style=wire] (475) to (477.center);
		\draw [style=wire] (476.center) to (473);
		\draw [style=wire] (477.center) to (474);
		\draw [style=wire] (478) to (475);
		\draw [style=wire] (470) to (479);
		\draw [style=wire] (480) to (478);
		\draw [style=wire] (479) to (481);
		\draw [style=wire] (482.center) to (483.center);
		\draw [style=wire] (483.center) to (485.center);
		\draw [style=wire] (485.center) to (484.center);
		\draw [style=wire] (484.center) to (482.center);
		\draw [style=wire] (486.center) to (487.center);
		\draw [style=wire] (487.center) to (489.center);
		\draw [style=wire] (489.center) to (488.center);
		\draw [style=wire] (488.center) to (486.center);
	\end{pgfonlayer}
\end{tikzpicture}
 \end{array}  \substack{=\\(\ref{string:!+02})} \begin{array}[c]{c}\begin{tikzpicture}
	\begin{pgfonlayer}{nodelayer}
		\node [style=component] (267) at (127, -0.75) {$\varepsilon$};
		\node [style=port] (268) at (127, 3.25) {};
		\node [style=port] (269) at (127, -1.75) {};
		\node [style=component] (354) at (127, 2.25) {$\eta$};
		\node [style=component] (363) at (127, 0.75) {$f+g$};
		\node [style=none] (364) at (126.25, 1.5) {};
		\node [style=none] (365) at (127.75, 1.5) {};
		\node [style=none] (366) at (126.25, 0) {};
		\node [style=none] (367) at (127.75, 0) {};
	\end{pgfonlayer}
	\begin{pgfonlayer}{edgelayer}
		\draw [style=wire] (267) to (269);
		\draw [style=wire] (268) to (354);
		\draw [style=wire] (364.center) to (365.center);
		\draw [style=wire] (365.center) to (367.center);
		\draw [style=wire] (367.center) to (366.center);
		\draw [style=wire] (366.center) to (364.center);
		\draw [style=wire] (363) to (267);
		\draw [style=wire] (363) to (354);
	\end{pgfonlayer}
\end{tikzpicture}
 \end{array} \substack{=\\(\ref{string:comonad-nat})} \begin{array}[c]{c}\begin{tikzpicture}
	\begin{pgfonlayer}{nodelayer}
		\node [style=port] (1594) at (265, 0) {};
		\node [style=port] (1595) at (265, -4) {};
		\node [style=component] (1596) at (265, -3) {$f+g$};
		\node [style=component] (1597) at (265, -1) {$\eta$};
		\node [style=component] (1598) at (265, -2) {$\varepsilon$};
	\end{pgfonlayer}
	\begin{pgfonlayer}{edgelayer}
		\draw [style=wire] (1596) to (1595);
		\draw [style=wire] (1594) to (1597);
		\draw [style=wire] (1597) to (1598);
		\draw [style=wire] (1598) to (1596);
	\end{pgfonlayer}
\end{tikzpicture}
 \end{array} \substack{=\\(\ref{string:coder})} \begin{array}[c]{c}\begin{tikzpicture}
	\begin{pgfonlayer}{nodelayer}
		\node [style=port] (427) at (86.5, 1.25) {};
		\node [style=port] (428) at (86.5, -0.75) {};
		\node [style=component] (429) at (86.5, 0.25) {$f+g$};
	\end{pgfonlayer}
	\begin{pgfonlayer}{edgelayer}
		\draw [style=wire] (427) to (429);
		\draw [style=wire] (429) to (428);
	\end{pgfonlayer}
\end{tikzpicture}
 \end{array}  
 \end{gather*}
 So we conclude that the additive enrichments are the same.  
 \end{proof}

\section{Hopf Coalgebra Modalities}\label{sec:hopf-coalg-mod}

Before we discuss differential categories, we take a slight detour to address the natural follow up question regarding Abelian group enrichment. It turns out that having negatives for the additive enrichment corresponds to asking that our monoidal/additive bialgebra modality comes equipped with an antipode which gives a natural Hopf monoid structure. In fact, it was shown in \cite[Prop 7.6]{lemay2019lifting} that an additive bialgebra modality has an antipode if and only if the base category has negatives. This makes sense since the additive enrichment is intrinsically linked to bialgebra convolution, and having inverses for bialgebra convolution corresponds precisely to having an antipode \cite[Pg.72]{sweedler1969hopf}, i.e., that the bialgebra is in fact a Hopf algebra. So in this section we introduce the notion of a \emph{monoidal Hopf coalgebra modality} in order to extend Thm \ref{thm:additive} to include negatives. Along the way, we also provide an alternative proof of \cite[Prop 7.6]{lemay2019lifting}, as well some slightly more general observations for convolution bialgebra modalities. 

So a \textbf{Hopf coalgebra modality}\footnote{We've elected to use the term Hopf \emph{coalegbra} modality instead of Hopf \emph{algebra} modality since the term algebra modality is used for the dual notion of a coalgebra modality. In particular, an algebra modality is a monad. So here we use the term Hopf coalgebra modality to emphasize that we have a comonad, hopefully avoiding any confusion or ambiguity.} on a symmetric monoidal category is an octuple $(\oc, \delta, \varepsilon, \Delta, \mathsf{e}, \nabla, \mathsf{u}, \mathsf{S})$ consisting of a bialgebra modality $(\oc, \delta, \varepsilon, \Delta, \mathsf{e}, \nabla, \mathsf{u})$ and a natural transformation $\mathsf{S}_A: \oc(A) \to \oc(A)$ such that the diagrams in Appendix \ref{sec:diagHopf} commute. The natural transformation $\mathsf{S}$ is called the \textbf{antipode} and is drawn as follows: 
\[ \begin{tikzpicture}
	\begin{pgfonlayer}{nodelayer}
		\node [style=object] (0) at (0, 2) {$\oc(A)$};
		\node [style=object] (1) at (0, 0) {$\oc(A)$};
		\node [style=component] (2) at (0, 1) {$\mathsf{S}$};
	\end{pgfonlayer}
	\begin{pgfonlayer}{edgelayer}
		\draw [style=wire] (0) to (2);
		\draw [style=wire] (2) to (1);
	\end{pgfonlayer}
\end{tikzpicture}
 \]
 and its naturality is drawn as: 
 \begin{align}
  \begin{array}[c]{c} \begin{tikzpicture}
	\begin{pgfonlayer}{nodelayer}
		\node [style=port] (16) at (-2.25, 0.5) {};
		\node [style=component] (17) at (-2.25, 1.25) {$\mathsf{S}$};
		\node [style=port] (18) at (-2.25, 3.75) {};
		\node [style=component] (19) at (-2.25, 2.5) {$f$};
		\node [style=none] (20) at (-2.75, 3) {};
		\node [style=none] (21) at (-1.75, 3) {};
		\node [style=none] (22) at (-2.75, 2) {};
		\node [style=none] (23) at (-1.75, 2) {};
	\end{pgfonlayer}
	\begin{pgfonlayer}{edgelayer}
		\draw [style=wire] (17) to (16);
		\draw [style=wire] (18) to (19);
		\draw [style=wire] (20.center) to (21.center);
		\draw [style=wire] (21.center) to (23.center);
		\draw [style=wire] (23.center) to (22.center);
		\draw [style=wire] (22.center) to (20.center);
		\draw [style=wire] (19) to (17);
	\end{pgfonlayer}
\end{tikzpicture}
\end{array} =  \begin{array}[c]{c}\begin{tikzpicture}
	\begin{pgfonlayer}{nodelayer}
		\node [style=port] (3) at (0, 3.75) {};
		\node [style=component] (4) at (0, 3) {$\mathsf{S}$};
		\node [style=port] (5) at (0, 0.5) {};
		\node [style=component] (7) at (0, 1.75) {$f$};
		\node [style=none] (8) at (-0.5, 2.25) {};
		\node [style=none] (9) at (0.5, 2.25) {};
		\node [style=none] (10) at (-0.5, 1.25) {};
		\node [style=none] (11) at (0.5, 1.25) {};
	\end{pgfonlayer}
	\begin{pgfonlayer}{edgelayer}
		\draw [style=wire] (4) to (3);
		\draw [style=wire] (5) to (7);
		\draw [style=wire] (8.center) to (9.center);
		\draw [style=wire] (9.center) to (11.center);
		\draw [style=wire] (11.center) to (10.center);
		\draw [style=wire] (10.center) to (8.center);
		\draw [style=wire] (4) to (7);
	\end{pgfonlayer}
\end{tikzpicture}
\end{array}
 \end{align}
 We then also ask that for each object $A$, $(\oc(A), \nabla_A, \mathsf{u}_A, \Delta_A, \mathsf{e}_A, \mathsf{S}_A)$ is a Hopf monoid (\ref{diag:Hopf}), which is drawn out as follows: 
\begin{align}\label{string:Hopf}
\begin{array}[c]{c}\begin{tikzpicture}
	\begin{pgfonlayer}{nodelayer}
		\node [style=component] (372) at (132, -4.25) {$\Delta$};
		\node [style=none] (373) at (131.25, -5) {};
		\node [style=none] (374) at (132.75, -5) {};
		\node [style=component] (375) at (132, -6.75) {$\nabla$};
		\node [style=none] (376) at (131.25, -6) {};
		\node [style=none] (377) at (132.75, -6) {};
		\node [style=component] (378) at (132.75, -5.5) {$\mathsf{S}$};
		\node [style=port] (379) at (132, -3.5) {};
		\node [style=port] (380) at (132, -7.5) {};
	\end{pgfonlayer}
	\begin{pgfonlayer}{edgelayer}
		\draw [style=wire] (372) to (373.center);
		\draw [style=wire] (372) to (374.center);
		\draw [style=wire] (375) to (376.center);
		\draw [style=wire] (375) to (377.center);
		\draw [style=wire] (379) to (372);
		\draw [style=wire] (375) to (380);
		\draw [style=wire] (373.center) to (376.center);
		\draw [style=wire] (377.center) to (378);
		\draw [style=wire] (374.center) to (378);
	\end{pgfonlayer}
\end{tikzpicture}
\end{array}
= \begin{array}[c]{c}\begin{tikzpicture}
	\begin{pgfonlayer}{nodelayer}
		\node [style=component] (368) at (130, -5) {$\mathsf{e}$};
		\node [style=port] (369) at (130, -3.5) {};
		\node [style=component] (370) at (130, -6) {$\mathsf{u}$};
		\node [style=port] (371) at (130, -7.5) {};
	\end{pgfonlayer}
	\begin{pgfonlayer}{edgelayer}
		\draw [style=wire] (369) to (368);
		\draw [style=wire] (371) to (370);
	\end{pgfonlayer}
\end{tikzpicture}
 \end{array} = \begin{array}[c]{c} \begin{tikzpicture}
	\begin{pgfonlayer}{nodelayer}
		\node [style=component] (372) at (132, -4.25) {$\Delta$};
		\node [style=none] (373) at (132.75, -5) {};
		\node [style=none] (374) at (131.25, -5) {};
		\node [style=component] (375) at (132, -6.75) {$\nabla$};
		\node [style=none] (376) at (132.75, -6) {};
		\node [style=none] (377) at (131.25, -6) {};
		\node [style=component] (378) at (131.25, -5.5) {$\mathsf{S}$};
		\node [style=port] (379) at (132, -3.5) {};
		\node [style=port] (380) at (132, -7.5) {};
	\end{pgfonlayer}
	\begin{pgfonlayer}{edgelayer}
		\draw [style=wire] (372) to (373.center);
		\draw [style=wire] (372) to (374.center);
		\draw [style=wire] (375) to (376.center);
		\draw [style=wire] (375) to (377.center);
		\draw [style=wire] (379) to (372);
		\draw [style=wire] (375) to (380);
		\draw [style=wire] (373.center) to (376.center);
		\draw [style=wire] (377.center) to (378);
		\draw [style=wire] (374.center) to (378);
	\end{pgfonlayer}
\end{tikzpicture}
 \end{array}
\end{align}

It is well known that for bimonoids, antipodes are unique \cite[Pg.71]{sweedler1969hopf}, so being a Hopf monoid is a property of a bimonoid. The same is true for bialgebra modalities. 

\begin{lem} If a bialgbera modality has an antipode, then it is unique. 
\end{lem}

It is also well known that for bicommutative Hopf monoids, the antipode is an isomorphism and its own inverse \cite[Prop 4.0.1.(6)]{sweedler1969hopf}, and that the antipode is a (co)monoid morphism \cite[Prop 4.0.1.(1-4)]{sweedler1969hopf}. 

\begin{lem}\label{lemma:Hopf-1} If $(\oc, \delta, \varepsilon, \Delta, \mathsf{e}, \nabla, \mathsf{u}, \mathsf{S})$ is a Hopf coalgebra modality, then the diagrams in Appendix (\ref{sec:lemma-Hopf-1}) commute, which are drawn as follows: 
\begin{align}\label{strings:S-lemma}
 \begin{array}[c]{c}\begin{tikzpicture}
	\begin{pgfonlayer}{nodelayer}
		\node [style=port] (16) at (-2.25, 0.5) {};
		\node [style=component] (17) at (-2.25, 1.5) {$\mathsf{S}$};
		\node [style=port] (18) at (-2.25, 3.75) {};
		\node [style=component] (19) at (-2.25, 2.75) {$\mathsf{S}$};
	\end{pgfonlayer}
	\begin{pgfonlayer}{edgelayer}
		\draw [style=wire] (17) to (16);
		\draw [style=wire] (18) to (19);
		\draw [style=wire] (19) to (17);
	\end{pgfonlayer}
\end{tikzpicture}
\end{array} =  \begin{array}[c]{c}\begin{tikzpicture}
	\begin{pgfonlayer}{nodelayer}
		\node [style=port] (16) at (-2.25, 0.5) {};
		\node [style=port] (18) at (-2.25, 3.75) {};
	\end{pgfonlayer}
	\begin{pgfonlayer}{edgelayer}
		\draw [style=wire] (18) to (16);
	\end{pgfonlayer}
\end{tikzpicture}
\end{array} && \begin{array}[c]{c}\begin{tikzpicture}
	\begin{pgfonlayer}{nodelayer}
		\node [style=component] (7) at (-3.5, 3.25) {$\mathsf{S}$};
		\node [style=port] (20) at (-4, 1.25) {};
		\node [style=component] (21) at (-3.5, 2.25) {$\Delta$};
		\node [style=none] (22) at (-4, 1.75) {};
		\node [style=none] (23) at (-3, 1.75) {};
		\node [style=port] (34) at (-3, 1.25) {};
		\node [style=port] (45) at (-3.5, 4) {};
	\end{pgfonlayer}
	\begin{pgfonlayer}{edgelayer}
		\draw [style=wire] (21) to (22.center);
		\draw [style=wire] (21) to (23.center);
		\draw [style=wire] (22.center) to (20);
		\draw [style=wire] (23.center) to (34);
		\draw [style=wire] (45) to (7);
		\draw [style=wire] (7) to (21);
	\end{pgfonlayer}
\end{tikzpicture}
 \end{array} = \begin{array}[c]{c}\begin{tikzpicture}
	\begin{pgfonlayer}{nodelayer}
		\node [style=component] (46) at (-1.5, 1.75) {$\mathsf{S}$};
		\node [style=port] (51) at (-1.5, 1) {};
		\node [style=component] (52) at (-1, 3) {$\Delta$};
		\node [style=none] (53) at (-1.5, 2.5) {};
		\node [style=none] (54) at (-0.5, 2.5) {};
		\node [style=port] (55) at (-0.5, 1) {};
		\node [style=component] (56) at (-0.5, 1.75) {$\mathsf{S}$};
		\node [style=port] (61) at (-1, 3.75) {};
	\end{pgfonlayer}
	\begin{pgfonlayer}{edgelayer}
		\draw [style=wire] (52) to (53.center);
		\draw [style=wire] (52) to (54.center);
		\draw [style=wire] (46) to (51);
		\draw [style=wire] (56) to (55);
		\draw [style=wire] (53.center) to (46);
		\draw [style=wire] (54.center) to (56);
		\draw [style=wire] (61) to (52);
	\end{pgfonlayer}
\end{tikzpicture}
 \end{array} && \begin{array}[c]{c}\begin{tikzpicture}
	\begin{pgfonlayer}{nodelayer}
		\node [style=component] (7) at (-3.5, 3.25) {$\mathsf{S}$};
		\node [style=component] (21) at (-3.5, 2.25) {$\mathsf{e}$};
		\node [style=port] (45) at (-3.5, 4) {};
	\end{pgfonlayer}
	\begin{pgfonlayer}{edgelayer}
		\draw [style=wire] (45) to (7);
		\draw [style=wire] (7) to (21);
	\end{pgfonlayer}
\end{tikzpicture}
 \end{array} = \begin{array}[c]{c}\begin{tikzpicture}
	\begin{pgfonlayer}{nodelayer}
		\node [style=component] (52) at (-0.75, 3.25) {$\mathsf{e}$};
		\node [style=port] (61) at (-0.75, 4.25) {};
	\end{pgfonlayer}
	\begin{pgfonlayer}{edgelayer}
		\draw [style=wire] (61) to (52);
	\end{pgfonlayer}
\end{tikzpicture}
 \end{array} && \begin{array}[c]{c}\begin{tikzpicture}
	\begin{pgfonlayer}{nodelayer}
		\node [style=component] (389) at (137.75, -3) {$\mathsf{S}$};
		\node [style=port] (390) at (137.25, -1) {};
		\node [style=component] (391) at (137.75, -2) {$\nabla$};
		\node [style=none] (392) at (137.25, -1.5) {};
		\node [style=none] (393) at (138.25, -1.5) {};
		\node [style=port] (394) at (138.25, -1) {};
		\node [style=port] (395) at (137.75, -3.75) {};
	\end{pgfonlayer}
	\begin{pgfonlayer}{edgelayer}
		\draw [style=wire] (391) to (392.center);
		\draw [style=wire] (391) to (393.center);
		\draw [style=wire] (392.center) to (390);
		\draw [style=wire] (393.center) to (394);
		\draw [style=wire] (395) to (389);
		\draw [style=wire] (389) to (391);
	\end{pgfonlayer}
\end{tikzpicture}
 \end{array} = \begin{array}[c]{c}\begin{tikzpicture}
	\begin{pgfonlayer}{nodelayer}
		\node [style=component] (396) at (139.25, 3) {$\mathsf{S}$};
		\node [style=port] (397) at (139.25, 3.75) {};
		\node [style=component] (398) at (139.75, 1.75) {$\nabla$};
		\node [style=none] (399) at (139.25, 2.25) {};
		\node [style=none] (400) at (140.25, 2.25) {};
		\node [style=port] (401) at (140.25, 3.75) {};
		\node [style=component] (402) at (140.25, 3) {$\mathsf{S}$};
		\node [style=port] (403) at (139.75, 1) {};
	\end{pgfonlayer}
	\begin{pgfonlayer}{edgelayer}
		\draw [style=wire] (398) to (399.center);
		\draw [style=wire] (398) to (400.center);
		\draw [style=wire] (396) to (397);
		\draw [style=wire] (402) to (401);
		\draw [style=wire] (399.center) to (396);
		\draw [style=wire] (400.center) to (402);
		\draw [style=wire] (403) to (398);
	\end{pgfonlayer}
\end{tikzpicture}
 \end{array} && \begin{array}[c]{c}\begin{tikzpicture}
	\begin{pgfonlayer}{nodelayer}
		\node [style=component] (404) at (143.75, 0) {$\mathsf{S}$};
		\node [style=component] (405) at (143.75, 1) {$\mathsf{u}$};
		\node [style=port] (406) at (143.75, -0.75) {};
	\end{pgfonlayer}
	\begin{pgfonlayer}{edgelayer}
		\draw [style=wire] (406) to (404);
		\draw [style=wire] (404) to (405);
	\end{pgfonlayer}
\end{tikzpicture}
 \end{array} = \begin{array}[c]{c}\begin{tikzpicture}
	\begin{pgfonlayer}{nodelayer}
		\node [style=component] (407) at (144.75, 4.25) {$\mathsf{u}$};
		\node [style=port] (408) at (144.75, 3.25) {};
	\end{pgfonlayer}
	\begin{pgfonlayer}{edgelayer}
		\draw [style=wire] (408) to (407);
	\end{pgfonlayer}
\end{tikzpicture}
 \end{array}
\end{align}
\end{lem}

A \textbf{monoidal Hopf coalgebra modality} is a decuple $(\oc, \delta, \varepsilon,  \Delta, \mathsf{e}, \mathsf{m}, \mathsf{m}_I, \nabla, \mathsf{u}, \mathsf{S})$ consisting of a monoidal bialgebra modality $(\oc, \delta, \varepsilon,  \Delta, \mathsf{e}, \mathsf{m}, \mathsf{m}_I, \nabla, \mathsf{u})$ and a Hopf coalgebra modality $(\oc, \delta, \varepsilon,  \Delta, \mathsf{e}, \nabla, \mathsf{u}, \mathsf{S})$, such that the diagrams in Appendix \ref{sec:Hopf-monoidal} commute. The first diagram says that $\mathsf{S}$ is a $\oc$-coalgebra morphism (\ref{diag:S-!coalg}), which is drawn as follows: 
\begin{align}\label{string:Hopf-coalg}
\begin{array}[c]{c}\begin{tikzpicture}
	\begin{pgfonlayer}{nodelayer}
		\node [style=port] (16) at (-2.25, 0.5) {};
		\node [style=component] (17) at (-2.25, 1.5) {$\delta$};
		\node [style=port] (18) at (-2.25, 3.75) {};
		\node [style=component] (19) at (-2.25, 2.75) {$\mathsf{S}$};
	\end{pgfonlayer}
	\begin{pgfonlayer}{edgelayer}
		\draw [style=wire] (17) to (16);
		\draw [style=wire] (18) to (19);
		\draw [style=wire] (19) to (17);
	\end{pgfonlayer}
\end{tikzpicture}
\end{array} =   \begin{array}[c]{c} \begin{tikzpicture}
	\begin{pgfonlayer}{nodelayer}
		\node [style=port] (3) at (0, 3.75) {};
		\node [style=component] (4) at (0, 3) {$\delta$};
		\node [style=port] (5) at (0, 0.5) {};
		\node [style=component] (7) at (0, 1.75) {$\mathsf{S}$};
		\node [style=none] (8) at (-0.5, 2.25) {};
		\node [style=none] (9) at (0.5, 2.25) {};
		\node [style=none] (10) at (-0.5, 1.25) {};
		\node [style=none] (11) at (0.5, 1.25) {};
	\end{pgfonlayer}
	\begin{pgfonlayer}{edgelayer}
		\draw [style=wire] (4) to (3);
		\draw [style=wire] (5) to (7);
		\draw [style=wire] (8.center) to (9.center);
		\draw [style=wire] (9.center) to (11.center);
		\draw [style=wire] (11.center) to (10.center);
		\draw [style=wire] (10.center) to (8.center);
		\draw [style=wire] (4) to (7);
	\end{pgfonlayer}
\end{tikzpicture}
\end{array}
\end{align}
The other required diagram states the compatibility between $\mathsf{m}$ and $\mathsf{S}$ (\ref{diag:Hopf-monoidal}), which is drawn out as follows: 
\begin{align}\label{string:Hopf-mon}
\begin{array}[c]{c}\begin{tikzpicture}
	\begin{pgfonlayer}{nodelayer}
		\node [style=port] (381) at (134.25, 1.5) {};
		\node [style=port] (382) at (133.75, 4.25) {};
		\node [style=component] (383) at (134.25, 2.5) {$\mathsf{m}$};
		\node [style=none] (384) at (133.75, 3) {};
		\node [style=none] (385) at (134.75, 3) {};
		\node [style=component] (386) at (134.75, 3.5) {$\mathsf{S}$};
		\node [style=port] (387) at (134.75, 4.25) {};
	\end{pgfonlayer}
	\begin{pgfonlayer}{edgelayer}
		\draw [style=wire] (381) to (383);
		\draw [style=wire] (383) to (384.center);
		\draw [style=wire] (383) to (385.center);
		\draw [style=wire] (384.center) to (382);
		\draw [style=wire] (385.center) to (386);
		\draw [style=wire] (387) to (386);
	\end{pgfonlayer}
\end{tikzpicture}
\end{array} = \begin{array}[c]{c} \begin{tikzpicture}
	\begin{pgfonlayer}{nodelayer}
		\node [style=port] (381) at (134.25, 1.5) {};
		\node [style=port] (382) at (134.75, 4.25) {};
		\node [style=component] (383) at (134.25, 3) {$\mathsf{m}$};
		\node [style=none] (384) at (134.75, 3.5) {};
		\node [style=none] (385) at (133.75, 3.5) {};
		\node [style=port] (387) at (133.75, 4.25) {};
		\node [style=component] (388) at (134.25, 2.25) {$\mathsf{S}$};
	\end{pgfonlayer}
	\begin{pgfonlayer}{edgelayer}
		\draw [style=wire] (383) to (384.center);
		\draw [style=wire] (383) to (385.center);
		\draw [style=wire] (384.center) to (382);
		\draw [style=wire] (387) to (385.center);
		\draw [style=wire] (383) to (388);
		\draw [style=wire] (388) to (381);
	\end{pgfonlayer}
\end{tikzpicture}
\end{array} = \begin{array}[c]{c}
\begin{tikzpicture}
	\begin{pgfonlayer}{nodelayer}
		\node [style=port] (381) at (134.25, 1.5) {};
		\node [style=port] (382) at (134.75, 4.25) {};
		\node [style=component] (383) at (134.25, 2.5) {$\mathsf{m}$};
		\node [style=none] (384) at (134.75, 3) {};
		\node [style=none] (385) at (133.75, 3) {};
		\node [style=component] (386) at (133.75, 3.5) {$\mathsf{S}$};
		\node [style=port] (387) at (133.75, 4.25) {};
	\end{pgfonlayer}
	\begin{pgfonlayer}{edgelayer}
		\draw [style=wire] (381) to (383);
		\draw [style=wire] (383) to (384.center);
		\draw [style=wire] (383) to (385.center);
		\draw [style=wire] (384.center) to (382);
		\draw [style=wire] (385.center) to (386);
		\draw [style=wire] (387) to (386);
	\end{pgfonlayer}
\end{tikzpicture}
\end{array}
\end{align}

We now discuss the relation between antipodes and our base additive category having \emph{negatives}. Following \cite{Blute2019,lemay2019lifting}, we say that an additive (symmetric monoidal) category $\mathbb{X}$ has \textbf{negatives} if each homset $\mathbb{X}(A,B)$ is in fact an Abelian group. Explicitly, this means that for every map $f: A \to B$, there is a map ${-f: A \to B}$ such that $f + (-f) = 0$. As is common, we write $f-g$ to express $f+ (-g)$. Of course, an additive (symmetric monoidal) category with negatives is equivalently a (symmetric monoidal) category enriched over Abelian groups. As such, this means that composition (and the monoidal product) preserves negatives: 
\begin{align}
f;(-g);h = -\left( f; g; h \right) && f;(g-h);k = f;g;k - f;h;k 
\end{align}
\begin{equation}\begin{gathered}\label{eq:neg-tensor}
f \otimes (-g) = -(f \otimes g) = (-f) \otimes g \\
f \otimes (g-h) = (f \otimes g) - (f \otimes h) \qquad (f-g) \otimes k  = (f \otimes k) - (g \otimes k) 
\end{gathered}\end{equation}
In particular, note that: 
\begin{align}\label{eq:-f}
f;-1_B = -f = -1_A;f
\end{align}
So an additive category having negatives is completely determined by whether identity morphisms have inverses with respect to the sum. In fact, an additive symmetric monoidal category having negatives is completely determined by if the identity for the monoidal unit has an inverse with respect to the sum \cite[Lemma 7.5]{lemay2019lifting}. Thus, this justifies that in string diagrams we draw: 
\begin{align}\label{string:-def2}
\begin{array}[c]{c}\begin{tikzpicture}
	\begin{pgfonlayer}{nodelayer}
		\node [style=port] (427) at (86.5, 1.25) {};
		\node [style=port] (428) at (86.5, -0.75) {};
		\node [style=component] (429) at (86.5, 0.25) {$-f$};
	\end{pgfonlayer}
	\begin{pgfonlayer}{edgelayer}
		\draw [style=wire] (427) to (429);
		\draw [style=wire] (429) to (428);
	\end{pgfonlayer}
\end{tikzpicture}
 \end{array} = - \begin{array}[c]{c}\begin{tikzpicture}
	\begin{pgfonlayer}{nodelayer}
		\node [style=port] (427) at (86.5, 1.25) {};
		\node [style=port] (428) at (86.5, -0.75) {};
		\node [style=component] (429) at (86.5, 0.25) {$f$};
	\end{pgfonlayer}
	\begin{pgfonlayer}{edgelayer}
		\draw [style=wire] (427) to (429);
		\draw [style=wire] (429) to (428);
	\end{pgfonlayer}
\end{tikzpicture}
 \end{array} &&
 \begin{array}[c]{c}\begin{tikzpicture}
	\begin{pgfonlayer}{nodelayer}
		\node [style=port] (427) at (86.5, 1.25) {};
		\node [style=port] (428) at (86.5, -0.75) {};
		\node [style=component] (429) at (86.5, 0.25) {$f-g$};
	\end{pgfonlayer}
	\begin{pgfonlayer}{edgelayer}
		\draw [style=wire] (427) to (429);
		\draw [style=wire] (429) to (428);
	\end{pgfonlayer}
\end{tikzpicture}
 \end{array} =  \begin{array}[c]{c}\begin{tikzpicture}
	\begin{pgfonlayer}{nodelayer}
		\node [style=port] (427) at (86.5, 1.25) {};
		\node [style=port] (428) at (86.5, -0.75) {};
		\node [style=component] (429) at (86.5, 0.25) {$f$};
	\end{pgfonlayer}
	\begin{pgfonlayer}{edgelayer}
		\draw [style=wire] (427) to (429);
		\draw [style=wire] (429) to (428);
	\end{pgfonlayer}
\end{tikzpicture}
 \end{array} -  \begin{array}[c]{c}\begin{tikzpicture}
	\begin{pgfonlayer}{nodelayer}
		\node [style=port] (427) at (86.5, 1.25) {};
		\node [style=port] (428) at (86.5, -0.75) {};
		\node [style=component] (429) at (86.5, 0.25) {$g$};
	\end{pgfonlayer}
	\begin{pgfonlayer}{edgelayer}
		\draw [style=wire] (427) to (429);
		\draw [style=wire] (429) to (428);
	\end{pgfonlayer}
\end{tikzpicture}
 \end{array}
\end{align}

We can now extend Thm \ref{thm:additive} and say that a monoidal Hopf coalgebra modality with a pre-codereliction induces an Abelian group enrichment. 

\begin{prop}\label{prop:neg} Let $(\oc, \delta, \varepsilon,  \Delta, \mathsf{e}, \mathsf{m}, \mathsf{m}_I, \nabla, \mathsf{u}, \mathsf{S})$ be a monoidal Hopf coalgebra modality with a pre-codereliction $\eta$ on a symmetric monoidal category $\mathbb{X}$. Then for the additive enrichment as defined in Thm \ref{thm:additive}, $\mathbb{X}$ is an additive symmetric monoidal category with negatives where for a map $f: A\to B$, $-f: A \to B$ is equal to following composite: 
\begin{align}
\xymatrixrowsep{1.75pc}\xymatrixcolsep{3pc}\xymatrix{A  \ar[r]^-{\eta_A} & \oc(A) \ar[r]^-{\mathsf{S}_A} & \oc(A) \ar[r]^-{\oc(f)} & \oc(B) \ar[r]^-{\varepsilon_B} & B    } 
\end{align}
which is drawn as follows: 
\begin{align}\label{string:-fdef}
- 
 
\end{align}
\end{prop}
\begin{proof} We need to show that $f-f=0$, or rather $f+(-f)=0$. To do so, similarly to the approach in the proof of Thm \ref{thm:additive}, we first need to show the following equality: 
\begin{align}\label{string:Sneg1}
%
\end{align}
So we compute: 
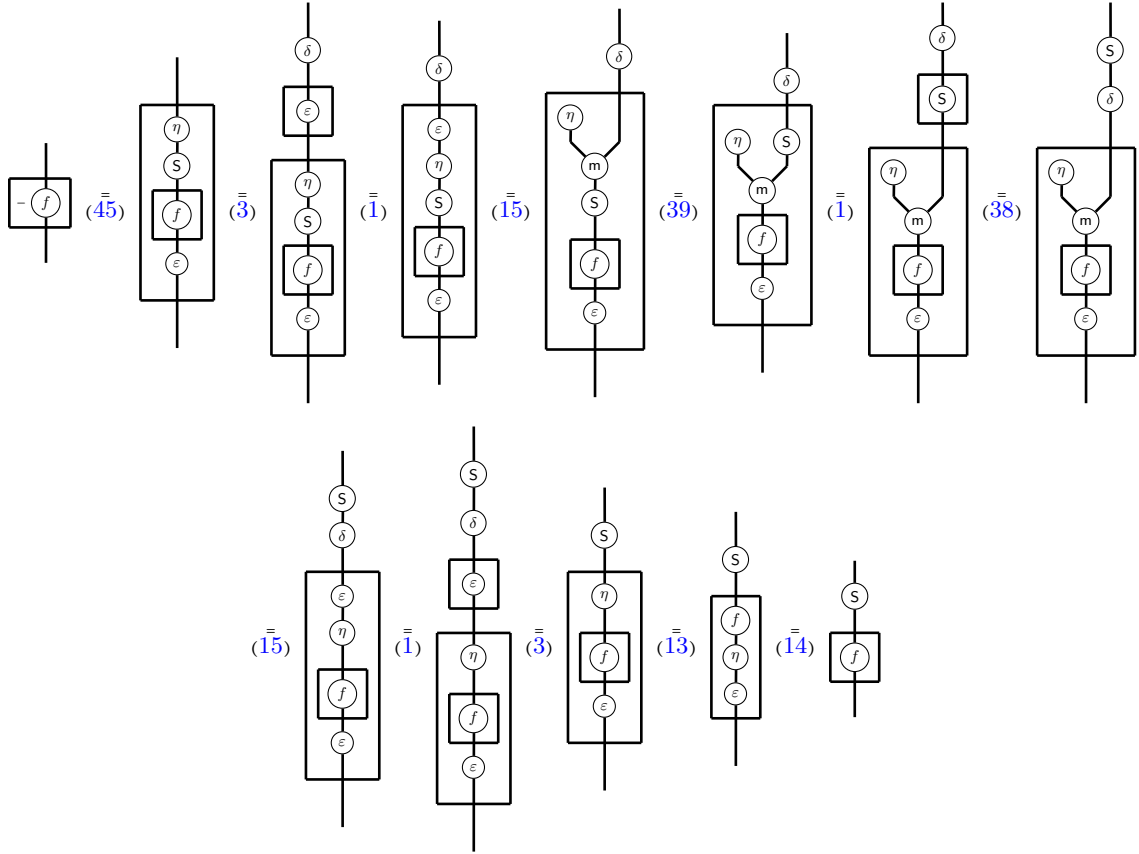
\begin{gather*}
%
\end{gather*}
So (\ref{string:Sneg1}) holds. Then using this identity and (\ref{string:!+0}), which recall we showed in the proof of Thm \ref{thm:additive}, we compute that: 
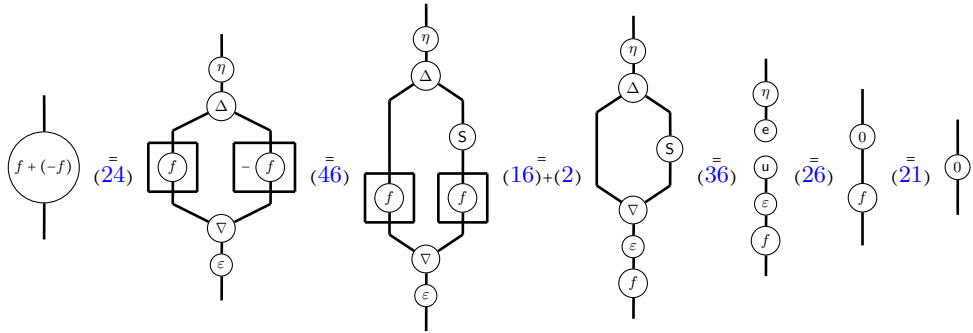
\begin{gather*}
 %
 
\end{gather*}
So we conclude that the base additive symmetric monoidal category has negatives. 
\end{proof}

Let us now revisit additive bialgebra modalities in additive symmetric monoidal categories with negatives. We first show that in the presence of negatives, convolution bialgebra modalities have antipodes. 

\begin{lem}\label{lemma:Hopf-2} Let $(\oc, \delta, \varepsilon,  \Delta, \mathsf{e}, \nabla, \mathsf{u})$ be a convolution bialgebra modality on an additive symmetric monoidal category with negatives. Then $(\oc, \delta, \varepsilon,  \Delta, \mathsf{e}, \nabla, \mathsf{u}, \mathsf{S})$ is a Hopf coalgebra modality where $\mathsf{S}_A := \oc(-1_A)$, which is drawn as follows: 
\begin{align}\label{string:Sdef}

\end{align}
Moreover, the diagrams in Appendix \ref{sec:lemma-Hopf-2} commute, which are drawn out as follows: 
\begin{align}\label{string:S+maps}
  %
 \end{align}
\end{lem}
\begin{proof} We first show that the equalities of (\ref{string:S+maps}) hold. The left most equality (which also implies naturality of $\mathsf{S}$) is immediate from functoriality and (\ref{eq:-f}), while the right most equality is immediate from naturality of $\varepsilon$. For the middle equalities, we use the convolution axiom: 
\begin{gather*}
%
 
\end{gather*}
Lastly, using the middle equality of (\ref{string:S+maps}), we can show that $\mathsf{S}$ is in fact an antipode\footnote{Note that by (co)commutativity, we need only show one of the equalities of (\ref{string:Hopf}).}:
\begin{gather*} %
\end{gather*}
So we conclude that $(\oc, \delta, \varepsilon,  \Delta, \mathsf{e}, \nabla, \mathsf{u}, \mathsf{S})$ is a Hopf coalgebra modality, as desired. 
\end{proof}

We now give an alternative proof of \cite[Prop 7.6]{lemay2019lifting}, which says that for an additive bialgebra modality the converse of Lemma \ref{lemma:Hopf-2} is true. 

\begin{lemC}[{\cite[Prop 7.6]{lemay2019lifting}}]\label{lemma:add-neg} Let $(\oc, \delta, \varepsilon,  \Delta, \mathsf{e}, \nabla, \mathsf{u})$ be an additive bialgebra modality on an additive symmetric monoidal category $\mathbb{X}$. Then $\mathbb{X}$ admits negatives if and only if $(\oc, \delta, \varepsilon,  \Delta, \mathsf{e}, \nabla, \mathsf{u})$ has an antipode $\mathsf{S}$.  
\end{lemC}
\begin{proof} The $\Rightarrow$ direction is Lemma \ref{lemma:Hopf-2}. For the $\Leftarrow$ direction, suppose that our additive bialgebra modality $(\oc, \delta, \varepsilon,  \Delta, \mathsf{e}, \nabla, \mathsf{u})$ has an antipode $\mathsf{S}$, so $(\oc, \delta, \varepsilon,  \Delta, \mathsf{e}, \nabla, \mathsf{u}, \mathsf{S})$ is a Hopf coalgebra modality. For a map $f$, define $-f$ as follows: 
\begin{align}
- \begin{array}[c]{c}\begin{tikzpicture}
	\begin{pgfonlayer}{nodelayer}
		\node [style=port] (427) at (86.5, 1.25) {};
		\node [style=port] (428) at (86.5, -0.75) {};
		\node [style=component] (429) at (86.5, 0.25) {$f$};
	\end{pgfonlayer}
	\begin{pgfonlayer}{edgelayer}
		\draw [style=wire] (427) to (429);
		\draw [style=wire] (429) to (428);
	\end{pgfonlayer}
\end{tikzpicture}
 \end{array} := \begin{array}[c]{c} \begin{tikzpicture}
	\begin{pgfonlayer}{nodelayer}
		\node [style=component] (404) at (144.25, 0) {$\mathsf{S}$};
		\node [style=component] (405) at (144.25, 0.75) {$\mathsf{m}$};
		\node [style=port] (436) at (145, 1.5) {};
		\node [style=port] (437) at (145, -1.5) {};
		\node [style=component] (438) at (145, 0) {$f$};
		\node [style=component] (439) at (144.25, -0.75) {$\varepsilon$};
	\end{pgfonlayer}
	\begin{pgfonlayer}{edgelayer}
		\draw [style=wire] (404) to (405);
		\draw [style=wire] (436) to (438);
		\draw [style=wire] (438) to (437);
		\draw [style=wire] (404) to (439);
	\end{pgfonlayer}
\end{tikzpicture}
 \end{array} 
\end{align}
We need to show that $f-f=0$. So we compute: 
\begin{gather*}
\begin{array}[c]{c}\begin{tikzpicture}
	\begin{pgfonlayer}{nodelayer}
		\node [style=port] (427) at (86.5, 1.25) {};
		\node [style=port] (428) at (86.5, -0.75) {};
		\node [style=component] (429) at (86.5, 0.25) {$f$};
	\end{pgfonlayer}
	\begin{pgfonlayer}{edgelayer}
		\draw [style=wire] (427) to (429);
		\draw [style=wire] (429) to (428);
	\end{pgfonlayer}
\end{tikzpicture}
 \end{array} -  \begin{array}[c]{c}\begin{tikzpicture}
	\begin{pgfonlayer}{nodelayer}
		\node [style=port] (427) at (86.5, 1.25) {};
		\node [style=port] (428) at (86.5, -0.75) {};
		\node [style=component] (429) at (86.5, 0.25) {$f$};
	\end{pgfonlayer}
	\begin{pgfonlayer}{edgelayer}
		\draw [style=wire] (427) to (429);
		\draw [style=wire] (429) to (428);
	\end{pgfonlayer}
\end{tikzpicture} \end{array} = \begin{array}[c]{c}\begin{tikzpicture}
	\begin{pgfonlayer}{nodelayer}
		\node [style=port] (427) at (86.5, 1.25) {};
		\node [style=port] (428) at (86.5, -0.75) {};
		\node [style=component] (429) at (86.5, 0.25) {$f$};
	\end{pgfonlayer}
	\begin{pgfonlayer}{edgelayer}
		\draw [style=wire] (427) to (429);
		\draw [style=wire] (429) to (428);
	\end{pgfonlayer}
\end{tikzpicture}
 \end{array} + \begin{array}[c]{c} \begin{tikzpicture}
	\begin{pgfonlayer}{nodelayer}
		\node [style=component] (404) at (144.25, 0) {$\mathsf{S}$};
		\node [style=component] (405) at (144.25, 0.75) {$\mathsf{m}$};
		\node [style=port] (436) at (145, 1.5) {};
		\node [style=port] (437) at (145, -1.5) {};
		\node [style=component] (438) at (145, 0) {$f$};
		\node [style=component] (439) at (144.25, -0.75) {$\varepsilon$};
	\end{pgfonlayer}
	\begin{pgfonlayer}{edgelayer}
		\draw [style=wire] (404) to (405);
		\draw [style=wire] (436) to (438);
		\draw [style=wire] (438) to (437);
		\draw [style=wire] (404) to (439);
	\end{pgfonlayer}
\end{tikzpicture}
 \end{array} \substack{=\\(\ref{string:comonad-monoidal})\\+(\ref{string:monoidalcomonoid})} \begin{array}[c]{c}\begin{tikzpicture}
	\begin{pgfonlayer}{nodelayer}
		\node [style=component] (443) at (147, 0.25) {$\mathsf{m}$};
		\node [style=port] (444) at (148.75, 1.5) {};
		\node [style=port] (445) at (148.75, -1.5) {};
		\node [style=component] (446) at (148.75, 0) {$f$};
		\node [style=component] (447) at (147, -0.5) {$\varepsilon$};
		\node [style=component] (448) at (148, 0.25) {$\mathsf{m}$};
		\node [style=component] (449) at (148, -0.5) {$\mathsf{e}$};
	\end{pgfonlayer}
	\begin{pgfonlayer}{edgelayer}
		\draw [style=wire] (444) to (446);
		\draw [style=wire] (446) to (445);
		\draw [style=wire] (448) to (449);
		\draw [style=wire] (443) to (447);
	\end{pgfonlayer}
\end{tikzpicture}
 \end{array}  + \begin{array}[c]{c}\begin{tikzpicture}
	\begin{pgfonlayer}{nodelayer}
		\node [style=component] (404) at (144.25, 0) {$\mathsf{S}$};
		\node [style=component] (405) at (144.25, 0.75) {$\mathsf{m}$};
		\node [style=port] (436) at (145, 1.5) {};
		\node [style=port] (437) at (145, -1.5) {};
		\node [style=component] (438) at (145, 0) {$f$};
		\node [style=component] (439) at (144.25, -0.75) {$\varepsilon$};
		\node [style=component] (440) at (143.25, 0.25) {$\mathsf{m}$};
		\node [style=component] (441) at (143.25, -0.5) {$\mathsf{e}$};
	\end{pgfonlayer}
	\begin{pgfonlayer}{edgelayer}
		\draw [style=wire] (404) to (405);
		\draw [style=wire] (436) to (438);
		\draw [style=wire] (438) to (437);
		\draw [style=wire] (404) to (439);
		\draw [style=wire] (440) to (441);
	\end{pgfonlayer}
\end{tikzpicture}
\end{array} \substack{=\\(\ref{strings:S-lemma})} \begin{array}[c]{c}\begin{tikzpicture}
	\begin{pgfonlayer}{nodelayer}
		\node [style=component] (442) at (148, 0) {$\mathsf{S}$};
		\node [style=component] (443) at (148, 0.75) {$\mathsf{m}$};
		\node [style=port] (444) at (148.75, 1.5) {};
		\node [style=port] (445) at (148.75, -1.5) {};
		\node [style=component] (446) at (148.75, 0) {$f$};
		\node [style=component] (447) at (148, -0.75) {$\mathsf{e}$};
		\node [style=component] (448) at (147, 0.25) {$\mathsf{m}$};
		\node [style=component] (449) at (147, -0.5) {$\varepsilon$};
	\end{pgfonlayer}
	\begin{pgfonlayer}{edgelayer}
		\draw [style=wire] (442) to (443);
		\draw [style=wire] (444) to (446);
		\draw [style=wire] (446) to (445);
		\draw [style=wire] (442) to (447);
		\draw [style=wire] (448) to (449);
	\end{pgfonlayer}
\end{tikzpicture}
 \end{array}  + \begin{array}[c]{c}\begin{tikzpicture}
	\begin{pgfonlayer}{nodelayer}
		\node [style=component] (404) at (144.25, 0) {$\mathsf{S}$};
		\node [style=component] (405) at (144.25, 0.75) {$\mathsf{m}$};
		\node [style=port] (436) at (145, 1.5) {};
		\node [style=port] (437) at (145, -1.5) {};
		\node [style=component] (438) at (145, 0) {$f$};
		\node [style=component] (439) at (144.25, -0.75) {$\varepsilon$};
		\node [style=component] (440) at (143.25, 0.25) {$\mathsf{m}$};
		\node [style=component] (441) at (143.25, -0.5) {$\mathsf{e}$};
	\end{pgfonlayer}
	\begin{pgfonlayer}{edgelayer}
		\draw [style=wire] (404) to (405);
		\draw [style=wire] (436) to (438);
		\draw [style=wire] (438) to (437);
		\draw [style=wire] (404) to (439);
		\draw [style=wire] (440) to (441);
	\end{pgfonlayer}
\end{tikzpicture}
\end{array} \\
\substack{=\\(\ref{string:ep-monoid})+(\ref{eq:add-tensor})} \begin{array}[c]{c} \begin{tikzpicture}
	\begin{pgfonlayer}{nodelayer}
		\node [style=component] (396) at (149, 4.5) {$\mathsf{m}$};
		\node [style=component] (398) at (149.5, 2.75) {$\nabla$};
		\node [style=none] (399) at (149, 3.25) {};
		\node [style=none] (400) at (150, 3.25) {};
		\node [style=component] (402) at (150, 3.75) {$\mathsf{S}$};
		\node [style=port] (444) at (151, 5) {};
		\node [style=port] (445) at (151, 1.5) {};
		\node [style=component] (446) at (151, 3) {$f$};
		\node [style=component] (450) at (150, 4.5) {$\mathsf{m}$};
		\node [style=component] (451) at (149.5, 2) {$\varepsilon$};
	\end{pgfonlayer}
	\begin{pgfonlayer}{edgelayer}
		\draw [style=wire] (398) to (399.center);
		\draw [style=wire] (398) to (400.center);
		\draw [style=wire] (399.center) to (396);
		\draw [style=wire] (400.center) to (402);
		\draw [style=wire] (444) to (446);
		\draw [style=wire] (446) to (445);
		\draw [style=wire] (450) to (402);
		\draw [style=wire] (398) to (451);
	\end{pgfonlayer}
\end{tikzpicture}
\end{array} \substack{=\\(\ref{string:monoidalcomonoid})} \begin{array}[c]{c} \begin{tikzpicture}
	\begin{pgfonlayer}{nodelayer}
		\node [style=component] (398) at (149.5, 2.75) {$\nabla$};
		\node [style=none] (399) at (149, 3.25) {};
		\node [style=none] (400) at (150, 3.25) {};
		\node [style=component] (402) at (150, 3.75) {$\mathsf{S}$};
		\node [style=port] (444) at (151, 6.25) {};
		\node [style=port] (445) at (151, 1.5) {};
		\node [style=component] (446) at (151, 3.75) {$f$};
		\node [style=component] (451) at (149.5, 2) {$\varepsilon$};
		\node [style=component] (452) at (149.5, 5.75) {$\mathsf{m}$};
		\node [style=component] (453) at (149.5, 5) {$\Delta$};
		\node [style=none] (454) at (149, 4.25) {};
		\node [style=none] (455) at (150, 4.25) {};
	\end{pgfonlayer}
	\begin{pgfonlayer}{edgelayer}
		\draw [style=wire] (398) to (399.center);
		\draw [style=wire] (398) to (400.center);
		\draw [style=wire] (400.center) to (402);
		\draw [style=wire] (444) to (446);
		\draw [style=wire] (446) to (445);
		\draw [style=wire] (398) to (451);
		\draw [style=wire] (453) to (454.center);
		\draw [style=wire] (453) to (455.center);
		\draw [style=wire] (452) to (453);
		\draw [style=wire] (455.center) to (402);
		\draw [style=wire] (454.center) to (399.center);
	\end{pgfonlayer}
\end{tikzpicture}
\end{array} \substack{=\\(\ref{string:Hopf})} \begin{array}[c]{c} \begin{tikzpicture}
	\begin{pgfonlayer}{nodelayer}
		\node [style=port] (444) at (151, 5.5) {};
		\node [style=port] (445) at (151, 1.5) {};
		\node [style=component] (446) at (151, 3.75) {$f$};
		\node [style=component] (451) at (149.5, 2) {$\varepsilon$};
		\node [style=component] (452) at (149.5, 5) {$\mathsf{m}$};
		\node [style=component] (453) at (149.5, 4) {$\mathsf{e}$};
		\node [style=component] (454) at (149.5, 3) {$\mathsf{u}$};
	\end{pgfonlayer}
	\begin{pgfonlayer}{edgelayer}
		\draw [style=wire] (444) to (446);
		\draw [style=wire] (446) to (445);
		\draw [style=wire] (452) to (453);
		\draw [style=wire] (454) to (451);
	\end{pgfonlayer}
\end{tikzpicture}
\end{array} \substack{=\\(\ref{string:monoidalcomonoid})+(\ref{string:ep-monoid})} \begin{array}[c]{c} \begin{tikzpicture}
	\begin{pgfonlayer}{nodelayer}
		\node [style=port] (455) at (152, -1) {};
		\node [style=port] (456) at (152, -3.5) {};
		\node [style=component] (457) at (152, -2.25) {$f$};
		\node [style=component] (458) at (151.25, -2.25) {$0$};
	\end{pgfonlayer}
	\begin{pgfonlayer}{edgelayer}
		\draw [style=wire] (455) to (457);
		\draw [style=wire] (457) to (456);
	\end{pgfonlayer}
\end{tikzpicture}
\end{array} \substack{=\\(\ref{eq:add-tensor})} \begin{array}[c]{c}\begin{tikzpicture}
	\begin{pgfonlayer}{nodelayer}
		\node [style=port] (427) at (86.5, 1.25) {};
		\node [style=port] (428) at (86.5, -0.75) {};
		\node [style=component] (429) at (86.5, 0.25) {$0$};
	\end{pgfonlayer}
	\begin{pgfonlayer}{edgelayer}
		\draw [style=wire] (427) to (429);
		\draw [style=wire] (429) to (428);
	\end{pgfonlayer}
\end{tikzpicture}
 \end{array} 
\end{gather*}
So we conclude that the base additive symmetric monoidal category has negatives. 
\end{proof}

As an immediate consequence, it follows that in the presence of negatives, any additive/monoidal bialgebra modality is a monoidal Hopf coalgebra modality. 

\begin{cor}\label{cor:mon-Hopf-S-2} Let $(\oc, \delta, \varepsilon,  \Delta, \mathsf{e}, \mathsf{m}, \mathsf{m}_I, \nabla, \mathsf{u})$ be a monoidal bialgebra modality on an additive symmetric monoidal category with negatives. Then $(\oc, \delta, \varepsilon,  \Delta, \mathsf{e}, \mathsf{m}, \mathsf{m}_I, \nabla, \mathsf{u}, \mathsf{S})$ is a monoidal Hopf coalgebra modality, where $\mathsf{S}_A := \oc(-1_A)$. 
\end{cor}
\begin{proof} By Lemma \ref{lemma:Hopf-2}, we know that $(\oc, \delta, \varepsilon,  \Delta, \mathsf{e}, \nabla, \mathsf{u}, \mathsf{S})$ is a Hopf coalgebra modality. Then since the antipode is of the form $\mathsf{S}_A := \oc(-1_A)$, (\ref{string:Hopf-coalg}) is precisely the naturality of $\delta$ (\ref{string:comonad-nat}), while (\ref{string:Hopf-mon}) follows from the naturality of $\mathsf{m}$ (\ref{string:monoidal-nat}) and (\ref{eq:neg-tensor}). Thus we conclude that $(\oc, \delta, \varepsilon,  \Delta, \mathsf{e}, \mathsf{m}, \mathsf{m}_I, \nabla, \mathsf{u}, \mathsf{S})$ is a monoidal Hopf coalgebra modality, as desired. 
\end{proof}

Moreover, since antipodes are unique, for any monoidal Hopf coalgebra modality on an additive symmetric monoidal category with negatives, its antipode must be given as in (\ref{string:Sdef}). 

\begin{cor}\label{cor:mon-Hopf-S-1} If $(\oc, \delta, \varepsilon,  \Delta, \mathsf{e}, \mathsf{m}, \mathsf{m}_I, \nabla, \mathsf{u}, \mathsf{S})$ is a monoidal Hopf coalgebra modality on an additive symmetric monoidal category with negatives, then $\mathsf{S}_A = \oc(-1_A)$. 
\end{cor}

\section{Differential Linear Categories Revisited}\label{sec:diff-lin}

In this final section, we bring everything together to provide a new axiomatization of a differential linear category. For an in-depth introduction to differential categories, we invite the reader to see \cite{blute2006differential,Blute2019,ehrhard2017introduction,garner2021cartesian,kerjean2023taylor}.  

It may be worthwhile to first recall the definition of a differential category in terms of a \emph{deriving transformation}, as was originally done in \cite{blute2006differential}. Borrowing terminology from \cite{garner2021cartesian}, a \textbf{differential modality} \cite[Def 4.3]{garner2021cartesian} on an additive symmetric monoidal category is a sextuple $(\oc, \delta, \varepsilon, \Delta, \mathsf{e}, \mathsf{d})$ consisting of a coalgebra modality $(\oc, \delta, \varepsilon, \Delta, \mathsf{e})$ and a natural transformation $\mathsf{d}_A: \oc(A) \otimes A \to \oc(A)$ such that the diagrams in Appendix \ref{sec:deriving} commute. The natural transformation $\mathsf{d}$ is called the \textbf{deriving transformation} and is drawn as follows\footnote{We note that we are drawing the deriving transformation differently then how it was drawn in \cite{Blute2019} for clarity.}: 
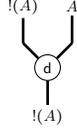
 \begin{align*}\begin{array}[c]{c} \begin{tikzpicture}
	\begin{pgfonlayer}{nodelayer}
		\node [style=object] (219) at (34.5, 1.75) {$\oc(A)$};
		\node [style=object] (220) at (34, 4) {$\oc(A)$};
		\node [style=component] (221) at (34.5, 2.75) {$\mathsf{d}$};
		\node [style=none] (222) at (34, 3.25) {};
		\node [style=none] (223) at (35, 3.25) {};
		\node [style=object] (224) at (35, 4) {$A$};
	\end{pgfonlayer}
	\begin{pgfonlayer}{edgelayer}
		\draw [style=wire] (219) to (221);
		\draw [style=wire] (221) to (222.center);
		\draw [style=wire] (221) to (223.center);
		\draw [style=wire] (223.center) to (224);
		\draw [style=wire] (222.center) to (220);
	\end{pgfonlayer}
\end{tikzpicture}
 \end{array} 
 \end{align*}
with naturality drawn as follows: 
 \begin{align}\label{string:d-nat}\begin{array}[c]{c} 
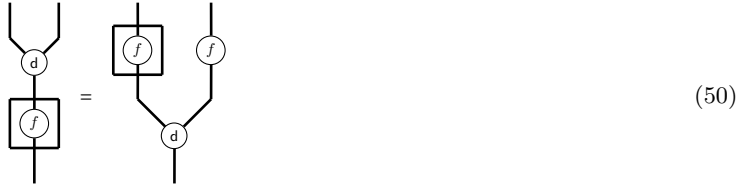
\begin{tikzpicture}
	\begin{pgfonlayer}{nodelayer}
		\node [style=component] (259) at (46.75, 1.75) {$f$};
		\node [style=none] (260) at (46.25, 1.25) {};
		\node [style=none] (261) at (47.25, 1.25) {};
		\node [style=none] (262) at (46.25, 2.25) {};
		\node [style=none] (263) at (47.25, 2.25) {};
		\node [style=port] (264) at (46.25, 4.25) {};
		\node [style=component] (265) at (46.75, 3) {$\mathsf{d}$};
		\node [style=none] (266) at (46.25, 3.5) {};
		\node [style=none] (267) at (47.25, 3.5) {};
		\node [style=port] (268) at (47.25, 4.25) {};
		\node [style=port] (269) at (46.75, 0.5) {};
	\end{pgfonlayer}
	\begin{pgfonlayer}{edgelayer}
		\draw [style=wire] (260.center) to (261.center);
		\draw [style=wire] (261.center) to (263.center);
		\draw [style=wire] (263.center) to (262.center);
		\draw [style=wire] (262.center) to (260.center);
		\draw [style=wire] (265) to (266.center);
		\draw [style=wire] (265) to (267.center);
		\draw [style=wire] (266.center) to (264);
		\draw [style=wire] (267.center) to (268);
		\draw [style=wire] (269) to (259);
		\draw [style=wire] (259) to (265);
	\end{pgfonlayer}
\end{tikzpicture}
 \end{array} = \begin{array}[c]{c}\begin{tikzpicture}
	\begin{pgfonlayer}{nodelayer}
		\node [style=component] (90) at (161.5, 3.25) {$f$};
		\node [style=none] (91) at (161, 2.75) {};
		\node [style=none] (92) at (162, 2.75) {};
		\node [style=none] (93) at (161, 3.75) {};
		\node [style=none] (94) at (162, 3.75) {};
		\node [style=port] (95) at (161.5, 4.25) {};
		\node [style=component] (96) at (162.25, 1.5) {$\mathsf{d}$};
		\node [style=none] (97) at (161.5, 2.25) {};
		\node [style=none] (98) at (163, 2.25) {};
		\node [style=port] (99) at (163, 4.25) {};
		\node [style=component] (100) at (163, 3.25) {$f$};
		\node [style=port] (105) at (162.25, 0.5) {};
	\end{pgfonlayer}
	\begin{pgfonlayer}{edgelayer}
		\draw [style=wire] (91.center) to (92.center);
		\draw [style=wire] (92.center) to (94.center);
		\draw [style=wire] (94.center) to (93.center);
		\draw [style=wire] (93.center) to (91.center);
		\draw [style=wire] (96) to (97.center);
		\draw [style=wire] (96) to (98.center);
		\draw [style=wire] (90) to (95);
		\draw [style=wire] (100) to (99);
		\draw [style=wire] (97.center) to (90);
		\draw [style=wire] (98.center) to (100);
		\draw [style=wire] (105) to (96);
	\end{pgfonlayer}
\end{tikzpicture}
 \end{array}
 \end{align}
Then the axioms of the deriving transformation\footnote{It is worth noting that the original definition of a deriving transformation in \cite[Def 2.5]{blute2006differential} only included \textbf{[D.1]} to \textbf{[D.4]} and not \textbf{[D.5]}. However, as explained in \cite[Sec 4]{Blute2019}, \textbf{[D.5]} is now commonly accepted as part of the definition of a deriving transformation.} (\ref{diag:deriving}) are drawn as follows: 
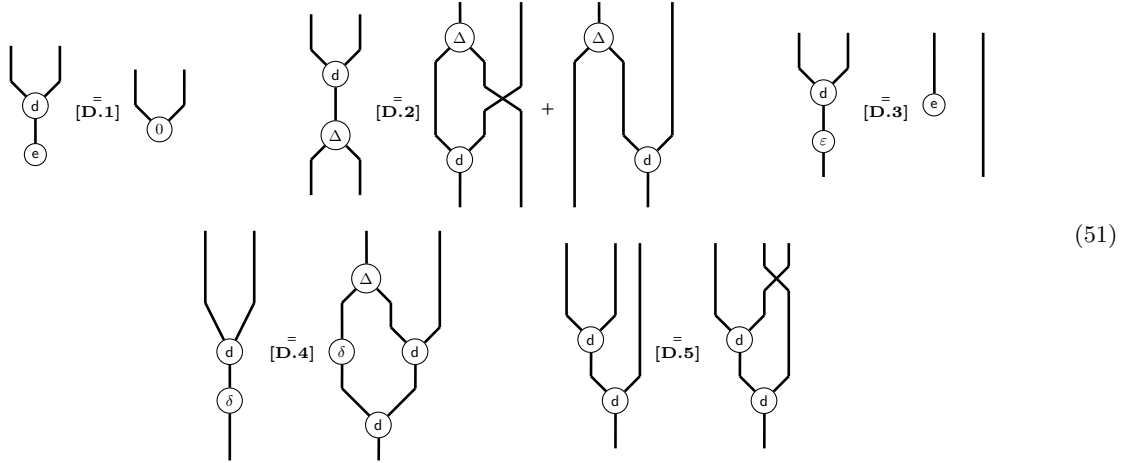
\begin{equation}\begin{gathered}
\begin{array}[c]{c} \begin{tikzpicture}
	\begin{pgfonlayer}{nodelayer}
		\node [style=component] (116) at (22.25, -0.25) {$\mathsf{e}$};
		\node [style=port] (117) at (21.75, 2) {};
		\node [style=component] (118) at (22.25, 0.75) {$\mathsf{d}$};
		\node [style=none] (119) at (21.75, 1.25) {};
		\node [style=none] (120) at (22.75, 1.25) {};
		\node [style=port] (121) at (22.75, 2) {};
	\end{pgfonlayer}
	\begin{pgfonlayer}{edgelayer}
		\draw [style=wire] (118) to (119.center);
		\draw [style=wire] (118) to (120.center);
		\draw [style=wire] (119.center) to (117);
		\draw [style=wire] (120.center) to (121);
		\draw [style=wire] (116) to (118);
	\end{pgfonlayer}
\end{tikzpicture}
 \end{array} \substack{=\\\text{\textbf{[D.1]}}} \begin{array}[c]{c} \begin{tikzpicture}
	\begin{pgfonlayer}{nodelayer}
		\node [style=port] (69) at (153.25, -8) {};
		\node [style=component] (70) at (153.75, -9.25) {$0$};
		\node [style=none] (71) at (153.25, -8.75) {};
		\node [style=none] (72) at (154.25, -8.75) {};
		\node [style=port] (73) at (154.25, -8) {};
	\end{pgfonlayer}
	\begin{pgfonlayer}{edgelayer}
		\draw [style=wire] (70) to (71.center);
		\draw [style=wire] (70) to (72.center);
		\draw [style=wire] (71.center) to (69);
		\draw [style=wire] (72.center) to (73);
	\end{pgfonlayer}
\end{tikzpicture}
 \end{array} \qquad \qquad \begin{array}[c]{c}\begin{tikzpicture}
	\begin{pgfonlayer}{nodelayer}
		\node [style=port] (29) at (17, -5.5) {};
		\node [style=component] (30) at (17.5, -4.25) {$\Delta$};
		\node [style=none] (31) at (17, -4.75) {};
		\node [style=none] (32) at (18, -4.75) {};
		\node [style=port] (33) at (18, -5.5) {};
		\node [style=port] (162) at (17, -1.75) {};
		\node [style=component] (163) at (17.5, -3) {$\mathsf{d}$};
		\node [style=none] (164) at (17, -2.5) {};
		\node [style=none] (165) at (18, -2.5) {};
		\node [style=port] (166) at (18, -1.75) {};
	\end{pgfonlayer}
	\begin{pgfonlayer}{edgelayer}
		\draw [style=wire] (30) to (31.center);
		\draw [style=wire] (30) to (32.center);
		\draw [style=wire] (32.center) to (33);
		\draw [style=wire] (31.center) to (29);
		\draw [style=wire] (163) to (164.center);
		\draw [style=wire] (163) to (165.center);
		\draw [style=wire] (164.center) to (162);
		\draw [style=wire] (165.center) to (166);
		\draw [style=wire] (163) to (30);
	\end{pgfonlayer}
\end{tikzpicture}
 \end{array} \substack{=\\\text{\textbf{[D.2]}}} \begin{array}[c]{c} \begin{tikzpicture}
	\begin{pgfonlayer}{nodelayer}
		\node [style=component] (106) at (164.5, -5.5) {$\Delta$};
		\node [style=none] (107) at (164, -6) {};
		\node [style=none] (108) at (165, -6) {};
		\node [style=port] (109) at (164.5, -4.75) {};
		\node [style=component] (110) at (164.5, -8) {$\mathsf{d}$};
		\node [style=none] (111) at (164, -7.5) {};
		\node [style=none] (112) at (165, -7.5) {};
		\node [style=port] (113) at (164.5, -9) {};
		\node [style=none] (114) at (165, -6.5) {};
		\node [style=none] (115) at (165, -7) {};
		\node [style=port] (116) at (165.75, -4.75) {};
		\node [style=port] (117) at (165.75, -9) {};
		\node [style=none] (118) at (165.75, -6.5) {};
		\node [style=none] (119) at (165.75, -7) {};
	\end{pgfonlayer}
	\begin{pgfonlayer}{edgelayer}
		\draw [style=wire] (106) to (107.center);
		\draw [style=wire] (106) to (108.center);
		\draw [style=wire] (106) to (109);
		\draw [style=wire] (110) to (111.center);
		\draw [style=wire] (110) to (112.center);
		\draw [style=wire] (113) to (110);
		\draw [style=wire] (107.center) to (111.center);
		\draw [style=wire] (112.center) to (115.center);
		\draw [style=wire] (108.center) to (114.center);
		\draw [style=wire] (116) to (118.center);
		\draw [style=wire] (118.center) to (115.center);
		\draw [style=wire] (114.center) to (119.center);
		\draw [style=wire] (119.center) to (117);
	\end{pgfonlayer}
\end{tikzpicture} \end{array} + \begin{array}[c]{c} \begin{tikzpicture}
	\begin{pgfonlayer}{nodelayer}
		\node [style=component] (106) at (164.25, -5.5) {$\Delta$};
		\node [style=none] (107) at (163.75, -6) {};
		\node [style=none] (108) at (164.75, -6) {};
		\node [style=port] (109) at (164.25, -4.75) {};
		\node [style=component] (110) at (165.25, -8) {$\mathsf{d}$};
		\node [style=none] (111) at (164.75, -7.5) {};
		\node [style=none] (112) at (165.75, -7.5) {};
		\node [style=port] (113) at (165.25, -9) {};
		\node [style=port] (116) at (165.75, -4.75) {};
		\node [style=port] (117) at (163.75, -9) {};
	\end{pgfonlayer}
	\begin{pgfonlayer}{edgelayer}
		\draw [style=wire] (106) to (107.center);
		\draw [style=wire] (106) to (108.center);
		\draw [style=wire] (106) to (109);
		\draw [style=wire] (110) to (111.center);
		\draw [style=wire] (110) to (112.center);
		\draw [style=wire] (113) to (110);
		\draw [style=wire] (116) to (112.center);
		\draw [style=wire] (108.center) to (111.center);
		\draw [style=wire] (107.center) to (117);
	\end{pgfonlayer}
\end{tikzpicture}
 \end{array} \qquad \qquad  \begin{array}[c]{c}\begin{tikzpicture}
	\begin{pgfonlayer}{nodelayer}
		\node [style=component] (0) at (22.25, -0.25) {$\varepsilon$};
		\node [style=port] (1) at (21.75, 2) {};
		\node [style=component] (2) at (22.25, 0.75) {$\mathsf{d}$};
		\node [style=none] (3) at (21.75, 1.25) {};
		\node [style=none] (4) at (22.75, 1.25) {};
		\node [style=port] (5) at (22.75, 2) {};
		\node [style=port] (6) at (22.25, -1) {};
	\end{pgfonlayer}
	\begin{pgfonlayer}{edgelayer}
		\draw [style=wire] (2) to (3.center);
		\draw [style=wire] (2) to (4.center);
		\draw [style=wire] (3.center) to (1);
		\draw [style=wire] (4.center) to (5);
		\draw [style=wire] (0) to (2);
		\draw [style=wire] (0) to (6);
	\end{pgfonlayer}
\end{tikzpicture}
 \end{array} \substack{=\\\text{\textbf{[D.3]}}} \begin{array}[c]{c}\begin{tikzpicture}
	\begin{pgfonlayer}{nodelayer}
		\node [style=port] (118) at (166.75, 2) {};
		\node [style=port] (119) at (167.75, 2) {};
		\node [style=component] (120) at (166.75, 0.5) {$\mathsf{e}$};
		\node [style=port] (122) at (167.75, -1) {};
	\end{pgfonlayer}
	\begin{pgfonlayer}{edgelayer}
		\draw [style=wire] (118) to (120);
		\draw [style=wire] (119) to (122);
	\end{pgfonlayer}
\end{tikzpicture}
\end{array} \\ 
\begin{array}[c]{c}\begin{tikzpicture}
	\begin{pgfonlayer}{nodelayer}
		\node [style=component] (259) at (58, -1) {$\delta$};
		\node [style=port] (264) at (57.5, 2.5) {};
		\node [style=component] (265) at (58, 0) {$\mathsf{d}$};
		\node [style=none] (266) at (57.5, 1) {};
		\node [style=none] (267) at (58.5, 1) {};
		\node [style=port] (268) at (58.5, 2.5) {};
		\node [style=port] (269) at (58, -2.25) {};
	\end{pgfonlayer}
	\begin{pgfonlayer}{edgelayer}
		\draw [style=wire] (265) to (266.center);
		\draw [style=wire] (265) to (267.center);
		\draw [style=wire] (266.center) to (264);
		\draw [style=wire] (267.center) to (268);
		\draw [style=wire] (269) to (259);
		\draw [style=wire] (259) to (265);
	\end{pgfonlayer}
\end{tikzpicture}
 \end{array} \substack{=\\\text{\textbf{[D.4]}}} \begin{array}[c]{c}\begin{tikzpicture}
	\begin{pgfonlayer}{nodelayer}
		\node [style=component] (123) at (171.25, -3.25) {$\delta$};
		\node [style=port] (124) at (171.75, -0.75) {};
		\node [style=port] (128) at (173.25, -0.75) {};
		\node [style=port] (129) at (172, -5.5) {};
		\node [style=component] (130) at (171.75, -1.75) {$\Delta$};
		\node [style=none] (131) at (171.25, -2.25) {};
		\node [style=none] (132) at (172.25, -2.25) {};
		\node [style=component] (133) at (172.75, -3.25) {$\mathsf{d}$};
		\node [style=none] (134) at (172.25, -2.75) {};
		\node [style=none] (135) at (173.25, -2.75) {};
		\node [style=component] (136) at (172, -4.75) {$\mathsf{d}$};
		\node [style=none] (137) at (171.25, -4) {};
		\node [style=none] (138) at (172.75, -4) {};
	\end{pgfonlayer}
	\begin{pgfonlayer}{edgelayer}
		\draw [style=wire] (130) to (131.center);
		\draw [style=wire] (130) to (132.center);
		\draw [style=wire] (124) to (130);
		\draw [style=wire] (131.center) to (123);
		\draw [style=wire] (133) to (134.center);
		\draw [style=wire] (133) to (135.center);
		\draw [style=wire] (132.center) to (134.center);
		\draw [style=wire] (128) to (135.center);
		\draw [style=wire] (136) to (137.center);
		\draw [style=wire] (136) to (138.center);
		\draw [style=wire] (123) to (137.center);
		\draw [style=wire] (133) to (138.center);
		\draw [style=wire] (136) to (129);
	\end{pgfonlayer}
\end{tikzpicture}
 \end{array} \qquad \qquad \begin{array}[c]{c}\begin{tikzpicture}
	\begin{pgfonlayer}{nodelayer}
		\node [style=port] (139) at (173, 1) {};
		\node [style=port] (140) at (173.5, 5.25) {};
		\node [style=component] (141) at (173, 2) {$\mathsf{d}$};
		\node [style=none] (142) at (173.5, 2.5) {};
		\node [style=none] (143) at (172.5, 2.5) {};
		\node [style=port] (144) at (173, 5.25) {};
		\node [style=component] (145) at (172.5, 3.25) {$\mathsf{d}$};
		\node [style=none] (146) at (173, 3.75) {};
		\node [style=none] (147) at (172, 3.75) {};
		\node [style=port] (148) at (172, 5.25) {};
	\end{pgfonlayer}
	\begin{pgfonlayer}{edgelayer}
		\draw [style=wire] (139) to (141);
		\draw [style=wire] (141) to (142.center);
		\draw [style=wire] (141) to (143.center);
		\draw [style=wire] (142.center) to (140);
		\draw [style=wire] (145) to (146.center);
		\draw [style=wire] (145) to (147.center);
		\draw [style=wire] (147.center) to (148);
		\draw [style=wire] (146.center) to (144);
		\draw [style=wire] (143.center) to (145);
	\end{pgfonlayer}
\end{tikzpicture}
 \end{array} \substack{=\\\text{\textbf{[D.5]}}} \begin{array}[c]{c} \begin{tikzpicture}
	\begin{pgfonlayer}{nodelayer}
		\node [style=port] (149) at (175.25, 1) {};
		\node [style=port] (150) at (175.75, 5.25) {};
		\node [style=component] (151) at (175.25, 2) {$\mathsf{d}$};
		\node [style=none] (152) at (175.75, 2.5) {};
		\node [style=none] (153) at (174.75, 2.5) {};
		\node [style=port] (154) at (175.25, 5.25) {};
		\node [style=component] (155) at (174.75, 3.25) {$\mathsf{d}$};
		\node [style=none] (156) at (175.25, 3.75) {};
		\node [style=none] (157) at (174.25, 3.75) {};
		\node [style=port] (158) at (174.25, 5.25) {};
		\node [style=none] (159) at (175.25, 4.75) {};
		\node [style=none] (160) at (175.75, 4.75) {};
		\node [style=none] (161) at (175.75, 4.25) {};
		\node [style=none] (162) at (175.25, 4.25) {};
	\end{pgfonlayer}
	\begin{pgfonlayer}{edgelayer}
		\draw [style=wire] (149) to (151);
		\draw [style=wire] (151) to (152.center);
		\draw [style=wire] (151) to (153.center);
		\draw [style=wire] (155) to (156.center);
		\draw [style=wire] (155) to (157.center);
		\draw [style=wire] (157.center) to (158);
		\draw [style=wire] (153.center) to (155);
		\draw [style=wire] (162.center) to (156.center);
		\draw [style=wire] (150) to (160.center);
		\draw [style=wire] (154) to (159.center);
		\draw [style=wire] (159.center) to (161.center);
		\draw [style=wire] (160.center) to (162.center);
		\draw [style=wire] (161.center) to (152.center);
	\end{pgfonlayer}
\end{tikzpicture}
 \end{array}
\end{gathered}\end{equation}
The axioms of a deriving transformation capture the fundamental identities from differential calculus. We call \textbf{[D.1]} the constant rule, \textbf{[D.2]} the Leibniz rule, \textbf{[D.2]} the linear rule, \textbf{[D.4]} the chain rule, and \textbf{[D.5]} the interchange rule. See \cite{blute2006differential,Blute2019} for more intuition on the axioms of a deriving transformation. Then a \textbf{differential category} \cite[Def 7]{Blute2019} is an additive symmetric monoidal category with a differential modality. 

A \textbf{monoidal differential modality} \cite[Def 4.3]{garner2021cartesian} on an additive symmetric monoidal category is an octuple $(\oc, \delta, \varepsilon, \Delta, \mathsf{e}, \mathsf{m}, \mathsf{m}_I, \mathsf{d})$ consisting of a monoidal coalgebra modality $(\oc, \delta, \varepsilon, \Delta, \mathsf{e}, \mathsf{m}, \mathsf{m}_I)$ and a differential modality $(\oc, \delta, \varepsilon, \Delta, \mathsf{e}, \mathsf{d})$. Now one might expect that the definition of a monoidal differential modality should also require some compatibility axioms between the deriving transformation and the monoidal structure, or even the induced bialgebra structure. However, one of the key results from \cite{Blute2019} is that these compatibilities come for free for a monoidal differential modality. Thus for a monoidal differential modality, the diagrams in Appendix \ref{sec:mon-deriving} commute \cite[Cor 7]{Blute2019}, which are drawn as follows: 
\begin{align}
\begin{array}[c]{c}\begin{tikzpicture}
	\begin{pgfonlayer}{nodelayer}
		\node [style=port] (227) at (38.25, -3) {};
		\node [style=port] (228) at (37.75, 0.5) {};
		\node [style=component] (229) at (38.25, -2) {$\nabla$};
		\node [style=none] (230) at (37.75, -1.5) {};
		\node [style=none] (231) at (38.75, -1.5) {};
		\node [style=port] (232) at (38.25, 0.5) {};
		\node [style=component] (233) at (38.75, -0.75) {$\mathsf{d}$};
		\node [style=none] (234) at (38.25, -0.25) {};
		\node [style=none] (235) at (39.25, -0.25) {};
		\node [style=port] (236) at (39.25, 0.5) {};
	\end{pgfonlayer}
	\begin{pgfonlayer}{edgelayer}
		\draw [style=wire] (227) to (229);
		\draw [style=wire] (229) to (230.center);
		\draw [style=wire] (229) to (231.center);
		\draw [style=wire] (230.center) to (228);
		\draw [style=wire] (233) to (234.center);
		\draw [style=wire] (233) to (235.center);
		\draw [style=wire] (235.center) to (236);
		\draw [style=wire] (234.center) to (232);
		\draw [style=wire] (231.center) to (233);
	\end{pgfonlayer}
\end{tikzpicture}
 \end{array} \substack{=\\\text{\textbf{[D.$\nabla$]}}} \begin{array}[c]{c}\begin{tikzpicture}
	\begin{pgfonlayer}{nodelayer}
		\node [style=port] (227) at (38.75, -3) {};
		\node [style=port] (228) at (39.25, 0.5) {};
		\node [style=component] (229) at (38.75, -2) {$\mathsf{d}$};
		\node [style=none] (230) at (39.25, -1.5) {};
		\node [style=none] (231) at (38.25, -1.5) {};
		\node [style=port] (232) at (38.75, 0.5) {};
		\node [style=component] (233) at (38.25, -0.75) {$\nabla$};
		\node [style=none] (234) at (38.75, -0.25) {};
		\node [style=none] (235) at (37.75, -0.25) {};
		\node [style=port] (236) at (37.75, 0.5) {};
	\end{pgfonlayer}
	\begin{pgfonlayer}{edgelayer}
		\draw [style=wire] (227) to (229);
		\draw [style=wire] (229) to (230.center);
		\draw [style=wire] (229) to (231.center);
		\draw [style=wire] (230.center) to (228);
		\draw [style=wire] (233) to (234.center);
		\draw [style=wire] (233) to (235.center);
		\draw [style=wire] (235.center) to (236);
		\draw [style=wire] (234.center) to (232);
		\draw [style=wire] (231.center) to (233);
	\end{pgfonlayer}
\end{tikzpicture}
 \end{array} && \begin{array}[c]{c}\begin{tikzpicture}
	\begin{pgfonlayer}{nodelayer}
		\node [style=port] (237) at (79, -7) {};
		\node [style=port] (238) at (78.5, -1.75) {};
		\node [style=component] (239) at (79, -6.25) {$\mathsf{m}$};
		\node [style=none] (240) at (78.5, -5.5) {};
		\node [style=none] (241) at (79.5, -5.5) {};
		\node [style=port] (367) at (80, -1.75) {};
		\node [style=component] (368) at (79.5, -4.25) {$\mathsf{d}$};
		\node [style=none] (369) at (80, -3.25) {};
		\node [style=none] (370) at (79, -3.25) {};
		\node [style=port] (371) at (79, -1.75) {};
	\end{pgfonlayer}
	\begin{pgfonlayer}{edgelayer}
		\draw [style=wire] (237) to (239);
		\draw [style=wire] (239) to (240.center);
		\draw [style=wire] (239) to (241.center);
		\draw [style=wire] (240.center) to (238);
		\draw [style=wire] (368) to (369.center);
		\draw [style=wire] (368) to (370.center);
		\draw [style=wire] (370.center) to (371);
		\draw [style=wire] (369.center) to (367);
		\draw [style=wire] (368) to (241.center);
	\end{pgfonlayer}
\end{tikzpicture}
 \end{array} \substack{=\\\text{\textbf{[D.m.r]}}} \begin{array}[c]{c}\begin{tikzpicture}
	\begin{pgfonlayer}{nodelayer}
		\node [style=component] (164) at (178, -6.5) {$\mathsf{d}$};
		\node [style=none] (165) at (179.25, -5.75) {};
		\node [style=none] (166) at (177.25, -5.75) {};
		\node [style=component] (167) at (177.25, -2.5) {$\Delta$};
		\node [style=none] (168) at (176.75, -3) {};
		\node [style=none] (169) at (177.75, -3) {};
		\node [style=port] (170) at (177.25, -1.75) {};
		\node [style=port] (171) at (179.25, -1.75) {};
		\node [style=none] (172) at (178.5, -3) {};
		\node [style=none] (173) at (179.25, -3) {};
		\node [style=component] (178) at (177.25, -5.25) {$\mathsf{m}$};
		\node [style=none] (179) at (176.75, -4.75) {};
		\node [style=none] (180) at (177.75, -4.75) {};
		\node [style=none] (182) at (178.5, -4.75) {};
		\node [style=none] (183) at (179.25, -4.25) {};
		\node [style=port] (184) at (178.5, -1.75) {};
		\node [style=port] (185) at (178, -7.25) {};
		\node [style=component] (186) at (177.75, -3.5) {$\varepsilon$};
		\node [style=none] (187) at (177.75, -4) {};
		\node [style=none] (188) at (178.5, -5.75) {};
		\node [style=none] (189) at (178.5, -4) {};
		\node [style=none] (190) at (177.75, -4.5) {};
		\node [style=none] (191) at (178.5, -4.5) {};
	\end{pgfonlayer}
	\begin{pgfonlayer}{edgelayer}
		\draw [style=wire] (164) to (165.center);
		\draw [style=wire] (164) to (166.center);
		\draw [style=wire] (167) to (168.center);
		\draw [style=wire] (167) to (169.center);
		\draw [style=wire] (170) to (167);
		\draw [style=wire] (178) to (179.center);
		\draw [style=wire] (178) to (180.center);
		\draw [style=wire] (168.center) to (179.center);
		\draw [style=wire] (173.center) to (183.center);
		\draw [style=wire] (178) to (166.center);
		\draw [style=wire] (171) to (173.center);
		\draw [style=wire] (164) to (185);
		\draw [style=wire] (169.center) to (186);
		\draw [style=wire] (186) to (187.center);
		\draw [style=wire] (182.center) to (188.center);
		\draw [style=wire] (188.center) to (164);
		\draw [style=wire] (183.center) to (165.center);
		\draw [style=wire] (184) to (172.center);
		\draw [style=wire] (172.center) to (189.center);
		\draw [style=wire] (187.center) to (191.center);
		\draw [style=wire] (191.center) to (182.center);
		\draw [style=wire] (189.center) to (190.center);
		\draw [style=wire] (190.center) to (180.center);
	\end{pgfonlayer}
\end{tikzpicture}
 \end{array} && \begin{array}[c]{c}\begin{tikzpicture}
	\begin{pgfonlayer}{nodelayer}
		\node [style=port] (237) at (79.5, -7) {};
		\node [style=port] (238) at (80, -1.75) {};
		\node [style=component] (239) at (79.5, -6.25) {$\mathsf{m}$};
		\node [style=none] (240) at (80, -5.5) {};
		\node [style=none] (241) at (79, -5.5) {};
		\node [style=port] (367) at (78.5, -1.75) {};
		\node [style=component] (368) at (79, -4.25) {$\mathsf{d}$};
		\node [style=none] (369) at (78.5, -3.25) {};
		\node [style=none] (370) at (79.5, -3.25) {};
		\node [style=port] (371) at (79.5, -1.75) {};
	\end{pgfonlayer}
	\begin{pgfonlayer}{edgelayer}
		\draw [style=wire] (237) to (239);
		\draw [style=wire] (239) to (240.center);
		\draw [style=wire] (239) to (241.center);
		\draw [style=wire] (240.center) to (238);
		\draw [style=wire] (368) to (369.center);
		\draw [style=wire] (368) to (370.center);
		\draw [style=wire] (370.center) to (371);
		\draw [style=wire] (369.center) to (367);
		\draw [style=wire] (368) to (241.center);
	\end{pgfonlayer}
\end{tikzpicture}
 \end{array} \substack{=\\\text{\textbf{[D.m.l]}}} \begin{array}[c]{c}\begin{tikzpicture}
	\begin{pgfonlayer}{nodelayer}
		\node [style=component] (164) at (178, -6.5) {$\mathsf{d}$};
		\node [style=none] (165) at (179.25, -5.75) {};
		\node [style=none] (166) at (177.25, -5.75) {};
		\node [style=component] (167) at (178.75, -2.5) {$\Delta$};
		\node [style=none] (168) at (178.25, -3) {};
		\node [style=none] (169) at (179.25, -3) {};
		\node [style=port] (170) at (178.75, -1.75) {};
		\node [style=port] (171) at (176.75, -1.75) {};
		\node [style=none] (172) at (177.75, -3) {};
		\node [style=none] (173) at (176.75, -3) {};
		\node [style=component] (178) at (177.25, -5.25) {$\mathsf{m}$};
		\node [style=none] (179) at (176.75, -4.75) {};
		\node [style=none] (180) at (177.75, -4.75) {};
		\node [style=none] (182) at (178.25, -3.5) {};
		\node [style=none] (183) at (176.75, -4.25) {};
		\node [style=port] (184) at (177.75, -1.75) {};
		\node [style=port] (185) at (178, -7.25) {};
		\node [style=component] (186) at (179.25, -4) {$\varepsilon$};
		\node [style=none] (187) at (179.25, -5.75) {};
		\node [style=none] (188) at (178.25, -5.75) {};
		\node [style=none] (189) at (177.75, -3.5) {};
		\node [style=none] (190) at (177.75, -4.25) {};
		\node [style=none] (191) at (178.25, -4.25) {};
	\end{pgfonlayer}
	\begin{pgfonlayer}{edgelayer}
		\draw [style=wire] (164) to (165.center);
		\draw [style=wire] (164) to (166.center);
		\draw [style=wire] (167) to (168.center);
		\draw [style=wire] (167) to (169.center);
		\draw [style=wire] (170) to (167);
		\draw [style=wire] (178) to (179.center);
		\draw [style=wire] (178) to (180.center);
		\draw [style=wire] (173.center) to (183.center);
		\draw [style=wire] (171) to (173.center);
		\draw [style=wire] (164) to (185);
		\draw [style=wire] (169.center) to (186);
		\draw [style=wire] (186) to (187.center);
		\draw [style=wire] (188.center) to (164);
		\draw [style=wire] (184) to (172.center);
		\draw [style=wire] (172.center) to (189.center);
		\draw [style=wire] (168.center) to (182.center);
		\draw [style=wire] (182.center) to (190.center);
		\draw [style=wire] (190.center) to (180.center);
		\draw [style=wire] (189.center) to (191.center);
		\draw [style=wire] (191.center) to (188.center);
		\draw [style=wire] (178) to (166.center);
		\draw [style=wire] (183.center) to (179.center);
	\end{pgfonlayer}
\end{tikzpicture}
 \end{array} 
\end{align}
where \textbf{[D.$\nabla$]} is called the $\nabla$-rule \cite[Def 8]{Blute2019} and \textbf{[D.m.r]} (resp. \textbf{[D.m.l]}) is called the right (resp. left) monoidal rule \cite[Sec 6]{Blute2019}. The monoidal rule was also called the strength rule in \cite[Def 4.2]{fiore2007differential}. We also note that in the string diagrams for the monoidal rule, the deriving transformation is of type $\mathsf{d}_{A \otimes B}: \oc(A \otimes B) \otimes A \otimes B \to \oc(A \otimes B)$, and is therefore drawn as taking in three inputs. Then a \textbf{differential linear category} \cite[Sec 6]{Blute2019} is an additive symmetric monoidal category with a monoidal differential modality. Examples of differential (linear) categories can be found in \cite[Sec 9]{Blute2019} and \cite[Ex 4.7]{garner2021cartesian}. In particular, every free exponential modality is a monoidal differential modality, and thus every additive Lafont category is a differential linear category \cite[Thm 21]{lemay2019lifting}. 

One of the fundamental theorems about differential linear categories is that the differential structure can also be axiomatized in terms of \emph{coderelictions}, which lines up more nicely with the Differential Linear Logic perspective \cite{ehrhard2017introduction}. Using the terminology of this paper, a \textbf{codereliction} \cite[Def 9]{Blute2019} for a monoidal coalgebra modality on an additive symmetric monoidal category is a pre-codereliction $\eta$ such that the diagrams from Appendix \ref{diag:coder} also commute, which are drawn as follows: 
\begin{align}
\begin{array}[c]{c}\begin{tikzpicture}
	\begin{pgfonlayer}{nodelayer}
		\node [style=component] (65) at (152.5, -0.25) {$\eta$};
		\node [style=component] (66) at (152.5, -1.25) {$\mathsf{e}$};
		\node [style=port] (67) at (152.5, 0.5) {};
	\end{pgfonlayer}
	\begin{pgfonlayer}{edgelayer}
		\draw [style=wire] (65) to (66);
		\draw [style=wire] (65) to (67);
	\end{pgfonlayer}
\end{tikzpicture}
 \end{array} \substack{=\\\text{\textbf{[cd.1]}}}  \begin{array}[c]{c}\begin{tikzpicture}
	\begin{pgfonlayer}{nodelayer}
		\node [style=component] (63) at (149.5, -1) {$0$};
		\node [style=port] (64) at (149.5, 0.5) {};
	\end{pgfonlayer}
	\begin{pgfonlayer}{edgelayer}
		\draw [style=wire] (63) to (64);
	\end{pgfonlayer}
\end{tikzpicture}
 \end{array} &&  \begin{array}[c]{c}\begin{tikzpicture}
	\begin{pgfonlayer}{nodelayer}
		\node [style=component] (68) at (153.75, -7) {$\eta$};
		\node [style=port] (69) at (153.25, -9.25) {};
		\node [style=component] (70) at (153.75, -8) {$\Delta$};
		\node [style=none] (71) at (153.25, -8.5) {};
		\node [style=none] (72) at (154.25, -8.5) {};
		\node [style=port] (73) at (154.25, -9.25) {};
		\node [style=port] (74) at (153.75, -6.25) {};
	\end{pgfonlayer}
	\begin{pgfonlayer}{edgelayer}
		\draw [style=wire] (70) to (71.center);
		\draw [style=wire] (70) to (72.center);
		\draw [style=wire] (71.center) to (69);
		\draw [style=wire] (72.center) to (73);
		\draw [style=wire] (68) to (70);
		\draw [style=wire] (68) to (74);
	\end{pgfonlayer}
\end{tikzpicture}
 \end{array} \substack{=\\\text{\textbf{[cd.2]}}} \begin{array}[c]{c}\begin{tikzpicture}
	\begin{pgfonlayer}{nodelayer}
		\node [style=port] (75) at (156.25, -1) {};
		\node [style=port] (76) at (155.25, -1) {};
		\node [style=component] (77) at (156.25, 0.5) {$\mathsf{u}$};
		\node [style=component] (78) at (155.25, 0.5) {$\eta$};
		\node [style=port] (79) at (155.25, 2) {};
	\end{pgfonlayer}
	\begin{pgfonlayer}{edgelayer}
		\draw [style=wire] (75) to (77);
		\draw [style=wire] (76) to (78);
		\draw [style=wire] (78) to (79);
	\end{pgfonlayer}
\end{tikzpicture}
 \end{array} + \begin{array}[c]{c}\begin{tikzpicture}
	\begin{pgfonlayer}{nodelayer}
		\node [style=port] (75) at (155.25, -1) {};
		\node [style=port] (76) at (156.25, -1) {};
		\node [style=component] (77) at (155.25, 0.5) {$\mathsf{u}$};
		\node [style=component] (78) at (156.25, 0.5) {$\eta$};
		\node [style=port] (79) at (156.25, 2) {};
	\end{pgfonlayer}
	\begin{pgfonlayer}{edgelayer}
		\draw [style=wire] (75) to (77);
		\draw [style=wire] (76) to (78);
		\draw [style=wire] (78) to (79);
	\end{pgfonlayer}
\end{tikzpicture}
 \end{array} &&  \begin{array}[c]{c}\begin{tikzpicture}
	\begin{pgfonlayer}{nodelayer}
		\node [style=port] (16) at (-2.25, 0.5) {};
		\node [style=component] (17) at (-2.25, 1.5) {$\delta$};
		\node [style=port] (18) at (-2.25, 3.75) {};
		\node [style=component] (19) at (-2.25, 2.75) {$\eta$};
	\end{pgfonlayer}
	\begin{pgfonlayer}{edgelayer}
		\draw [style=wire] (17) to (16);
		\draw [style=wire] (18) to (19);
		\draw [style=wire] (19) to (17);
	\end{pgfonlayer}
\end{tikzpicture}
\end{array} \substack{=\\\text{\textbf{[cd.4]}}}  \begin{array}[c]{c}\begin{tikzpicture}
	\begin{pgfonlayer}{nodelayer}
		\node [style=component] (80) at (159, -4.25) {$\delta$};
		\node [style=component] (82) at (160, -4.25) {$\eta$};
		\node [style=port] (83) at (160, -2.75) {};
		\node [style=component] (84) at (159.5, -5.25) {$\nabla$};
		\node [style=none] (85) at (159, -4.75) {};
		\node [style=none] (86) at (160, -4.75) {};
		\node [style=port] (87) at (159.5, -6) {};
		\node [style=component] (88) at (159, -3.5) {$\mathsf{u}$};
		\node [style=component] (89) at (160, -3.5) {$\eta$};
	\end{pgfonlayer}
	\begin{pgfonlayer}{edgelayer}
		\draw [style=wire] (84) to (85.center);
		\draw [style=wire] (84) to (86.center);
		\draw [style=wire] (80) to (85.center);
		\draw [style=wire] (82) to (86.center);
		\draw [style=wire] (84) to (87);
		\draw [style=wire] (88) to (80);
		\draw [style=wire] (89) to (82);
		\draw [style=wire] (89) to (83);
	\end{pgfonlayer}
\end{tikzpicture}\end{array} 
\end{align}
For a codereliction, the axiom \textbf{[cd.1]} is called the constant rule, \textbf{[cd.2]} is called the Leibniz (or product) rule, and \textbf{[cd.4]} is called the chain rule\footnote{In \cite[Sec 5]{Blute2019}, this version of \textbf{[cd.4]} was called the \emph{alternative} chain rule. However here we follow Fiore's approach \cite{fiore2007differential} and the more recent convention of calling \textbf{[cd.4]} the chain rule for the codereliction.}. 

We now briefly review the bijective correspondence between coderelictions and deriving transformations, that is, how we can build one from the other, as given in Appendix \ref{sec:d-eta}. So let $(\oc, \delta, \varepsilon, \Delta, \mathsf{e}, \mathsf{m}, \mathsf{m}_I)$ be a monoidal coalgebra modality on an additive symmetric monoidal category, with induced additive bialgebra modality $(\oc, \delta, \varepsilon, \Delta, \mathsf{e}, \nabla, \mathsf{u})$. Then given a codereliction $\eta$, define $\mathsf{d}$ as follows: 
\begin{align}
\begin{array}[c]{c}\begin{tikzpicture}
	\begin{pgfonlayer}{nodelayer}
		\node [style=port] (192) at (158.75, 3.75) {};
		\node [style=component] (193) at (159.5, 1.5) {$\mathsf{d}$};
		\node [style=none] (194) at (158.75, 2.25) {};
		\node [style=none] (195) at (160.25, 2.25) {};
		\node [style=port] (196) at (160.25, 3.75) {};
		\node [style=port] (198) at (159.5, 0.5) {};
	\end{pgfonlayer}
	\begin{pgfonlayer}{edgelayer}
		\draw [style=wire] (193) to (194.center);
		\draw [style=wire] (193) to (195.center);
		\draw [style=wire] (198) to (193);
		\draw [style=wire] (192) to (194.center);
		\draw [style=wire] (196) to (195.center);
	\end{pgfonlayer}
\end{tikzpicture}
\end{array} := 
\begin{array}[c]{c} \begin{tikzpicture}
	\begin{pgfonlayer}{nodelayer}
		\node [style=port] (95) at (161.5, 3.75) {};
		\node [style=component] (96) at (162.25, 1.5) {$\nabla$};
		\node [style=none] (97) at (161.5, 2.25) {};
		\node [style=none] (98) at (163, 2.25) {};
		\node [style=port] (99) at (163, 3.75) {};
		\node [style=component] (100) at (163, 3) {$\eta$};
		\node [style=port] (105) at (162.25, 0.5) {};
	\end{pgfonlayer}
	\begin{pgfonlayer}{edgelayer}
		\draw [style=wire] (96) to (97.center);
		\draw [style=wire] (96) to (98.center);
		\draw [style=wire] (100) to (99);
		\draw [style=wire] (98.center) to (100);
		\draw [style=wire] (105) to (96);
		\draw [style=wire] (95) to (97.center);
	\end{pgfonlayer}
\end{tikzpicture}
\end{array} 
\end{align}
Then $\mathsf{d}$ is a deriving transformation \cite[Thm 3]{Blute2019}. Conversely, given instead a deriving transformation $\mathsf{d}$, define $\eta$ as follows: 
\begin{align}
\begin{array}[c]{c} \begin{tikzpicture}
	\begin{pgfonlayer}{nodelayer}
		\node [style=port] (76) at (156.25, -1) {};
		\node [style=component] (78) at (156.25, 0.5) {$\eta$};
		\node [style=port] (79) at (156.25, 2) {};
	\end{pgfonlayer}
	\begin{pgfonlayer}{edgelayer}
		\draw [style=wire] (76) to (78);
		\draw [style=wire] (78) to (79);
	\end{pgfonlayer}
\end{tikzpicture}
\end{array} := 
\begin{array}[c]{c} \begin{tikzpicture}
	\begin{pgfonlayer}{nodelayer}
		\node [style=port] (95) at (163, 3.75) {};
		\node [style=component] (96) at (162.25, 1.5) {$\mathsf{d}$};
		\node [style=none] (97) at (163, 2.25) {};
		\node [style=none] (98) at (161.5, 2.25) {};
		\node [style=component] (100) at (161.5, 3) {$\mathsf{u}$};
		\node [style=port] (105) at (162.25, 0.5) {};
	\end{pgfonlayer}
	\begin{pgfonlayer}{edgelayer}
		\draw [style=wire] (96) to (97.center);
		\draw [style=wire] (96) to (98.center);
		\draw [style=wire] (98.center) to (100);
		\draw [style=wire] (105) to (96);
		\draw [style=wire] (95) to (97.center);
	\end{pgfonlayer}
\end{tikzpicture}
\end{array} 
\end{align}
Then $\eta$ is a codereliciton \cite[Thm 3]{Blute2019}. Moreover, these constructions are inverses of each other. As such for monoidal coalgebra modalities, there is a bijective correspondence between coderelictions and deriving transformations \cite[Thm 4]{Blute2019}. Thus, a differential linear category is equivalently an additive symmetric monoidal category with a monoidal coalgebra modality that comes equipped with a codereliction.

Now it turns out that two of the axioms of a codereliction are in fact redundant, that is, provable. Indeed, in an additive enriched setting, the constant rule \textbf{[cd.1]} is provable from simply naturality \cite[Lemma 6]{Blute2019}, while the product rule \textbf{[cd.2]} is provable from naturality and the linear rule \textbf{[cd.3]}. Therefore, a codereliction is equivalently simply a natural transformation which satisfies the linear rule \textbf{[cd.3]}, the chain rule \textbf{[cd.4]}, and either monoidal rule \textbf{[cd.m.r/l]} \cite[Cor 5]{Blute2019}. In other words, a codereliction is a pre-codereliction which also satisfies the chain rule \textbf{[cd.4]}. However, observe that the three remaining axioms only require the monoidal bialgebra modality structural maps and no longer seems to require sums or zeros. Thus we can technically define a codereliction for a monoidal bialgebra modality on a symmetric monoidal category to be a pre-codereliction on the underlying monoidal coalgebra modality which also satisfies the chain rule \textbf{[cd.4]}. However, it follows from Thm \ref{thm:additive} that this gives us an equivalent characterization of a differential linear category without including additive enrichment in the axiomatization. So bringing all of this together, we obtain: 

\begin{thm}\label{thm:diff-lin} A differential linear category is equivalently a symmetric monoidal category with a monoidal bialgebra modality equipped with a codereliction. 
\end{thm}

It may also be of interest to ask how, in the presence of negatives, the deriving transformation/codereliction interacts with the antipode. 

\begin{prop}\label{prop:der-neg} For a monoidal differential modality $(\oc, \delta, \varepsilon, \Delta, \mathsf{e}, \mathsf{m}, \mathsf{m}_I, \mathsf{d})$, with induced codereliction $\eta$, on an additive symmetirc monoidal category with negatives, for the induced antipode $\mathsf{S}$, the diagrams in Appendix \ref{prop:der-neg} commute, which are drawn as follows:
\begin{align}
\begin{array}[c]{c}\begin{tikzpicture}
	\begin{pgfonlayer}{nodelayer}
		\node [style=component] (362) at (208.5, 1.5) {$\mathsf{S}$};
		\node [style=port] (363) at (208, 3.75) {};
		\node [style=component] (364) at (208.5, 2.5) {$\mathsf{d}$};
		\node [style=none] (365) at (208, 3) {};
		\node [style=none] (366) at (209, 3) {};
		\node [style=port] (367) at (209, 3.75) {};
		\node [style=port] (368) at (208.5, 0.75) {};
	\end{pgfonlayer}
	\begin{pgfonlayer}{edgelayer}
		\draw [style=wire] (364) to (365.center);
		\draw [style=wire] (364) to (366.center);
		\draw [style=wire] (365.center) to (363);
		\draw [style=wire] (366.center) to (367);
		\draw [style=wire] (362) to (364);
		\draw [style=wire] (362) to (368);
	\end{pgfonlayer}
\end{tikzpicture}
 \end{array} = - \begin{array}[c]{c}\begin{tikzpicture}
	\begin{pgfonlayer}{nodelayer}
		\node [style=port] (213) at (205.75, 3.75) {};
		\node [style=port] (356) at (206.75, 3.75) {};
		\node [style=component] (357) at (206.25, 1.75) {$\mathsf{d}$};
		\node [style=none] (358) at (206.75, 2.5) {};
		\node [style=none] (359) at (205.75, 2.5) {};
		\node [style=component] (360) at (205.75, 3) {$\mathsf{S}$};
		\node [style=port] (361) at (206.25, 0.75) {};
	\end{pgfonlayer}
	\begin{pgfonlayer}{edgelayer}
		\draw [style=wire] (357) to (358.center);
		\draw [style=wire] (357) to (359.center);
		\draw [style=wire] (359.center) to (360);
		\draw [style=wire] (361) to (357);
		\draw [style=wire] (356) to (358.center);
		\draw [style=wire] (213) to (360);
	\end{pgfonlayer}
\end{tikzpicture}
 \end{array} && \begin{array}[c]{c}\begin{tikzpicture}
	\begin{pgfonlayer}{nodelayer}
		\node [style=port] (16) at (-2.25, 0.5) {};
		\node [style=component] (17) at (-2.25, 1.5) {$\mathsf{S}$};
		\node [style=port] (18) at (-2.25, 3.75) {};
		\node [style=component] (19) at (-2.25, 2.75) {$\eta$};
	\end{pgfonlayer}
	\begin{pgfonlayer}{edgelayer}
		\draw [style=wire] (17) to (16);
		\draw [style=wire] (18) to (19);
		\draw [style=wire] (19) to (17);
	\end{pgfonlayer}
\end{tikzpicture}
\end{array} = - \begin{array}[c]{c}\begin{tikzpicture}
	\begin{pgfonlayer}{nodelayer}
		\node [style=port] (427) at (86.5, 1.25) {};
		\node [style=port] (428) at (86.5, -0.75) {};
		\node [style=component] (429) at (86.5, 0.25) {$\eta$};
	\end{pgfonlayer}
	\begin{pgfonlayer}{edgelayer}
		\draw [style=wire] (427) to (429);
		\draw [style=wire] (429) to (428);
	\end{pgfonlayer}
\end{tikzpicture}
 \end{array}
\end{align}
\end{prop}
\begin{proof} Recall from Cor \ref{cor:mon-Hopf-S-2} that the antipode is of the form $\mathsf{S}_A := \oc(-1_A)$. Then the equality on the left follows from naturality of $\mathsf{d}$ (\ref{string:d-nat}) and (\ref{eq:neg-tensor}), while the equality on the right follows from naturality of $\eta$ (\ref{string:coder-nat}). 
\end{proof}

We conclude this paper by discussing uniqueness of coderelictions. In \cite[Thm 21]{lemay2021coderelictions}, it was shown that for free exponential modalities, coderelictions were unique. Here however, we can now extend this and state that for any monoidal coalgebra modality, coderelictions are unique. 

\begin{thm}\label{thm:coder-unique} For a monoidal coalgebra modality (equivalently an additive bialgebra modality or a monoidal bialgebra modality) on additive symmetric monoidal category, if a codereliction exists, then it is unique.
\end{thm}
\begin{proof} In Prop \ref{prop:precoderunique} we showed that a pre-codereliciton, if it exists, was unique. Thus since a codereliction is a pre-codereliction by definition, it follows that if one exists, it must be unique. 
\end{proof}

 Moreover, since there is a bijective correspondence between coderelictions and deriving transformations, we also have uniqueness of deriving transformation for monoidal coalgebra modalities. 

\begin{cor} For a monoidal coalgebra modality on additive symmetric monoidal category, if a deriving transformation exists, then it is unique.
\end{cor}

Thus we are properly justified in saying that in categorical models of Differential Linear Logic, there really is only one way to differentiate smooth maps. However, it is important to stress that the above corollary does not necessarily extend to mere coalgebra modalities. Indeed, it is still an open question whether or not deriving transformations are unique for non-monoidal coalgebra modalities.

\section*{Acknowledgment}
  \noindent The author would like to thank Thomas Ehrhard, Richard Blute, and Robin Cockett for useful discussions and their support of this research. This material is based upon work supported by the AFOSR under award number FA9550-24-1-0008.


\bibliographystyle{alphaurl}      
\bibliography{reference}   

\appendix
\section{Commutative Diagrams}\label{sec:appendix}

\subsection{Diagrams for Comonads}\label{sec:diagramsforcomonad}
\begin{equation}\begin{gathered}\label{diag:comonad}\xymatrixrowsep{1.75pc}\xymatrixcolsep{5pc}\xymatrix{ 
        \oc(A)  \ar[r]^-{\delta_A} \ar[d]_-{\delta_A} \ar@{=}[dr]^-{}& \oc \oc(A) \ar[d]^-{\varepsilon_{\oc(A)}}  & \oc(A)  \ar[r]^-{\delta_A}\ar[d]_-{\delta_A} & \oc\oc(A)  \ar[d]^-{\delta_{\oc(A)}}\\
        \oc \oc(A) \ar[r]_-{\oc(\varepsilon_A)} & \oc(A)  & \oc \oc(A) \ar[r]_-{\oc(\delta_A)} & \oc \oc \oc(A)} \end{gathered}\end{equation}

\subsection{Diagrams for Coalgebra Modalities}\label{sec:diagramsforcoalgmod}

\begin{equation}\begin{gathered}\label{diag:comonoid}\xymatrixrowsep{1.75pc}\xymatrixcolsep{3pc}\xymatrix{\oc(A)  \ar[r]^-{\Delta_A} \ar[d]_-{\Delta_A} & \oc(A)  \otimes  \oc(A) \ar[d]^-{\Delta_A \otimes 1_{\oc(A)}} & \oc(A) \ar@{=}[dr] \ar[r]^-{\Delta_A}  \ar[d]_-{\Delta_A}   & \oc(A)  \otimes  \oc(A) \ar[d]^-{1_{\oc(A)} \otimes \mathsf{e}_A}  \\
      \oc(A) \otimes  \oc(A) \ar[r]_-{1_{\oc(A)} \otimes \Delta_A} & \oc(A)  \otimes  \oc(A)  \otimes  \oc(A)& \oc(A) \otimes  \oc(A) \ar[r]_-{\mathsf{e}_A \otimes 1_{\oc(A)}}  & \oc(A) } \\
      \xymatrixrowsep{1.75pc}\xymatrixcolsep{5pc}\xymatrix{ \oc(A)  \ar[r]^-{\Delta_A}  \ar[dr]_-{\Delta_A} & \oc(A)  \otimes  \oc(A)  \ar[d]^-{\sigma_{\oc(A),\oc(A)}}\\
 & \oc(A)  \otimes  \oc(A) } \end{gathered}\end{equation} 
\begin{equation}\begin{gathered}\label{diag:deltacomonoid}\xymatrixrowsep{1.75pc}\xymatrixcolsep{5pc}\xymatrix{\oc(A) \ar[r]^-{\delta_A}  \ar[d]_-{\Delta_A} & \oc\oc(A) \ar[d]^-{\Delta_{\oc(A)}} & \oc(A) \ar[r]^-{\delta_A} \ar[dr]_-{\mathsf{e}_A} & \oc\oc(A) \ar[d]^-{\mathsf{e}_{\oc(A)}}   \\
      \oc(A) \otimes \oc(A) \ar[r]_-{\delta_A \otimes \delta_A} & \oc \oc(A) \otimes \oc \oc(A) & & I }  \end{gathered}\end{equation}

\subsection{Diagrams for Symmetric Monoidal Endofunctor}\label{sec:diagramsforsmc}

\begin{equation}\begin{gathered}\label{diag:smendo} \xymatrixrowsep{1.75pc}\xymatrixcolsep{3pc}\xymatrix{\oc(A)  \! \otimes \! \oc(B)  \! \otimes \! \oc(C)  \ar[r]^-{\mathsf{m}_{A,B} \otimes 1_{\oc(C)}} \ar[d]_-{1_{\oc(A)} \otimes \mathsf{m}_{B,C}} & \oc(A \otimes B) \! \otimes \! \oc(C) \ar[d]^-{\mathsf{m}_{A \otimes B, C}}& \oc(A)  \ar@{=}[dr] \ar[d]_-{\mathsf{m}_I \otimes 1_{\oc(A)}} \ar[r]^-{1_{\oc(A)} \otimes \mathsf{m}_I} & \oc(A)  \otimes  \oc(I) \ar[d]^-{\mathsf{m}_{A,I}} \\  
      \oc(A)  \otimes  \oc(B  \otimes C) \ar[r]_-{\mathsf{m}_{A,B \otimes C}} & \oc(A \otimes B \otimes C) & \oc(I) \otimes  \oc(A) \ar[r]_-{\mathsf{m}_{I,A}} & \oc(A)  } \\
      \xymatrixrowsep{1.75pc}\xymatrixcolsep{5pc}\xymatrix{ \oc(A)  \otimes  \oc(B)  \ar[r]^-{\sigma_{\oc(A), \oc(B)}}  \ar[d]_-{\mathsf{m}_{A,B}} & \oc(B)  \otimes  \oc(A) \ar[d]^-{\mathsf{m}_{B,A}}\\  
 \oc(A \otimes B) \ar[r]_-{\oc(\sigma_{A,B})} &  \oc(B \otimes A)  } \end{gathered}\end{equation}

 \subsection{Diagrams for Monoidal Coalgebra Modality}\label{sec:diagramsformcoalg}
 \begin{equation}\begin{gathered}\label{diag:monoidalcomonad}\xymatrixrowsep{1.75pc}\xymatrixcolsep{4pc}\xymatrix{\oc(A)  \otimes  \oc(B)  \ar[d]_-{\delta_A \otimes \delta_B} \ar[r]^-{\mathsf{m}_{A,B}}  & \oc(A  \otimes  B) \ar[dd]^-{\delta_{A \otimes B}} & I \ar[r]^-{ \mathsf{m}_I} \ar[d]_-{\mathsf{m}_I} & \oc(I) \ar[d]^-{\delta_{I}}   \\
  \oc\oc(A)  \otimes  \oc\oc(B)  \ar[d]_-{\mathsf{m}_{\oc(A),\oc(B)}}  & &\oc(I) \ar[r]_-{\oc(\mathsf{m}_I)}& \oc \oc(I)   \\
  \oc (\oc(A)  \otimes  \oc(B))  \ar[r]_-{\oc(\mathsf{m}_{A,B})} & \oc\oc(A  \otimes  B)
  } \\
  \xymatrixrowsep{1.75pc}\xymatrixcolsep{5pc}\xymatrix{ \oc(A)  \otimes  \oc(B) \ar[dr]_-{\varepsilon_A \otimes \varepsilon_B} \ar[r]^-{\mathsf{m}_{A,B}} & \oc(A  \otimes  B) \ar[d]^-{\varepsilon_{A \otimes B}} & I \ar@{=}[dr]^-{} \ar[r]^-{ \mathsf{m}_I} & \oc(I) \ar[d]^-{\varepsilon_A} \\
 & A  \otimes  B  & &  I  
  } \end{gathered}\end{equation}
  \begin{equation}\begin{gathered}\label{diag:monoidalcomonoid}\xymatrixrowsep{1.75pc}\xymatrixcolsep{3pc}\xymatrix{\oc(A)  \otimes  \oc(B)  \ar[r]^-{\Delta_A \otimes \Delta_B} \ar[dd]_-{\mathsf{m}_{A,B}}  & \oc(A)  \otimes  \oc(A)  \otimes  \oc(B)  \otimes  \oc(B) \ar[d]^-{1_{\oc(A)} \otimes \sigma_{\oc(A),\oc(B)} \otimes 1_{\oc(B)}}  &  I \ar[r]^-{\mathsf{m}_I} \ar[dr]_-{\mathsf{m}_I \otimes \mathsf{m}_I} & \oc(I) \ar[d]^-{\Delta_I}   \\
   &\oc(A)  \otimes  \oc(B) \otimes \oc(A)  \otimes  \oc(B) \ar[d]^-{\mathsf{m}_{A,B} \otimes \mathsf{m}_{A,B}}  & & \oc(I) \otimes  \oc(I) \\
 \oc(A  \otimes  B) \ar[r]_-{\Delta_{A \otimes B}} & \oc(A  \otimes  B) \otimes \oc(A  \otimes  B)  } \\
 \xymatrixrowsep{1.75pc}\xymatrixcolsep{5pc}\xymatrix{ \oc(A)  \otimes  \oc(B) \ar[dr]_-{\mathsf{e}_A \otimes \mathsf{e}_B} \ar[r]^-{\mathsf{m}_{A,B}} & \oc(A  \otimes  B) \ar[d]^-{\mathsf{e}_{A \otimes B}} &   I \ar@{=}[dr]^-{} \ar[r]^-{\mathsf{m}_I} & \oc(I) \ar[d]^-{\mathsf{e}_I}  \\
   & I  & & I }\end{gathered}\end{equation}

      \begin{equation}\begin{gathered} \label{diag:!coalgcomonoid}\xymatrixrowsep{1.75pc}\xymatrixcolsep{2.5pc}\xymatrix{\oc(A) \ar[d]_-{\Delta_A} \ar[rr]^-{\delta_A} & & \oc \oc(A) \ar[d]^-{\oc(\Delta_A)} & \oc(A) \ar[d]_-{\mathsf{e}_A} \ar[r]^-{\delta_A} & \oc \oc(A) \ar[d]^-{\oc(\mathsf{e}_A)} \\
    \oc(A) \otimes \oc(A) \ar[r]_-{\delta_A \otimes \delta_A} & \oc \oc(A) \otimes \oc \oc(A) \ar[r]_-{\mathsf{m}_{\oc(A),\oc(A)}} & \oc(\oc(A) \! \otimes \! \oc(A)) &   I \ar[r]_-{\mathsf{m}_I} & \oc(I) 
  }  \end{gathered}\end{equation}

  \subsection{Diagrams for Pre-Coderelictions}\label{sec:precoder}

   \begin{equation}\begin{gathered}\label{diag:precoder}
  \xymatrixrowsep{1.75pc}\xymatrixcolsep{2.5pc}\xymatrix{A  \ar@{=}[dr]_-{} \ar[r]^-{\eta_A} & \oc(A)\ar[d]^-{\varepsilon_A} & \oc(A) \! \otimes \! B  \ar[d]_-{\varepsilon_A \otimes 1_B} \ar[r]^-{1_{\oc(A)}\otimes \eta_B} & \oc(A) \! \otimes \!  \oc (B) \ar[d]^-{\mathsf{m}_{A,B}} &  A \! \otimes \! \oc(B)  \ar[d]_-{1_A \otimes \varepsilon_B} \ar[r]^-{\eta_A \otimes 1_{\oc(B)}} & \oc(A) \! \otimes \! \oc (B) \ar[d]^-{\mathsf{m}_{A,B}}  \\
  & A & A \! \otimes \! B \ar[r]_-{\eta_{A \otimes B}} & \oc(A \! \otimes \! B) & A \! \otimes \! B \ar[r]_-{\eta_{A \otimes B}} & \oc(A \! \otimes \! B)  }
\end{gathered}\end{equation}

\subsection{Diagrams for Lemma \ref{lemma:precoder}}\label{sec:lemmaprecoder}

 \begin{equation}\begin{gathered}
  \xymatrixrowsep{1.75pc}\xymatrixcolsep{5pc}\xymatrix{\oc(A)  \ar[rr]^-{\eta_I \otimes 1_{\oc(A)}}  \ar[dr]^-{\varepsilon_A}  \ar[dd]_-{1_{\oc(A)} \otimes \eta_I}  && \oc(I) \otimes \oc(A) \ar[dd]^-{\mathsf{m}_{I,A}} \\
  & A  \ar[dr]^-{\eta_A} \\ 
  \oc(A) \otimes \oc(I)  \ar[rr]_-{\mathsf{m}_{A,I}} && \oc(A)  } 
\end{gathered}\end{equation}
 \begin{equation}\begin{gathered}
  \xymatrixrowsep{1.75pc}\xymatrixcolsep{5pc}\xymatrix{A \otimes B  \ar[r]^-{\eta_A \otimes \eta_B}  \ar[dr]_-{\eta_{A \otimes B}} & \oc(A) \otimes \oc(B) \ar[d]^-{\mathsf{m}_{A,B}} \\
  & \oc(A \otimes B)  } 
\end{gathered}\end{equation}

\subsection{Diagrams for Bialgebra Modalities} \label{sec:diagbialgmod}
\begin{equation}\begin{gathered}\label{diag:monoid} \xymatrixrowsep{2pc}\xymatrixcolsep{3pc}\xymatrix{ \oc(A)  \otimes  \oc(A)  \otimes  \oc(A)  \ar[r]^-{\nabla_A \otimes 1_{\oc(A)}} \ar[d]_-{1_{\oc(A)} \otimes \nabla_A} & \oc(A)  \otimes  \oc(A) \ar[d]^-{\nabla_A} & \oc(A) \ar@{=}[dr] \ar[r]^-{1_{\oc(A)} \otimes \mathsf{u}_A}  \ar[d]_-{\mathsf{u}_A \otimes 1_{\oc(A)}}   & \oc(A)  \otimes  \oc(A) \ar[d]^-{\nabla_A}  \\
      \oc(A) \otimes  \oc(A) \ar[r]_-{\nabla_A} & \oc(A)& \oc(A) \otimes  \oc(A) \ar[r]_-{\nabla_A}  & \oc(A) } \\
      \xymatrixrowsep{1.75pc}\xymatrixcolsep{5pc}\xymatrix{ \oc(A)  \otimes  \oc(A)  \ar[r]^-{\sigma_{\oc(A),\oc(A)}} \ar[dr]_-{\nabla_A}   &  \oc(A)  \otimes  \oc(A) \ar[d]^-{\nabla_A} \\
  & \oc(A) }\end{gathered}\end{equation}
 \begin{equation}\begin{gathered}\label{diag:bimonoid}\xymatrixrowsep{2pc}\xymatrixcolsep{3pc}\xymatrix{\oc(A)  \otimes  \oc(A)  \ar[r]^-{\Delta_A \otimes \Delta_A} \ar[dd]_-{\nabla_A}  & \oc(A)  \otimes  \oc(A)  \otimes  \oc(A)  \otimes  \oc(A) \ar[d]^-{1_{\oc(A)} \otimes \sigma_{\oc(A),\oc(A)} \otimes 1_{\oc(A)}}  &  I \ar[r]^-{\mathsf{u}_A} \ar[dr]_-{\mathsf{u}_A \otimes \mathsf{u}_A} & \oc(A) \ar[d]^-{\Delta_A}   \\
   &\oc(A)  \otimes  \oc(A) \otimes \oc(A)  \otimes  \oc(A) \ar[d]^-{\nabla_A \otimes \nabla_A}  & & \oc(A) \otimes  \oc(A) \\
 \oc(A) \ar[r]_-{\Delta_A} & \oc(A) \otimes \oc(A)  } \\
 \xymatrixrowsep{1.75pc}\xymatrixcolsep{5pc}\xymatrix{ \oc(A)  \otimes  \oc(A) \ar[dr]_-{\mathsf{e}_A \otimes \mathsf{e}_A} \ar[r]^-{\nabla_A} & \oc(A) \ar[d]^-{\mathsf{e}_A} &   I \ar@{=}[dr]^-{} \ar[r]^-{\mathsf{u}_A} & \oc(A) \ar[d]^-{\mathsf{e}_A}  \\
   & I  & & I }\end{gathered}\end{equation}

   \subsection{Diagrams for Monoidal Bialgebra Modalities} \label{sec:diagmonbialgmod}

    \begin{equation}\begin{gathered}\label{diag:nablamonspecial}
  \xymatrixrowsep{1.75pc}\xymatrixcolsep{3.5pc}\xymatrix{\oc(A) \! \otimes \! \oc(B) \otimes \! \oc(B)  \ar[d]|-{\Delta_A \otimes 1_{\oc(B)} \otimes 1_{\oc(B)}} \ar[rr]^-{1_{\oc(A)} \otimes \nabla_B} && \oc(A) \! \otimes \! \oc (B) \ar[dd]^-{\mathsf{m}_{A,B}}  \\
\oc(A) \! \otimes \! \oc(A) \otimes \! \oc(B) \! \otimes \! \oc(B) \ar[d]|-{1_{\oc(A)} \otimes \sigma_{\oc(A),\oc(B)} \otimes 1_{\oc(B)}}  \\  
\oc(A) \!  \otimes \! \oc(B) \!\otimes\! \oc(A)  \!\otimes \! \oc(B) \ar[r]_-{\mathsf{m}_{A,B} \otimes \mathsf{m}_{A,B}} &  \oc(A \!\otimes\! B) \otimes  \oc(A \!\otimes\! B)  \ar[r]_-{\nabla_{A,B}} & \oc(A \otimes B)  } \\ 
  \xymatrixrowsep{1.75pc}\xymatrixcolsep{5pc}\xymatrix{\oc(A) \! \otimes  \!\oc(A) \otimes\! \oc(B)  \ar[d]|-{1_{\oc(A)} \otimes 1_{\oc(A)} \otimes \Delta_B} \ar[rr]^-{ \nabla_A  \otimes 1_{\oc(A)}} && \oc(A) \!\otimes\! \oc (B) \ar[dd]^-{\mathsf{m}_{A,B}}  \\
\oc(A) \!\otimes\! \oc(A) \!\otimes\! \oc(B) \!\otimes \!\oc(B) \ar[d]|-{1_{\oc(A)} \otimes \sigma_{\oc(A),\oc(B)} \otimes 1_{\oc(B)}}  \\  
\oc(A) \! \otimes \! \oc(B)\! \otimes \!\oc(A) \! \otimes  \!\oc(B) \ar[r]_-{\mathsf{m}_{A,B} \otimes \mathsf{m}_{A,B}} &  \oc(A \!\otimes\! B) \! \otimes \! \oc(A \!\otimes\! B)  \ar[r]_-{\nabla_{A,B}} & \oc(A \!\otimes\! B)  } \\ 
  \xymatrixrowsep{1.75pc}\xymatrixcolsep{5pc}\xymatrix{\oc(A)  \ar[d]_-{\mathsf{e}_A} \ar[r]^-{1_{\oc(A)} \otimes \mathsf{u}_B} & \oc(A) \otimes \oc (B) \ar[d]^-{\mathsf{m}_{A,B}}  & \oc(B)  \ar[d]_-{\mathsf{e}_B} \ar[r]^-{\mathsf{u}_A \otimes 1_{\oc(B)}} & \oc(A) \otimes \oc (B) \ar[d]^-{\mathsf{m}_{A,B}}  \\
I \ar[r]_-{\mathsf{u}_{A \otimes B}} & \oc(A \otimes B) & I \ar[r]_-{\mathsf{u}_{A \otimes B}} & \oc(A \otimes B)  } 
\end{gathered}\end{equation}

 \begin{equation}\begin{gathered} \label{diag:!coalgmonoid}\xymatrixrowsep{1.75pc}\xymatrixcolsep{2.5pc}\xymatrix{\oc(A) \otimes \oc(A) \ar[d]_-{\nabla_A}  \ar[r]^-{\delta_A \otimes \delta_A}  & \oc \oc(A) \otimes \oc \oc(A) \ar[r]^-{\mathsf{m}_{\oc(A),\oc(A)}} & \oc(\oc(A) \otimes \oc(A)) \ar[d]^-{\oc(\nabla_A)} & I  \ar[d]_-{\mathsf{u}_A} \ar[r]^-{\mathsf{m}_I} & \oc(I) \ar[d]^-{\oc(\mathsf{u}_A)} \\
    \oc(A) \ar[rr]_-{\delta_A} &  & \oc \oc(A)  &   \oc(A) \ar[r]_-{\delta_A} & \oc\oc(A) 
  }  \end{gathered}\end{equation}

   \subsection{Diagrams for Pre-Additive Bialgebra Modalities} \label{sec:diagpreaddbialgmod}

   \begin{equation}\begin{gathered} \label{diag:ep-monoid}\xymatrixrowsep{1.75pc}\xymatrixcolsep{5pc}\xymatrix{ \oc(A)  \otimes  \oc(A)   \ar[r]^-{\nabla_A} \ar[dr]_-{\substack{(\varepsilon_A \otimes \mathsf{e}_A) \\ + (\mathsf{e}_A \otimes \varepsilon_A)}} & \oc(A) \ar[d]^-{\varepsilon_A} & I  \ar[r]^-{\mathsf{u}_A}  \ar[dr]_-{0} &  \oc(A)  \ar[d]^-{\varepsilon_A} \\
     & A & & A }\end{gathered}\end{equation}
     
   \subsection{Diagrams for Convolution Bialgebra Modalities} \label{sec:diagaddbialgmod}

     \begin{equation}\begin{gathered} \label{diag:add-bialg}\xymatrixrowsep{1.75pc}\xymatrixcolsep{5pc}\xymatrix{ \oc(A)   \ar[r]^-{\oc(f+g)} \ar[d]_-{\Delta_A} & \oc(A) & \oc(A) \ar[rr]^-{\oc(0)}  \ar[dr]_-{\mathsf{e}_A} &&  \oc(B)   \\
   \oc(A) \otimes \oc(A) \ar[r]_-{\oc(f) \otimes \oc(g)} & \oc(B) \ar[u]_-{\nabla_A} & & I \ar[ur]_-{\mathsf{u}_B} }\end{gathered}\end{equation}

   \subsection{Diagrams for Correspondence between Monoidal Coalgebra Modalities to an Additive Bialgebra Modalities}\label{sec:moncoalg-addbialg}

   \begin{align}
\begin{array}[c]{c} \nabla_A 
 \end{array} := \left( \begin{array}[c]{c} \xymatrixrowsep{1.75pc}\xymatrixcolsep{7pc}\xymatrix{ \oc(A) \otimes \oc(A)  \ar[r]^-{\delta_A \otimes \delta_A}  & \oc\oc(A) \otimes \oc\oc(A) \ar[r]^-{\mathsf{m}_{\oc(A),\oc(A)}} & \\
 \oc\left( \oc(A) \otimes \oc(A) \right)  \ar[r]^-{\oc\left( (\varepsilon_A \otimes \mathsf{e}_A) + (\mathsf{e}_A \otimes \varepsilon_A) \right)} & \oc(A) }
 \end{array}\right) 
\end{align} 
\begin{align}
\begin{array}[c]{c} \mathsf{u}_A 
 \end{array} := \begin{array}[c]{c} \xymatrixrowsep{1.75pc}\xymatrixcolsep{5pc}\xymatrix{ I  \ar[r]^-{\mathsf{m}_I}  & \oc(I)  \ar[r]^-{\oc(0)} &\oc(A) }
 \end{array}
\end{align}

  \begin{align}
\begin{array}[c]{c} \mathsf{m}_{A,B} 
 \end{array} := \left( \begin{array}[c]{c} \xymatrixrowsep{1.75pc}\xymatrixcolsep{7pc}\xymatrix{ \oc(A) \otimes \oc(B)  \ar[r]^-{\delta_A \otimes \delta_B}  &\\ 
 \oc\oc(A) \otimes \oc\oc(B) \ar[r]^-{\oc(1_{\oc(A)} \otimes \mathsf{u}_B) \otimes \oc(\mathsf{u}_A \otimes 1_{\oc(A)})} & \\
 \oc\left(\oc(A) \otimes \oc(B) \right) \otimes  \oc\left(\oc(A) \otimes \oc(B) \right) \ar[r]^-{\nabla_{\oc(A) \otimes \oc(B)}} & \\
\oc\left(\oc(A) \otimes \oc(B) \right) \ar[r]^-{\delta_{\oc(A) \otimes \oc(B)}} & \\
\oc\oc\left(\oc(A) \otimes \oc(B) \right)  \ar[r]^-{\oc\left(\Delta_{\oc(A) \otimes \oc(B)} \right)} & \\
\oc\left( \oc\left(\oc(A) \otimes \oc(B) \right) \otimes \oc\left(\oc(A) \otimes \oc(B) \right) \right) \ar[r]^-{\oc\left( \oc(\varepsilon_A \otimes \mathsf{e}_B) \otimes \oc(\mathsf{e}_A \otimes \varepsilon_B) \right)} & \\
\oc(\oc(A) \otimes \oc(B))  \ar[r]^-{\oc\left( \varepsilon_A \otimes \varepsilon_B \right)} & \oc(A \otimes B)  }
 \end{array} \right) 
\end{align}
\begin{align}
\begin{array}[c]{c} \mathsf{m}_I 
 \end{array} := \begin{array}[c]{c} \xymatrixrowsep{1.75pc}\xymatrixcolsep{5pc}\xymatrix{ I  \ar[r]^-{\mathsf{u}_I}  & \oc(I)  \ar[r]^-{\delta_I} &\oc\oc(I) \ar[r]^-{\oc(\mathsf{e}_I)} & \oc(I)  }
 \end{array}
\end{align}

\subsection{Diagram for Hopf Coalgebra Modalities} \label{sec:diagHopf}
\begin{equation}\begin{gathered}\label{diag:Hopf} \xymatrixrowsep{1.75pc}\xymatrixcolsep{3pc}\xymatrix{ & \oc(A) \otimes \oc(A) \ar[rr]^-{1_{\oc(A)} \otimes \mathsf{S}_A} && \oc(A) \otimes \oc(A) \ar[dr]^-{\nabla_A} \\
\oc(A) \ar[ur]^-{\Delta_A}  \ar[dr]_-{\Delta_A} \ar[rr]^-{\mathsf{e}_A} & & I \ar[rr]^-{\mathsf{u}_A} && \oc(A) \\ 
& \oc(A) \otimes \oc(A) \ar[rr]_-{\mathsf{S}_A \otimes 1_{\oc(A)}} && \oc(A) \otimes \oc(A) \ar[ur]_-{\nabla_A} }\end{gathered}\end{equation}

\subsection{Diagrams for Lemma \ref{lemma:Hopf-1}} \label{sec:lemma-Hopf-1}

\begin{equation}\begin{gathered}\label{diag:S-iso} \xymatrixrowsep{1.75pc}\xymatrixcolsep{5pc}\xymatrix{ \oc(A)  \ar[r]^-{\mathsf{S}_A} \ar@{=}[dr]_-{} & \oc(A) \ar[d]^-{\mathsf{S}_A} \\
& \oc(A) }\end{gathered}\end{equation}

\begin{equation}\begin{gathered}\label{diag:S-comonoid} \xymatrixrowsep{1.75pc}\xymatrixcolsep{5pc}\xymatrix{ \oc(A)  \ar[r]^-{\mathsf{S}_A} \ar[d]_-{\Delta_A} & \oc(A) \ar[d]^-{\Delta_A} & \oc(A)  \ar[r]^-{\mathsf{S}_A} \ar[dr]_-{\mathsf{e}_A} & \oc(A) \ar[d]^-{\mathsf{e}_A} \\ 
\oc(A) \otimes \oc(A) \ar[r]_-{\mathsf{S}_A \otimes \mathsf{S}_A} &  \oc(A) \otimes \oc(A) & & I  }\end{gathered}\end{equation}

\begin{equation}\begin{gathered}\label{diag:S-monoid} \xymatrixrowsep{1.75pc}\xymatrixcolsep{5pc}\xymatrix{ \oc(A) \otimes \oc(A) \ar[r]^-{\mathsf{S}_A \otimes \mathsf{S}_A} \ar[d]_-{\nabla_A} & \oc(A) \otimes \oc(A) \ar[d]^-{\nabla_A} & I \ar[r]^-{\mathsf{u}_A} \ar[dr]_-{\mathsf{u}_A} & \oc(A) \ar[d]^-{\mathsf{S}_A} \\ 
  \oc(A) \ar[r]_-{\mathsf{S}_A} & \oc(A) & & \oc(A)  }\end{gathered}\end{equation}

\subsection{Diagrams for Monoidal Hopf Coalgebra Modalities} \label{sec:Hopf-monoidal}

\begin{equation}\begin{gathered}\label{diag:S-!coalg} \xymatrixrowsep{1.75pc}\xymatrixcolsep{5pc}\xymatrix{ \oc(A)  \ar[r]^-{\mathsf{S}_A} \ar[d]_-{\delta_A} & \oc(A) \ar[d]^-{\delta_A} \\
\oc\oc(A) \ar[r]_-{\oc(\mathsf{S}_A)} & \oc\oc(A) }\end{gathered}\end{equation}

\begin{equation}\begin{gathered}\label{diag:Hopf-monoidal} \xymatrixrowsep{1.75pc}\xymatrixcolsep{5pc}\xymatrix{ & \oc(A) \otimes \oc(B) \ar[dr]^-{\mathsf{m}_{A,B}} \\
\oc(A) \otimes \oc(B) \ar[ur]^-{1_{\oc(A)} \otimes \mathsf{S}_B}  \ar[dr]_-{\mathsf{S}_A \otimes 1_{\oc(B)}} \ar[r]^-{\mathsf{m}_{A,B}}  & \oc(A \otimes B) \ar[r]^-{\mathsf{S}_{A \otimes B}} & \oc(A \otimes B) \\ 
& \oc(A) \otimes \oc(B)  \ar[ur]_-{\mathsf{m}_{A,B}} }\end{gathered}\end{equation}

\subsection{Diagrams for Lemma \ref{lemma:Hopf-2}} \label{sec:lemma-Hopf-2}

\begin{equation}\begin{gathered}\label{diag:-f} \xymatrixrowsep{1.75pc}\xymatrixcolsep{2.5pc}\xymatrix{ & & \oc(A) \otimes \oc(A) \ar[r]^-{1_{\oc(A)} \otimes \mathsf{S}_A}
\oc(A) & \oc(A) \otimes \oc(A)  \ar[r]^-{\oc(f) \otimes \oc(g)} & \oc(B) \otimes \oc(B) \ar[d]^-{\nabla_B} \\
\oc(A) \ar[r]^-{\mathsf{S}_A} \ar[d]_-{\oc(f)} \ar[dr]^-{\oc(-f)} & \oc(A) \ar[d]^-{\oc(f)} & \oc(A) \ar[rr]^-{\oc(f-g)} \ar[u]^-{\Delta_A} \ar[d]_-{\Delta_A}  &&  \oc(B) \\
\oc(B) \ar[r]_-{\mathsf{S}_B} & \oc(B) & \oc(A) \otimes \oc(A) \ar[r]_-{\oc(f) \otimes \oc(g)} & \oc(B) \otimes \oc(B)\ar[r]_-{1_{\oc(B)} \otimes \mathsf{S}_B} & \oc(B) \otimes \oc(B) \ar[u]_-{\nabla_B}  }\end{gathered}\end{equation}

\begin{equation}\begin{gathered}\label{diag:S-epsilon} \xymatrixrowsep{1.75pc}\xymatrixcolsep{5pc}\xymatrix{ \oc(A)  \ar[r]^-{\mathsf{S}_A} \ar[dr]_-{-\varepsilon_A} & \oc(A) \ar[d]^-{\varepsilon_A} \\
& A }\end{gathered}\end{equation}

 \subsection{Diagrams for Differential Modalities}\label{sec:deriving}

    \begin{equation}\begin{gathered}\label{diag:deriving}
  \xymatrixrowsep{1.75pc}\xymatrixcolsep{3pc}\xymatrix{\oc(A) \otimes A  \ar[dr]_-{0} \ar[r]^-{\mathsf{d}_A} & \oc(A) \ar[d]^-{\mathsf{e}_A} &   \\
  & I } \\
   \xymatrixrowsep{1.75pc}\xymatrixcolsep{5pc}\xymatrix{ \oc(A) \otimes A \ar[rrr]^-{\mathsf{d}_A} \ar[d]_-{\Delta_A \otimes 1_A} &&& \oc(A) \ar[d]^-{\Delta_A} \\ 
  \oc(A) \otimes \oc(A) \otimes A \ar[rrr]_-{(1_{\oc(A)} \otimes \mathsf{d}_A) + (1_{\oc(A)} \otimes \sigma_{\oc(A),A});(\mathsf{d}_A \otimes 1_{\oc(A)}) } &&& \oc(A) \otimes \oc(A) } \\
   \xymatrixrowsep{1.75pc}\xymatrixcolsep{3pc}\xymatrix{ \oc(A) \otimes A \ar[r]^-{\mathsf{d}_A} \ar[dr]_-{\mathsf{e}_A \otimes 1_A} & \oc(A)  \ar[d]^-{\varepsilon_A}   \\
& A } \\ 
  \xymatrixrowsep{1.75pc}\xymatrixcolsep{3pc}\xymatrix{ \oc(A) \otimes A \ar[rr]^-{\mathsf{d}_A} \ar[d]_-{\Delta_A \otimes 1_A} && \oc(A) \ar[d]^-{\delta_A} \\ 
\oc(A) \otimes \oc(A) \otimes A \ar[r]_-{\delta_A \otimes \mathsf{d}_A} & \oc\oc(A) \otimes \oc(A) \ar[r]_-{\mathsf{d}_{\oc(A)}} & \oc(A) } \\ 
  \xymatrixrowsep{1.75pc}\xymatrixcolsep{3pc}\xymatrix{ \oc(A) \otimes A \otimes A \ar[d]_-{1_{\oc(A)} \otimes \sigma_{A,A}} \ar[rr]^-{\mathsf{d}_A \otimes 1_A} && \oc(A) \otimes A \ar[d]^-{\mathsf{d}_A} \\
  \oc(A) \otimes A \otimes A \ar[r]_-{\mathsf{d}_A \otimes 1_A} & \oc(A) \otimes A \ar[r]_-{\mathsf{d}_A} & \oc(A) 
}
\end{gathered}\end{equation}

\subsection{Diagrams for Monoidal Differential Modalities}\label{sec:mon-deriving}

  \begin{equation}\begin{gathered}\label{diag:diff-mon}
  \xymatrixrowsep{2pc}\xymatrixcolsep{5pc}\xymatrix{\oc(A)  \otimes  \oc(B) \otimes B \ar[d]|-{\Delta_A \otimes 1_{\oc(B)} \otimes 1_{B}} \ar[rr]^-{1_{\oc(A)} \otimes \mathsf{d}_B} && \oc(A) \otimes \oc (B) \ar[ddd]^-{\mathsf{m}_{A,B}}  \\
\oc(A) \otimes \oc(A) \otimes \oc(B) \otimes B \ar[d]|-{1_{\oc(A)} \otimes \varepsilon_A \otimes 1_{\oc(B)} \otimes 1_B}  \\ 
\oc(A) \otimes A \otimes \oc(B) \otimes B \ar[d]|-{1_{\oc(A)} \otimes \sigma_{A,\oc(B)} \otimes 1_{B}}  \\ 
\oc(A)  \otimes  \oc(B) \otimes A  \otimes B \ar[r]_-{\mathsf{m}_{A,B} \otimes 1_A \otimes 1_B} &  \oc(A \otimes B) \otimes  A \otimes B  \ar[r]_-{\mathsf{d}_{A,B}} & \oc(A \otimes B)  } \\ 
   \xymatrixrowsep{2pc}\xymatrixcolsep{5pc}\xymatrix{\oc(A) \otimes A  \otimes  \oc(B) \ar[d]|-{1_{\oc(A)} \otimes 1_A \otimes \Delta_B} \ar[rr]^-{\mathsf{d}_A \otimes 1_{\oc(B)}} && \oc(A) \otimes \oc (B) \ar[ddd]^-{\mathsf{m}_{A,B}}  \\
\oc(A) \otimes A \otimes \oc(B) \otimes \oc(B) \ar[d]|-{1_{\oc(A)} \otimes 1_A \otimes 1_{\oc(B)} \otimes \varepsilon_B}  \\ 
\oc(A) \otimes A \otimes \oc(B) \otimes B \ar[d]|-{1_{\oc(A)} \otimes \sigma_{A,\oc(B)} \otimes 1_{B}}  \\ 
\oc(A)  \otimes  \oc(B) \otimes A  \otimes B \ar[r]_-{\mathsf{m}_{A,B} \otimes 1_A \otimes 1_B} &  \oc(A \otimes B) \otimes  A \otimes B  \ar[r]_-{\mathsf{d}_{A,B}} & \oc(A \otimes B)  } \\ 
  \xymatrixrowsep{1.75pc}\xymatrixcolsep{5pc}\xymatrix{\oc(A) \otimes \oc(A) \otimes A  \ar[d]_-{\nabla_A \otimes 1_A} \ar[r]^-{1_{\oc(A)} \otimes \mathsf{d}_A}  & \oc(A) \otimes \oc (A) \ar[d]^-{\nabla}   \\
\oc(A) \otimes A \ar[r]_-{\mathsf{d}_{A}} & \oc(A)  } 
\end{gathered}\end{equation}

\subsection{Diagrams for Correspondence between Coderelictions and Deriving Transformations}\label{sec:d-eta}

 \begin{align}
\begin{array}[c]{c} \mathsf{d}_A 
 \end{array} := \begin{array}[c]{c} \xymatrixrowsep{1.75pc}\xymatrixcolsep{5pc}\xymatrix{ \oc(A) \otimes A \ar[r]^-{1_{\oc(A)} \otimes \eta_A}  & \oc(A) \otimes \oc(A) \ar[r]^-{\nabla_A} & \oc(A) }
 \end{array}
\end{align}
\begin{align}
\begin{array}[c]{c} \eta_A 
 \end{array} := \begin{array}[c]{c} \xymatrixrowsep{1.75pc}\xymatrixcolsep{5pc}\xymatrix{A \ar[r]^-{\mathsf{u}_A \otimes 1_A}  & \oc(A) \otimes A \ar[r]^-{\mathsf{d}_A} & \oc(A) }
 \end{array}
\end{align}

 \subsection{Diagrams for Coderelictions}\label{sec:coder}

   \begin{equation}\begin{gathered}\label{diag:coder}
  \xymatrixrowsep{1.75pc}\xymatrixcolsep{5pc}\xymatrix{A  \ar[dr]_-{0} \ar[r]^-{\mathsf{e}_A} & \oc(A)\ar[d]^-{\eta_A} & A \ar[dr]_-{\substack{(\eta_A \otimes \mathsf{u}_A) \\ + (\mathsf{u}_A \otimes \eta_A)} } \ar[r]^-{\eta_A} & \oc(A) \ar[d]^-{\Delta_A}   \\
  & I & & \oc(A) \otimes \oc(A)  } \\
   \xymatrixrowsep{1.75pc}\xymatrixcolsep{4pc}\xymatrix{ A \ar[d]_-{\mathsf{u}_A \otimes \eta_A} \ar[rr]^-{\eta_A} && \oc(A) \ar[d]^-{\delta_A} \\
  \oc(A) \otimes \oc(A) \ar[r]_-{\delta_A \otimes \eta_{\oc(A)}} & \oc\oc(A) \otimes \oc\oc(A)  \ar[r]_-{\nabla_{\oc(A)}} &  \oc\oc(A)
  }
\end{gathered}\end{equation}

 \subsection{Diagrams for Proposition \ref{prop:der-neg}}\label{sec:der-neg}

\begin{equation}\begin{gathered}\label{diag:S-der} \xymatrixrowsep{1.75pc}\xymatrixcolsep{5pc}\xymatrix{ \oc(A) \otimes A \ar[r]^-{\mathsf{d}_A} \ar[d]_-{\mathsf{S}_A \otimes 1_A} & \oc(A) \ar[d]^-{\mathsf{S}_A} & A  \ar[r]^-{\eta_A} \ar[dr]_-{-\eta_A} & \oc(A) \ar[d]^-{\mathsf{S}_A} \\
\oc(A) \otimes A \ar[r]_-{-\mathsf{d}_A} & \oc(A) & & \oc(A) }\end{gathered}\end{equation}

\end{document}